# A DEGENERATION FORMULA OF GW-INVARIANTS

JUN LI

ABSTRACT. This is the sequel to the paper [Li]. In this paper, we construct the virtual moduli cycles of the degeneration of the moduli of stable morphisms constructed in [Li]. We also construct the virtual moduli cycles of the moduli of relative stable morphisms of a pair of a smooth divisor in a smooth variety. Based on these, we prove a degeneration formula of the Gromov-Witten invariants.

## 0. INTRODUCTION

This is the second part of the project initiated in [Li].

Like Donaldson invariants of 4-manifolds, Gromov-Witten invariants are intersection theories on the moduli spaces of stable morphisms to varieties or symplectic manifolds. However, for Gromov-Witten invariants one needs to use the virtual intersection theories to define the Gromov-Witten invariants. Namely the intersection theories defined based on the virtual moduli cycles. Such cycles were first constructed by Tian and the author [LT1, LT2] for algebraic varieties, and an alternative construction was achieved by Behrend-Fantechi [Beh, BF]. For Gromov-Witten invariants of general symplectic manifolds see [Ru1, RT, FO, LT3, Ru2, Si1] and the equivalence of these constructions see [LT4, Si2, KKP]

The goal of this project is to investigate the structure of the Gromov-Witten invariants when $X$ degenerates to a singular scheme. Specifically, we will prove a degeneration formula of the Gromov-Witten invariants in algebraic geometry. This is the analogous of the Donaldson-Floer theory in gauge theory. Such a degeneration theory was investigated by several groups using analysis [EGH, IP1, IP2, LR, Tia].

Here is the situation we will study in this paper. Let $W \to C$ be a family of projective schemes over a smooth curve $C$ with $0 \in C$ a distinguished point so that the fibers $W_t$ over $t \neq 0 \in C$ are smooth varieties and the special fiber $W_0$ has two smooth irreducible components $Y_1$ and $Y_2$ intersecting transversally along a connected smooth divisor $D \in W_0$. We will call $Y_i^{\mathrm{rel}} \triangleq (Y_i, D_i)$, where $D_i = D \subset Y_i$, the relative pairs after decomposing $W_0$. In this paper, we will construct the Gromov-Witten invariants of $W_t$ for all $t \in C$; We will construct the relative Gromov-Witten invariants of any relative pair $Z^{\mathrm{rel}} \triangleq (Z, D)$ of a smooth projective variety and a smooth divisor; We will prove that the Gromov-Witten invariants of $W_t$ are locally constant for $t \in C$; Finally we will prove a degeneration formula relating the Gromov-Witten invariants of $W_0$ with the relative Gromov-Witten invariants of $Y_1^{\mathrm{rel}}$ and $Y_2^{\mathrm{rel}}$, one in cycle form and the other in numerical form.

*Date*: October 1, 2001.

Supported partially by NSF grants.

The first part of this paper is devoted to construct the Gromov-Witten invariants of the singular scheme $W_0$. To achieve this, we need to construct a new relative moduli space of stable morphisms to $W/C$. This is achieved by refining the notion of stable morphisms, worked out in the first part of this project [Li]. Recall that in [Li] we constructed a stack $\mathfrak{W}$ that includes all expanded degenerations of $W$. We then introduced the notion of pre-stable, pre-deformable and stable morphisms to $\mathfrak{W}$. For simplicity, we denote by $\Gamma$ the triple consisting of the genus, the number of marked points and the degree of the stable morphisms. The main result of [Li] is the existence of the moduli space $\mathfrak{M}(\mathfrak{W}, \Gamma)$ of stable morphisms to $\mathfrak{W}$ of topological type $\Gamma$ as a Deligne-Mumford stack, separated and proper over $C$. Applying parallel construction to a relative pair $Z^{\mathrm{rel}}$ we constructed the stack of expanded relative pairs $\mathfrak{Z}^{\mathrm{rel}}$. We then defined the notion of relative stable morphisms to $\mathfrak{Z}^{\mathrm{rel}}$ and showed that the moduli $\mathfrak{M}(\mathfrak{Z}^{\mathrm{rel}}, \Gamma)$ of relative stable morphisms to $\mathfrak{Z}^{\mathrm{rel}}$ with topological type[1] $\Gamma$ is also a separated and proper Deligne-Mumford stack.

In this paper, we first construct the standard obstruction theory of $\mathfrak{M}(\mathfrak{W}, \Gamma)$, $\mathfrak{M}(\mathfrak{W}_0, \Gamma) = \mathfrak{M}(\mathfrak{W}, \Gamma) \times_C 0$ and $\mathfrak{M}(\mathfrak{Z}^{\mathrm{rel}}, \Gamma)$. We then show that they are all perfect, thus allow us to construct their respective virtual moduli cycles. Based on the virtual moduli cycle $[\mathfrak{M}(\mathfrak{W}_0, \Gamma)]^{\mathrm{virt}}$, we define the Gromov-Witten invariants of $W_0$

$$\Psi_\Gamma^{W_0} : H^*(W_0)^{\times k} \times H^*(\mathfrak{M}_{g,k}) \longrightarrow \mathbb{Q}$$

in the standard way, where $\Gamma = (g, k, A)$. Similarly, for $Z^{\mathrm{rel}}$ we define its relative Gromov-Witten invariants

$$\Psi_\Gamma^{Z^{\mathrm{rel}}} : H^*(Z)^{\times k} \times H^*(\mathfrak{M}_{\Gamma^\circ}) \longrightarrow H_*(D^r)$$

to be

$$\Psi_\Gamma^{Z^{\mathrm{rel}}}(\alpha, \beta) = \mathbf{q}_*\left(\mathrm{ev}^*(\alpha) \cup \pi_\Gamma^*(\beta) \left[\mathfrak{M}(\mathfrak{Z}^{\mathrm{rel}}, \Gamma)\right]^{\mathrm{virt}}\right) \in H_*(D^r).$$

Here $k$ (resp. $r$) is the number of ordinary (resp. distinguished) marked points of the domain curves, ev is the evaluation morphism associated to the ordinary marked points and $\mathbf{q} : \mathfrak{M}(\mathfrak{Z}^{\mathrm{rel}}, \Gamma) \to D^r$ is the evaluation associated to the distinguished marked points. The symbol $\Gamma$ is the topological type of the relative stable morphisms and $\pi_\Gamma : \mathfrak{M}(\mathfrak{Z}^{\mathrm{rel}}, \Gamma) \to \mathfrak{M}_{\Gamma^\circ}$ is the forgetful morphism, where $\mathfrak{M}_{\Gamma^\circ}$ is the moduli of stable nodal curves[2] whose topology is given by the data in $\Gamma$.

The invariants $\Phi_\Gamma^{W_0}$ and $\Phi_{\Gamma'}^{Z^{\mathrm{rel}}}$ have the expected properties. For instance, $\Phi_\Gamma^{W_t}$ is locally constant for $t \in C$, and $\Phi_{\Gamma'}^{Z^{\mathrm{rel}}}$ is invariant under any smooth deformation of $Z^{\mathrm{rel}}$.

The second part of this paper is to derive a degeneration formula of the Gromov-Witten invariants associated to the (degeneration) family $W$. As explained in [Li], we expect to have a formula relating the Gromov-Witten invariants of $W_0$, and hence of $W_t$, to the relative Gromov-Witten invariants of $Y_i^{\mathrm{rel}}$. In this paper, we

---

[1] Here the topological type of a relative stable morphism $f \colon X \to Z^{\mathrm{rel}}$ is given by a weighted graph $\Gamma$. $\Gamma$ consists of vertices, ordered legs and ordered roots, where each leg and root has one end attached to one vertex. Each vertex represents a connected component of the domain $X$, each leg (resp. root) represents a marked point (resp. distinguished marked point) on $X$ and all distinguished marked points are mapped to $D$ under $f$. The vertices and roots of $\Gamma$ are assigned weights, representing the degree of $f$ along the associated component of $X$ and the contact order of $f$ along the normal direction of $D \subset Z$. Here naturally we insist that $f$ maps all distinguished points to $D$.

[2] By which we mean nodal curves having no vector fields, not necessary connected.

prove the following degeneration formula

$$\Phi_\Gamma^{W_t}(\alpha(t), \beta) = \sum_{\eta \in \bar\Omega} \frac{\mathbf{m}(\eta)}{|\mathrm{Eq}(\eta)|} \sum_{j \in K_\eta} \left[ \Psi_{\Gamma_1}^{Y_1^{\mathrm{rel}}}(j_1^*\alpha(0), \beta_{\eta,1,j}) \bullet \Psi_{\Gamma_2}^{Y_2^{\mathrm{rel}}}(j_2^*\alpha(0), \beta_{\eta,2,j}) \right]_0.$$

We will explain the notation momentarily. We call the above the degeneration formula in numerical form. There is a parallel degeneration formula in cycle form

$$[\mathfrak{M}(\mathfrak{W}_0, \Gamma)]^{\mathrm{virt}} = \sum_{\eta \in \bar\Omega} \frac{\mathbf{m}(\eta)}{|\mathrm{Eq}(\eta)|} \Phi_{\eta*} \Delta^! \left( [\mathfrak{M}(\mathfrak{Y}_1^{\mathrm{rel}}, \Gamma_1)]^{\mathrm{virt}} \times [\mathfrak{M}(\mathfrak{Y}_2^{\mathrm{rel}}, \Gamma_2)]^{\mathrm{virt}} \right).$$

The degeneration formula in numerical form is an easy consequence of the degeneration formula in cycle form.

We now explain the notation in the degeneration formulas. Let $\Gamma_1$ be a topological type of a relative stable morphism to $Y_1^{\mathrm{rel}}$. Then by evaluating on the distinguished marked points of the domains of the relative stable morphisms, we obtain the evaluation morphism

$$\mathbf{q}_1 : \mathfrak{M}(\mathfrak{Y}_1^{\mathrm{rel}}, \Gamma_1) \longrightarrow D^r,$$

where $r$ is the number of *distinguished marked points*. Now let $\Gamma_1$ and $\Gamma_2$ be so that they have identical number of roots. (Note that roots are associated to the distinguished points of the domains of relative stable morphisms.) Then we have a pair of evaluation morphisms $\mathbf{q}_1$ and $\mathbf{q}_2$ and thus can form the Cartesian product

$$\mathfrak{M}(\mathfrak{Y}_1^{\mathrm{rel}}, \Gamma_1) \times_{D^r} \mathfrak{M}(\mathfrak{Y}_2^{\mathrm{rel}}, \Gamma_2).$$

Let $D^r \to D^r \times D^r$ be the diagonal morphism. The virtual moduli cycle of the above Cartesian product is

$$\Delta^! \left( [\mathfrak{M}(\mathfrak{Y}_1^{\mathrm{rel}}, \Gamma_1)]^{\mathrm{virt}} \times [\mathfrak{M}(\mathfrak{Y}_2^{\mathrm{rel}}, \Gamma_2)]^{\mathrm{virt}} \right).$$

The set $\bar\Omega$ is the set of equivalence classes of all admissible triples $(\Gamma_1, \Gamma_2, I)$. Here $\Gamma_1$ and $\Gamma_2$ are two weighted graphs associated to two relative stable morphisms. They are required to have identical number of roots and the weights of the $i$-th roots of $\Gamma_1$ and $\Gamma_2$ are identical for all $i$. Hence if $(f_1, f_2)$ be an element in $\mathfrak{M}(\mathfrak{Y}_1^{\mathrm{rel}}, \Gamma_1) \times_{D^r} \mathfrak{M}(\mathfrak{Y}_2^{\mathrm{rel}}, \Gamma_2)$ with $X_i$ the domain of $f_i$, then we can glue the $i$-th distinguished marked points (associated to the $i$-th root) of $X_1$ with the $i$-th distinguished marked point of $X_2$ for all $i$ to obtain a morphism $f : X \to W_0$. As part of the requirement on the triple $(\Gamma_1, \Gamma_2, I)$ the curve $X$ is connected of arithmetic genus $g$ and $f$ is a stable morphism topological type $\Gamma$. Lastly, $I$ is a rule concerning the ordering of the union of the ordinary marked points of $X_1$ and $X_2$. Note that in case $\eta$ has $r$ roots, then any permutation $\sigma \in S_r$ defines a new element $\eta^\sigma$ by reordering the roots of $\eta$ according to $\sigma$. For $\eta_1, \eta_2 \in \Omega$, we say $\eta_1 \sim \eta_2$ if $\eta_1 = \eta_2^\sigma$ for some $\sigma$. The set $\bar\Omega$ is the set of equivalence classes $\Omega / \sim$. The notation $\mathrm{Eq}(\eta)$ appeared in the degeneration formula is the set of all $\sigma \in S_r$ so that $\eta = \eta^\sigma$.

Let $\eta \in \bar\Omega$, then the above construction associates to every pair $(f_1, f_2)$ a stable morphism $f$ in $\mathfrak{M}(\mathfrak{W}, \Gamma)$. This defines a morphism

$$\Phi_\eta : \mathfrak{M}(\mathfrak{Y}_1^{\mathrm{rel}}, \Gamma_1) \times_{D^r} \mathfrak{M}(\mathfrak{Y}_2^{\mathrm{rel}}, \Gamma_2) \longrightarrow \mathfrak{M}(\mathfrak{W}_0, \Gamma).$$

A Lemma in [Li] asserts that $\Phi_\eta$ is a local immersion and the degree of $\Phi_\eta$, as morphism to $\mathrm{Im}(\Phi_\eta)$, is $|\mathrm{Eq}(\eta)|$. Then the degeneration formula in cycle form asserts that the union of

$$\Phi_{\eta*}\Delta^!\big([\mathfrak{M}(\mathfrak{Y}_1^{\mathrm{rel}},\Gamma_1)]^{\mathrm{virt}}\times[\mathfrak{M}(\mathfrak{Y}_2^{\mathrm{rel}},\Gamma_2)]^{\mathrm{virt}}$$

with multiplicity $\mathfrak{m}(\eta)/|\operatorname{Eq}(\eta)|$ is the virtual moduli cycle of $\mathfrak{M}(\mathfrak{W}_0,\Gamma)$. Here $\mathfrak{m}(\eta)$ is the product of the weights of the roots of $\eta$ (or of the $\Gamma_1$ in $\eta$).

We now explain briefly the strategy to prove the degeneration formula. Let $(\mathbf{1},\mathbf{t})$ be the pair of the trivial line bundle on $\mathfrak{M}(\mathfrak{W},\Gamma)$ and the pull back of a section $t\in\Gamma(\mathcal{O}_C)$ so that $t^{-1}(0)$ is the origin $0\in C$. Then the virtual moduli cycle $[\mathfrak{M}(\mathfrak{W}_0,\Gamma)]^{\mathrm{virt}}$ is the intersection

$$[\mathfrak{M}(\mathfrak{W},\Gamma)]^{\mathrm{virt}}\cap\mathbf{t}^{-1}(0)\in A_*\mathfrak{M}(\mathfrak{W}_0,\Gamma).$$

Using the notion of localized top Chern class, this is $c_1(\mathbf{1},\mathbf{t})[\mathfrak{M}(\mathfrak{W},\Gamma)]^{\mathrm{virt}}$.

It turns out that to each $\eta\in\bar{\Omega}$ there associates a pair $(\mathbf{L}_\eta,\mathbf{t}_\eta)$ of a line bundle $\mathbf{L}_\eta$ on $\mathfrak{M}(\mathfrak{W},\Gamma)$ and a section $\mathbf{t}_\eta\in\Gamma(\mathbf{L}_\eta)$ so that $\mathbf{1}\cong\otimes_{\eta\in\Omega}\mathbf{L}_\eta$ as line bundles and under this isomorphism $\mathbf{t}=\Pi_{\eta\in\Omega}\mathbf{t}_\eta$. This says that $\mathfrak{M}(\mathfrak{W}_0,\Gamma)$ "virtually" is a union of normal crossing divisors, each associated to an $\eta\in\Omega$ and is defined by the vanishing of $\mathbf{t}_\eta$. This can be seen as follows: Let $f\in\mathfrak{M}(\mathfrak{W}_0,\Gamma)$ be a general point, say represented by $f\colon X\to W_0$. Then $f_1=f|_{X_1}$ with $X_1=f^{-1}(Y_1)$ defines a relative stable morphism to $Y_1^{\mathrm{rel}}$. Similarly we have the induced relative stable morphism $f_2\colon X_2\to Y_2^{\mathrm{rel}}$. Let $\Gamma_1$ and $\Gamma_2$ be the topological types of $f_1$ and $f_2$, respectively. The fact that $f$ can be reconstructed from the pair $(f_1,f_2)$ provides us a triple $\eta=(\Gamma_1,\Gamma_2,I)$, which belongs to $\Omega$. The general points $f\in\mathfrak{M}(\mathfrak{W}_0,\Gamma)$ that share identical $\eta$ defines a closed subset in $\mathfrak{M}(\mathfrak{W}_0,\Gamma)$. This set is homeomorphic to $\mathbf{t}_\eta^{-1}(0)$. The miracle is that such closed subset carries a natural closed subscheme structure, and is in fact defined by the vanishing of a Cartier divisor $(\mathbf{L}_\eta,\mathbf{t}_\eta)$ and the tensor product of all such $(\mathbf{L}_\eta,\mathbf{t}_\eta)$ is the $(\mathbf{1},\mathbf{t})$ that defines the moduli stack $\mathfrak{M}(\mathfrak{W}_0,\Gamma)$. Then by applying the known identity of the localized top Chern class, we have

$$[\mathfrak{M}(\mathfrak{W}_0,\Gamma)]^{\mathrm{virt}}=c_1(\mathbf{1},\mathbf{t})[\mathfrak{M}(\mathfrak{W},\Gamma)]^{\mathrm{virt}}=\sum_{\eta\in\Omega}c_1(\mathbf{L}_\eta,\mathbf{t}_\eta)[\mathfrak{M}(\mathfrak{W},\Gamma)]^{\mathrm{virt}}.$$

To prove the degeneration formula in cycle form, we need to show that

$$c_1(\mathbf{L}_\eta,\mathbf{t}_\eta)[\mathfrak{M}(\mathfrak{W},\Gamma)]^{\mathrm{virt}}=\frac{\mathfrak{m}(\eta)}{|\operatorname{Eq}(\eta)|}\Phi_{\eta*}\Delta^!\big([\mathfrak{M}(\mathfrak{Y}_1^{\mathrm{rel}},\Gamma_1)]^{\mathrm{virt}}\times[\mathfrak{M}(\mathfrak{Y}_2^{\mathrm{rel}},\Gamma_2)]^{\mathrm{virt}}\big).$$

This is done as follows: First we will show that the vanishing locus $\mathbf{t}_\eta^{-1}(0)\subset\mathfrak{M}(\mathfrak{W},\Gamma)$ is homeomorphic to the image stack of $\Phi_\eta$. More than that, we will show that the cycle $c_1(\mathbf{L}_\eta,\mathbf{t}_\eta)[\mathfrak{M}(\mathfrak{W},\Gamma)]^{\mathrm{virt}}$ is a multiple of the virtual moduli cycle of $\mathfrak{M}(\mathfrak{Y}_1^{\mathrm{rel}},\Gamma_1)\times_{D^r}\mathfrak{M}(\mathfrak{Y}_2^{\mathrm{rel}},\Gamma_2)$, endowed with the obstruction theory induced by the Cartesian product. This leads to the formula above. Here the multiplicity is $\mathfrak{m}(\eta)/|\operatorname{Eq}(\eta)|$, where both the numerator and the denominator are as mentioned before.

Finally, we explain the notations in the degeneration formula. This formula is an immediate corollary of the degeneration formula in cycle form. The only new symbols are $j_{i*}$, $K_\eta$ and $\beta_{\eta,i,j}$. First $j_i\colon Y_i\to W_0$ is the inclusion and hence $j_i^*\alpha(0)$ is the pull back cohomology. Secondly, given any $\eta=(\Gamma_1,\Gamma_2,I)\in\Omega$, we can form the moduli space of stable curves (not necessary connected but with no vector fields) whose topology are determined by $\Gamma_i$. We denote such moduli space by $\mathfrak{M}_{\Gamma_i^0}$. For any pair $(C_1,C_2)\in\mathfrak{M}_{\Gamma_1^0}\times\mathfrak{M}_{\Gamma_2^0}$ we can glue $C_1$ with $C_2$ pairwise along all pairs of the $i$-th distinguished nodes. This defines a morphism $G_\eta : \mathfrak{M}_{\Gamma_1^0} \times \mathfrak{M}_{\Gamma_2^0} \longrightarrow \mathfrak{M}_{g,k}$. The terms $\beta_{\eta,i,j}$ are terms appeared in the Kunneth decomposition,

$$G_\eta^*(\beta) = \sum_{j \in K_\eta} \beta_{\eta,1,j} \boxplus \beta_{\eta,2,j},$$

assuming it exists.

As mentioned in the introduction of [Li], the construction of the moduli stacks $\mathfrak{M}(\mathfrak{W}, \Gamma)$ and $\mathfrak{M}(\mathfrak{Z}^{\mathrm{rel}}, \Gamma)$, and the derivation of the degeneration formula in this paper will be useful in studying several problems in algebraic geometry, some related to mathematical physics. Some of these will be addressed in the future research.

The degeneration formula of Gromov-Witten invariants for symplectic manifolds has been pursued by several groups prior this work. (In symplectic setting this is when a smooth symplectic manifold degenerates to a union of two smooth symplectic manifolds intersecting transversally, called the symplectic sum.) In [Tia], Tian studied the symplectic sums for semi-positive symplectic manifolds and showed how to derive the decomposition formula of the Gromov-Witten invariants in this setting. Later, A.Li-Ruan [LR] worked out a version of the relative Gromov-Witten invariants and the degeneration formula of Gromov-Witten invariants for general symplectic manifolds and symplectic sum. The degeneration formula in numerical form proved in this paper is analogous to the degeneration formula in [LR]. A parallel theory was developed by Ionel-Parker around the same time [IP1, IP2, IP3]. Their formula works for more general cases and is largely analogous to ours. It contains a correction term, which is expected to be trivial when the symplectic sum is along a (holomorphic) divisor. The SFT theory of Eliashberg-Givental-Hofer [EGH] is a very general theory part of which can be interpreted as research along this line. The degeneration formula in cycle form proved here is new.

This paper consists of five sections. In the first section, we work out the obstruction theory of $\mathfrak{M}(\mathfrak{W}, \Gamma)$, $\mathfrak{M}(\mathfrak{Z}^{\mathrm{rel}}, \Gamma)$ and other related moduli stacks. The main result of this section is that the standard obstruction theories of these moduli stacks are perfect. The second section is devoted to construct the virtual moduli cycles of these moduli stacks. We present a modified construction of virtual moduli cycles which allow one to construct such cycles without assuming the existence of global locally free sheaves that resolve the obstruction sheaves, as assumed in [LT2] (and also assumed in [BF] but recently removed in [Kr2]). The Gromov-Witten invariants of $W_0$ and the relative Gromov-Witten invariants of $Z^{\mathrm{rel}}$ are constructed in this section. In section three, we constructed the line bundles with sections $(\mathbf{L}_\eta, \mathbf{t}_\eta)$ mentioned in the introduction and in section four we demonstrated how to derive the degeneration formula, assuming a series of key Lemmas. The last section is devoted to the proof of these key Lemmas. In Appendix, we give a expression of the obstruction space of a closed point in $\mathfrak{M}(\mathfrak{W}, \Gamma)$ and in $\mathfrak{M}(\mathfrak{Z}^{\mathrm{rel}}, \Gamma)$ in terms of some known cohomology groups.

## 1. Deformation theory of log-morphisms

In this section we will first recall the notion of morphisms between schemes with log-structures (in short log-morphism). We will then show that the notion of pre-deformable morphisms introduced in [Li] is a special class of log-morphisms. In the end, we will work out the deformation theory of log-morphisms.

## 1.1. Pre-deformable morphisms and log-morphisms.

Let $f : \mathcal{X} \to W[n]$ be a flat family of pre-deformable morphisms over $S$ with $g : S \to \mathbf{A}^{n+1}$ the morphism underlying $f$ (cf. [Li, §2]). In this subsection, we will give $\mathcal{X}$ and $W[n]$ canonical log structures and show that $f$ is a morphism between schemes with log structures. We will then use the sheaf of log-differentials to describe the deformation of pre-deformable morphisms. All materials concerning schemes with log structures are drawn from the papers of Kato [Ka1, Ka2].

We first recall the notion of logarithmic structures (log structures) on schemes, following [Ka1, Ka2]. Let $X$ be any scheme with $\mathcal{O}_X$ its structure sheaf. We view $\mathcal{O}_X$ as a sheaf of monoids under multiplication.

**Definition 1.1.** *A pre-log structure on $X$ is a homomorphism $\alpha : \mathcal{M} \to \mathcal{O}_X$ of sheaves of monoids, where $\mathcal{M}$ is a sheaf of commutative monoids on the étale site $X_{et}$ of $X$. The pre-log structure $\alpha : \mathcal{M} \to \mathcal{O}_X$ is said to be a log structure if $\alpha$ induces an isomorphism $\alpha : \alpha^{-1}(\mathcal{O}_X^{\times}) \to \mathcal{O}_X^{\times}$, where $\mathcal{O}_X^{\times}$ is the subsheaf of invertible elements in $\mathcal{O}_X$.*

Given a pre-log structure, one can construct canonically an associated log structure $\alpha^a : \mathcal{M}^a \to \mathcal{O}_X$, where $\mathcal{M}^a = (\mathcal{M} \oplus \mathcal{O}_X^{\times})/\alpha^{-1}(\mathcal{O}_X^{\times})$ (cf. [Ka2, §2]). A morphism $(X, \mathcal{M}) \to (Y, \mathcal{N})$ of schemes with log structures is defined to be a pair $(f, h)$ of a morphism of schemes $f : X \to Y$ and a homomorphism $h : f^{-1}(\mathcal{N}) \to \mathcal{M}$ that satisfy the obvious commutativity condition: The composite $f^{-1}(\mathcal{N}) \to \mathcal{M} \to \mathcal{O}_X$ is identical to $f^{-1}(\mathcal{N}) \to f^{-1}(\mathcal{O}_Y) \to \mathcal{O}_X$. For convenience, given a scheme $(X, \mathcal{M})$ with local log structure, we shall abbreviate it to $X^{\dagger}$ with the local log structure $\mathcal{M}$ implicitly understood. Accordingly, we will abbreviate a morphism $(X, \mathcal{M}) \to (Y, \mathcal{N})$ between schemes with local log structures by $f : X^{\dagger} \to Y^{\dagger}$.

A typical example of a scheme with log structure is the log structure of a pair $(X, D)$ of a smooth scheme $X$ and a divisor $D \subset X$ with normal crossing singularities (cf. [Ka1, (1.5)]). Let $\pi_n : W[n] \to \mathbf{A}^{n+1}$ be the expanded family constructed in [Li, §2] associated to a family $W \to C$ and an étale $C \to \mathbf{A}^1$. Then $W[n] \times_{\mathbf{A}^1} 0 \subset W[n]$ is a divisor with normal crossing singularities having $n+2$ irreducible components. Here $\mathbf{A}^{n+1} \to \mathbf{A}^1$ is the morphism $(z_1, \cdots, z_{n+1}) \mapsto z_1 \cdots z_{n+1}$. We let $W[n]^{\dagger}$ and $\mathbf{A}^{n+1\dagger}$ be $W[n]$ and $\mathbf{A}^{n+1}$ endowed with the log-structures induced by the pairs

$$W[n] \times_{\mathbf{A}^1} 0 \subset W[n] \quad \text{and} \quad \mathbf{A}^{n+1} \times_{\mathbf{A}^1} 0 \subset \mathbf{A}^{n+1},$$

respectively.

Let $X^{\dagger}$ be a log-scheme. A chart of $X^{\dagger}$ consists of an étale neighborhood $U \subset X$, a constant sheaf of monoids $P$ and a homomorphism $P \to \mathcal{O}_U$ so that the associated log structure $P^a$ of $P$ is isomorphic to the log structure on $U$. We now describe the charts of $W[n]^{\dagger}$. Recall that the projection $\pi_n$ is smooth away from $n+1$ disjoint smooth codimension 2 subvarieties $\mathbf{D}_1, \cdots, \mathbf{D}_{n+1}$, indexed so that the component $\mathbf{D}_l$ surjects onto the $l$-th coordinate hyperplane $H_l \subset \mathbf{A}^{n+1}$. Let $y \in \mathbf{D}_l$ be any point. A chart of $W[n]$ along $y$ consists of a pair $(\mathcal{W}, \psi)$, where $\mathcal{W}$ is a Zariski open subset of $y \in W[n]$ and $\psi$ is a smooth morphism

$$(1.1) \qquad \psi : \mathcal{W} \longrightarrow \operatorname{Spec} \Bbbk[w_1, w_2] \otimes_{\Bbbk[t_l]} \Gamma(\mathbf{A}^{n+1}) \triangleq \Theta_l$$

so that the canonical projections $\mathcal{W} \to \mathbf{A}^{n+1}$ is the composite of $\psi$ with the projection $\Theta_l \to \mathbf{A}^{n+1}$. Here the two projections in the fiber product are defined via $t_l \mapsto w_1 w_2$ and by viewing $t_l$ as the $l$-th standard coordinate variable of $\mathbf{A}^{n+1}$.

Let $\mathbb{N}^2 \to \mathcal{O}_\mathcal{W}$ be the homomorphism of monoids[3] defined by $e_i \mapsto \psi^*(w_i)$ and let $\mathbb{N}^{n+1} \to \mathcal{O}_{\mathbf{A}^{n+1}}$ be defined via $e_l \mapsto t_l$. We then form the product $\mathbb{N}^2 \times_{\mathbb{N}_l} \mathbb{N}^{n+1}$ where $\mathbb{N}_l \equiv \mathbb{N}$, $\mathbb{N}_l \to \mathbb{N}^2$ is defined by $e \mapsto e_1 + e_2$ while $\mathbb{N}_l \to \mathbb{N}^{n+1}$ is the inclusion as the $l$-th copy in $\mathbb{N}^{n+1}$. Because of the relation $t_l = w_1 w_2$, the homomorphism

$$(1.2) \qquad \mathbb{N}^2 \times_{\mathbb{N}_l} \mathbb{N}^{n+1} \longrightarrow \mathcal{O}_\mathcal{W}$$

induced by $\mathbb{N}^2 \to \mathcal{O}_\mathcal{W}$ and $\mathbb{N}^{n+1} \to \mathcal{O}_{\mathbf{A}^{n+1}} \to \mathcal{O}_\mathcal{W}$ defines a pre-logarithmic structure on $\mathcal{W}$. This defines a chart of $W[n]^\dagger$ near $y \in \mathbf{D}_l$. Similarly, the log-structure of $\mathbf{A}^{n+1}$ is given by the homomorphism $\mathbb{N}^{n+1} \to \mathcal{O}_{\mathbf{A}^{n+1}}$ via $e_l \mapsto t_l$.

Now let $f: \mathcal{X} \to W[n]$ be a flat family of pre-deformable morphisms over $S$, defined in [Li, §2]. Out next step is to show that we can endow $\mathcal{X}/S$ with a log structure so that $f$ becomes a morphism between schemes with log structures. We will use the following convention throughout this paper. We call $p: U \to V$ an étale neighborhood of the family $\mathcal{X}/S$ if $U$ and $V$ are étale neighborhoods of $\mathcal{X}$ and $S$, respectively, such that the diagram

$$(1.3) \qquad \begin{array}{ccc} U & \longrightarrow & \mathcal{X} \\ p \downarrow & & \pi \downarrow \\ V & \longrightarrow & S. \end{array}$$

is commutative.

**Definition 1.2.** *Let $\pi: \mathcal{X} \to S$ be a flat family of nodal curves and let $x \in \mathcal{X}$ be a node of the fibers of the family $\mathcal{X}/S$. A chart of the nodes of $\mathcal{X}/S$ near $x$ consists of an affine étale neighborhood $p: \mathcal{U} \to \mathcal{V}$ of $\mathcal{X}/S$ near $x$, two regular functions $z_1$ and $z_2 \in \Gamma(\mathcal{O}_\mathcal{U})$ and a regular function $s \in \Gamma(\mathcal{O}_\mathcal{V})$ satisfying $z_1 z_2 = p^*(s)$ such that the homomorphism*

$$(1.4) \qquad \phi: \Bbbk[z_1, z_2] \otimes_{\Bbbk[s]} \Gamma(\mathcal{O}_\mathcal{V}) \longrightarrow \Gamma(\mathcal{O}_\mathcal{U})$$

*is an étale homomorphism, that $\{z_1 = z_2 = 0\} \subset \mathcal{U}$ is connected and that the induced homomorphism*

$$\hat{\phi}: \big(\Bbbk[z_1, z_2] \otimes_{\Bbbk[s]} \Gamma(\mathcal{O}_\mathcal{V})\big)\hat{\ } \longrightarrow \Gamma(\mathcal{O}_\mathcal{U})\hat{\ }$$

*is an isomorphism.*

A few remarks are in order. First $\phi$ is defined by viewing $z_1$, $z_2$ and $s$ on one hand as formal variables and on the other hand as regular functions. This should cause no confusion since their roles are clear from the context. Also, the homomorphism $\Bbbk[s] \to \Bbbk[z_1, z_2]$ is defined by $s \mapsto z_1 z_2$. As a convention, in this section we will use $\hat{L}$ and $\hat{h}$ (or $L\hat{\ }$ and $h\hat{\ }$) to denote the $I$-adic completion of the ring $L$ and the image of $h \in L$ under the homomorphism $L \to \hat{L}$, assuming $I = (z_1, z_2)$ is an ideal of $L$.

**Lemma 1.3.** *Let $\pi: \mathcal{X} \to S$ be a flat family of nodal curves and let $x \in \mathcal{X}$ be a node in $\mathcal{X}_y$, where $y = \pi(x)$. Then there exists a chart of the nodes of $\mathcal{X}/S$ near $x$.*

*Proof.* The proof is straightforward and will be omitted. $\qquad \square$

We next recall the notion of morphisms of pure contact. We let $(\mathcal{U}/\mathcal{V}, \phi)$ be a chart of a node of $\mathcal{X}/S$, as in Definition 1.2. For convenience we let $A = \Gamma(\mathcal{O}_\mathcal{V})$ and $R = \Gamma(\mathcal{O}_\mathcal{U})$. As before we will view $z_1$ and $z_2$ (resp. $s$) as regular functions of $\mathcal{U}$ (resp. $\mathcal{V}$), thus are elements in $R$ (resp. $A$). We let $\mathcal{S}$ (resp $\mathcal{T}$) be the set of all $a \in R$ (resp $a \in A$) so that their residue classes $\bar{a} \in R/(z_1, z_2)$ (resp. $\bar{a} \in A/(s)$) are units. The set $\mathcal{S}$ and $\mathcal{T}$ form multiplicative systems in $R$ and $A$ respectively. We let $R_\mathcal{S}$ and $A_\mathcal{T}$ be the localization of the respective rings. As before, we let $\Bbbk[t] \to \Bbbk[w_1, w_2]$ be defined by $t \mapsto w_1 w_2$. We let $\hat{A}$ be the $(s)$-adic completion of $A$.

---

[3]We will denote by $\mathbb{N}$ the additive monoid of non-negative integers, by $\mathbb{N}^k$ the direct product monoid with standard generators $e_1, \cdots, e_k$.

**Definition 1.4.** *Let $\Bbbk[t] \to A$ be a homomorphism and $\varphi \colon \Bbbk[w_1, w_2] \to R$ be a $\Bbbk[t]$-homomorphism. We say $\varphi$ is formally of pure contact of order $m$ if the induced homomorphism $\hat{\varphi} \colon \Bbbk[w_1, w_2] \to \big(\Bbbk[\![z_1, z_2]\!] \otimes_{\Bbbk[\![s]\!]} A\big)\hat{\ }$ is of pure contact of order $m$, as defined in Definition 2.3 in* [Li]. *Namely, there are units $\hat{h}_1, \hat{h}_2 \in (\Bbbk[\![z_1, z_2]\!] \times_{\Bbbk[\![s]\!]} A)\hat{\ }$ satisfying $\hat{h}_1 \hat{h}_2 \in \hat{A}$ so that $\hat{\varphi}(w_i) = \hat{h}_i z_i^m$, possibly after exchanging $w_1$ and $w_2$. We say $\varphi$ is of pure contact of order $m$ if there are units $h_1, h_2 \in R_\mathcal{S}$ satisfying $h_1 h_2 \in A_\mathcal{T}$ such that $\varphi(w_i) = z_i^m h_i$ in $R_\mathcal{S}$, possibly after exchanging $w_1$ and $w_2$.*

We have the following facts whose proof can be found in the Appendix.

**Lemma 1.5.** *The notion of pure contact is independent of the choice of the charts of the nodes of $\mathcal{U}/\mathcal{V}$.*

**Lemma 1.6.** *Let the notation be as in Definition 1.4. Then $\varphi$ is of pure contact if and only if it is formally of pure contact.*

Let us continue to denote by $f$ a pre-deformable morphism $f \colon \mathcal{X} \to W[n]$ over $S$. To define log-structure on $\mathcal{X}$ we need to introduce charts of $f$. Let $x \in \mathcal{X}$ be any closed point so that $f(x) \notin \mathbf{D}$. Here $\mathbf{D}$ is the union of all $\mathbf{D}_l$. A chart of $f$ near $x$ is a triple $(\mathcal{U}_\alpha/\mathcal{V}_\alpha, \mathcal{W}_\alpha, f_\alpha)$ of an étale neighborhood $\mathcal{U}_\alpha/\mathcal{V}_\alpha$ of $x \in \mathcal{X}$, a Zariski neighborhood $\mathcal{W}_\alpha \subset W[n] - \mathbf{D}$ and $f_\alpha = f|_{\mathcal{U}_\alpha}$ so that $f_\alpha(\mathcal{U}) \subset \mathcal{W}_\alpha$. Charts of this type will be called of the first kind. Next let $x \in \mathcal{X}$ be a distinguished node of $\mathcal{X}/S$, namely $f(x) \in \mathbf{D}_{l_\alpha}$ for some $l_\alpha$. We pick a chart $(\mathcal{W}_\alpha, \psi_\alpha)$ of $f(x) \in W[n]$ with $\psi_\alpha$ as in (1.1). Because $f$ is pre-deformable, by Lemma 1.6 we can find a chart of the nodes of $\mathcal{X}/S$ near $x$, say given by $(\mathcal{U}_\alpha/\mathcal{V}_\alpha, \phi_\alpha)$ as in Definition 1.2, so that $\mathcal{U}_\alpha \to \mathcal{X} \to W[n]$ lifts to $f_\alpha \colon \mathcal{U}_\alpha \to \mathcal{W}_\alpha$, and further there are units $h_{\alpha,1}, h_{\alpha,2} \in \Gamma(\mathcal{O}_{\mathcal{U}_\alpha})^\times$, elements $g_{\alpha,1}, \cdots, g_{\alpha,n+1} \in \Gamma(\mathcal{O}_{\mathcal{V}_\alpha})$ and an integer $m_\alpha$ so that

$$(1.5) \qquad h_{\alpha,1} h_{\alpha,2} \in \Gamma(\mathcal{O}_{\mathcal{V}_\alpha}), \quad f_\alpha^*(w_{\alpha,i}) = z_{\alpha,i}^{m_\alpha} h_{\alpha,i} \quad \text{and} \quad f_\alpha^*(t_j) = g_{\alpha,j}$$

[4] for $i = 1, 2$ and $j = 1, \cdots, n+1$. Here we comment that we will not distinguish $w_{\alpha,i}$ (resp. $z_{\alpha,i}$) with $\psi_\alpha^*(w_{\alpha,i})$ (resp. $\phi_\alpha^*(z_{\alpha,i})$). In other words, we will view $w_{\alpha,i} \in \mathcal{O}_{\mathcal{W}_\alpha}$ (resp. $z_{\alpha,i} \in \mathcal{O}_{\mathcal{U}_\alpha}$). Similarly we will view $t_l \in \mathcal{O}_{W[n]}$ via the pull back $\mathcal{O}_{\mathbf{A}^{n+1}} \to \mathcal{O}_{W[n]}$. Note that since $w_{\alpha,1} w_{\alpha,2} = t_{l_\alpha}$ and $z_{\alpha,1} z_{\alpha,2} = s_\alpha$, we must have

$$(1.6) \qquad\qquad g_{\alpha,l_\alpha} = s_\alpha^{m_\alpha}(h_{\alpha,1} h_{\alpha,2}).$$

We will call such triplet $(\mathcal{U}_\alpha/\mathcal{V}_\alpha, \mathcal{W}_\alpha, f_\alpha)$ with $(\phi_\alpha, \psi_\alpha)$ understood a chart of $f$ of the second kind.

**Simplification 1.7.** *In case $(\mathcal{U}_\alpha/\mathcal{V}_\alpha, \mathcal{W}_\alpha, f_\alpha)$ is a chart of the second kind, for simplicity we assume $h_{\alpha,1} \equiv h_{\alpha,2} \equiv 1$. This is possible possibly after an étale base change of $\mathcal{U}_\alpha/\mathcal{V}_\alpha$.*

---

[4] By this we mean it is in the image of the pull back homomorphism $\Gamma(\mathcal{O}_{\mathcal{V}_\alpha}) \to \Gamma(\mathcal{O}_{\mathcal{U}_\alpha})$.

We now cover $\mathcal{X}$ by charts of $f$, say $(\mathcal{U}_\alpha / \mathcal{V}_\alpha, \mathcal{W}_\alpha, f_\alpha)$ indexed by $\alpha \in \Lambda$, of the first or the second kinds satisfying the Simplification 1.7. Let $\alpha$ be a chart of the second kind with $\phi_\alpha$ and $\psi_\alpha$ understood. We then let $M_\alpha^0 = \mathbb{N}^2$ (resp. $N_\alpha^0 = \mathbb{N}$) and let $M_\alpha^0 \to \mathcal{O}_{\mathcal{U}_\alpha}$ (resp. $N_\alpha \to \mathcal{O}_{\mathcal{V}_\alpha}$) be the pre-log structure defined by $e_{\alpha,i} \mapsto z_{\alpha,i}$ (resp. $e_\alpha \mapsto s_\alpha$). Note that $N_\alpha^0 \to M_\alpha^0$ defined by $e_\alpha \mapsto e_{\alpha,1} + e_{\alpha,2}$ makes the projection $\mathcal{U}_\alpha \to \mathcal{V}_\alpha$ a morphism between schemes with (their respective associated) log structures. We now define the desired log structures on $S$. Let $\xi \in S$ be any closed point and let $\mathcal{X}_\xi$ be the fiber of $\mathcal{X}$ over $\xi$ with $\Sigma = f^{-1}(\mathbf{D}) \cap \mathcal{X}_\xi$ be the set of distinguished nodes of $\mathcal{X}_\xi$. We then let $\Lambda_0$ be a collection of $\alpha \in \Lambda$ so that $\{\mathcal{U}_\alpha\}_{\alpha \in \Lambda_0}$ covers a neighborhood of $\mathcal{X}_\xi \subset \mathcal{X}$. We let $K_l$ be those $\alpha$ so that $f^{-1}(\mathbf{D}_l) \cap \mathcal{U}_\alpha \ne \emptyset$. We let $K$ be the union of $K_1, \cdots, K_{n+1}$. By eliminating redundant $\alpha$ from $\Lambda_0$ we can assume that each node of $\mathcal{X}_\xi \cap f^{-1}(\mathbf{D}_l)$ is covered by at most one $\alpha \in K_l$. We then pick an étale neighborhood $\mathcal{V}$ of $\xi \in S$ so that $\mathcal{V} \to S$ factor through $\mathcal{V}_\alpha \to S$ for all $\alpha \in \Lambda_0$. For each $l \in [n+1]$ we let $\bar{N}_l = \oplus_{\alpha \in K_l} N_\alpha^0$ and let $N_l$ be the quotient (monoid) of $\bar{N}_l$ by the relations $m_\alpha e_\alpha = m_\beta e_\beta$ for all pairs $\alpha, \beta \in K_l$. In case $K_l = \emptyset$, we agree $N_l = \bar{N}_l = \mathbb{N}$. Note that for $\alpha \in K_l$ the homomorphism $\bar{N}_l \to \mathcal{O}_{\mathcal{V}}$ defined by $e_\alpha \mapsto s_\alpha$ and the pull back $\mathcal{O}_{\mathcal{V}_\alpha} \to \mathcal{O}_{\mathcal{V}}$ descends to a homomorphism $N_l \to \mathcal{O}_{\mathcal{V}}$, because of the relations (1.5) and (1.6). When $K_l = \emptyset$, we define $N_l \to \mathcal{O}_{\mathcal{V}}$ via $e \mapsto g^*(t_l)$. We define the pre-log structures $N_{\mathcal{V}}$ and $\bar{N}_{\mathcal{V}}$ on $\mathcal{V}$ to be the ones given by the direct sums

$$(1.7) \qquad N_{\mathcal{V}} \triangleq \oplus_{l=1}^{n+1} N_l \longrightarrow \mathcal{O}_{\mathcal{V}} \quad \text{and} \quad \bar{N}_{\mathcal{V}} \triangleq \oplus_{l=1}^{n+1} \bar{N}_l \longrightarrow \mathcal{O}_{\mathcal{V}}.$$

We denote by $\mathcal{N}_{\mathcal{V}}$ and $\bar{\mathcal{N}}_{\mathcal{V}}$ the associated log structures on $\mathcal{V}$.

We next define the desired pre-log structure on $\mathcal{U}_\alpha$. By replacing $\mathcal{U}_\alpha$ by $\mathcal{U}_\alpha \times_{\mathcal{V}_\alpha} \mathcal{V}$, we can assume $\mathcal{V}_\alpha = \mathcal{V}$. Accordingly we let $N_\alpha$ and $\bar{N}_\alpha \to \mathcal{O}_{\mathcal{V}}$ be that induced by $N_{\mathcal{V}}$ and $\bar{N}_{\mathcal{V}}$. Now let $l$ be so that $\alpha \in K_l$. Recall that $\mathcal{U}_\alpha$ has a pre-log structure given by $M_\alpha^0 \to \mathcal{O}_{\mathcal{U}_\alpha}$. Let $N_\alpha^0 \to M_\alpha^0$ be as before (defined by $e_\alpha \mapsto e_{\alpha,1} + e_{\alpha,2}$). Then we have the obvious homomorphism

$$(1.8) \qquad \bar{M}_\alpha = M_\alpha^0 \times_{N_\alpha^0} \bar{N}_l \oplus \left( \oplus_{l' \ne l} \bar{N}_{l'} \right) \longrightarrow \mathcal{O}_{\mathcal{U}_\alpha}$$

and

$$(1.9) \qquad M_\alpha = M_\alpha^0 \times_{N_\alpha^0} N_l \oplus \left( \oplus_{l' \ne l} N_{l'} \right) \longrightarrow \mathcal{O}_{\mathcal{U}_\alpha}.$$

Here $M_\alpha^0 \times_{N_\alpha^0} \bar{N}_l \to \mathcal{O}_{\mathcal{U}_\alpha}$ is induced by $M_\alpha^0 \to \mathcal{O}_{\mathcal{U}_\alpha}$ and $\bar{N}_l \to \mathcal{O}_{\mathcal{V}_\alpha} \to \mathcal{O}_{\mathcal{U}_\alpha}$ while $\bar{N}_{l'} \to \mathcal{O}_{\mathcal{U}_\alpha}$ is the composite of $N_{l'} \to \mathcal{O}_{\mathcal{V}_\alpha} \to \mathcal{O}_{\mathcal{U}_\alpha}$. The arrow in (1.9) is defined similarly. They define two pre-log structures on $\mathcal{U}_\alpha$. We let $\mathcal{M}_\alpha \to \mathcal{O}_{\mathcal{U}_\alpha}$ and $\bar{\mathcal{M}}_\alpha \to \mathcal{O}_{\mathcal{U}_\alpha}$ be the associated log-structures. Note that the obvious $N_\alpha \to M_\alpha$ and $\bar{M}_\alpha \to \bar{N}_\alpha$ make the projection $\mathcal{U}_\alpha \to \mathcal{V}_\alpha$ a morphism between schemes with respective log-structures.

**Proposition 1.8.** *The log structures $(\mathcal{U}_\alpha, \mathcal{M}_\alpha)$ and $(\mathcal{V}_\alpha, \mathcal{N}_\alpha)$ patch together to form log structures $\mathcal{M}$ on $\mathcal{X}$ and $\mathcal{N}$ on $S$. The same conclusion holds for $(\mathcal{U}_\alpha, \bar{\mathcal{M}}_\alpha)$ and $(\mathcal{V}_\alpha, \bar{\mathcal{N}}_\alpha)$. The collection of homomorphisms $N_\alpha \to M_\alpha$ makes $\mathcal{X}^\dagger \equiv (\mathcal{X}, \mathcal{M})$ a log scheme over $S^\dagger \equiv (S, \mathcal{N})$. Further, the morphism $f$ is naturally a morphism between schemes with log structures $\mathcal{X}^\dagger / S^\dagger \to W[n]^\dagger / \mathbf{A}^{n+1\dagger}$.*

*Proof.* The fact that the so defined log structures on $\mathcal{U}_\alpha$ and $\mathcal{V}_\alpha$ patch together to form log structures on $\mathcal{X}$ and $S$ is obvious. We now study the morphism $f$. We first investigate the morphism $g : S \to \mathbf{A}^{n+1}$ underlying $f$. Let $\xi \in S$ be any closed point with $(\mathcal{V}, N_{\mathcal{V}})$ a chart of $S^\dagger$, constructed before. Recall that the log structure on $\mathbf{A}^{n+1\dagger}$ is given by the pre-log structure $\mathbb{N}^{n+1} \to \Gamma(\mathbf{A}^{n+1})$ via $e_l \mapsto t_l$. To show $g|_{\mathcal{V}} : \mathcal{V} \to \mathbf{A}^{n+1}$ is a morphism between schemes with log structures we need to define a homomorphism $\mathbb{N}^{n+1} \to N_{\mathcal{V}}$ that satisfies the obvious compatibility condition. By definition $N_{\mathcal{V}} = \oplus_l N_l$ and $N_l$ is $\mathbb{N}$ in case $K_l = \emptyset$ and is the quotient of $\bar{N}_l = \oplus_{\alpha \in K_l} N_\alpha^0$ otherwise. In the first case we define $\mathbb{N} \to N_l$ to be the unique isomorphism. In the later case we let $\mathbb{N} \to N_l$ be induced by $1 \mapsto m_\alpha e_\alpha \in N_l$ for some $\alpha \in K_l$. By the definition of $N_l$ such definition is independent of $\alpha \in K_l$. The desired homomorphism $\mathbb{N}^{n+1} \to N_{\mathcal{V}}$ is the direct sum of these $n + 1$ copies $\mathbb{N} \to N_l$. It is direct to check that this defines $\mathcal{V} \to \mathbf{A}^{n+1}$ a morphism between schemes with log structures. By working over a covering of $S$, this concludes that $g$ is a morphism between $S^\dagger \to \mathbf{A}^{n+1\dagger}$.

The definition making $f$ a morphism between $\mathcal{X}^\dagger$ and $W[n]^\dagger$ and compatible to $g : S^\dagger \to \mathbf{A}^{n+1\dagger}$ is similar, relying on the relations (1.5) and (1.6) and the assumption that all $h_{\alpha,i} \equiv 1$. This completes the proof of the Proposition. □

We conclude this subsection by commenting the equivalence of the deformations of pre-deformable morphisms and the deformations of log-morphisms.

**Definition 1.9.** *1. A log-extension of $(\mathcal{V}_\alpha, \mathcal{N}_\alpha)$ by $I$ consists of an extension $\tilde{\mathcal{V}}_\alpha$ of $\mathcal{V}_\alpha$ by $I$ as schemes, and an extension $N_\alpha \to \mathcal{O}_{\tilde{\mathcal{V}}_\alpha}$ of $N_\alpha \to \mathcal{O}_{\mathcal{V}_\alpha}$. We denote such extension by $(\tilde{\mathcal{V}}_\alpha, \tilde{\mathcal{N}}_\alpha)$.*
*2. A flat log-extension of $(\mathcal{U}_\alpha / \mathcal{V}_\alpha, \mathcal{M}_\alpha / \mathcal{N}_\alpha)$ by $I$ consists of an extension $(\tilde{\mathcal{V}}_\alpha, \tilde{\mathcal{N}}_\alpha)$ of $(\mathcal{V}_\alpha, \mathcal{N}_\alpha)$ by $I$, a flat extension $\tilde{\mathcal{U}}_\alpha \to \tilde{\mathcal{V}}_\alpha$ of $\mathcal{U}_\alpha \to \mathcal{V}_\alpha$ and an extension $M_\alpha \to \mathcal{O}_{\tilde{\mathcal{U}}_\alpha}$ of the pre-log structure $M_\alpha \to \mathcal{O}_{\mathcal{U}_\alpha}$ of which the following holds: a. The projection $\tilde{\mathcal{U}}_\alpha \to \tilde{\mathcal{V}}_\alpha$ is a log-morphism under the given $N_\alpha \to M_\alpha$; b. Away from the distinguished nodes of $\mathcal{U}_\alpha$ the pre-log structure $M_\alpha \to \mathcal{O}_{\tilde{\mathcal{U}}_\alpha}$ is the pull back of $N_\alpha \to \mathcal{O}_{\tilde{\mathcal{V}}_\alpha}$; c. Near the distinguished nodes in $\mathcal{U}_\alpha$ the projection $\tilde{\mathcal{U}}_\alpha \to \tilde{\mathcal{V}}_\alpha$ is log-smooth.*

We have the following Lemma which says that extending $f$ as pre-deformable morphism is equivalent to extending $f$ as log-morphisms.

**Lemma 1.10.** *Let $f : \mathcal{X}/S \to W[n]/\mathbf{A}^{n+1}$ be as before with the canonical log structures understood. Let $\tilde{S} \supset S$ be an extension of $S$. Suppose $\tilde{S}^\dagger$ is a log-extension of $S$, $\tilde{\mathcal{X}}^\dagger \to \tilde{S}^\dagger$ is a flat log-extension of $\mathcal{X}/S$ and $\tilde{f} : \tilde{\mathcal{X}}^\dagger / \tilde{S}^\dagger \to W[n]^\dagger / \mathbf{A}^{n+1\dagger}$ is an extension of $f$ as log-morphism. Then $\tilde{f}$ is a pre-deformable extension of $f$ and the log structures on $\tilde{\mathcal{X}}/\tilde{S}$ induced by $\tilde{f}$ coincide with the log structure of $\tilde{\mathcal{X}}^\dagger / \tilde{S}^\dagger$.*

*Proof.* The proof is straightforward. Since we do not need this Lemma in this paper, we shall omit the proof here. □

## 1.2. Deformation of pre-deformable morphisms.

The goal of this subsection is to work out the deformation theory of pre-deformable morphisms. Based on the equivalence Lemma, it is natural to work out the deformation theory of pre-deformable morphisms in the frame work of log-morphisms. However, the deformation theory of log-morphisms worked out in [Ka1, Ka2] deal with the question on how to extend families over $S^\dagger$ to $\tilde{S}^\dagger$. In our situation, the log structure on $S$ for a pre-deformable family $f : \mathcal{X}/S \to W[n]/\mathbf{A}^{n+1}$ relies on the morphism $f$. Hence for any extension of base $S \subset \tilde{S}$ the extension of log structure to $\tilde{S}$ is part of the extension problem. Hence the deformation theory of Kato [Ka2] can not be applied directly. In this subsection, we will work with the deformation of pre-deformable directly.

We begin with recalling some basic notion in deformation theory. Our treatment follows [Art, LT2]. Let $A$ be an $\bar{A}$-algebra. We let $\mathfrak{Tri}_{A/\bar{A}}$ be the category whose objects are all triples $(B, I, \varphi)$ where $B$ are $\bar{A}$-algebras, $I$ are ideals of $B$ such that $I^2 = 0$ and $\varphi$ are $\bar{A}$-homomorphisms $A \to B/I$. Let $\xi = (B, I, \varphi)$ and $\xi' = (B', I', \varphi')$ be two objects in $\mathfrak{Tri}_{A/\bar{A}}$. A morphism from $\xi$ to $\xi'$ consists of an $\bar{A}$-homomorphism $r \colon B \to B'$ so that $r(I) \subset I'$ and $\varphi' = \varphi \circ r_0$ where $r_0 \colon B/I \to B'/I'$ is the induced homomorphism. We let $\mathfrak{Mod}_A$ be the category whose objects are pairs $(B, I)$ where $B$ are $A$-algebras and $I$ are $B$-modules. Morphisms from $(B, I)$ to $(B', I')$ are pairs $(r, \tilde{r})$ where $r \colon B \to B'$ are $A$-homomorphisms and $\tilde{r}$ are $B$-homomorphisms $I \to I'$. We let $\mathfrak{Mod}_A^*$ be the category whose objects are $(v, B, I)$ where $(B, I) \in Ob(\mathfrak{Mod}_A)$ and $v \in I$. Morphisms from $(B, I, v)$ to $(B', I', v')$ are pairs $(r, \tilde{r})$ as in $\mathfrak{Mod}_A$ so that $\tilde{r}(v) = v'$.

We define $\mathfrak{Def}_{A/\bar{A}} \colon \mathfrak{Tri}_{A/\bar{A}} \to (\text{Sets})$ be the functor that associates to any $\xi = (B, I, \varphi)$ the set of all $\bar{A}$-homomorphisms $A \to B$ extending $\varphi \colon A \to B/I$. (In case $\bar{A}$ is understood, we will omit the subscript $\bar{A}$.) It is known that under some mild conditions [Shl] this set admits $B/I$-module structures. Once this is the case, then after fixing a reference element $a \in \mathfrak{Def}_A(\xi)$ we can give $\mathfrak{Def}_A(\xi)$ a natural $B/I$-module structure. In particular, if $B = B/I * I$, which is the trivial extension of $B/I$ by the module $I$, then $\mathfrak{Def}_A(\xi)$ contains the trivial extension $B \to A$ induced by $B = B/I * I \to B/I \to A$. Using this as the reference element, we obtain a canonical module structure on $\mathfrak{Def}_A(\xi)$. Note that there is a natural functor $\mathfrak{Mod}_A \to \mathfrak{Tri}_{A/\bar{A}}$ that sends $(B, I)$ to $(B * I, I, \varphi)$ where $\varphi \colon A \to B$ is the obvious homomorphism. We let $\mathfrak{Def}_A^1 \colon \mathfrak{Mod}_A \to \mathfrak{Mod}_A$ be the composite of $\mathfrak{Def}_A$ with $\mathfrak{Mod}_A \to \mathfrak{Tri}_{A/\bar{A}}$. We call $\mathfrak{Def}_A^1$ the functor of the first order deformations.

Now let $E^\bullet$ be any complex of $A$-modules. Then for any integer $i$ the $i$-th cohomology of $E^\bullet$ defines a functor $\mathfrak{h}^i(E^\bullet) \colon \mathfrak{Mod}_A \to \mathfrak{Mod}_A$ via $(B, I) \mapsto (B, h^i(E^\bullet \otimes_A I))$.

**Definition 1.11.** *Let $S = \operatorname{Spec} A$ be an affine scheme over $T = \operatorname{Spec} \bar{A}$. A perfect obstruction theory of $S/T$ consists of a two term complex of finitely generated free $A$-modules $E^\bullet = [E^1 \to E^2]$ (indexed at $[1, 2]$) and an obstruction assignment $\mathfrak{ob}$ taking value in the second cohomology of $E^\bullet$ of which the following hold:*

*1. The functor $\mathfrak{Def}^1$ is isomorphic to the functor $\mathfrak{h}^1(E^\bullet)$.*

*2. For any triple $(B, I, \varphi) \in Ob(\mathfrak{Tri}_{A/\bar{A}})$, the element*

$$\mathfrak{ob}(B, I, \varphi) \in \mathfrak{h}^2(E^\bullet)(I) = h^2(E^\bullet \otimes_A I)$$

*is the obstruction class to extending $\varphi \colon A \to B/I$ to an $\bar{A}$-homomorphism $A \to B$.*

*3. The obstruction assignment*

$$(B, I, \varphi) \mapsto (B, h^2(E^\bullet \otimes_A I), \mathfrak{ob}(B, I, \varphi))$$

*is a functor from $\mathfrak{Tri}_{A/\bar{A}}$ to $\mathfrak{Mod}_A^*$. Namely it satisfies the base change property.*

A few remarks are in order here. First, in [LT2], we only considered the case where $T = \operatorname{Spec} \Bbbk$. Here we need to study the relative case for the proof of the degeneration formula. When $T = \operatorname{Spec} \Bbbk$, we will omit $T$ from the notation. When $T$ is non-trivial, we will call the obstruction theory so defined the relative obstruction theory. Secondly, when we restrict to the subcategory of all triples $(B, I, \varphi)$ so that $B$ are Artin local rings with residue fields $\Bbbk$, then the above data is the obstruction theory to deforming $\varphi(\mathrm{Spec}\,\Bbbk)$ in $S$. Thirdly, since $\mathrm{Hom}_A(\Omega_{A/\bar{A}} \otimes_A B/I, I)$ is canonically isomorphic to $\mathfrak{Def}^1_{A/\bar{A}}(I)$, we have $\ker\{E^{2\vee} \to E^{1\vee}\} \cong \Omega_{A/\bar{A}}$. Lastly, in the definition we can replace $E^\bullet$ by its associated complex of sheaves of $\mathcal{O}_{\mathrm{Spec}\,B}$-modules and modify the wording accordingly. This is convenient if we work with obstruction of schemes or stacks. We will call $h^2(E^\bullet)$ the obstruction module and call its associated sheaf (of $\mathcal{O}_{\mathrm{Spec}\,B}$-modules) obstruction sheaf.

We now investigate the deformation theory of pre-deformable morphisms. Let $\Gamma$ be the data $(g, k, b)$ representing the genus, the number of marked points and the degree of the maps. We let $\mathfrak{M}(W[n], \Gamma)^{\mathrm{st}}$ be the moduli space of stable morphisms to $W[n]$ of prescribed topological type that are also pre-deformable as morphisms to the family $W[n]/C[n]$ and are stable as morphisms to the stack $\mathfrak{W}$. As argued in section 2 in [Li], it is a Deligne-Mumford stack and it comes with a tautological morphism

$$\mathfrak{M}(W[n], \Gamma)^{\mathrm{st}} \longrightarrow \mathfrak{M}(\mathfrak{W}, \Gamma).$$

In the remainder of this section, we will cover $\mathfrak{M}(W[n], \Gamma)^{\mathrm{st}}$ by affine étale charts and construct canonical obstruction theory of each of these charts. The obstruction theory of $\mathfrak{M}(\mathfrak{W}, \Gamma)$ will be the descent of the obstruction theory of $\mathfrak{M}(W[n], \Gamma)^{\mathrm{st}}$. Since this study is local, during our study we are free to shrink an open chart $S \to \mathfrak{M}(W[n], \Gamma)^{\mathrm{st}}$ if necessary.

In the remainder of this subsection we fix an affine étale chart $S \to \mathfrak{M}(W[n], \Gamma)^{\mathrm{st}}$ with $f: \mathcal{X} \to W[n]$ its universal family. We let $S = \mathrm{Spec}\,A$ and let $\mathcal{D} \subset \mathcal{X}$ be the divisor of the union of all marked sections of $\mathcal{X}/S$. We fix a collection of charts $(\mathcal{U}_\alpha/\mathcal{V}_\alpha, f_\alpha)$ of $f$ that covers $\mathcal{X}$ satisfying the simplification assumption. In case $\alpha$ is a chart of the second kind, we will reserve the symbols $\mathcal{W}_\alpha, z_{\alpha,i}, w_{\alpha,i}, \phi_\alpha, \psi_\alpha, m_\alpha$ and $l_\alpha$ for the data associated to the chart $\alpha$. Our first task is to show that locally there are no obstruction to extending pre-deformable morphisms $f_\alpha$.

We begin with the notion of flat extensions of an étale neighborhood $\mathcal{U}_\alpha/\mathcal{V}_\alpha$. Let $I$ be an $A$-module. We say $\tilde{\mathcal{V}}_\alpha$ is an extension of $\mathcal{V}_\alpha$ by $I$ if $\mathcal{V}_\alpha \subset \tilde{\mathcal{V}}_\alpha$ is a subscheme with the ideal sheaf $\tilde{\mathcal{I}}$ of $\mathcal{V}_\alpha \subset \tilde{\mathcal{V}}_\alpha$ satisfying $(\tilde{\mathcal{I}})^2 = 0$ and the module $\Gamma(\tilde{\mathcal{I}})$ is isomorphic to $I$ as $\Gamma(\mathcal{O}_\mathcal{V})$-modules. We say $\tilde{\mathcal{U}}_\alpha/\tilde{\mathcal{V}}_\alpha$ is a small extension of $\mathcal{U}_\alpha/\mathcal{V}_\alpha$ by $I$ if $\tilde{\mathcal{V}}_\alpha$ is an extension of $\mathcal{V}_\alpha$ by $I$ and $\tilde{\mathcal{U}}_\alpha$ is a flat extension of $\mathcal{U}_\alpha$ over $\tilde{\mathcal{V}}_\alpha$. Now let $\tilde{\mathcal{U}}_\alpha/\tilde{\mathcal{V}}_\alpha$ be an extension of $\mathcal{U}_\alpha/\mathcal{V}_\alpha$ by $I$. We first consider the case where $\alpha$ is a chart of the second kind. By Lemma 1.6, we can assume that the $\phi_\alpha$ of $\mathcal{U}_\alpha/\mathcal{V}_\alpha$ in (1.4) is extended to $\tilde{\phi}_\alpha$ of $\tilde{\mathcal{U}}_\alpha/\tilde{\mathcal{V}}_\alpha$. We say an extension $\tilde{f}_\alpha: \tilde{\mathcal{U}}_\alpha \to W[n]$ of $f_\alpha: \mathcal{U}_\alpha \to W[n]$ is a pre-deformable extension if $\tilde{\mathcal{U}}_\alpha \to W[n] \to \mathbf{A}^{n+1}$ factor through $\tilde{\mathcal{V}}_\alpha \to \mathbf{A}^{n+1}$ and if the family $\tilde{f}_\alpha: \tilde{\mathcal{U}}_\alpha \to W[n]$ is pre-deformable along $\mathbf{D}_{l_\alpha}$.

We now define the space $\mathrm{Hom}_{\mathcal{U}_\alpha}(f^*\Omega_{W[n]}, I)^\dagger$ that will parameterize all such extensions. The group $\mathrm{Hom}_{\mathcal{U}_\alpha}(f^*\Omega_{W[n]}, I)^\dagger$ is the set of triples

$$(1.10) \qquad (\varphi, \eta_1, \eta_2) \in \mathrm{Hom}_{\mathcal{U}_\alpha}(f^*\Omega_{W[n]}, \mathcal{I}_\alpha) \oplus \mathcal{I}_\alpha^{\oplus 2}, \quad \mathcal{I}_\alpha = I \otimes_A \mathcal{O}_{\mathcal{U}_\alpha}$$

that obey the following condition:

$$(1.11) \qquad \varphi(f^*dw_{\alpha,i}) = f^*(w_{\alpha,i}) \cdot \eta_i, \quad \varphi(f^*dt_l) \in I_\alpha \quad \text{and} \quad \eta_1 + \eta_2 \in I_\alpha$$

for $i = 1, 2$ and $l = 1, \cdots, n+1$. Here $I_\alpha = I \otimes_A \mathcal{O}_{\mathcal{V}_\alpha}$. Since $f_\alpha(\mathcal{U}_\alpha) \subset \mathcal{W}_\alpha$ and $w_{\alpha,i} \in \Gamma(\mathcal{O}_{\mathcal{W}_\alpha})$, $f^*dw_{\alpha,i} \in f^*\Omega_{W[n]} \otimes_{\mathcal{O}_\mathcal{X}} \mathcal{O}_{\mathcal{U}_\alpha}$ and hence $\varphi(f^*dw_{\alpha,i})$ makes sense. Note also that because $t_{l_\alpha} = w_{\alpha,1}w_{\alpha,2}$, from (1.11) we have

$$(1.12) \qquad s_\alpha^{m_\alpha}(\eta_1 + \eta_2) = \varphi(f^*dt_{l_\alpha}).$$

Clearly, $\mathrm{Hom}_{\mathcal{U}_\alpha}(f^*\Omega_{W[n]}, I)^\dagger$ is an $A$-module and is $A$-flat if $I$ is $A$-flat.

When $\alpha$ is chart of the first kind, we define $\mathrm{Hom}_{\mathcal{U}_\alpha}(f^*\Omega_{W[n]}, I)^\dagger$ be the subgroup of $\varphi \in \mathrm{Hom}_{\mathcal{U}_\alpha}(f^*\Omega_{W[n]}, \mathcal{I})$ such that $\varphi(dt_l) \in I_\alpha$ for all $l$. Note that in case $\alpha$ and $\beta$ are two charts of $f$ and $p_\alpha : \mathcal{U}_{\alpha\beta} \triangleq \mathcal{U}_\alpha \times_\mathcal{X} \mathcal{U}_\beta \to \mathcal{U}_\alpha$ is the projection. Then there is a canonical restriction $A$-homomorphism

$$p_\alpha^* : \mathrm{Hom}_{\mathcal{U}_\alpha}(f^*\Omega_{W[n]}, I)^\dagger \longrightarrow \mathrm{Hom}_{\mathcal{U}_{\alpha\beta}}(f^*\Omega_{W[n]}, I)^\dagger.$$

We now state and prove the following local deformation Lemma.

**Lemma 1.12.** *Let $\tilde{\mathcal{U}}_\alpha/\tilde{\mathcal{V}}_\alpha$ be an extension of $\mathcal{U}_\alpha/\mathcal{V}_\alpha$ by $I$. Then $f_\alpha$ automatically extends to a pre-deformable morphism $\tilde{f}_\alpha : \tilde{\mathcal{U}}_\alpha \to W[n]$. Further, after fixing one such extension, say $\tilde{f}'_\alpha$, the space of all such extensions is canonically isomorphic to the space $\mathrm{Hom}_{\mathcal{U}_\alpha}(f^*\Omega_{W[n]}, I)^\dagger$.*

*Proof.* We will consider the case where $\alpha$ is a chart of the second kind. The other case is simpler. First, we extend $\phi_\alpha$ to a parameterization $\tilde{\phi}_\alpha$ of $\tilde{\mathcal{U}}_\alpha/\tilde{\mathcal{V}}_\alpha$ with $\tilde{z}_{\alpha,i}$ and $\tilde{s}_\alpha$ the corresponding extensions of $z_{\alpha,i}$ and $s_\alpha$, respectively. We consider the composite $\psi_\alpha \circ f_\alpha : \mathcal{U}_\alpha \to \Theta_{l_\alpha}$, where $\Theta_{l_\alpha}$ is defined in (1.1). Let $h_{\alpha,i} \equiv 1 \in \Gamma(\mathcal{O}_{\mathcal{U}_\alpha})$ and $g_{\alpha,l} \in \Gamma(\mathcal{O}_{\mathcal{V}_\alpha})$ be part of the definition of the chart of $f_\alpha$, as in (1.5) and (1.6). Since both $\mathcal{U}_\alpha$ and $\mathcal{V}_\alpha$ are affine, we can extend $h_{\alpha,i}$ and $g_{\alpha,l}$ to $\tilde{h}_\alpha \in \Gamma(\mathcal{O}_{\tilde{\mathcal{U}}_\alpha})$ and $\tilde{g}_{\alpha,l} \in \Gamma(\mathcal{O}_{\tilde{\mathcal{V}}_\alpha})$, respectively, so that $\tilde{h}_{\alpha,1}\tilde{h}_{\alpha,2} \in \Gamma(\mathcal{O}_{\tilde{\mathcal{V}}_\alpha})$ and $\tilde{g}_{\alpha,l_\alpha} = \tilde{s}_\alpha^{m_\alpha}(\tilde{h}_{\alpha,1}\tilde{h}_{\alpha,2})$. We then define $F : \tilde{\mathcal{U}}_\alpha \to \Theta_{l_\alpha}$ by $w_{\alpha,i} \mapsto \tilde{z}_{\alpha,i}^{m_\alpha}\tilde{h}_{\alpha,i}$ for $i = 1,2$ and $t_l \mapsto \tilde{g}_{\alpha,l}$ for $l = 1,\cdots, n+1$. Since $\mathcal{W}_\alpha \to \Theta_{l_\alpha}$ is smooth, we can lift $F$ to an extension $\tilde{f}_\alpha : \tilde{\mathcal{U}}_\alpha \to W[n]$ of $f_\alpha$. Namely, $\tilde{f}_\alpha$ is an extension of $f_\alpha$ such that $\psi_\alpha \circ \tilde{f}_\alpha = F$. The morphism $\tilde{f}_\alpha$ is a desired extension.

Now let $\tilde{f}'_\alpha$ be a fixed pre-deformable extension of $f_\alpha$ to $\tilde{\mathcal{U}}_\alpha/\tilde{\mathcal{V}}_\alpha$. Let $\tilde{f}_\alpha$ be any pre-deformable extension of $f_\alpha$. As morphisms from $\tilde{\mathcal{U}}_\alpha$ to $\Theta_{l_\alpha}$, $\psi_\alpha \circ \tilde{f}_\alpha$ is defined by

$$(1.13) \qquad (\psi_\alpha \circ \tilde{f}_\alpha)^*(w_{\alpha,i}) = \tilde{z}_{\alpha,i}^{m_\alpha}\tilde{h}_{\alpha,i} \quad \text{and} \quad (\psi_\alpha \circ \tilde{f}_\alpha)^*(t_l) = \tilde{g}_{\alpha,l}$$

and $\psi_\alpha \circ \tilde{f}'_\alpha$ is defined by

$$(1.14) \qquad (\psi_\alpha \circ \tilde{f}'_\alpha)^*(w_{\alpha,i}) = \tilde{z}_{\alpha,i}^{m_\alpha}\tilde{h}'_{\alpha,i} \quad \text{and} \quad (\psi_\alpha \circ \tilde{f}'_\alpha)^*(t_l) = \tilde{g}'_{\alpha,l}.$$

We let $\eta_i = \tilde{h}_{\alpha,i} - \tilde{h}'_{\alpha,i} \in \mathcal{I}_\alpha$. Then $\eta_1 + \eta_2 \in I_\alpha$ since $\tilde{h}_{\alpha,1}\tilde{h}_{\alpha,2}$ and $\tilde{h}'_{\alpha,1}\tilde{h}'_{\alpha,2} \in \Gamma(\mathcal{O}_{\tilde{\mathcal{V}}_\alpha})$. Now let $\varphi \in \mathrm{Hom}_{\mathcal{U}_\alpha}(f_\alpha^*\Omega_\mathcal{W}, \mathcal{I})$ be defined by the difference[5] $\mathbf{d}(\tilde{f}_\alpha - \tilde{f}'_\alpha)$. It can be easily checked, based on (1.13) and (1.14), that $(\varphi, \eta_1, \eta_2)$ is in $\mathrm{Hom}_{\mathcal{U}_\alpha}(f^*\Omega_{W[n]}, I)^\dagger$.

It remains to check that this correspondence is one-one and onto, which is straightforward. This completes the proof of the Lemma. $\blacksquare$

**Remark 1.13.** *In the remainder part of this paper, we will call $(\varphi, \eta_1, \eta_2)$ the log-difference of $\tilde{f}_\alpha$ and $\tilde{f}'_\alpha$, denoted by $\mathbf{d}^\dagger(\tilde{f}_\alpha - \tilde{f}'_\alpha)$. Furthermore, if $\tilde{f}_1$, $\tilde{f}_2$ and $\tilde{f}_3$*

---

[5]Let $\iota : A' \to A$ be a small ring extension with $I = \ker\{\iota\}$. We call $A' \to A$ a small extension if $I^2 = 0$. Let $B$ be any ring. Let $f_1, f_2 : B \to A$ be two ring homomorphisms so that $\iota \circ f_1 = \iota \circ f_2$. We define the difference of $f_1$ and $f_2$ to be the map $\mathbf{d}(f_1 - f_2) : \Omega_B \to I$ defined by $b \otimes 1 - 1 \otimes b \mapsto f_1(b) - f_2(b) \in I$. It is an $A$-homomorphism $\Omega_B \otimes_B A \to I$. Note that once $f_1$ is fixed, then $f_2$ is uniquely determined by $\mathbf{d}(f_1 - f_2)$. Conversely, any $\varphi \in \mathrm{Hom}_A(\Omega_B \otimes_B A, I)$ defines a unique homomorphism $f_2 : B \to A$ so that $\mathbf{d}(f_1 - f_2) = \varphi$.

*are three pre-deformable extensions of $f_\alpha$ to $\tilde{\mathcal{U}}_\alpha$, then $\mathbf{d}^\dagger(\tilde{f}_3 - \tilde{f}_1) = \mathbf{d}^\dagger(\tilde{f}_3 - \tilde{f}_2) + \mathbf{d}^\dagger(\tilde{f}_2 - \tilde{f}_1)$.*

Before we proceed, we introduce more convention which we will follow throughout this paper. We let $\iota_\alpha : \mathcal{U}_\alpha \to \mathcal{X}$ and $\jmath_\alpha : \mathcal{V}_\alpha \to S$ be the tautological projections. As a convention, for each $\alpha \in \Lambda$ we let $A_\alpha = \Gamma(\mathcal{V}_\alpha)$ with $A \to A_\alpha$ the tautological homomorphism, for $\alpha$ and $\beta \in \Lambda$ we will denote by $\mathcal{U}_{\alpha\beta}$ the product $\mathcal{U}_\alpha \times_\mathcal{X} \mathcal{U}_\beta$ and by $p_\alpha$ and $p_\beta$ the projections $\mathcal{U}_{\alpha\beta} \to \mathcal{U}_\alpha$ and to $\mathcal{U}_\beta$ respectively. The same convention applies to multi-indices in the obvious way. We let $\pi_\alpha : \mathcal{U}_\alpha \to \mathcal{V}_\alpha$ be the projection. For sheaf of $\mathcal{O}_\mathcal{X}$-modules $\mathcal{A}$ (resp. $A$-module $I$; divisor $\mathcal{D}$) we will denote by $\mathcal{A}_\alpha$ (resp. $I_\alpha$, resp. $\mathcal{D}_\alpha$) the pull back sheaf $\iota_\alpha^* \mathcal{A}$ (resp. $I \otimes_A A_\alpha$; resp. $\mathcal{D} \times_\mathcal{X} \mathcal{U}_\alpha$). Also, for $A$-module $I$, we will use $\mathcal{I}$ (resp. $\mathcal{I}_\alpha$) to denote the sheaf of $\mathcal{O}_\mathcal{X}$-modules $\mathcal{O}_\mathcal{X} \otimes_A I$ (resp. $\mathcal{O}_{\mathcal{U}_\alpha} \otimes_A I$). We denote by $d_\alpha : \mathcal{F}_\alpha^0 \to \mathcal{F}_\alpha^2$ the pull back of $d : F^0 \to F^1$ in (1.15). For sheaves of $\mathcal{O}_\mathcal{X}$-modules $\mathcal{A}$ and $\mathcal{B}$ we agree that $\mathrm{Hom}_{\mathcal{U}_\alpha}(\mathcal{A}, \mathcal{B}) = \mathrm{Hom}_{\mathcal{U}_\alpha}(\mathcal{A}_\alpha, \mathcal{B}_\alpha)$ and for $A$-modules $I$ and $J$ we agree $\mathrm{Hom}_{A_\alpha}(I, J) = \mathrm{Hom}_{A_\alpha}(I_\alpha, J_\alpha)$.

We now study the deformation of $f : \mathcal{X}/S \to W[n]/\mathbf{A}^{n+1}$. By the deformation theory of nodal curves with sections, there is a complex of free $A$-modules

$$(1.15) \qquad\qquad F^\bullet = [F^0 \xrightarrow{\ d\ } F^1]$$

so that for any $A$-module $I$,

$$(1.16) \qquad\qquad \mathrm{Ext}_\mathcal{X}^i(\Omega_{\mathcal{X}/S}(\mathcal{D}), \mathcal{I}) = h^i(F^\bullet \otimes_A I).$$

Now let $\tilde{S} \triangleq \mathrm{Spec}\, A * F^{1\vee}$ [6] with $S \subset \tilde{S}$ be the immersion induced by the obvious projection $A * F^{1\vee} \to A$. Let $\mathbf{1} \in F^1 \otimes_A F^{1\vee}$ be the identity element and let $[\mathbf{1}] \in \mathrm{Ext}_\mathcal{X}^1(\Omega_{\mathcal{X}/S}(\mathcal{D}), F^{1\vee})$ be the associated element. The element $[\mathbf{1}]$ defines a family $\tilde{\mathcal{X}}/\tilde{S}$ extending the family $\mathcal{X}/S$, using (1.16). It has the following properties:

*Firstly, let $I$ be any $A$-module and $a \in F^1 \otimes_A I$ be any element. Let $T \triangleq \mathrm{Spec}\, A * I$ and $\mathcal{X}_T/T$ be the extension of $\mathcal{X}$ defined by the cohomology class $[a] \in \mathrm{Ext}_\mathcal{X}^1(\Omega_{\mathcal{X}/S}(\mathcal{D}), I)$ of $a \in F^1 \otimes_A I$. Then $\mathcal{X}_T$ is the pull back of $\tilde{\mathcal{X}}$ under the morphism $T \to \tilde{S}$ defined by $A * F^{1\vee} \to A * I$ via $(x, y) \mapsto (x, a(y))$. Namely we have isomorphism $\tilde{\mathcal{X}} \times_{\tilde{S}} T \cong \mathcal{X}_T$.*

*Secondly, let*

$$\tilde{T} \triangleq \mathrm{Spec}\, A * F^{0\vee} \to \tilde{S} = \mathrm{Spec}\, A * F^{1\vee}$$

*be the morphism defined by $A * F^{1\vee} \to A * F^{0\vee}$ via $(x, y) \mapsto (x, d^\vee(y))$. Then we have isomorphism over $\tilde{T}$*

$$(1.17) \qquad\qquad \tilde{\mathcal{X}} \times_{\tilde{S}} \tilde{T} \cong \mathcal{X} \times_S \tilde{T},$$

*where the projection $\tilde{T} \to S$ is defined by the obvious inclusion $A \hookrightarrow A * F^{0\vee}$.*

As mentioned, $\mathbf{1} \in F^1 \otimes_A F^{1\vee}$ defines an extension $\tilde{\mathcal{X}}/\tilde{S}$, which we fix from now on. We let $\tilde{\mathcal{U}}_\alpha$ (resp. $\tilde{\mathcal{V}}_\alpha$) be the étale neighborhood of $\tilde{\mathcal{X}}$ (resp. $\tilde{S}$) that is the minimal extension of $\mathcal{U}_\alpha$ (resp. $\mathcal{V}_\alpha$) [7]. Then $\{\tilde{\mathcal{U}}_\alpha/\tilde{\mathcal{V}}_\alpha\}$ forms a covering of $\tilde{\mathcal{X}}/\tilde{S}$. By Lemma 1.12, $f_\alpha : \mathcal{U}_\alpha \to \mathcal{W}_\alpha$ can be extended to a pre-deformable morphism $\tilde{\mathcal{U}}_\alpha \to W[n]$. For each $\alpha$ we pick one such extension $\zeta_\alpha : \tilde{\mathcal{U}}_\alpha \to W[n]$ once and for all. Now let $\alpha, \beta \in \Lambda$ be a pair so that $\mathcal{U}_{\alpha\beta} \neq \emptyset$. Let $\tilde{\mathcal{U}}_{\alpha\beta} = \tilde{\mathcal{U}}_\alpha \times_{\tilde{\mathcal{X}}} \tilde{\mathcal{U}}_\beta$ and let $\tilde{p}_\alpha \colon \tilde{\mathcal{U}}_{\alpha\beta} \to \tilde{\mathcal{U}}_\alpha$ be the projection, following our convention. Since both $\zeta_\alpha \circ \tilde{p}_\alpha$ and $\zeta_\beta \circ \tilde{p}_\beta$ are pre-deformable extensions of $f_{\alpha\beta} \colon \mathcal{U}_{\alpha\beta} \to W[n]$, by Lemma 1.12 their difference defines an element

$$\zeta_{\alpha\beta} \triangleq \mathbf{d}^\dagger(\zeta_\beta \circ \tilde{p}_\beta - \zeta_\alpha \circ \tilde{p}_\alpha) \in \operatorname{Hom}_{\mathcal{U}_{\alpha\beta}}(f^*\Omega_{W[n]}, F^{1\vee})^\dagger.$$

---

[6] In this section we will denote by $A * I$ the trivial ring extension of $A$ by the $A$-module $I$.

[7] We say $\tilde{\mathcal{U}}_\alpha$ is a minimal extension of the étale neighborhood $\mathcal{U}_\alpha$ if $\tilde{\mathcal{U}}_\alpha \to \tilde{\mathcal{X}}$ is an étale neighborhood and $\tilde{\mathcal{U}}_\alpha \times_{\tilde{\mathcal{X}}} \mathcal{X} \cong \mathcal{U}_\alpha$.

This defines a homomorphism

(1.18) $$\zeta_{\alpha\beta}(\cdot) \colon F^1 \longrightarrow \operatorname{Hom}_{\mathcal{U}_{\alpha\beta}}(f^*\Omega_{W[n]}, A)^\dagger.$$

Now to each $\alpha$ we construct a homomorphism

(1.19) $$\zeta_\alpha(\cdot) \colon F^0 \longrightarrow \operatorname{Hom}_{\mathcal{U}_\alpha}(f^*\Omega_{W[n]}, A)^\dagger.$$

Let $\mathbf{1} \in F^0 \otimes_A F^{0\vee}$ be the identity element. Then $d(\mathbf{1}) \in F^1 \otimes_A F^{0\vee}$ ($d$ is the differential in the complex $F^\bullet$) defines a homomorphism $d(\mathbf{1}) \colon F^{1\vee} \to F^{0\vee}$ which induces a morphism

$$\tau \colon T \triangleq \operatorname{Spec} A * F^{0\vee} \to \operatorname{Spec} A * F^{1\vee} = \tilde{S}.$$

Let $\mathcal{X}_T/T$ be the pull back family of $\tilde{\mathcal{X}}/\tilde{S}$ via $\tau$ with $q \colon \mathcal{X}_T \to \tilde{\mathcal{X}}$ the induced the projection. By the second property of the complex $F^\bullet$, $\mathcal{X}_T$ is isomorphic to $\mathcal{X} \times_S T$, where $T \to S$ is induced by the obvious inclusion $A \to A * F^{0\vee}$. We let $q_0 \colon \mathcal{X}_T \to \mathcal{X}$ be the induced projection. Clearly, $U_\alpha = \tilde{\mathcal{U}}_\alpha \times_{\tilde{S}} T$ is canonically isomorphic to $\mathcal{U}_\alpha \times_S T$ under the isomorphism $\mathcal{X}_T \cong \mathcal{X} \times_S T$. Then restricting to $U_\alpha$, both $\zeta_\alpha \circ q$ and $f \circ q_0$ are pre-deformable extensions of $f_\alpha$. By Lemma 1.12, the difference

$$\mathbf{d}^\dagger(\zeta_\alpha \circ q - f \circ q_0) \in \operatorname{Hom}_{\mathcal{U}_\alpha}(f^*\Omega_{W[n]}, F^{0\vee})^\dagger$$

It defines the desired homomorphism $\zeta_\alpha(\cdot)$ in (1.19).

**Lemma 1.14.** *Both $\zeta_\alpha(\cdot)$ and $\zeta_{\alpha\beta}(\cdot)$ are homomorphisms of $A$-modules. Further, for any $a \in F^0$ we have*

$$-p_\alpha^*\zeta_\alpha(a) + p_\beta^*\zeta_\beta(a) = \zeta_{\alpha\beta}(d(a)) \in \operatorname{Hom}_{\mathcal{U}_{\alpha\beta}}(f^*\Omega_{W[n]^\dagger}, A)^\dagger,$$

*where $p_\alpha \colon \mathcal{U}_{\alpha\beta} \to \mathcal{U}_\alpha$ is the projection and $p_\alpha^*\zeta_\alpha(\cdot)$ is the pull back homomorphism.*

*Proof.* The proof is straightforward and will be omitted. □

In the next part, we will construct the complex that will be a part of the perfect obstruction theory of $S \subset \mathfrak{M}(W[n], \Gamma)^{\mathrm{st}}$. Let $I$ be any $A$-module. We let $\mathbf{D}(I)^\bullet$ be the Čech complex

$$\mathbf{D}(I)^\bullet = \mathbf{C}^\bullet(\Lambda, \mathcal{H}om(f^*\Omega_{W[n]}, I)^\dagger)$$

associating to the covering $\Lambda$, where

$$\Gamma(\mathcal{U}_{\alpha_1 \cdots \alpha_k}, \mathcal{H}om(f^*\Omega_{W[n]}, I)^\dagger) = \operatorname{Hom}_{\mathcal{U}_{\alpha_1 \cdots \alpha_k}}(f^*\Omega_{W[n]}, I)^\dagger$$

with $\partial \colon \mathbf{D}(I)^\bullet \to \mathbf{D}(I)^{\bullet+1}$ the coboundary differential in the Čech complex. We let

$$\delta_0 \colon F^0 \otimes_A I \longrightarrow \mathbf{C}^0(\Lambda, \mathcal{H}om(f^*\Omega_{W[n]}, I)^\dagger)$$

be defined by $\delta(a)_\alpha = \zeta_\alpha(a)$ and let

$$\delta_1 \colon F^1 \otimes_A I \longrightarrow \mathbf{C}^1(\Lambda, \mathcal{H}om(f^*\Omega_{W[n]}, I)^\dagger)$$

be defined by $\delta_1(b)_{\alpha\beta} = \zeta_{\alpha\beta}(b)$. By Lemma 1.14, $\delta_i$ are homomorphisms of $A$-modules. We now show that this defines a homomorphism of complexes

$$\delta \colon F^\bullet \otimes_A I \longrightarrow \mathbf{D}(I)^\bullet.$$

To prove this, we need to check that $\delta_1 \circ d = \partial \circ \delta_0$ on $F^0 \otimes_A I$ and $\partial \circ \delta_1 = 0$ on $F^1 \otimes_A I$. We will check $\partial \circ \delta = 0$ on $F^1 \otimes_A I$ and leave the other to the readers. Let $b \in F^1 \otimes_A I$ be any element. By definition,

$$(\partial \circ \delta_1)(b)_{\alpha\beta\gamma} = \delta_1(b)_{\alpha\beta} - \delta_1(b)_{\alpha\gamma} + \delta_1(b)_{\beta\gamma} = \zeta_{\alpha\beta}(b) - \zeta_{\alpha\gamma}(b) + \zeta_{\beta\gamma}(b)$$

as elements in $\mathrm{Hom}_{\mathcal{U}_{\alpha\beta\gamma}}(f^*\Omega_{W[n]}, I)^{\dagger}$, where the summation is taken after pulling back each term in the summation to this module in the obvious way. But this vanishes because of the relation $\zeta_{\beta\gamma} - \zeta_{\alpha\gamma} + \zeta_{\alpha\beta} = 0$, following Remark 1.13. This shows that $\partial \circ \delta_1 = 0$ on $F^1 \otimes_A I$.

In the end, we define the complex $\mathbf{E}(I)^\bullet$ by

$$(1.20) \qquad \mathbf{E}(I)^k = \underset{i+j=k}{\oplus} (F^i \otimes_A I \oplus \mathbf{D}(I)^{j-1})$$

with the differential $d_{\mathbf{E}} \colon \mathbf{E}(I)^k \to \mathbf{E}(I)^{k+1}$ defined so that its restriction to $F^0 \otimes_A I$ and $F^1 \otimes_A I$ are $d \oplus -\delta_0$ and $\delta_1$, respectively, and its restriction to $\mathbf{D}(I)^i$ is $\partial$. When $I = A$, we will abbreviate $\mathbf{D}(A)^\bullet$ and $\mathbf{E}(A)^\bullet$ to $\mathbf{D}^\bullet$ and $\mathbf{E}^\bullet$ respectively.

**Lemma 1.15.** *We assume $\Lambda$ is a sufficiently fine cover of $f$. Then for any $A$-module $I$, $\mathbf{E}(I)^\bullet$ is canonically isomorphic to $\mathbf{E}^\bullet \otimes_A I$ as complexes. Further, the complex $\mathbf{E}^\bullet$ is a complex of flat $A$-modules.*

*Proof.* The proof is straightforward and will be omitted.                                   □

We remark that the complex $\mathbf{E}^\bullet$ just constructed depends on the choice of the atlas $\Lambda$. To emphasize this dependence we shall denote it by $\mathbf{E}_\Lambda^\bullet$. Let $\Lambda'$ be any atlas of $f$ that is a refinement of $\Lambda$ with associated complex $\mathbf{E}_{\Lambda'}^\bullet$. Following [Mil, III.2], we can define a homomorphism of complexes $\mathbf{E}_\Lambda^\bullet \to \mathbf{E}_{\Lambda'}^\bullet$ which induces a map of cohomologies

$$\rho(\Lambda, \Lambda') \colon h^i(\mathbf{E}_\Lambda^\bullet \otimes_A I) \longrightarrow h^i(\mathbf{E}_{\Lambda'}^\bullet \otimes_A I)$$

and then form the direct limit $\lim_{\to} h^i(\mathbf{E}_\Lambda \otimes_A I)$ taken over all charts of $f$. Note that this limit is a functor from $\mathfrak{Mod}_A$ to $\mathfrak{Mod}_A$, denoted by $\mathfrak{h}^i(\mathbf{E}^\bullet)$. We now assume $\Lambda$ is fine enough so that $\mathfrak{h}^i(\mathbf{E}^\bullet)(I) = h^i(\mathbf{E}_\Lambda^\bullet \otimes_A I)$ for all $A$-module $I$. We fix such a $\Lambda$ once and for all, and abbreviate the resulting complex $\mathbf{E}_\Lambda^\bullet$ to $\mathbf{E}^\bullet$.

We now prove the main result of this section.

**Proposition 1.16.** *Let $\mathfrak{Def}_A^1$ be the functor of the first order deformations of morphisms to $S = \mathrm{Spec}\, A$, which is naturally a functor from $A$-modules to $A$-modules. Then $\mathfrak{Def}_A^1$ is naturally isomorphic to the functor $\mathfrak{h}^1(\mathbf{E}^\bullet)$.*

*Proof.* Let $\xi = (B, I)$ be any object in $\mathfrak{Mod}_A$. We first show that there is a canonical isomorphism $\mathfrak{Def}_A^1(\xi) \cong h^1(\mathbf{E}^\bullet \otimes_A I)$.

Let $T = \mathrm{Spec}\, B$ and $\tilde{T} = \mathrm{Spec}\, B * I$, the trivial extension by $I$. Since $B$ is an $A$-algebra, there is a tautological morphism $T \to S$. Let $x \in \mathfrak{Def}_S^1(\xi)$ [8] be any element, associated to an extension $\tilde{T} \to S$ of $T \to S$, and let $f_T \colon \mathcal{X}_T \to W[n]$ (resp. $f_{\tilde{T}} \colon \mathcal{X}_{\tilde{T}} \to W[n]$) be the pull back family of $f$ via $T \to S$ (resp $\tilde{T} \to S$). We let $\mathcal{D}_T \subset \mathcal{X}_T$ and $\mathcal{D}_{\tilde{T}} \subset \mathcal{X}_{\tilde{T}}$ be the associated divisors of marked points and let $\mathcal{I} = \mathcal{O}_{\mathcal{X}_T} \otimes_B I$. First of all, since $\mathcal{X}_{\tilde{T}}$ is a flat extension of $\mathcal{X}_T$ to $\tilde{T}$, it associates to a unique element in

$$(1.21) \qquad \mathrm{Ext}^1_{\mathcal{X}_T}(\Omega_{\mathcal{X}_T/T}(\mathcal{D}_T), \mathcal{I}) \cong h^1(F^\bullet \otimes_A I).$$

---

[8] We use the convention $\mathfrak{Def}_S = \mathfrak{Def}_A$ when $S = \mathrm{Spec}\, A$.

Hence it is represented by $[a] \in h^1(F^\bullet \otimes_A I)$ for an $a \in F^1 \otimes_A I$. Now let $\tilde{S}$ be the trivial extension of $S$ by $F^{1\vee}$ and let $\tilde{\mathcal{X}}/\tilde{S}$ be the family defined before (1.17). Then $a$ defines a morphism $\varphi_a \colon \tilde{T} \to \tilde{S}$ via the homomorphism $A \oplus F^{1\vee} \to B \oplus I$ defined by $(x, y) \mapsto (\bar{x}, a(y))$, where $\bar{x}$ is the image of $x$ in $B$. By the first property of (1.16), $a$ defines an isomorphism $\gamma_a \colon \mathcal{X}_{\tilde{T}} \cong \tilde{\mathcal{X}} \times_{\tilde{S}} \tilde{T}$, where the projection $\tilde{T} \to \tilde{S}$ is via $\varphi_a$. We let $\tilde{\mathcal{U}}_\alpha = \tilde{\mathcal{U}}_\alpha \times_{\tilde{S}} \tilde{T}$ and $\tilde{V}_\alpha = \tilde{\mathcal{V}}_\alpha \times_{\tilde{S}} \tilde{T}$. Then each $\tilde{U}_\alpha / \tilde{V}_\alpha$ extends to a chart $(\tilde{U}_\alpha / \tilde{V}_\alpha, f_{\tilde{T}, \alpha})$ of $f_{\tilde{T}}$, where $f_{\tilde{T}, \alpha} \triangleq f_{\tilde{T}}|_{\tilde{U}}$. We let $\tilde{\rho}_\alpha \colon \tilde{U}_\alpha \to \tilde{\mathcal{U}}_\alpha$ be the projection. Then over each $\tilde{U}_\alpha$ we have two pre-deformable extensions of $f_T$: One is $f_{\tilde{T}, \alpha}$ and the other is the composite $\zeta_\alpha \circ \tilde{\rho}_\alpha$. Let

$$b_\alpha = \mathbf{d}^\dagger(f_{\tilde{T}, \alpha} - \zeta_\alpha \circ \tilde{\rho}_\alpha) \in \mathrm{Hom}_{\mathcal{U}_\alpha}(f_T^* \Omega_{W[n]}, I)^\dagger$$

and let $b = \{b_\alpha\} \in \mathbf{D}^0 \otimes_A I$. We claim that $(a, b) \in \mathbf{E}^1 \otimes_A I$ is in the kernel of $d_{\mathbf{E}}$. For this, we only need to check that

$$d_{\mathbf{E}}(a, b)_{\alpha\beta} = \delta(a)_{\alpha\beta} + \partial(b)_{\alpha\beta} = \zeta_{\alpha\beta}(a) + (b_\beta - b_\alpha) = \zeta_{\alpha\beta}(a) - \mathbf{d}^\dagger(\zeta_\beta \circ \tilde{\rho}_\beta - \zeta_\alpha \circ \tilde{\rho}_\alpha)$$

vanishes for all pairs $(\alpha, \beta)$. But this follows immediately from the definition of $\zeta_{\alpha\beta}(\cdot)$. This shows that $(a, b)$ defines a cohomology class $[(a, b)] \in h^1(\mathbf{E}^\bullet \otimes_A I)$.

Next, we show that $[(a, b)]$ is independent of the choices of $a$ and the isomorphisms $\gamma_a \colon \mathcal{X}_{\tilde{T}} \cong \tilde{\mathcal{X}} \times_{\tilde{S}} \tilde{T}$. Let $a' \in F^1 \times_A I$ be another element so that $[a'] = [a]$ in (1.21). Then $a - a' = d_I(c)$ for a $c \in F^0 \otimes_A I$, where $d_I \colon F^0 \otimes_A I \to F^1 \otimes_A I$ is induced by $d$ in (1.15). Now let

$$\varphi_{a'} \colon \tilde{T} \to \tilde{S}, \quad \gamma_{a'} \colon \mathcal{X}_{\tilde{T}} \cong \tilde{\mathcal{X}} \times_{\tilde{S}} \tilde{T}, \quad \tilde{\rho}'_\alpha \colon \tilde{U}_\alpha \to \tilde{\mathcal{U}}_\alpha, \quad b'_\alpha = \mathbf{d}^\dagger(f_{\tilde{T}, \alpha} - \zeta_\alpha \circ \tilde{\rho}'_\alpha)$$

be objects defined similarly with $a$ replaced by $a'$. Let $b' = \{b'_\alpha\}$. We claim that

$$d_{\mathbf{E}}(c) = (a, b) - (a', b') \in \mathbf{E}^1 \otimes_A I.$$

Once this is established then $[(a, b)] = [(a', b')] \in h^1(\mathbf{E}^\bullet \otimes_A I)$, which shows that $[(a, b)]$ only depends on the class $x \in \mathfrak{Def}^1_S(v)$. This way we obtain a map

$$\mathbf{T}(\xi) \colon \mathfrak{Def}^1_S(v) \longrightarrow h^1(E^\bullet \otimes_A I).$$

Now we prove the claim. Since

$$d_{\mathbf{E}}(c) = (\partial(c), -\delta(c)) = (d_I(c), -\delta(c))$$

and $d_I(c) = a - a'$, it suffices to show that $\delta(c)_\alpha = \zeta(c)_\alpha$ is identical to

$$-(b - b')_\alpha = -(b_\alpha - b'_\alpha) = -\mathbf{d}^\dagger(f_{\tilde{T}, \alpha} - \zeta_\alpha \circ \tilde{\rho}_\alpha) + \mathbf{d}^\dagger(f_{\tilde{T}, \alpha} - \zeta_\alpha \circ \tilde{\rho}'_\alpha) = \mathbf{d}^\dagger(\zeta_\alpha \circ \tilde{\rho}_\alpha - \zeta_\alpha \circ \tilde{\rho}'_\alpha).$$

The proof that $\zeta_\alpha(c) = \mathbf{d}^\dagger(\zeta_\alpha \circ \tilde{\rho}_\alpha - \zeta_\alpha \circ \tilde{\rho}'_\alpha)$ is routine and will be omitted.

We now show that the map $\mathbf{T}(\xi)$ is one-one and onto. We give an outline of the proof since it is standard. We first show that it is one-one. Let $x \in \mathfrak{Def}^1_S(\xi)$ be any element so that $\mathbf{T}(\xi)(x) = 0 \in h^1(\mathbf{E}^\bullet \otimes_A I)$. We let $(a, b)$ be the pair constructed associated to the family $f_{\tilde{T}} \colon \mathcal{X}_{\tilde{T}} \to W[n]$ following the previous discussion. Since $[(a, b)] = 0$, there is a $c \in F^0 \otimes_A I$ so that $(a, b) = d_{\mathbf{E}}(c)$. This implies at first that $a = d_I(c)$. Hence $\mathcal{X}_{\tilde{T}} \cong \mathcal{X}_T \times_T \tilde{T}$ under the obvious projection $\tilde{T} \to T$ via $B \to B \oplus I$. Let $\rho \colon \mathcal{X}_{\tilde{T}} \to \mathcal{X}_T$ be the projection. Because $\mathbf{T}(\xi)(x)$ is well-defined, $\mathbf{T}(\xi)(x)$ is also represented by $(0, b')$, where $b' = \{b'_\alpha\}$ and

$$b'_\alpha = \mathbf{d}^\dagger(f_{\tilde{T}, \alpha} - f_{T, \alpha} \circ \rho_\alpha) \in \mathrm{Hom}_{U_\alpha}(f_T^* \Omega_{W[n]}, I)^\dagger.$$

Again, since $\mathbf{T}(\xi)(x) = 0$, there must be a $c \in F^0 \otimes_A I$ so that $d_I(c) = 0$ and $\delta(c) = b'$. Hence $c$ lifts to an element in $\mathrm{Ext}^0_{\mathcal{X}_T}(\Omega_{\mathcal{X}_T/T}(\mathcal{D}_T), \mathcal{I})$, which defines a new isomorphism $\gamma: \mathcal{X}_{\tilde{T}} \cong \mathcal{X}_T \times_T \tilde{T}$. Further, $\delta(c) = b'$ implies that under this new isomorphism $f_{\tilde{T}}$ is the constant extension of $f_T$. Namely, $f_{\tilde{T}}$ is the pull back of $f_T$ via the projection $\tilde{T} \to T$. This proves that $x = 0$ in $\mathfrak{Def}^{\dagger}_S(\xi)$ and hence $\mathbf{T}(\xi)$ is one-one.

The onto part is similar. Since this argument is standard in deformation theory, we will omit it here. In the end, we need to check that $\mathbf{T}(\xi)$ is a homomorphism of modules and that $\mathbf{T}$ is an isomorphism of functors. This is straightforward and will be omitted.                                                                                                    □

**Corollary 1.17.** *Let $\xi = (B, I, \varphi_0)$ be any object in $\mathfrak{Tri}_S$. Suppose $\mathfrak{Def}_S(\xi) \neq \emptyset$, then $\mathfrak{Def}_S(\xi)$ is isomorphic to the set $h^1(\mathbf{E}^{\bullet} \otimes_A I)$.*

*Proof.* Since $S$ is an étale chart of $\mathfrak{M}(W[n], \Gamma)^{\mathrm{st}}$, which is a Deligne-Mumford stack, the standard fact in deformation theory shows that once an extension $\varphi: \operatorname{Spec} B \to S$ of $\varphi_0$ is fixed, then the space of all such extensions is canonically isomorphic to $\mathfrak{Def}^1_S(\xi_0) \cong h^1(\mathbf{E}^{\bullet} \otimes_A I)$, where $\xi_0 = (B/I * I, I, \varphi_0)$. This proves the corollary.        □

**Proposition 1.18.** *There is a natural obstruction theory to deformation of the family of pre-deformable $f: \mathcal{X} \to W[n]$ over $S$ that takes values in $\mathfrak{h}^2(\mathbf{E}^{\bullet})$.*

*Proof.* The construction of the obstruction class is standard, as shown in [LT1] for smooth targets. Let $\xi = (B, I, \varphi)$ be an object in $\mathfrak{Tri}_S$. Let $T = \operatorname{Spec} B/I$ and let $f_T: \mathcal{X}_T \to W[n]$ with $\mathcal{D}_T \subset \mathcal{X}_T$ be the divisor of marked points be the pull back of $f$ via $T \to S$. Then to find an extension $\tilde{r}: \tilde{T} \to S$ of $T \to S$ is equivalent to extend the family $f_T$ to a family of pre-deformable morphisms over $\tilde{T}$. We now construct the obstruction class to this extension problem.

Since deformation of pointed curves is unobstructed, we can extend $\mathcal{X}_T$ to a family over $\tilde{T}$. By our choice of $F^{\bullet}$, such extension can be realized by an extension $\tilde{T} \to \tilde{S}$ of $T \to S$. We let $\mathcal{X}_{\tilde{T}}$ be the pull back family $\tilde{\mathcal{X}} \times_{\tilde{S}} \tilde{T}$. Let $\tilde{U}_{\alpha} = \tilde{\mathcal{U}}_{\alpha} \times_{\tilde{S}} \tilde{T}$, $\tilde{V}_{\alpha} = \tilde{\mathcal{V}}_{\alpha} \times_{\tilde{S}} \tilde{T}$, $U_{\alpha} = \tilde{U}_{\alpha} \times_{\mathcal{X}_T} \mathcal{X}_T$ and let $V_{\alpha} = \tilde{V}_{\alpha} \times_{\mathcal{X}_{\tilde{T}}} \mathcal{X}_T$. Also we extend $U_{\alpha}/V_{\alpha}$ to a chart $(U_{\alpha}/V_{\alpha}, f_{T,\alpha}, \mathcal{W}_{\alpha})$ of $f_T$. This way $\tilde{U}_{\alpha}/\tilde{V}_{\alpha}$ is a minimal extension of $U_{\alpha}/V_{\alpha}$ to $\mathcal{X}_{\tilde{T}}/\tilde{T}$. Next we let $\tilde{p}_{\alpha}: \tilde{U}_{\alpha} \to \tilde{\mathcal{U}}_{\alpha}$ be the projection. For any $\alpha$ we pick a family of pre-deformable extension $\tilde{h}_{\alpha}: \tilde{U}_{\alpha} \to W[n]$ of $f_{T,\alpha}$. We let $\tilde{p}_{\alpha}: \tilde{U}_{\alpha\beta} \to \tilde{U}_{\alpha}$ be the projection. We then let

$$(1.22) \qquad b_{\alpha\beta} = -\mathbf{d}^{\dagger}(\tilde{h}_{\beta} \circ \tilde{p}_{\beta} - \tilde{h}_{\alpha} \circ \tilde{p}_{\alpha}) \in \operatorname{Hom}_{U_{\alpha\beta}}(f_T^* \Omega_{W[n]}, I)^{\dagger}.$$

We let $b = \{b_{\alpha\beta}\}$, which belongs to $\mathbf{D}^1 \otimes_A I \subset \mathbf{E}^2 \otimes_A I$. It follows from the Remark 1.13 that $b$ is a cocycle, and thus defines a cohomology class $[b]$ in $h^2(\mathbf{E}^{\bullet} \otimes_A I)$. The technical part of the proof is to check that the cohomology class $[b]$ is independent of the choices of $\tilde{T} \to \tilde{S}$ and $\tilde{h}_{\alpha}$. The argument for this is straightforward though tedious, and will be omitted.

We now show that it is an obstruction class to extending $f_T$ to families of pre-deformable morphisms over $\tilde{T}$. First, if such extensions do exist, then we can choose $\mathcal{X}_{\tilde{T}}$ and $\tilde{h}_{\alpha}$ be data coming from one of such extensions. Then the corresponding $b' = 0$ as cocycle and thus $[b] = [b'] = 0$. This shows that $[b] = 0$ whenever extensions of $f_T$ exist. Now assume $[b] = 0$. We first look at the case where the cycle $b \in \mathbf{D}^1 \otimes_A I$ is a coboundary in $\mathbf{D}^{\bullet} \otimes_A I$. Namely, there is a $c \in \mathbf{D}^0 \otimes_A I$ so that $b = \partial(c)$. Let $c = \{c_{\alpha}\}$ with $c_{\alpha} \in \operatorname{Hom}_{U_{\alpha}}(f_T^* \Omega_{W[n]^{\dagger}}, I)^{\dagger}$. Then by Lemma 1.12, we can find pre-deformable extension $f_{\tilde{T},\alpha} : \tilde{U}_\alpha \to W[n]$ of $f_{T,\alpha}$ so that $\mathbf{d}^\dagger(f_{\tilde{T},\alpha} - \tilde{h}_\alpha) = c_\alpha$. Then over $\tilde{U}_{\alpha\beta}$, the difference of the pull backs $f_{\tilde{T},\beta} \circ \tilde{p}_\beta$ and $f_{\tilde{T},\alpha} \circ \tilde{p}_\alpha$ is

$$\mathbf{d}^\dagger(f_{\tilde{T},\beta} \circ \tilde{p}_\beta - f_{\tilde{T},\alpha} \circ \tilde{p}_\alpha) = \mathbf{d}^\dagger(f_{\tilde{T},\beta} \circ \tilde{p}_\beta - \tilde{h}_\beta \circ \tilde{p}_\beta) - \mathbf{d}^\dagger(f_{\tilde{T},\alpha} \circ \tilde{p}_\alpha - \tilde{h}_\alpha \circ \tilde{p}_\alpha) +$$
$$+ \mathbf{d}^\dagger(\tilde{h}_\beta \circ \tilde{p}_\beta - \tilde{h}_\alpha \circ \tilde{p}_\alpha) = (c_\beta - c_\alpha) - b_{\alpha\beta} = 0.$$

Hence $\{f_{\tilde{T},\alpha}\}$ patch together to form a desired extension $\mathcal{X}_{\tilde{T}} \to W[n]$ of $f_T$.

In general, assume $[b] = 0$, then there is an $a \in F^0 \otimes_A I$ and $c \in \mathbf{D}^0 \otimes_A I$ so that

$$b = d_{\mathbf{E}}(a, c) = \delta(a) + \partial(c).$$

Now let $\tilde{r}' : \tilde{T} \to \tilde{S}$ be a new extension of $T \to S$ so that $\mathbf{d}(\tilde{r}' - \tilde{r}) = a$, where $\tilde{r} : \tilde{T} \to \tilde{S}$ is the morphism used to construct the cycle $b$. One checks that if one uses the new extension $\tilde{\mathcal{X}} \times_{\tilde{S}} \tilde{T}$, where $\tilde{T} \to \tilde{S}$ is via $\tilde{r}'$, to construct a similar cycle $b' \in \mathbf{D}^1 \otimes_A I$, then $b'$ is a coboundary in $\mathbf{D}^\bullet \otimes_A I$. This reduces the situation to the previous case studied, and hence confirms that a pre-deformable extension of $f_T$ over $\tilde{T}$ can be found. This shows that $[b]$ is an obstruction class to extending $f_T$ to families of pre-deformable morphisms over $\tilde{T}$, or equivalently the obstruction class to extending $T \to S$ to $\tilde{T} \to S$.

We define

$$\mathfrak{ob}(B, I, \varphi) = [b] \in h^2(\mathbf{E}^\bullet \otimes_A I).$$

Again it is direct to check that this assignment defines a functor from $\mathfrak{Tri}_S$ to $\mathfrak{Mod}_S^*$. This completes the proof of the Proposition. $\qquad\square$

We now summarize the result of this section in the following theorem. We need a vanishing Lemma whose proof will be provided in Proposition 5.1.

**Lemma 1.19.** *For sufficiently fine $\Lambda$, we have $h^i(\mathbf{D}^\bullet \otimes_A I) = 0$ for any $A$-module $I$ and $i \geq 2$.*

**Theorem 1.20.** *Let $S$ be an affine chart of the moduli stack $\mathfrak{M}(W[n], \Gamma)^{st}$. Then the obstruction theory just defined is a perfect obstruction theory of $S$.*

*Proof.* Let $f : \mathcal{X} \to W[n]$ be the universal family over $S$. If suffices to show that there is a complex of finitely generated free $A$-modules $E^\bullet = [E^1 \to E^2]$ so that it is quasi-isomorphic to $\mathbf{E}^\bullet$, where $\mathbf{E}^\bullet$ is the complex associated to a sufficiently find atlas $\Lambda$ of $f$. Since $h^i(\mathbf{D}^\bullet \otimes_A I) = 0$ for $i \geq 2$ and any $I$, $h^i(\mathbf{E}^\bullet \otimes_A I) = 0$ for $i \geq 3$ and any $I$. Hence there is a bounded subcomplex $\tilde{E}^\bullet$ of flat $A$-modules so that it is quasi-isomorphic to $\mathbf{E}^\bullet$. Then we can apply the standard technique [Har, III.12] to find a bounded subcomplex $E^\bullet$ of finitely generated free $A$-modules that is quasi-isomorphic to $\mathbf{E}^\bullet$. Finally, since $h^i(\mathbf{E} \otimes_A I) \neq 0$ only for $i = 1$ and 2, we can choose $E^\bullet$ to be of the form $[E^1 \to E^2]$. This proves the Theorem. $\qquad\square$

In the Appendix, we will express the cohomology $h^\bullet(\mathbf{E}^\bullet)$ in terms of some know cohomologies.

1.3. **Obstruction to deforming relative stable morphisms.** We will follow the notation developed in [Li, section 4] concerning relative stable morphisms. Let $(Z, D)$ be a polarized relative pair and let $Z[n]^{\text{rel}}$ be the expanded relative pair constructed there. Recall that $Z[n]^{\text{rel}}$ consists of a proper variety $Z[n]$ over $\mathbf{A}^n$ and a smooth divisor $D[n] \subset Z[n]$ that is isomorphic to $D \times \mathbf{A}^n$ under the projection $Z[n] \to Z \times \mathbf{A}^n$. The pair $(Z[n], D[n])$ over $\mathbf{A}^n$ also admits an equivariant $G[n]$ action ($= G_m^{\times n}$) whose action on $\mathbf{A}^n$ is the standard one[9]. The fibers of $Z[n]/\mathbf{A}^n$ has at most normal crossing singularities and the singular locus of all the fibers of $Z[n]/\mathbf{A}^n$ is a disjoint union of smooth varieties $\mathbf{B}_1, \cdots, \mathbf{B}_n$, indexed so that $B_l$ surjects to the $l$-th coordinate hyperplane $H_l \subset \mathbf{A}^n$.

In [Li] we used admissible graph to describe the topological type of relative stable morphisms to $(Z[n], D[n])$. Recall that a weighted graph $\Gamma$ (introduced in [Li]) consists of a collection of vertices $V_\Gamma$, an ordered collection of weighted roots $R_\Gamma$ and an ordered collection of legs $L_\Gamma$ (the later two are half edges with one ends attached to vertices and the other ends left free) plus two weight functions $g, b : V_\Gamma \to \mathbb{Z}_{\geq 0}$ and a multiplicity assignment $\mu : R_\Gamma \to \mathbb{Z}^+$. We require $\Gamma$ to be relatively connected in the sense that either $\Gamma$ is connected or each of its vertex has at least one root attached to it. Recall that a relative morphism to $Z[n]^{\mathrm{rel}}$ of type $\Gamma$ consists of a pointed complete nodal curve $X$ and a morphism $f : X \to Z[n]$ that has the following property: First the marked points of $X$ are labeled by the roots and the legs of $\Gamma$. For consistency, we will reserve the notation $p_i \in X$ for $i = 1, \cdots, k$ and $q_j \in X$ for $j = 1, \cdots, r$ to denote the marked points associated to the ordered legs and roots of $\Gamma$. we will use $\mu_1, \cdots, \mu_r$ to denote the weights of the respective roots of $\Gamma$; Secondly, the connected components of $X$ are labeled by $a \in V_\Gamma$ and the arithmetic genus of the component $X_a$ is $g(a)$. Further, in case a root or a leg is attached to a vertex $a$ then its associated marked point must lie in the connected component $X_a$; Thirdly, restricting to each connected component $X_a$ the morphism $f|_{X_a}$ with all the marked points in $X_a$ is an ordinary stable morphism to $Z[n]$ of degree $b(a)$ (using the polarization on $Z$ chosen implicitly); Lastly, as divisor $f^{-1}(D[n]) = \mu_1 q_1 + \cdots + \mu_r q_r$. Here we implicitly assume the former is a proper divisor.

We recall the notion of relative stable morphisms to $\mathfrak{Z}^{\mathrm{rel}}$ defined in [Li]. Let $f : X \to Z[n]$ be a relative morphism of type $\Gamma$, as described. Recall fibers of $Z[n]/\mathbf{A}^n$ have normal crossing singularities along $\mathbf{B}_1, \cdots, \mathbf{B}_n$. We say $f$ is pre-deformable if it is pre-deformable along all $\mathbf{B}_i$, as defined in [Li, Section 2]. We say $f$ is stable as a relative morphism to $\mathfrak{Z}^{\mathrm{rel}}$ if $f$ is pre-deformable and $\mathrm{Aut}_{\mathfrak{Z}}(f)$ is finite. Here $\mathrm{Aut}_{\mathfrak{Z}}(f)$ is the group of all pairs $(a, b)$ where $a : X \to X$ are automorphisms and $b \in G[n]$ so that $f \circ a = f^b$. (Here we view $b$ as an automorphism $b : Z[n] \to Z[n]$ using the $G[n]$ action on $Z[n]$ and $f^b$ is the composite of $f$ with $b$.) As was proved in [Li], the moduli of all relative morphisms to $Z[n]^{\mathrm{rel}}$ of type $\Gamma$ that are stable as morphisms to $\mathfrak{Z}^{\mathrm{rel}}$ form a Deligne-Mumford stack. We will denote this stack by $\mathfrak{M}(Z[n]^{\mathrm{rel}}, \Gamma)^{st}$. The goal of this subsection is to describe the obstruction theory of this moduli stack.

Similar to the case $f : \mathcal{X} \to W[n]$, for any family of relative stable morphisms $f : \mathcal{X} \to Z[n]$ over $S$ there is a canonical log structure on $\mathcal{X}/S$ and on $Z[n]/\mathbf{A}^n$ that makes $f$ a morphism between log schemes. First, the log structure on $\mathbf{A}^n$ (resp. $Z[n]$) is given by the divisor $\cup_{l=1}^n H_l \subset \mathbf{A}^n$ (resp. $Z[n] \times_{\mathbf{A}^1} 0 \cup D[n] \subset Z[n]$). As to the log structure on $\mathcal{X}/S$, we first note that if we let $W/\mathbf{A}^1$ be $Z[1]/\mathbf{A}^1$, then $Z[n]/\mathbf{A}^n = W[n-1]/\mathbf{A}^n$ and $f$ is a family of pre-deformable morphisms to $W[n-1]$. We let $(\mathcal{X}, \mathcal{M}')$ and $(S, \mathcal{N})$ be the associated log structures of $f : \mathcal{X} \to W[n-1]$. On the other hand, let $r$ be the number of distinguished marked points of $f$. Then

---

In short, $Z[1]$ is the blowing up of $Z \times \mathbf{A}^1$ along $D \times 0$ and $D[1]$ is the proper transform of $D \times \mathbf{A}^1$. $Z[2]$ is the blowing up of $Z[1] \times \mathbf{A}^1$ along $D[1] \times \mathbf{A}^1$, etc. The $G[n]$ action is the unique lifting of its standard action on $Z \times \mathbf{A}^n$.

locally the log structure of $f^{-1}(D[n]) \subset \mathcal{X}$ is given by a pre-log structure $\mathbb{N}^r \to \mathcal{O}_U$ for open $U \subset \mathcal{X}$. The desired log structure of $\mathcal{X}$ over $U$ is given by the associated log structure of the pre-log structure $\mathcal{M}'|_U \oplus \mathbb{N}^r \to \mathcal{O}_U$. It is obvious that this gives $\mathcal{X}$ a log structure $\mathcal{M}$, making it a log-scheme over $(S, \mathcal{N})$ and making $f$ a log-morphism between $\mathcal{X}^\dagger/S^\dagger \to Z[n]^\dagger/\mathbf{A}^{n\dagger}$.

Let $S = \operatorname{Spec} A$ be an affine chart of $\mathfrak{M}(Z[n]^{\mathrm{rel}}, \Gamma)^{st}$ with $f : \mathcal{X} \to Z[n]$ the universal family and $q_i$, $p_j : S \to \mathcal{X}$ its marked sections. We cover $f$ by charts of the first or the second kind. Let these charts be $(\mathcal{U}_\alpha/\mathcal{V}_\alpha, f_\alpha, \mathcal{Z}_\alpha)$ indexed by $\Lambda$ as defined before (1.5) with $W[n]$ (resp. $\mathbf{D}_l$; resp. $n+1$) replaced by $Z[n]$ (resp. $\mathbf{B}_l$; resp. $n$). We let $\mathcal{D} \subset \mathcal{X}$ be the divisor of the locus of all marked sections of $\mathcal{X}/S$. We now fix a complex of finite rank free $A$-modules $F^\bullet = [F^0 \xrightarrow{d} F^1]$ so that

$$h^\bullet(F^\bullet \otimes_A I) = \operatorname{Ext}^\bullet_{\mathcal{X}}(\Omega_{\mathcal{X}/S}(\mathcal{D}), \mathcal{I}), \quad \mathcal{I} = I \otimes_{\mathcal{O}_S} \mathcal{O}_{\mathcal{X}}$$

for all $A$-module $I$. We form the group $\operatorname{Hom}_{\mathcal{U}_\alpha}(f^*\Omega_{Z[n]}, I)^\dagger$ as follows: In case $\mathcal{U}_\alpha/\mathcal{V}_\alpha$ is a chart away from $f^{-1}(D[n])$, this group is $\operatorname{Hom}_{\mathcal{U}_\alpha}(f^*\Omega_{W[n-1]}, I)^\dagger$ with $W[n-1]/\mathbf{A}^n = Z[n]/\mathbf{A}^n$ understood. Now let $\mathcal{U}_\alpha/\mathcal{V}_\alpha$ be a chart of some points in $f^{-1}(D[n])$ that is away from $f^{-1}(\mathbf{B})$. Then by shrinking $\mathcal{U}_\alpha/\mathcal{V}_\alpha$ and $\mathcal{Z}_\alpha$ if necessary, we can assume that there is a section $z_\alpha \in \Gamma(\mathcal{O}_{\mathcal{Z}_\alpha})$ so that $z_\alpha = 0$ is the divisor $\mathcal{Z}_\alpha \cap D[n]$. We define $\operatorname{Hom}_{\mathcal{U}_\alpha}(f^*\Omega_{Z[n]}, I)^\dagger$ be the subgroup of

$$(\varphi, \eta) \in \operatorname{Hom}_{\mathcal{U}_\alpha}(f^*\Omega_{Z[n]}, \mathcal{I}_\alpha) \oplus \mathcal{I}_\alpha$$

so that

$$\varphi(f^*dz_\alpha) = f_\alpha^*(z_\alpha) \cdot \eta \quad \text{and} \quad \varphi(dt_l) \in I_\alpha, \quad \forall l.$$

Using the $A$-modules $\operatorname{Hom}_{\mathcal{U}_\alpha}(f^*\Omega_{Z[n]}, I)^\dagger$, we can form the complex $(\mathbf{D}(I)^\bullet, \partial)$, just as we did for $W[n]$ before (1.20). As before, the element $\mathbf{1} \in F^1 \otimes_A F^{1\vee}$ defines a flat extension of $\mathcal{X}/S$ to $\tilde{\mathcal{X}}/\tilde{S}$, where $\tilde{S} \triangleq \operatorname{Spec} A * F^{1\vee}$, with extended sections $\tilde{q}_i$, $\tilde{p}_j : \tilde{S} \to \tilde{\mathcal{X}}$. We let $\tilde{\mathcal{U}}_\alpha/\tilde{\mathcal{V}}_\alpha$ be the minimal extension of $\mathcal{U}_\alpha/\mathcal{V}_\alpha$ as an étale neighborhood of $\tilde{\mathcal{X}}/\tilde{S}$. We then pick $\zeta_\alpha : \tilde{\mathcal{U}}_\alpha \to Z[n]$ that is an extension of $f_\alpha : \mathcal{U}_\alpha \to Z[n]$ so that $\zeta_\alpha$ is pre-deformable and

$$(1.23) \qquad \zeta_\alpha^{-1}(D[n]) = \sum_{j=1}^r \mu_j \tilde{q}_j(\tilde{\mathcal{V}}_\alpha).$$

To construct the corresponding complex $\mathbf{E}(I)^\bullet$ we need two homomorphisms

$$(1.24) \qquad \zeta_\alpha(\cdot) : F^0 \longrightarrow \operatorname{Hom}_{\mathcal{U}_\alpha}(f^*\Omega_{Z[n]}, A)^\dagger$$

and

$$(1.25) \qquad \zeta_{\alpha\beta}(\cdot) : F^1 \longrightarrow \operatorname{Hom}_{\mathcal{U}_{\alpha\beta}}(f^*\Omega_{Z[n]}, A)^\dagger.$$

First, as before we argue that the difference of $\zeta_\alpha$ and $\zeta_\beta$ over $\tilde{\mathcal{U}}_{\alpha\beta}$ canonically defines an element $\zeta_{\alpha\beta} \in \operatorname{Hom}_{\mathcal{U}_{\alpha\beta}}(f^*\Omega_{Z[n]}, F^{1\vee})^\dagger$, which naturally defines a homomorphism as required in (1.25). Here the log differential, namely $f^*(dz_\alpha)/f^*(z_\alpha) \mapsto \eta$, appear near $D[n]$ because of the constraint (1.23). The construction of (1.24) is similar. Namely, locally over $\operatorname{Spec} A * F^{0\vee}$ there are two extensions of $f$. One is given by the pull back of $\zeta_\alpha$ and the other is given by the pull back of $f$. Their difference then gives rise to the homomorphism (1.24).

Once all such data are constructed, we then go ahead to form the homomorphism of complexes $F^\bullet \otimes_A I \to \mathbf{D}(I)^\bullet$, form a new complex $\mathbf{E}(I)^\bullet$ and check that there is a complex of finite rank free $A$-modules $E^\bullet = [E^1 \to E^2]$ so that it is quasi-isomorphic to $\mathbf{E}^\bullet$, parallel to the argument in the previous subsection.

We now state the main result of this section.

**Lemma 1.21.** *Let the notation be as before and let $\Lambda$ be a sufficiently fine covering of $f$ by charts of the first or the second kinds. Then the complex $\mathbf{E}^\bullet$ is a complex of flat $A$-modules. Further, for any $A$-module $I$ we have $\mathbf{E}(I)^\bullet = \mathbf{E}^\bullet \otimes_A I$ and $h^i(\mathbf{E}^\bullet \otimes I) = 0$ for $i \neq 1, 2$. In particular, there is a complex $E^\bullet = [E^1 \to E^2]$ of finitely generated free $A$-modules so that it is quasi-isomorphic to $\mathbf{E}^\bullet$.*

**Theorem 1.22.** *Let the notation be as in the previous Lemma and let $S = \operatorname{Spec} A$ be the affine chart of $\mathfrak{M}(Z[n]^{rel}, \Gamma)^{st}$ as before. Then there is a perfect obstruction theory of $S$ taking value in the complex $E^\bullet$. In particular, the functor of the first order deformations $\mathfrak{Def}_S^1$ is isomorphic to $\mathfrak{h}^1(E^\bullet)$ and there is an obstruction assignment $\mathfrak{ob}$ taking values in $\mathfrak{h}^2(E^\bullet)$ that satisfies the required base change property.*

*Proof.* The proof of the Lemma and the Theorem are parallel to that of Lemmas 1.15, 1.19 and the Theorem 1.20. The only new ingredient is about preserving the divisor $f^{-1}(D[n]) = \sum \mu_j q_j(S)$. Since $q_j(S) \subset \mathcal{X}$ is a divisor smooth over $S$ and $D[n] \subset Z[n]$ is a smooth divisor, that the deformation of morphisms preserving this relation is given by the sheaf of log-differentials is well-known, for example see [Kol]. Since the proof is routine and parallel to what we did before, we shall omit it. This completes the proof of the Theorem. $\qquad\square$

## 2. Gromov-Witten invariants

In this section, we will define the virtual moduli cycle of $\mathfrak{M}(\mathfrak{W}, \Gamma)$, $\mathfrak{M}(\mathfrak{W}_t, \Gamma)$ and $\mathfrak{M}(\mathfrak{Z}^{rel}, \Gamma)$, thus defining the Gromov-Witten invariants of the family $W$, of the singular variety $W_0$ and the relative Gromov-Witten invariants of the pair $(Z, D)$. In the next section, we will prove the decomposition (degeneration) formula relating the Gromov-Witten invariants of $W_t$ to the relative Gromov-Witten invariants of the pairs $(Y_1, D_1)$ and $(Y_2, D_2)$.

2.1. **Perfect obstruction theories of $\mathfrak{M}(\mathfrak{W}, \Gamma)$ and $\mathfrak{M}(\mathfrak{Z}^{\mathbf{rel}}, \Gamma)$.** Recall that the construction of the virtual cycles of moduli stacks is based on the choice of their perfect obstruction theories. In this section, we will show that the perfect obstruction theories constructed in the previous section naturally induce perfect obstruction theories of $\mathfrak{M}(\mathfrak{W}, \Gamma)$ and $\mathfrak{M}(\mathfrak{Z}^{rel}, \Gamma)$.

Let $\mathbf{M}$ be a proper Deligne-Mumford stack with an atlas $\Lambda$ consisting of finitely many affine étale morphisms $\iota_\alpha \colon S_\alpha = \operatorname{Spec} A_\alpha \to \mathbf{M}$. We first recall the definition of a perfect obstruction theory of $\mathbf{M}$.

**Definition 2.1.** *A perfect obstruction theory of $\mathbf{M}$ (over the atlas $\Lambda$) consists of a choice of perfect obstruction theory $(E_\alpha^\bullet, \mathfrak{ob}_\alpha)$ of $S_\alpha$ for each $\alpha \in \Lambda$ so that they satisfy the following compatibility condition: First, let $\mathcal{Ob}_\alpha$ be the obstruction sheaf (i.e. $= \operatorname{Coker}\{E_\alpha^1 \to E_\alpha^2\}$), then the collection $\{\mathcal{Ob}_\alpha\}_{\alpha \in \Lambda}$ descends to a (global) sheaf of $\mathcal{O}_\mathbf{M}$-modules $\mathcal{Ob}_\mathbf{M}$. Secondly, the obstruction assignments $\mathfrak{ob}_\alpha$ and $\mathfrak{ob}_\beta$ are identical when pulled back to $S_{\alpha\beta}$, using the given isomorphisms.*

Note that $\{\mathcal{Ob}_\alpha\}$ descends means that over $S_{\alpha\beta}$ the pull back of $\mathcal{Ob}_\alpha$ and of $\mathcal{Ob}_\beta$ are isomorphic and that such isomorphisms satisfy the cocycle condition on $S_{\alpha\beta\gamma}$.

Our immediate goal is to show that the obstruction theories of $\mathfrak{M}(W[n], \Gamma)^{st}$ naturally induces a perfect obstruction theory of $\mathfrak{M}(\mathfrak{W}, \Gamma)$.

**Theorem 2.2.** *There is a natural perfect obstruction theory of $\mathfrak{M}(\mathfrak{W}, \Gamma)$ induced by the perfect obstruction theories of $\mathfrak{M}(W[n], \Gamma)^{st}$ constructed in the previous section.*

*Proof.* Let $S$ be an affine chart of $\mathfrak{M}(\mathfrak{W}, \Gamma)$. Without lose of generality, we can assume $S$ is one of the chart constructed in the proof of Theorem 3.10 in [Li]. Namely, there is a chart $\bar{S} \subset \mathfrak{M}(W[n], \Gamma)^{st}$ for some $n$ so that $S$ is a closed subscheme of $\bar{S}$ and $S \to \mathfrak{M}(\mathfrak{W}, \Gamma)$ is induced by $S \to \bar{S}$ and $\mathfrak{M}(W[n], \Gamma)^{st} \to \mathfrak{M}(\mathfrak{W}, \Gamma)$. We let $S = \operatorname{Spec} A$ and $\bar{S} = \operatorname{Spec} \bar{A}$. We then let $\bar{E}^{\bullet} = [\bar{E}^1 \to \bar{E}^2]$ be the complex of $\bar{A}$-modules provided by Theorem 1.20 for the chart $\bar{S}$.

We begin with the functor of the first order deformations in $S$. As argued in the proof of [Li, Theorem 3.10], there is a neighborhood $U$ of $S \times \{e\} \subset S \times G[n]$ so that the morphism $S \times G[n] \to \mathfrak{M}(W[n], \Gamma)^{st}$ induced by the $G[n]$-action lifts to an étale $\phi \colon U \to \bar{S}$. Since $U \to \bar{S}$ is étale, each vector $v \in T_e G[n]$ defines a first order deformation of the inclusion $S \to \bar{S}$, and hence an element $\tilde{v} \in h^1(\bar{E}^{\bullet} \otimes_{\bar{A}} A)$. This induces a homomorphism

$$(2.1) \qquad T_e G[n] \otimes_{\Bbbk} A \longrightarrow h^1(\bar{E}^{\bullet} \otimes_{\bar{A}} A) \longrightarrow \bar{E}^1 \otimes_{\bar{A}} A.$$

Since elements in $\mathfrak{M}(W[n], \Gamma)^{st}$ associate to stable morphisms to $\mathfrak{W}$, at each closed point $p \in S$ the homomorphism $T_e G[n] \to T_p \bar{S}/T_p S$ induced by the group action is injective. Hence the cokernel of (2.1) is also free. Now let $E^1$ be the cokernel of (2.1) and let $E^2 = \bar{E}^2 \otimes_{\bar{A}} A$. Let $L \subset \bar{E}^1 \otimes_{\bar{A}} A$ be the image module of (2.1). Since $L$ lies in the image of the second arrow (above), the image of $L$ in $E^2$ is trivial. Hence $\bar{E}^1 \otimes_{\bar{A}} A \to \bar{E}^2 \otimes_{\bar{A}} A$ lifts to $E^1 \to E^2$.

We next show that the natural obstruction theory of $S$ takes values in the cohomology theory of $E^{\bullet}$. First, since (2.1) has free cokernel, $E^{\bullet}$ is a two-term complex of finitely generated free $A$-modules. Secondly, that the functor of the first order deformations in $S$ is given by the functor $\mathfrak{h}^1(E^{\bullet})$ is obvious since the morphism

$$(2.2) \qquad S \times G[n] \supset U \xrightarrow{\phi} \bar{S}$$

induced by the group action is étale near $S \times \{e\}$. Finally, let $(B, I, \varphi) \in Ob(\mathfrak{Tri}_S)$, then it is also an object in $\mathfrak{Tri}_{\bar{S}}$ and thus has an obstruction class

$$\mathfrak{ob}(B, I, \varphi) \in h^2(\bar{E}^{\bullet} \otimes_{\bar{A}} I) = h^2(E^{\bullet} \otimes_A I)$$

to extending $\varphi \colon \operatorname{Spec} B/I \to S$ to $\operatorname{Spec} B \to \bar{S}$. Because $G[n]$ is smooth and $\phi \colon U \to \bar{S}$ (in (2.2)) is étale, $\mathfrak{ob}((B, I, \varphi))$ is also an obstruction class to extending $\varphi$ to $\operatorname{Spec} B \to S$.

Now let $S_{\alpha}$ be charts of $\mathfrak{M}(\mathfrak{W}, \Gamma)$ with $E_{\alpha}^{\bullet}$ their complexes that are part of their obstruction theories. Then it is direct to check that the collection $\{h^2(E_{\alpha}^{\bullet})\}$ form a sheaf over $\mathfrak{M}(\mathfrak{W}, \Gamma)$, and the obstruction assignments $\mathfrak{ob}_{\alpha}$ are compatible. This completes the proof of the Theorem. $\qquad \square$

We now state the theorem concerning the obstruction theory of $\mathfrak{M}(\mathfrak{Z}^{rel}, \Gamma)$.

**Theorem 2.3.** *The perfect obstruction theory of $\mathfrak{M}(Z[n]^{rel}, \Gamma)^{st}$ constructed in the previous section naturally induces a perfect obstruction theory of $\mathfrak{M}(\mathfrak{Z}^{rel}, \Gamma)$.*

*Proof.* We will omit the proof here since it is exactly the same as the proof of the previous Theorem. $\qquad \square$

The next issue is about the obstruction theory of the substack $\mathfrak{M}(\mathfrak{W}_t, \Gamma)$ defined by the fiber product

$$\mathfrak{M}(\mathfrak{W}_t, \Gamma) = \mathfrak{M}(\mathfrak{W}, \Gamma) \times_C t,$$

where $t \in C$ is a closed point. Clearly, when $t \neq 0$ the stack $\mathfrak{M}(\mathfrak{W}_t, \Gamma)$ is naturally isomorphic to the module stack of stable morphisms to $W_t$ of topological type $\Gamma$, which itself admits a natural obstruction theory as worked out in [LT1]. The obstruction theory of $\mathfrak{M}(\mathfrak{W}_0, \Gamma)$ deserves more attention since it was not known before.

We now study the obstruction theory of $\mathfrak{M}(\mathfrak{W}_0, \Gamma)$. Let $W_0[n] = W[n] \times_C 0$. As in the case of $\mathfrak{M}(\mathfrak{W}, \Gamma)$, we only need to work out the obstruction theory of

$$\mathfrak{M}(W_0[n], \Gamma)^{\mathrm{st}} \triangleq \mathfrak{M}(W[n], \Gamma)^{\mathrm{st}} \times_C 0.$$

Let $S = \operatorname{Spec} A$ be an affine chart of $\mathfrak{M}(W[n], \Gamma)^{\mathrm{st}}$. Then $S_0 = S \times_C 0$ is an affine chart of $\mathfrak{M}(W_0[n], \Gamma)^{\mathrm{st}}$. We let $A_0$ be the quotient ring of $A$ so that $S_0 = \operatorname{Spec} A_0$. As before we let $f \colon \mathcal{X} \to W[n]$ be the universal family over $S$ and let $f_0 \colon \mathcal{X}_0 \to W_0[n]$ be the restriction of $f$ to $S_0 \subset S$. We fix a sufficiently fine covering $(\mathcal{U}_\alpha / \mathcal{V}_\alpha, f_\alpha, \mathcal{W}_\alpha)$ of $f$ indexed by $\Lambda$ and let $\mathbf{E}^\bullet$ be the associated complex of $A$-modules constructed in section 1.1. We let $\underline{\mathcal{V}_0}$ be the coverings $\{\mathcal{V}_{0,\alpha}\}$ of $S_0$ with $\mathcal{V}_{0,\alpha} = \mathcal{V}_\alpha \times_C 0$. Similarly we let $\mathcal{U}_{0,\alpha} = \mathcal{U}_\alpha \times_C 0$ and let $f_{0,\alpha} \colon \mathcal{U}_{0,\alpha} \to W_0[n]$ be the restriction of $f$ to $\mathcal{U}_{0,\alpha}$. For any $A_0$-module $I$, we define

$$(2.3) \qquad \Gamma_{\mathcal{V}_{0,\alpha_1 \cdots \alpha_m}}(I) \triangleq I \otimes_{A_0} \Gamma(\mathcal{O}_{\mathcal{V}_{0,\alpha_1 \cdots \alpha_m}}).$$

Using these, we can form a Čech complex $\mathbf{C}^\bullet(\underline{\mathcal{V}_0}, A_0)$ with the standard coboundary operation. Let $\mathbf{E}_0^\bullet = \mathbf{E}^\bullet \otimes_A A_0$. We next construct a homomorphism

$$\bar{\delta}_i \colon \mathbf{E}_0^i \to \mathbf{C}^{i-1}(\underline{\mathcal{V}_0}, A_0) \otimes_{\Bbbk} T_0 C.$$

Here we understand $\mathbf{C}^{-1}(\underline{\mathcal{V}_0}, A_0) = 0$. Let $\xi \in \mathbf{E}_0^1$ be any element. We write $\xi = (a, b)$ with $a \in F^1 \otimes_A A_0$ and $b \in \mathbf{C}^0(\Lambda, \mathcal{H}om(f^* \Omega_{W[n]}, A_0)^\dagger)$, as in section 1. Then by the construction in section 1.1, to each $\alpha \in \Lambda$ the element $a$ defines an extension of $\mathcal{U}_{0,\alpha} / \mathcal{V}_{0,\alpha}$ to $\tilde{\mathcal{U}}_{0,\alpha} / \tilde{\mathcal{V}}_{0,\alpha}$ by $A_{0,\alpha} (= \Gamma_{\mathcal{V}_{0,\alpha}}(A_0))$. The extension $\zeta_\alpha$ chosen before (1.19) induces an extension $\tilde{\zeta}_{0,\alpha} \colon \tilde{\mathcal{U}}_{0,\alpha} \to W[n]$ of $f_{0,\alpha} \colon \mathcal{U}_{0,\alpha} \to W_0[n]$. We let $\pi_n \colon W[n] \to C$ be the tautological projection. We now consider $\pi_n \circ \tilde{\zeta}_{0,\alpha} \colon \tilde{\mathcal{U}}_{0,\alpha} \to C$. Since $\pi_n \circ \tilde{\zeta}_{0,\alpha}|_{\mathcal{U}_{0,\alpha}}$ factor through $0 \in C$,

$$\mathbf{d}(\pi_n \circ \tilde{\zeta}_{0,\alpha} - 0) \in \Gamma_{\mathcal{U}_{0,\alpha}}(\mathcal{A}_{0,\alpha}) \otimes_{\Bbbk} T_0 C, \quad \mathcal{A}_{0,\alpha} = A_0 \otimes_{A_0} \mathcal{O}_{\mathcal{U}_{0,\alpha}}.$$

Because $\tilde{\zeta}_{0,\alpha}$ is a pre-deformable extension, the above element lies in $A_{0,\alpha} \otimes_{\Bbbk} T_0 C$. We define $\bar{\delta}(a)_\alpha$ to be this element. As to $\bar{\delta}(b_\alpha)$, since $b_\alpha \in \Gamma(\mathcal{U}_\alpha, \mathcal{H}om(f^* \Omega_{W[n]}, A_0)^\dagger)$, $b_\alpha$ induces a homomorphism $(\pi_n \circ f)^* T_0^\vee C \to \mathcal{A}_{0,\alpha}$. Again this is an element in $A_{0,\alpha} \otimes_{\Bbbk} T_0 C$. We define $\bar{\delta}_1(b_\alpha)$ to be this element. Clearly, this construction carries over to the case of multi-indices. This defines a map (a homomorphism) of complexes

$$\bar{\delta} \colon \mathbf{E}_0^\bullet \longrightarrow \mathbf{C}^{\bullet-1}(\underline{\mathcal{V}_0}, A_0) \otimes_{\Bbbk} T_0 C.$$

It is direct to check that this is a homomorphism of complexes. We let $\mathbf{F}^\bullet$ be the associated complex defined by $\mathbf{F}^i = \mathbf{E}_0^i \oplus \mathbf{C}^{i-2}(\underline{\mathcal{V}_0}, A_0) \otimes_{\Bbbk} T_0 C$ whose differential is the obvious induced one. Then we have a short exact sequence of complexes

$$0 \Longrightarrow \mathbf{C}^{\bullet-1}(\underline{\mathcal{V}_0}, A_0) \otimes_{\Bbbk} T_0 C \Longrightarrow \mathbf{F}^\bullet \Longrightarrow \mathbf{E}^\bullet \otimes_A A_0 \Longrightarrow 0$$

which induces a long exact sequence of cohomologies (for any $A_0$-module $I$)

$$0 \longrightarrow h^1(\mathbf{F}^\bullet \otimes_{A_0} I) \longrightarrow h^1(\mathbf{E}_0^\bullet \otimes_{A_0} I) \longrightarrow I \otimes_{\Bbbk} T_0 C \longrightarrow$$
$$\longrightarrow h^2(\mathbf{F}^\bullet \otimes_{A_0} I) \longrightarrow h^2(\mathbf{E}_0^\bullet \otimes_{A_0} I) \longrightarrow 0.$$

Here we used the fact that since $\underline{\mathcal{V}_0}$ is an étale covering of $S_0$, $h^j(\mathbf{C}^\bullet(\underline{\mathcal{V}_0}, I)) = I$ when $j = 0$ and vanishes when $j \geq 1$.

Since terms in $\mathbf{C}^\bullet(\underline{\mathcal{V}_0}, A_0)$ are flat $A_0$-modules, we can pick a complex of finitely generated free $A_0$-modules $F^\bullet = [F^1 \to F^2]$ so that $F^\bullet$ is quasi-isomorphic to $\mathbf{F}^\bullet$.

**Proposition 2.4.** *The chart $S_0$ admits a natural perfect obstruction theory taking values in the cohomology of the complex $F^\bullet$.*

*Proof.* We need to check that the functor of the first order deformations $\mathfrak{Def}^1_{S_0}$ is isomorphic to the functor $\mathfrak{h}^1(F^\bullet)$ and that there is an obstruction assignment taking values in $\mathfrak{h}^2(F^\bullet)$ that satisfies the required base change property.

The fact that the functor $\mathfrak{Def}^1_{S_0}$ is isomorphic to $\mathfrak{h}^1(F^\bullet) \equiv \mathfrak{h}^1(\mathbf{F}^\bullet_0)$ follows directly from the definition and will be omitted. Now we construct the obstruction assignment. Let $(B, I, \varphi)$ be any object in $\mathfrak{Tri}_{S_0}$. Let $T = \operatorname{Spec} B/I$ and $f_T: \mathcal{X}_T \to W_0[n]$ be the pull back family under $\varphi: T \to S_0$. Let $\mathcal{U}_{T,\alpha}/\mathcal{V}_{T,\alpha}$ be the pull back of $\mathcal{U}_{0,\alpha}/\mathcal{V}_{0,\alpha}$ and let $f_{T,\alpha}: \mathcal{U}_{T,\alpha} \to W_0[n]$ be the restriction of $f_T$ to $\mathcal{U}_{T,\alpha}$. Recall that in constructing the obstruction class to extending $T \to S_0$ to $\tilde{T} = \operatorname{Spec} B \to S$, we first extend $\mathcal{X}_T/T$ to $\mathcal{X}_{\tilde{T}}/\tilde{T}$ and extend $f_{T,\alpha}: \mathcal{U}_{T,\alpha} \to W_0[n]$ to $f_{\tilde{T},\alpha}: \mathcal{U}_{\tilde{T},\alpha} \to W[n]$ where $\mathcal{U}_{\tilde{T},\alpha}/\mathcal{V}_{\tilde{T},\alpha}$ is the minimal extension of $\mathcal{U}_{T,\alpha}/\mathcal{V}_{T,\alpha}$ in $\mathcal{X}_{\tilde{T}}/\tilde{T}$. We then use the difference of $f_{\tilde{T},\alpha}$ and $f_{\tilde{T},\beta}$ to build a cocycle $a \in \mathbf{E} \otimes_{A_0} I$. Let $\pi_n: W[n] \to C$ be the tautological projection as before. Since $\pi_n \circ f_{T,\alpha}$ factor through $0 \in C$, $\pi_n \circ f_{\tilde{T},\alpha} \in \mathcal{I}_\alpha$. Further because $f_{\tilde{T},\alpha}$ is a pre-deformable extension, it lies in $I_\alpha$. Hence the collection $\{\pi_n \circ f_{\tilde{T},\alpha}\}$ defines a cochain $c \in \mathbf{C}^0(\underline{\mathcal{V}_0}, I)$. It is routine to check that the pair $(a, c) \in \mathbf{F}^2 \otimes_{A_0} I$ is closed, and hence defines a cohomology class $[(a, c)] \in h^2(\mathbf{F}^\bullet \otimes_{A_0} I)$. Further, it is routine to check that this class is independent of the choice of the extensions $f_{\tilde{T},\alpha}$, and that it is an obstruction class to extending $\varphi: T \to S_0$ to $\tilde{T} \to S_0$. We define $\mathfrak{ob}_0$ be the assignment that assigns $(B, I, \varphi) \in Ob(\mathfrak{Tri}_{S_0})$ to this class in $h^2(\mathbf{F}^\bullet \otimes_{A_0} I) = h^2(F^\bullet \otimes_{A_0} I)$. $\square$

**Theorem 2.5.** *The obstruction theories of the charts $S_0$ so defined induce a perfect obstruction theory of $\mathfrak{M}(\mathfrak{W}_0, \Gamma)$.*

*Proof.* The only thing needs to be checked is that the sheaves $h^2(\mathbf{F}^\bullet)$ over the charts $S_0$ of $\mathfrak{M}(\mathfrak{W}_0, \Gamma)$ just constructed descends to a global sheaf over $\mathfrak{M}(\mathfrak{W}_0, \Gamma)$ and that the obstruction assignments over these charts are mutually compatible. This is routine and will be omitted. $\square$

We comment that so far all the results concerning $\mathfrak{W}_0$, including its construction, are based on the existence of the smoothing $W$ of $W_0$. It is not difficult to see that we can construct $W_0[n]$ from $W_0$ directly, without reference to $W$. Therefore, we can define $\mathfrak{W}_0$ and the moduli stack $\mathfrak{M}(\mathfrak{W}_0, \Gamma)$ directly without assuming the existence of $W$. The construction of the perfect obstruction theory of $\mathfrak{M}(\mathfrak{W}_0, \Gamma)$ without using $W$ is a little tricky, but should be doable. Since we will not use this in this paper, we will content with assuming the existence of a smoothing $W$ of $W_0$.

2.2. **Gromov-Witten invariants.** The goal of this subsection is to construct the virtual moduli cycles of the moduli stacks $\mathfrak{M}(\mathfrak{W}, \Gamma)$, $\mathfrak{M}(\mathfrak{W}_t, \Gamma)$ and $\mathfrak{M}(\mathfrak{Z}^{\mathrm{rel}}, \Gamma)$ and to define their respective Gromov-Witten invariants.

Currently, there are two constructions of virtual moduli cycles in algebraic geometry. One is the original construction by Tian and the author. They assumed that the moduli space admits a perfect obstruction theory. They then constructed a global cone that function as a virtual normal cone. Such cone was constructed using the (algebraic) Kuranishi maps of the obstruction theory [LT1, LT2]. In their construction they made a technical assumption that there is a global vector bundle on the moduli space that makes the obstruction sheaf its quotient. With this vector bundle, the cone becomes a subcone of this vector bundle, and the virtual moduli cycle is the intersection of this cone with the zero section of this vector bundle, using Gysin map. The alternative construction of Behrend and Fantechi [Beh, BF] works along a parallel line. They constructed a similar cone as an Artin stack, assuming the moduli space admits a perfect obstruction theory. They then obtain a cone cycle in a vector bundle by assuming the existence of a global vector bundle, as in the original construction of Tian and the author. These two constructions yield identical cycles [KKP]. Recently, by working out the intersection theory on Artin stacks, Kresh [Kr2] showed that one can construct the virtual moduli cycle without relying on the existence of a global vector bundle as mentioned, thus removing this technical condition. This makes the construction of virtual moduli cycles more versatile. After seeing Kresh's work, we realized that by applying a simple trick we can remove the technical condition of the existence of such vector bundles in our construction of the virtual moduli cycles. In the following, we will present this modified construction.

We begin with the general situations. Let $\mathbf{M}$ be a proper and separated DM-stack. It follows from the definition that $\mathbf{M}$ can be covered by a quasi-projective scheme. We let $\mathcal{O}b$ be a sheaf of $\mathcal{O}_{\mathbf{M}}$-modules. We assume that there is a finite collection of schemes $S_\alpha$ and smooth morphisms $\rho_\alpha : S_\alpha \to \mathbf{M}$, indexed by a set $\Lambda$, so that the collection of images $\rho_\alpha(S_\alpha) \subset \mathbf{M}$ form an open covering of $\mathbf{M}$. We next assume that to each $\alpha \in \Lambda$ there is a locally free sheaf of $\mathcal{O}_{S_\alpha}$-modules $\mathcal{E}_\alpha$, a surjective homomorphism of sheaves $\mathcal{E}_\alpha \to \rho_\alpha^* \mathcal{O}b$ and a cone cycle $C_\alpha \in Z_* \mathrm{Vect}(\mathcal{E}_\alpha)$ that satisfy the following *cycle consistency criteria*. Here we denote by $\mathrm{Vect}(\mathcal{E}_\alpha)$ the vector bundle over $S_\alpha$ so that its sheaf of sections is $\mathcal{E}_\alpha$. In this paper, by abuse of notation we will view a vector bundle as its total space. We first fix a few notations before we state the criteria. Let $p \in \mathbf{M}$ be any closed point. We pick an étale morphism $\varphi : (X, \bar{p}) \to (\mathbf{M}, p)$ and let $\hat{X}_p$ be the formal completion of $X$ along $\bar{p}$. We let $G_p$ be the automorphism group of $p \in \mathbf{M}$. Note that $G_p$ acts naturally on $\hat{X}_p$ and up to $G_p$ the scheme $\hat{X}_p$ is canonical. We next let $V_{\bar{p}} = \varphi^* \mathcal{O}b \otimes_{\mathcal{O}_X} \Bbbk_{\bar{p}}$ and let $V_{\hat{X}_p}$ be the vector bundle $V_{\bar{p}} \times \hat{X}_p$ over $\hat{X}_p$. Again $G_p$ acts on $V_p$ and $V_{\hat{X}_p}$, and up to $G_p$ they are canonical.

**Cycle consistency criteria**. *We say the collection $\mathcal{C} = \{(S_\alpha, \mathcal{E}_\alpha, C_\alpha)\}_\Lambda$ satisfies the cycle consistency criteria at $p \in \mathbf{M}$ if there is a cycle $C_p \in Z_* V_{\hat{X}_p}$ invariant under $G_p$ of which the following hold. Let $\alpha \in \Lambda$ be any index, let $S_{\alpha,p} = S_\alpha \times_{\mathbf{M}} p$ and let $\hat{S}_\alpha = S_\alpha \times_{\mathbf{M}} \hat{X}_p$. We let $\mathrm{pr}_i$ be the i-th projection of the product $\hat{S}_\alpha = S_\alpha \times_{\mathbf{M}} \hat{X}_p$. Then there is a surjective homomorphism of vector bundles*

$$\Phi_1 : \mathrm{Vect}(\mathcal{E}_\alpha) \times_{S_\alpha} \hat{S}_\alpha \longrightarrow V_{\hat{X}_p} \times_{\hat{X}_p} \hat{S}_\alpha$$

*extending the canonical (composite) homomorphism*

$$\mathrm{Vect}(\mathcal{E}_\alpha) \times_{S_\alpha} S_{\alpha,p} \longrightarrow \mathrm{Vect}(\mathrm{pr}_1^* \rho_\alpha^* \mathcal{O}b|_{S_{\alpha,p}}) \equiv \mathrm{Vect}(\mathrm{pr}_2^* \varphi^* \mathcal{O}b|_{S_{\alpha,p}}) \equiv V_{\bar{p}} \times S_{\alpha,p}{}^{10}$$

---

[10] For sheaves $\mathcal{F}$ of $\mathcal{O}_Z$-modules and closed subscheme $X \subset Z$ we use $\mathcal{F}|_X$ to denote $\mathcal{F} \otimes_{\mathcal{O}_Z} \mathcal{O}_X$.

*so that $\Phi_1^* C_{\tilde\xi} = \Phi_2^* C_\alpha$. Here the first arrow above is induced by $\mathcal{E}_\alpha \to p_\alpha^* \mathcal{O}b$, $\Phi_2$ is the tautological flat morphism $\mathrm{Vect}(\mathcal{E}_\alpha) \times_{S_\alpha} \hat{S}_\alpha \to \mathrm{Vect}(\mathcal{E}_\alpha)$ and $\Phi_i^*$ are the flat pull back homomorphism of cycles.*

We will call the collection $\hat{\mathcal{C}} \triangleq \{C_p \subset V_{\hat{X}_p}\}_{p \in \mathbf{M}}$ satisfying the above critereia the infinitesimal models of the collection $\mathcal{C}$. Accordingly we will call $\mathcal{C}$ a local model of $\hat{\mathcal{C}}$. In the following, we say the collection $\mathcal{C}$ is consistent if there is a $\hat{\mathcal{C}}$ as above that satisfy the above criteria. Conversely, given $\hat{\mathcal{C}}$, we say it can be algebraized if there is a $\mathcal{C}$ so that they satisfy the above criteria. Note that once the infinitesimal models exist, then the property of the local model $\mathcal{C}$ is completely determined by the infinitesimal models. This is the key to many of the results concerning virtual moduli cycles.

Given the collection $\mathcal{C} = \{(S_\alpha, \mathcal{E}_\alpha, C_\alpha)\}_\Lambda$ over $(\mathbf{M}, \mathcal{O}b)$ we now construct a canonical cycle $[\mathcal{C}] \in A_* \mathbf{M}$ as follows. For each $\alpha$ we let $\Xi_\alpha$ be the set of irreducible components of $C_\alpha$. For $a \in \Xi_\alpha$ we denote by $N_a$ the irreducible variety (component) in $C_\alpha$ associated to $a$ and let $m_a$ be the multiplicity of $N_a$ in $C_\alpha$. Then we have

$$(2.4) \qquad C_\alpha = \sum_{a \in \Xi_\alpha} m_a \, N_a \in Z_* \mathrm{Vect}(\mathcal{E}_\alpha).$$

For any $a \in \Xi_\alpha$ we define the base stack of $a$ to be the minimal closed integral substack $\mathbf{Y}_a \subset \mathbf{M}$ so that the natural $N_a \to \mathbf{M}$ factor through $\mathbf{Y}_a \subset \mathbf{M}$. We let $j_a : \mathbf{Y}_a^0 \to \mathbf{Y}_a$ be the (maximal) dense open substack so that the pull back sheaf $j_a^* \mathcal{O}b$ is locally free. Then $\mathbf{F}_a^0 \triangleq \mathrm{Vect}(j_a^* \mathcal{O}b)$ is a vector bundle stack over $\mathbf{Y}_a^0$. Further, the natural morphism

$$(2.5) \qquad \eta_a : \mathrm{Vect}(\mathcal{E}_\alpha)|_{\rho_\alpha^{-1}(\mathbf{Y}_a^0)} \longrightarrow \mathbf{F}_a^0 \equiv \mathrm{Vect}(j_a^* \mathcal{O}b)$$

induced by $\mathcal{E}_\alpha \to \rho_\alpha^* \mathcal{O}b$ is flat. We let $\mathbf{N}_a^0 \subset \mathbf{F}_a^0$ be the image stack of $N_a|_{\rho_\alpha^{-1}(\mathbf{Y}_a^0)}$ under $\eta_a$ with the reduced stack structure. Clearly, by cycle consistency criteria the flat pull back $\eta_a^* \mathbf{N}_a^0$ contains $N_a|_{\rho_\alpha^{-1}(\mathbf{Y}_a^0)}$ as one of its irrducible components. In the following, we will call $\mathbf{N}_a^0 \subset \mathbf{F}_a^0$ the intrinsic representative of $a$. Note that since we choose $\mathbf{Y}_a^0$ to be the maximal possible open substack of $\mathbf{Y}_a$ so that $j_a^* \mathcal{O}b$ is locally free, the open $\mathbf{Y}_a^0 \subset \mathbf{Y}_a$ and the representative $\mathbf{N}_a^0 \subset \mathbf{F}_a^0$ only depend on $a$.

Now let $b \in \Xi_\beta$ be any element with $\mathbf{Y}_b$ and $\mathbf{N}_b^0 \subset \mathbf{F}_b^0$ its base stack and intrinsic representative. We say $a \sim b$ if $\mathbf{Y}_a = \mathbf{Y}_b$ and $\mathbf{N}_a^0 = \mathbf{N}_b^0$ in $\mathbf{F}_a^0 \equiv \mathbf{F}_b^0$. This defines an equivalence relation $\sim$ on $\cup_{\alpha \in \Lambda} \Xi_\alpha$. We define $\Xi = (\bigcup_{\alpha \in \Lambda} \Xi_\alpha)/\sim$. Again by the cycle consistency criteria whenever $a \sim b$ then $m_a = m_b$. Hence each $a \in \Xi$ has an associated multiplicity $m_a$, a base substack $\mathbf{Y}_a$ and an intrinsic representative $\mathbf{N}_a^0 \subset \mathbf{F}_a^0$ over an open substack $\mathbf{Y}_a^0 \subset \mathbf{Y}_a$.

Assuming there is a global locally free sheaf $\mathcal{E}$ on $\mathbf{M}$ making $\mathcal{O}b$ its quotient sheaf, then over each $\mathbf{Y}_a^0$ we have a flat projection $\mathrm{Vect}(\mathcal{E})|_{\mathbf{Y}_a^0} \to \mathbf{F}_a^0$. We let $\mathbf{N}_a \subset \mathrm{Vect}(\mathcal{E})$ be the closure of the pull back of $\mathbf{N}_a^0 \subset \mathbf{F}_a^0$ under this projection. The associated cycle $[\mathcal{C}]$ is then defined to be

$$(2.6) \qquad [\mathcal{C}] = \sum_{a \in \Xi} m_a 0^!_{\mathrm{Vect}(\mathcal{E})}[\mathbf{N}_a] \in A_* \mathbf{M},$$

where $0^!_{\mathrm{Vect}(\mathcal{E})}$ is the Gysin map of the 0-section of $\mathrm{Vect}(\mathcal{E})$. This is essentially the original construction of Tian and the author.

We now back to the general situation (without assuming the existence of such $\mathcal{E}$). We need to define a map $\xi : \Xi \to A_*\mathbf{M}$ so that $\xi(a)$ is the cycle $0^!_{\mathrm{Vect}(\mathcal{E})}[\mathbf{N}_a]$ should a global $\mathcal{E}$ exist. Let $a \in \Xi$ be any element. Since $\mathbf{M}$ can be covered by a quasi-projective scheme, there is a normal projective variety $Y_a$ and a generically finite surjective morphism $\varphi_a : Y_a \to \mathbf{Y}_a$. By abuse of notation, we also view $\varphi_a$ as the composite of $Y_a \to \mathbf{Y}_a$ with $\mathbf{Y}_a \to \mathbf{M}$. Since $Y_a$ is projective, there is a locally free sheaf of $\mathcal{O}_{Y_a}$-modules $\mathcal{F}_a$ so that $\varphi_a^*\mathcal{O}b$ is a quotient sheaf of $\mathcal{F}_a$. We denote by $F_a$ the vector bundle $\mathrm{Vect}(\mathcal{F}_a)$ over $Y_a$. Let $Y_a^0 \subset Y_a$ be a dense open subset so that $Y_a^0 \to \mathbf{Y}_a^0$ is étale. Then the morphism $F_a|_{Y_a^0} \to \mathbf{F}_a^0$ induced by $\mathcal{F}_a \to \varphi_a^*\mathcal{O}b$ is a flat morphism. We let $N_a$ be the closure in $F_a$ of the flat pull-back of $\mathbf{N}_a^0 \subset \mathbf{F}_a^0$. Note that $N_a$ only depends on $Y_a$ and $\mathcal{F}_a \to \varphi_a^*\mathcal{O}b$. The cycle $N_a \subset F_a$ will be called a representative of $a \in \Xi$. With $N_a \subset F_a$ chosen, we define

$$(2.7) \qquad \xi(a) = \deg(\varphi_a)^{-1} \varphi_{a*} 0^!_{F_a}[N_a],$$

where $0_{F_a}$ is the zero section of $F_a$, $0^!_{F_a}$ is the Gysin homomorphism $Z_*F_a \to A_*Y_a$ (of the zero section of $0_{F_a}$) and $\varphi_{a*}$ is the push-forward homomorphism of cycles. The degree $\deg(\varphi_a)$ is the degree of the morphism $\varphi_a : Y_a \to \mathbf{Y}_a$ defined in [Vis]. Finally, we define

$$(2.8) \qquad [\mathcal{C}] = \sum_{a \in \Xi} m_a\, \xi(a) \in A_*\mathbf{M}.$$

Note that this construction coincides with that in (2.6) in case a global locally free sheaf $\mathcal{E}$ exists.

We will call this construction the basic construction and call $[\mathcal{C}]$ the associated cycle of the collection $\mathcal{C}$.

**Lemma 2.6.** *Let the notation be as before. Then $\xi(a)$ is independent of the choice of $Y_a$ and $\mathcal{F}_a$.*

*Proof.* Let $\varphi_{a,1} : Y_{a,1} \to \mathbf{Y}_a$ and $\varphi_{a,2} : Y_{a,2} \to \mathbf{Y}_a$ be two normal varieties and generically finite dominant morphisms and let $N_{a,1} \subset F_{a,1}$ and $N_{a,2} \subset F_{a,2}$ be the respective choices of the representatives of $a$ over $Y_{a,1}$ and $Y_{a,2}$. To prove the Lemma it suffices to show that

$$\deg(\varphi_{a,1})^{-1} \varphi_{a,1*} 0^!_{F_{a,1}}[N_{a,1}] = \deg(\varphi_{a,2})^{-1} \varphi_{a,2*} 0^!_{F_{a,2}}[N_{a,2}].$$

We let $Y_a$ be the normalization of an irreducible component of $Y_{a,1} \times_{\mathbf{Y}_a} Y_{a,2}$ that is dominant and generically finite over $\mathbf{Y}_a$. We let $p_i : Y_a \to Y_{a,i}$ be the projection induced by the $i$-th projection of $Y_{a,1} \times_{\mathbf{Y}_a} Y_{a,2}$. We pick a locally free sheaf $\mathcal{F}_a$ on $Y_a$ and surjective homomorphism $\mathcal{F}_a \to p_i^*\mathcal{F}_{a,i}$ so that the diagram

$$
\begin{array}{ccc}
\mathcal{F}_a & \longrightarrow & p_1^*\mathcal{F}_{a,1} \\
\downarrow & & \downarrow \\
p_2^*\mathcal{F}_{a,2} & \longrightarrow & p_1^*\varphi_{a,1}^*\mathcal{O}b \equiv p_2^*\varphi_{a,2}^*\mathcal{O}b
\end{array}
$$

is commutative. Now let $F_a = \mathrm{Vect}(\mathcal{F}_a)$, let $U \subset Y_a$ be a dense open subset so that the projections $U \to \mathbf{Y}_a$, $U \to Y_{a,1}$ and $U \to Y_{a,2}$ are flat. Then the flat pull-back of $\mathbf{N}_a^0 \subset \mathbf{F}_a^0$ via the induced $F_a|_U \to \mathbf{N}_a^0$ is identical to the flat pull back of $N_{a,i} \subset F_{a,i}$ under the flat morphism $F_a|_U \to F_{a,i}$. Further, it is direct to check that $\deg(\varphi_a) = \deg(\varphi_{a,i})\deg(p_i)$. Hence if we let $\tilde{N}_{a,i} \subset F_{a,i} \times_{Y_{a,i}} Y_a \triangleq p_i^*F_{a,i}$ be the closure of the flat pull back of $N_{a,i}|_{U_i}$ (under $F_{a,i} \times_{Y_{a,i}} U \to F_{a,i}$), then for $i = 1$ and $2$,

$$\varphi_{a,i*} 0^!_{F_{a,i}} [N_{a,i}] = \deg(p_i)^{-1} \varphi_{a*} 0^!_{p_i^* F_{a,i}} [\bar{N}_{a,i}] = \deg(p_i)^{-1} \varphi_{a*} 0^!_{F_a} [N_a].$$

This, combined with the identity about the degrees, proves the Lemma. $\qquad \blacksquare$

We now show how to apply this construction to construct the virtual moduli cycle of a Deligne-Mumford stack $\mathbf{M}$ endowed with a perfect obstruction theory. Let $\{S_\alpha\}_{\alpha \in \Lambda}$ be an atlas of $\mathbf{M}$ and let $\{(\mathcal{E}_\alpha^\bullet, \mathfrak{ob}_\alpha)\}_\Lambda$ be the data associated to the perfect obstruction theory of $\mathbf{M}$ as in Definition 2.1. Here $\mathcal{E}_\alpha^\bullet = [\mathcal{E}_\alpha^1 \to \mathcal{E}_\alpha^2]$ is a complex of finite rank locally free sheaves of $\mathcal{O}_{S_\alpha}$-modules. We let $\mathcal{O}b_\mathbf{M}$ be the sheaf of $\mathcal{O}_\mathbf{M}$-modules that is the descent of $\mathrm{Coker}\{\mathcal{E}_\alpha^1 \to \mathcal{E}_\alpha^2\}$. Following [LT2, Section 3], to each $\alpha$ we can construct a canonical cone cycle $C_\alpha \subset \mathrm{Vect}(\mathcal{E}_\alpha^2)$, using the relative Kuranishi-maps constructed from the perfect obstruction theory of $S_\alpha$. We recall the main results in [LT2, Section 3]. Let $p \in \mathbf{M}$ be any closed point and let $\bar{p} \in S_\alpha$ be a lift of $p$ for some $\alpha \in \Lambda$. We let $T_p = h^1(\mathcal{E}_\alpha^\bullet \otimes_{\mathcal{O}_{S_\alpha}} \Bbbk_{\bar{p}})$ and $O_p = h^2(\mathcal{E}_\alpha^\bullet \otimes_{\mathcal{O}_{S_\alpha}} \Bbbk_{\bar{p}})$. Then $T_p$ is the Zariski tangent space of $S_\alpha$ at $\bar{p}$ and the obstruction theory $\mathfrak{ob}_\alpha$ induces an obstruction theory to deformation of $\bar{p}$ in $S_\alpha$ taking values in $O_p$. Note that $G_p$ naturally acts on $T_p$ and $O_p$, and both $T_p$ and $O_p$ are canonical (depending only on $p$) up to the symmetry $G_p$. Let $X_p = \mathrm{Spec}\, \Bbbk[\![T_p^\vee]\!]$ and let $f_{\bar{p}} \in \Bbbk[\![T_p^\vee]\!] \otimes_\Bbbk O_p$ be a Kuranishi map of the obstruction theory to deforming $\bar{p}$ in $S_\alpha$. We let $\hat{X}_p = \mathrm{Spec}\, \Bbbk[\![T_p^\vee]\!]/(f_{\bar{p}})$. Then the formal completion of $S_\alpha$ along $\bar{p}$, denoted by $\hat{S}_\alpha$, is canonically isomorphic to $\hat{X}_p$. As before, we let $V_{\hat{X}_p}$ be the total scheme of the vector bundle $O_p \times \hat{X}_p$ over $\hat{X}_p$. Following the recipe in constructing the virtual moduli cycle, we consider the normal cone $C_p$ to $\hat{X}_p$ in $X_p$. The cone $C_p$ is naturally a subcone of $V_{\hat{X}_p}$. Though the Kuranishi map $f_{\bar{p}}$ is not unique and also depends on the choice of the lifting $\bar{p} \in S_\alpha$, the cone $C_p$ is unique in the following sense.

**Lemma 2.7.** *Let the notation be as before. Let $0 \in \hat{X}_p$ be the unique closed point. Then the subscheme $C_p \times_{\hat{X}_p} 0 \subset V_{\hat{X}_p} \times_{\hat{X}_p} 0$ is $G_p$-invariant and only depends on $p \in \mathbf{M}$. In particular, the supports and the multiplicities of the irreducible components of $C_p$ are independent of the choice of $\bar{p}$ and $f_{\bar{p}}$.*

The collection $\{C_p \subset V_{\hat{X}_p}\}_{p \in \mathbf{M}}$ form an infinitesimal models. By constructing relative Kuranishi family, Tian and the author established the following existence result [LT2, Section 3].

**Lemma 2.8.** *Let $\mathcal{F}$ be a locally free sheaf of $\mathcal{O}_{S_\alpha}$-modules and $\eta : \mathcal{F} \to \mathcal{O}b_\alpha$ be a surjective sheaf homomorphism. Then there is a unique cycle $C_\mathcal{F} \in Z_* \mathrm{Vect}(\mathcal{F})$ that satisfies the cycle consistency criteria for the pair $(\mathbf{M}, \mathcal{O}b_\mathbf{M})$ and the infinitesimal models $\{C_p \subset V_{\hat{X}_p}\}_{p \in \mathbf{M}}$.*

Applying Lemma 2.8 to $\rho_\alpha : S_\alpha \to \mathbf{M}$ and the quotient homomorphism $\mathcal{E}_\alpha^2 \to \rho_\alpha^* \mathcal{O}b_\mathbf{M}$ we obtain a cycle $C_\alpha \in Z_* \mathrm{Vect}(\mathcal{E}_\alpha^2)$. Since the cycles $C_p$ are unique, Lemma 2.7 and 2.8 combined shows that the collection $\{(S_\alpha, \mathcal{E}_\alpha^2, C_\alpha)\}_\Lambda$ satisfies the cycle consistency criteria for $(\mathbf{M}, \mathcal{O}b_\mathbf{M})$. We define the virtual moduli cycle of $\mathbf{M}$ (with the given perfect obstruction theory) to be the associated cycle of this collection, and denoted by $[\mathbf{M}]^{\mathrm{virt}}$.

In the end, by applying this construction to the stacks $\mathfrak{M}(\mathfrak{W}, \Gamma)$, $\mathfrak{M}(\mathfrak{W}_0, \Gamma)$ and $\mathfrak{M}(\mathfrak{Z}^{\mathrm{rel}}, \Gamma)$ we obtain moduli cycles $[\mathfrak{M}(\mathfrak{W}, \Gamma)]^{\mathrm{virt}}$, $[\mathfrak{M}(\mathfrak{W}_0, \Gamma)]^{\mathrm{virt}}$ and $[\mathfrak{M}(\mathfrak{Z}^{\mathrm{rel}}, \Gamma)]^{\mathrm{virt}}$. As in [LT2], for $t \in C$ we define the GW-invariant of $W_t$ to be the homomorphism

$$\Psi_\Gamma^{W_t} : H^*(W_t)^{\times k} \times H^*(\mathfrak{M}_{g,k}) \longrightarrow H_0(\mathrm{pt}) \cong \mathbb{Q}$$

defined by

$$\Psi_\Gamma^{W_t}(\alpha, \beta) = \mathbf{q}_{*0}\big(\mathrm{ev}^*(\alpha) \cup \pi_{g,k}^*(\beta)\, [\mathfrak{M}(\mathfrak{W}_t, \Gamma)]^{\mathrm{virt}}\big),$$

where $\pi_{g,k}$ and ev are the forgetful and evaluation morphisms, $g$ and $k$ are the genus and the number of marked points of the topological type $\Gamma$, $\mathbf{q}\colon \mathfrak{M}(\mathfrak{W}_t, \Gamma) \to \{t\} \subset C$ is the projection and $\mathbf{q}_{*0}$ is the push-forward $A_*\mathfrak{M}(\mathfrak{W}_t, \Gamma) \to H_0(pt)$ at degree 0. (Here we use $H_*$ to denote the ordinary homology theory in case the ground field is $\mathbb{C}$. Otherwise one can use Chow rings to define the GW-invariants.)

The Gromov-Witten invariants of $\mathfrak{W}$ is the homomorphism

$$\Psi_\Gamma^{W/C} : H_C^0(R^*\pi_*\mathbb{Q}_W)^{\times k} \times H^*(\mathfrak{M}_{g,k}) \longrightarrow H_2^{\mathrm{BM}}(C) \cong \mathbb{Q}$$

defined via a similar formula with $\mathbf{q}_{*0}$ replaced by $\mathbf{q}_{*1} : A_*\mathfrak{M}(\mathfrak{W}, \Gamma) \to H_2^{BM}(C)$. Here $\mathbb{Q}_W$ is the sheaf of locally constant sections on $W$ taking values in $\mathbb{Q}$, $\pi\colon W \to C$ is the tautological projection and $H_2^{BM}$ is the Borel-Moore homology of the open complex curve $C$.

Now we define the relative Gromov-Witten invariants of $Z^{\mathrm{rel}} = (Z, D)$. Let $\Gamma$ be an admissible weighted graph defined in [Li]. As mentioned before, it determines the topological type of relative morphisms to $\mathfrak{Z}^{\mathrm{rel}}$. We let $\Gamma^o$ be the sub-data (in $\Gamma$) relating to the domain curves. (Namely, the connected components, the genus and both kinds of marked points of the domain curves.) Let $\mathfrak{M}_{\Gamma^o}$ be the moduli space of stable curves with topological types $\Gamma^o$. Here a curve $C$ of topological type $\Gamma^o$ is stable if all connected components of $C$ are stable pointed curves. Clearly, $\mathfrak{M}_{\Gamma^o}$ is a Deligne-Mumford stack. As in the ordinary case, there is a forgetful morphism $\pi_{\Gamma^o} : \mathfrak{M}(\mathfrak{Y}^{\mathrm{rel}}, \Gamma) \longrightarrow \mathfrak{M}_{\Gamma^o}$. We define the relative GW-invariants to be the homomorphism

$$\Psi_\Gamma^{Z^{\mathrm{rel}}} : H^*(Z)^{\times k} \times H^*(\mathfrak{M}_{\Gamma^o}) \longrightarrow H_*(D^r)$$

(recall $k$ and $r$ are the numbers of legs and roots of $\Gamma$) defined by

$$\Psi_\Gamma^{Z^{\mathrm{rel}}}(\alpha, \beta) = \mathbf{q}_*\big(\mathrm{ev}^*(\alpha) \cup \pi_{\Gamma^o}^*(\beta)\, [\mathfrak{M}(\mathfrak{Z}^{\mathrm{rel}}, \Gamma)]^{\mathrm{virt}}\big) \in H_*(D^r).$$

Here

$$(2.9) \qquad \mathbf{q} : \mathfrak{M}(\mathfrak{Z}^{\mathrm{rel}}, \Gamma) \to D^r \quad \text{and} \quad \mathrm{ev} : \mathfrak{M}(\mathfrak{Z}^{\mathrm{rel}}, \Gamma) \longrightarrow Z^k$$

are the morphisms defined by evaluating on the distinguished and the ordinary marked sections respectively.

## 3. Degenerations of Gromov-Witten invariants

In this and the next section, we will prove the degeneration formula of the Gromov-Witten invariants of the family $W/C$ stated in the introduction of this paper. We will state the first version of the degeneration formula in the first subsection. We will state the reduction Lemmas in subsection 6.2. The proof of these Lemmas will be postponed to the next section.

3.1. **The first version of the degeneration formula.** The first step to prove the degeneration formula is to express $\mathfrak{M}(\mathfrak{W}_0, \Gamma)$ as a union of Cartier-divisors in $\mathfrak{M}(\mathfrak{W}, \Gamma)$. Here we fix a $\Gamma = (g, k, b)$ once and for all. We first define the notion of Cartier-divisor of an algebraic stack $\mathbf{M}$.

**Definition 3.1.** *Let $\mathbf{M}$ be an algebraic stack. A Cartier-divisor (in short C-divisor) on $\mathbf{M}$ is a pair $(\mathbf{L}, \mathbf{s})$ where $\mathbf{L}$ is a line bundle on $\mathbf{M}$ and $\mathbf{s}$ is a section of $\mathbf{L}$. An isomorphism between $(\mathbf{L}, \mathbf{s})$ and $(\mathbf{L}', \mathbf{s}')$ consists of an isomorphism $\mathbf{L} \cong \mathbf{L}'$ so that $\mathbf{s} \equiv \mathbf{s}'$ under this isomorphism.*

We comment that a C-divisor $(L, s)$ over a scheme is a pseudo-divisor (defined in [Ful]) via $(L, Z, s)$ where $Z = X - s^{-1}(0)$. Note that when $(L, s)$ and $(L', s')$ are two C-divisors, then $(L, s) \otimes (L', s') \triangleq (L \otimes L', ss')$ is also a C-divisor.

We now let $\Omega$ be the set of all admissible triples defined in [Li, Section 4] and reviewed in the introduction. Recall that $\eta = (\Gamma_1, \Gamma_2, I) \in \Omega$ is an admissible triple if $\Gamma_1$ and $\Gamma_2$ are two weighted graphs of identical numbers of roots and $I$ is an order preserving inclusion $I: [k_1] \to [k]$ where $k_i$ is the number of legs of $\Gamma_i$ and $k = k_1 + k_2$. Let $C_1$ and $C_2$ be two curves of topological types $\Gamma_1$ and $\Gamma_2$, respectively. Then we can identify the $i$-th distinguished marked point $q_{1,i} \in C_1$ with the $i$-th distinguished marked point $q_{2,i} \in C_2$ for all $i$ to obtain a new curve $C \triangleq C_1 \sqcup C_2$. It has $k$ marked points, ordered according to $I$. As part of the definition of $\Omega$, we require that the multiplicities of the $i$-th roots of $\Gamma_1$ and $\Gamma_2$ are identical, that $C$ is connected of genus $g$ and $b = \sum_{x \in V(\Gamma_1) \cup V(\Gamma_2)} b(x)$ ($b(x)$ accounts for the degree of the stable morphism along the connected component labeled by $x$). We will call $C$ the gluing of $C_1$ and $C_2$ along distinguished marked points. This gluing construction can be applied to a pair of families of curves. Hence for each $\eta = (\Gamma_1, \Gamma_2, I) \in \Omega$, we have a closed local immersion of stacks

$$(3.1) \qquad \Phi_\eta : \mathfrak{M}(\mathfrak{Y}_1^{\mathrm{rel}}, \Gamma_1) \times_{D^r} \mathfrak{M}(\mathfrak{Y}_2^{\mathrm{rel}}, \Gamma_2) \longrightarrow \mathfrak{M}(\mathfrak{W}, \Gamma),$$

defined in [Li, (4.8)]. Here $\mathfrak{M}(\mathfrak{Y}_i^{\mathrm{rel}}, \Gamma_i) \to D^r$ is the evaluation morphism $\mathbf{q}$ in (2.9) and the morphism (3.1) is defined by sending any pair $((f_1, \mathcal{X}_1), (f_2, \mathcal{X}_2))$ to the family $(f_1 \sqcup f_2, \mathcal{X}_1 \sqcup \mathcal{X}_2)$. Now let $\eta = (\Gamma_1, \Gamma_2, I) \in \Omega$. Following [Li], we define $\mathfrak{M}(\mathfrak{Y}_1^{\mathrm{rel}} \sqcup \mathfrak{Y}_2^{\mathrm{rel}}, \eta)$ be the image stack of (3.1). As was shown in [Li, Section 4],

$$(3.2) \qquad \mathfrak{M}(\mathfrak{Y}_1^{\mathrm{rel}}, \Gamma_1) \times_{D^r} \mathfrak{M}(\mathfrak{Y}_1^{\mathrm{rel}}, \Gamma_1) \to \mathfrak{M}(\mathfrak{Y}_1^{\mathrm{rel}} \sqcup \mathfrak{Y}_2^{\mathrm{rel}}, \eta)$$

is finite, étale of pure degree $|\mathrm{Eq}(\eta)|$. Here by abuse of notion we also use $\Phi_\eta$ to denote this induced morphism.

In this subsection, to each $\eta \in \Omega$ we will define a C-divisor $(\mathbf{L}_\eta, \mathbf{s}_\eta)$ on $\mathfrak{M}(\mathfrak{W}, \Gamma)$ so that the vanishing locus (as topological space) of $\mathbf{s}_\eta$ is $\mathfrak{M}(\mathfrak{Y}_1^{\mathrm{rel}} \sqcup \mathfrak{Y}_2^{\mathrm{rel}}, \eta)$. We begin with the study of line bundles on $\mathbf{A}^{n+1}$. We continue to use the convention introduced in [Li, Section 1] concerning subsets of $\mathbf{A}^{n+1}$. For $l \in [n+1]$ we will still denote by $\mathbf{H}_l \subset \mathbf{A}^{n+1}$ the $l$-th coordinate hyperplane of $\mathbf{A}^{n+1}$. We define $L_l$ be the line bundle on $\mathbf{A}^{n+1}$ so that its sheaf of sections is $\mathcal{O}_{\mathbf{A}^{n+1}}(\mathbf{H}_l)$ and define $s_l$ be the section of $L_l$ that is the constant section $1 \in \mathcal{O}_{\mathbf{A}^{n+1}}$ under the canonical inclusion $\mathcal{O}_{\mathbf{A}^{n+1}} \subset \mathcal{O}_{\mathbf{A}^{n+1}}(\mathbf{H}_l) \equiv \mathcal{O}_{\mathbf{A}^{n+1}}(L_l)$. Recall that $\mathbf{A}^{n+1}$ is a $G[n]$-subscheme ($G[n] \triangleq GL(1)^{\times n}$) as defined in [Li, Section 1]. Since $\mathbf{H}_l \subset \mathbf{A}^{n+1}$ is invariant under the $G[n]$-action, there is a unique $G[n]$-linearization on $L_l$ so that the section $s_l$ is $G[n]$-invariant. We fix such a linearization. Now let $J: [m+1] \to [n+1]$ be an order preserving embedding. Following the convention in [Li, Section 1], $J$ defines a standard embedding[11] $\gamma_J : \mathbf{A}^{m+1} \to \mathbf{A}^{n+1}$ and hence defines a pull back C-divisor $\gamma_J^*(L_l, s_l)$ on $\mathbf{A}^{m+1}$. There are two possibilities: One is when $l \neq \mathrm{Im}(J)$. Then $\mathrm{Im}(\gamma_J) \cap \mathbf{H}_l = \emptyset$ and hence there is a canonical isomorphism[12] $\gamma_J^*(L_l, s_l) \cong (\mathbf{1}_{\mathbf{A}^{m+1}}, 1)$. The other case is when $J(l') = l$ for some $l' \in [m+1]$, in which case we have $\gamma_J^*(L_l, s_l) \cong (L_{l'}, s_{l'})$.

Now let $J : [n_1 + 1] \to [n_2 + 1]$ be an order preserving embedding. Let $S$ be any scheme, $\tau : S \to \mathbf{A}^{n_1+1}$ and $\rho : S \to G[n_2]$ be two morphisms with $\gamma_J : \mathbf{A}^{n_1+1} \to \mathbf{A}^{n_2+1}$ the standard embedding. As in [Li], we define $(\gamma_J \circ \tau)^\rho : S \to \mathbf{A}^{n_2+1}$ be the morphisms induced by $\gamma_J \circ \tau : S \to \mathbf{A}^{n_2+1}$ and the $G[n_2]$-action on $\mathbf{A}^{n_2+1}$ via $\rho$.

**Lemma 3.2.** *Let $J$, $\tau$ and $\rho$ be as before. In case $l_2 = J(l_1)$ then we have a natural isomorphism $\big((\gamma_J \circ \tau)^\rho\big)^*(L_{l_2}, s_{l_2}) = \tau^*(L_{l_1}, s_{l_1})$. In case $l_2 \notin \mathrm{Im}(J)$ the same identity holds with $(L_{l_1}, s_{l_1})$ replaced by $(\mathbf{1}, 1)$.*

*Proof.* We have the canonical isomorphism $(\gamma_J \circ \tau)^*(L_{l_2}, s_{l_2}) \cong \tau^*(L_{l_1}, s_{l_1})$. The required isomorphism is then induced by the canonical isomorphism

$$\big((\gamma_J \circ \tau)^\rho\big)^*(L_{l_2}, s_{l_2}) \cong (\gamma_J \circ \tau)^*(L_{l_2}, s_{l_2})$$

induced by the $G[n_2]$-linearization on $(L_{l_2}, s_{l_2})$. This proves the Lemma. ∎

We now construct the required C-divisor $(\mathbf{L}_\eta, \mathbf{s}_\eta)$ on $\mathfrak{M}(\mathfrak{W}, \Gamma)$. Let $S \to \mathfrak{M}(\mathfrak{W}, \Gamma)$ be any chart with $f : \mathcal{X} \to \mathcal{W}$ its universal family. Without loss of generality, we can assume that $\mathcal{W} = W[n] \times_{C[n]} S$ via a $\tau : S \to C[n]$. We let

(3.3)               $S_\eta = S \times_{\mathfrak{M}(\mathfrak{W}, \Gamma)} \mathfrak{M}(\mathfrak{Y}_1^{\mathrm{rel}} \sqcup \mathfrak{Y}_2^{\mathrm{rel}}, \eta).$

In case $S_\eta = \emptyset$, we define $(\mathbf{L}_\eta, \mathbf{s}_\eta)|_S$ be $(\mathbf{1}_S, 1)$. When $S_\eta \neq \emptyset$, we consider the tautological projection $\rho_\eta : S_\eta \to \mathfrak{M}(\mathfrak{Y}_1^{\mathrm{rel}} \sqcup \mathfrak{Y}_2^{\mathrm{rel}}, \eta)$ and the composite $\mathbf{p}_i \circ \rho_\eta : S_\eta \to \mathfrak{M}(\mathfrak{Y}_i^{\mathrm{rel}}, \Gamma_i)$. Here

$$\mathbf{p}_i : \mathfrak{M}(\mathfrak{Y}_1^{\mathrm{rel}} \sqcup \mathfrak{Y}_2^{\mathrm{rel}}, \eta) \longrightarrow \mathfrak{M}(\mathfrak{Y}_i^{\mathrm{rel}}, \Gamma_i)$$

is the $i$-th projection, which exists if we replace $S$ by an étale cover of $S$. By shrinking $S$ if necessary, we can assume that the pull back of the universal family of $\mathfrak{M}(\mathfrak{Y}_i^{\mathrm{rel}}, \Gamma_i)$ to $S_\eta$ via the morphism $\mathbf{p}_i \circ \rho_\eta$ is given by a family $f_i : \mathcal{X}_i \to \mathcal{Y}_i$, where $\mathcal{Y}_i$ is an effective relative pair in $\mathfrak{Y}_i^{\mathrm{rel}}(S_\eta)$ associated to a morphism $\tau_i : S_\eta \to \mathbf{A}^{n_i}$. Following the discussion leading to the proof of Proposition [Li, Prop. 4.12], the tautological family over $S_\eta$ of the morphism $\Phi_\eta|_{S_\eta} : S_\eta \to \mathfrak{M}(\mathfrak{W}, \Gamma)$ is then represented by the family

$$f_1 \sqcup f_2 : \mathcal{X}_1 \sqcup \mathcal{X}_2 \longrightarrow \mathcal{Y}_1^o \sqcup \mathcal{Y}_2$$

with $\mathcal{Y}_1^o \sqcup \mathcal{Y}_2 \in \mathfrak{W}(S_\eta)$ given by the morphism $\tau_\eta : S_\eta \to C[n]$, where $n = n_1 + n_2$, defined in [Li, (4.4)][13]. By definition $f_1 \sqcup f_2$ is isomorphic to the restriction of $f$ to the family over $S_\eta$, denoted by $f|_{S_\eta}$. Namely there are isomorphisms shown below

---

[11]For instance in case $J : [2] \to [4]$ is defined by $J(1) = 1$ and $J(2) = 3$, then $\gamma_J : \mathbf{A}^2 \to \mathbf{A}^4$ is defined by $\gamma_J(t_1, t_2) = (t_1, 1, t_2, 1)$.

[12]We use bold $\mathbf{1}_X$ with subscription $X$ to denote the trivial line bundle on $X$.

[13]$\mathcal{Y}_1^o \sqcup \mathcal{Y}_2$ is the result of gluing the distinguished divisors of $\mathcal{Y}_1$ and $\mathcal{Y}_2$ in the obvious way.

that make the following diagram commutative

$$(3.4) \qquad \begin{array}{ccc} \mathcal{X}_1 \sqcup \mathcal{X}_2 & \xrightarrow{f_1 \sqcup f_2} & \mathcal{Y}_1^o \sqcup \mathcal{Y}_2 \\ \cong \downarrow & & \cong \downarrow \\ \mathcal{X} \times_S S_\eta & \xrightarrow{f|_{S_\eta}} & \mathcal{W} \times_S S_\eta. \end{array}$$

Now let $\mathbf{D}_1, \cdots, \mathbf{D}_{n+1}$ be the $n+1$ components of the singular locus of the fibers of $W[n]$ over $C[n]$. For any closed $z \in S_\eta$ there is an integer $l_z \in [n+1]$ so that the images of the distinguished divisors $\mathcal{D}_{1,z} \subset \mathcal{Y}_{1,z}$ and of $\mathcal{D}_{2,z} \subset \mathcal{Y}_{2,z}$ under the obvious morphism

$$\mathcal{D}_{1,z} \cong \mathcal{D}_{2,z} \subset \mathcal{Y}_1^o \sqcup \mathcal{Y}_2 \xrightarrow{\cong} \mathcal{W} \times_S S_\eta \longrightarrow W[n]$$

lie in $\mathbf{D}_{l_z}$. Clearly, $l_z$ is locally constant on $S_\eta$. Hence by shrinking $S$ if necessary we can assume that it is constant on $S_\eta$, say is $l_\eta \in \mathbb{Z}$. In the following we will call $f|_{S_\eta} = f_1 \sqcup f_2$ the $\eta$-decomposition of $f$ and call the divisor $\mathbf{D}_{l_\eta} \subset W[n]$ the locus where the $\eta$-decomposition of $f|_{S_\eta}$ takes place.

**Definition 3.3.** *Let $S \to \mathfrak{M}(\mathfrak{W}, \Gamma)$ be a chart with $f : \mathcal{X} \to \mathcal{W}$ the universal family, where $\mathcal{W} = W[n] \times_{C[n]} S$. We say $S$ is $\eta$-admissible if there is an integer $l \in [n+1]$ so that the tautological $S_\eta \to C[n]$, where $S_\eta$ is defined in (3.3), factor through $C[n] \times_{\mathbf{A}^{n+1}} H_l \subset C[n]$ and that the divisor $\mathbf{D}_l \subset W[n]$ (or the locus $\mathcal{X} \times_{W[n]} \mathbf{D}_l \subset \mathcal{X}$) is where the $\eta$-decomposition of $f|_{S_\eta}$ takes place.*

Clearly for each $\eta \in \Omega$ we can find an atlas $\Lambda$ of $\mathfrak{M}(\mathfrak{W}, \Gamma)$ so that all its charts are $\eta$-admissible. Let $\Lambda$ be such an atlas and let $S_\alpha$ be any chart in this atlas. We let $f_\alpha : \mathcal{X}_\alpha \to W[n_\alpha] \times_{C[n_\alpha]} S_\alpha$ be the universal family with $\tau_\alpha : S_\alpha \to C[n_\alpha]$ the tautological morphism. We let $l_\alpha$ be the integer so that $\mathbf{D}_{l_\alpha} \subset W[n_\alpha]$ is where the $\eta$-decomposition of $f_\alpha|_{S_{\alpha,\eta}}$ takes place. We then define the C-divisor $(\mathbf{L}_{\eta,\alpha}, \mathbf{s}_{\eta,\alpha})$ on $S_\alpha$ to be the pull back of the C-divisor $(L_{l_\alpha}, s_{l_\alpha})$ on $\mathbf{A}^{n_\alpha+1}$ via $S \xrightarrow{\tau_\alpha} C[n_\alpha] \to \mathbf{A}^{n_\alpha+1}$. Because the chart $S_\alpha$ is $\eta$-admissible, the vanishing locus of $\mathbf{s}_{\eta,\alpha}$ is exactly $S_\alpha \cap \mathfrak{M}(\mathfrak{Y}_1^{\mathrm{rel}} \sqcup \mathfrak{Y}_2^{\mathrm{rel}}, \eta)$.

**Lemma 3.4.** *The collection $(\mathbf{L}_{\eta,\alpha}, \mathbf{s}_{\eta,\alpha})_{\alpha \in \Lambda}$ forms a C-divisor on $\mathfrak{M}(\mathfrak{W}, \Gamma)$.*

*Proof.* Let $S_\alpha$ and $S_\beta$ be two charts in $\Lambda$. We consider $S_{\alpha\beta} = S_\alpha \times_{\mathfrak{M}(\mathfrak{W}, \Gamma)} S_\beta$ with its projections $\rho_\alpha$, $\rho_\beta$. Let $f_\alpha$ and $f_\beta$ be the universal families over $S_\alpha$ and $S_\beta$ and let $\rho_\alpha^*(f_\alpha)$ and $\rho_\beta^*(f_\beta)$ be the pull back families. We let the isomorphism of the families $\rho_\alpha^*(f_\alpha)$ and $\rho_\beta^*(f_\beta)$ be given by the diagram

$$(3.5) \qquad \begin{array}{ccc} \mathcal{X}_\alpha \times_{S_\alpha} S_{\alpha\beta} & \xrightarrow{f_\alpha|_{S_{\alpha\beta}}} & \mathcal{W}_\alpha \times_{S_\alpha} S_{\alpha\beta} \\ \varphi_1 \downarrow \cong & & \varphi_2 \downarrow \cong \\ \mathcal{X}_\beta \times_{S_\beta} S_{\alpha\beta} & \xrightarrow{f_\beta|_{S_{\alpha\beta}}} & \mathcal{W}_\beta \times_{S_\beta} S_{\alpha\beta}. \end{array}$$

We distinguish two cases: The first is when $S_{\alpha\beta}$ is disjoint from $\mathfrak{M}(\mathfrak{Y}_1^{\mathrm{rel}} \sqcup \mathfrak{Y}_2^{\mathrm{rel}}, \eta)$. Then the pull backs of $(\mathbf{L}_{\eta,\alpha}, \mathbf{s}_{\eta,\alpha})$ and $(\mathbf{L}_{\eta,\beta}, \mathbf{s}_{\eta,\beta})$ to $S_{\alpha\beta}$ are canonically isomorphic to the trivial C-divisor $(\mathbf{1}_{S_{\alpha\beta}}, 1)$, hence they are naturally isomorphic to each other. The other case is when $S_{\alpha\beta} \cap \mathfrak{M}(\mathfrak{Y}_1^{\mathrm{rel}} \sqcup \mathfrak{Y}_2^{\mathrm{rel}}, \eta) \neq \emptyset$. Since $S_\alpha$ and $S_\beta$ are $\eta$-admissible, they have the associated morphisms $\tau_\alpha : S_\alpha \to C[n_\alpha]$ and $\tau_\beta : S_\beta \to C[n_\beta]$ and the associated integers $l_\alpha$ and $l_\beta$ respectively. Because of the isomorphisms in (3.5), the subscheme $\mathbf{D}_{l_\alpha} \times_{C[n_\alpha]} S_{\alpha\beta} \subset \mathcal{W}_\alpha \times_{S_\alpha} S_{\alpha\beta}$ is isomorphic to $\mathbf{D}_{l_\beta} \times_{C[n_\beta]} S_{\alpha\beta} \subset \mathcal{W}_\beta \times_{S_\beta} S_{\alpha\beta}$ under $\varphi_2$. Now let $T$ be any open subset of $S_{\alpha\beta}$ so that the restriction of $\varphi_2$ to $\mathcal{W}_\alpha \times_{S_\alpha} T \cong \mathcal{W}_\beta \times_{S_\beta} T$ is induced by a sequence of effective arrows[14]. Then by Lemma 3.2, the restriction to $T$ of the pull back $\rho_\alpha^*(\mathbf{L}_{\eta,\alpha}, \mathbf{s}_{\eta,\alpha})$ is canonically isomorphic to the restriction to $T$ of $\rho_\beta^*(\mathbf{L}_{\eta,\beta}, \mathbf{s}_{\eta,\beta})$. By [Li, Lemma 1.8], we can cover $S_{\alpha\beta}$ by such $T$'s. Further, applying Lemma 3.2 again we see immediately that the isomorphisms $\rho_\alpha^*(\mathbf{L}_{\eta,\alpha}, \mathbf{s}_{\eta,\alpha})|_T \cong \rho_\beta^*(\mathbf{L}_{\eta,\beta}, \mathbf{s}_{\eta,\beta})|_T$ patch together to form an isomorphism

$$\varphi_{\beta\alpha} : \rho_\alpha^*(\mathbf{L}_{\eta,\alpha}, \mathbf{s}_{\eta,\alpha}) \overset{\cong}{\longrightarrow} \rho_\beta^*(\mathbf{L}_{\eta,\beta}, \mathbf{s}_{\eta,\beta}).$$

Since the isomorphism $\varphi_{\alpha\beta}$ canonically depends on the isomorphism (3.5), over $S_{\alpha\beta\gamma}$ we have $\varphi_{\alpha\beta} \circ \varphi_{\beta\gamma} = \varphi_{\alpha\gamma}$. Thus the collection $(\mathbf{L}_{\eta,\alpha}, \mathbf{s}_{\eta,\alpha})_{\alpha \in \Lambda}$ coupled with the isomorphisms $\varphi_{\alpha\beta}$ defines a C-divisor on $\mathfrak{M}(\mathfrak{W}, \Gamma)$. We denote the resulting C-divisor by $(\mathbf{L}_\eta, \mathbf{s}_\eta)$. $\qquad\blacksquare$

We now indicate the relation between the C-divisor $(\mathbf{L}_\eta, \mathbf{s}_\eta)$ and the canonical log structure on $\mathfrak{M}(\mathfrak{W}, \Gamma)$. Recall that the log structure defined in section 1.1 defines a canonical log structure on $\mathfrak{M}(W[n], \Gamma)^{\mathrm{st}}$, which is $G[n]$-equivariant and thus descends to a canonical log structure $\mathcal{N}$ on $\mathfrak{M}(\mathfrak{W}, \Gamma)$. The line bundle $\mathbf{L}_\eta$ with the section $\mathbf{s}_\eta$ defines also a log structure $\mathcal{L}_\eta$ on $\mathfrak{M}(\mathfrak{W}, \Gamma)$. It follows from the construction of $N_\alpha$ (see (1.7)) and the $(\mathbf{L}_\eta, \mathbf{s}_\eta)$ that the identity map is a log morphism $(\mathfrak{M}(\mathfrak{W}, \Gamma), \mathcal{L}_\eta) \to (\mathfrak{M}(\mathfrak{W}, \Gamma), \mathcal{N})$. Using the chart $\mathcal{V}$ in (1.7), this is given by the homomorphism of pre-log structure $\mathbb{N} \to N_l \subset N_\mathcal{V}$ defined by $1 \mapsto m_\alpha e_\alpha$ for $\alpha \in K_l$.

Associating to each closed point $t \in C$, considered as an effective divisor in $C$, we have a C-divisor $(L_t, r_t)$ such that $\mathcal{O}_C(L_t) = \mathcal{O}_C(t)$ and that $r_t$ is the section induced by the constant section $1 \in \Gamma(C, \mathcal{O}_C)$ together with the natural homomorphism $\mathcal{O}_C(L_t) \to \mathcal{O}_C(t)$. We let $(\mathbf{L}_t, \mathbf{r}_t)$ be the C-divisor on $\mathfrak{M}(\mathfrak{W}, \Gamma)$ that is the pull back of $(L_t, r_t)$ via the tautological projection $\mathfrak{M}(\mathfrak{W}, \Gamma) \to C$. When $t = 0 \in C$ we denote the corresponding C-divisor by $(\mathbf{L}_0, \mathbf{r}_0)$.

Let $\eta = (\Gamma_1, \Gamma_2, I) \in \Omega$ be any triple having $r$ roots. Then for any permutation $\sigma \in S_r$ we have a new triple $\eta^\sigma = (\Gamma_1^\sigma, \Gamma_2^\sigma, I)$ that is derived by permuting the order of the $r$-ordered roots of $\Gamma_1$ and $\Gamma_2$ according to $\sigma$. We say two elements $\eta_1, \eta_2 \in \Omega$ are similar if there is a permutation $\sigma$ so that $\eta_2 = \eta_1^\sigma$. We let $\bar{\Omega}$ be the set of similar classes of $\Omega$. Note that when $\eta_1$ is similar to $\eta_2$, then the C-divisors $(\mathbf{L}_{\eta_1}, \mathbf{s}_{\eta_1}) \equiv (\mathbf{L}_{\eta_2}, \mathbf{s}_{\eta_2})$.

**Proposition 3.5.** *The tensor product of the C-divisors $\{(\mathbf{L}_\eta, \mathbf{s}_\eta) \mid \eta \in \bar{\Omega}\}$ is isomorphic to $(\mathbf{L}_0, \mathbf{r}_0)$.*

*Proof.* Let $\Lambda$ be an atlas of $\mathfrak{M}(\mathfrak{W}, \Gamma)$ so that all its charts are $\eta$-admissible for all $\eta \in \bar{\Omega}$. Since $\bar{\Omega}$ is a finite set, such atlas does exist. Now let $S_\alpha$ be any chart in $\Lambda$. We let $\bar{\Omega}_\alpha$ be those triples $\eta \in \bar{\Omega}$ so that $\mathrm{Im}(\Phi_\eta) \cap S_\alpha \neq \emptyset$. Recall that $(\mathbf{L}_{\eta,\alpha}, \mathbf{s}_{\eta,\alpha}) \cong (\mathbf{1}_{S_\alpha}, 1)$ canonically when $\eta \notin \bar{\Omega}_\alpha$. Therefore

$$\otimes_{\eta \in \bar{\Omega}}(\mathbf{L}_\eta, \mathbf{s}_\eta)|_{S_\alpha} \equiv \otimes_{\eta \in \bar{\Omega}_\alpha}(\mathbf{L}_{\eta,\alpha}, \mathbf{s}_{\eta,\alpha}).$$

Now let $f_\alpha : \mathcal{X}_\alpha \to \mathcal{W}_\alpha$ be the universal family over $S_\alpha$ and let $\tau_\alpha : S_\alpha \to C[n_\alpha]$ be so that $\mathcal{W}_\alpha = \tau_\alpha^* W[n_\alpha]$. To each $\eta \in \bar{\Omega}_\alpha$ we let $l_\eta$ be the integer defined in Definition 3.3. The assignment $\eta \mapsto l_\eta$ defines a function $\bar{\Omega}_\alpha \to [n_\alpha + 1]$. Because of [Li, Lemma 4.13], this assignment is one-to-one. Now let $K_\alpha \subset [n_\alpha + 1]$ be the image set of this assignment and let $\bar\tau_\alpha \colon S_\alpha \to \mathbf{A}^{n_\alpha+1}$ be the composition of $\tau_\alpha$ with the projection $C[n_\alpha] \to \mathbf{A}^{n_\alpha+1}$. Clearly, if $l \in [n_\alpha + 1] - K_\alpha$, then $\bar\tau_\alpha(S_\alpha) \cap \mathbf{H}_l = \emptyset$ and hence $\bar\tau_\alpha^*(L_l, s_l) \equiv (\mathbf{1}_{S_\alpha}, 1)^{15}$. Otherwise, $l = l_\eta$ for a unique $\eta \in \tilde\Omega_\alpha$ and then $\bar\tau_\alpha^*(L_l, s_l) \equiv (\mathbf{L}_{\eta,\alpha}, \mathbf{s}_{\eta,\alpha})$. Therefore, we have canonical isomorphisms

$$\otimes_{\eta \in \tilde\Omega_\alpha} (\mathbf{L}_{\eta,\alpha}, \mathbf{s}_{\eta,\alpha})|_{S_\alpha} \cong \otimes_{l \in K_\alpha} \bar\tau_\alpha^*(L_l, s_l) \cong \bar\tau_\alpha^*\big(\otimes_{l \in K_\alpha}(L_l, s_l)\big) \cong (\mathbf{L}_0, \mathbf{r}_0)|_{S_\alpha}.$$

---

[14]For the definition of effective arrows please see [Li, Section 1].

Because the above isomorphisms are canonical, they are compatible over $S_{\alpha\beta}$ and hence define an isomorphism of C-divisors as required by the Proposition. $\qquad\square$

We now derive the first version of the degeneration formula. We first recall the notion of localized top Chern class of a vector bundle with a section. Let $E$ be a rank $m$ vector bundle over a scheme $X$ and $s$ a section of $E$. The localized top Chern class of $(E, s)$ is the homomorphism

$$c_m(E, s) \colon A_* X \to A_{*-m} s^{-1}(0)$$

defined in [Ful] as follows: Let $Z$ be any variety in $X$. We take the normal cone $N_{s^{-1}(0) \cap Z/Z}$ to $s^{-1}(0) \cap Z$ in $Z$ and then define $c_m(E, s)([Z]) = 0_E^!(N_{s^{-1}(0) \cap Z/Z})$, where $0_E^!$ is the Gysin map of the zero section $0_E$ of $E|_{s^{-1}(0)}$. This defines a homomorphism of the group of cycles. This construction can be extended to the case where $X$ is an algebraic stack [Vis] with $A_* X$ understood to be the cycle group with rational coefficients.

We consider the moduli stack $\mathfrak{M}(\mathfrak{W}, \Gamma)$ and its Gromov-Witten invariants

$$\Psi_\Gamma^{W/C} \colon H_C^0(R^*\pi_*\mathbb{Q}_W)^{\times k} \times H^*(\mathfrak{M}_{g,k}) \longrightarrow H_2^{\mathrm{BM}}(C) \cong \mathbb{Q}.$$

Now let $\xi \in C$ be any closed point and let $H_C^0(R^*\pi_*\mathbb{Q}_W) \to H^*(W_\xi)$, denoted by $\alpha \mapsto \alpha(\xi)$, be induced by $W_\xi \to W$. We let $H_2^{\mathrm{BM}}(C) \to \mathbb{Q}$ be the Gysin homomorphism defined by intersecting with the divisor $\xi \in C$. Then we have a commutative diagram

(3.6)
$$
\begin{array}{ccc}
H_C^0(R^*\pi_*\mathbb{Q}_W)^{\times k} \times H^*(\mathfrak{M}_{g,k}) & \xrightarrow{\ \Psi_\Gamma^{W/C}\ } & H_2^{\mathrm{BM}}(C) \cong \mathbb{Q} \\
\downarrow & & \downarrow \\
H^*(W_\xi)^{\times k} \times H^*(\mathfrak{M}_{g,k}) & \xrightarrow{\ \Psi_\Gamma^{W_\xi}\ } & \mathbb{Q}.
\end{array}
$$

This was proved in [LT1] for $\xi \neq 0$ except that there we used the existence of a global vector bundle in defining the GW-invariants. This will be proved later in this paper again. We let $\Gamma = (g, k, b)$ be the triple as before and let $\alpha \in H^0(R^*\pi_*\mathbb{Q}_W)^{\times k}$ and $\beta \in H^*(\mathfrak{M}_{g,k})$ be any classes. As before, we let $ev_\xi \colon \mathfrak{M}(W_\xi, \Gamma) \to W_\xi^{\times k}$ be the evaluation morphism by the ordinary marked points of the stable morphisms, and let $\pi_{g,k} \colon \mathfrak{M}(\mathfrak{W}, \Gamma) \to \mathfrak{M}_{\mathfrak{g},k}$ be the forgetful map.

**Theorem 3.6.** *For any closed $\xi \in C$, we have*

$$\Psi_\Gamma^{W_\xi}(\alpha(\xi), \beta) = \mathbf{q}_{*0}\Big(\sum_{\eta \in \Omega}\big(ev_0^*(\alpha(0)) \cup \pi_{g,k}^*(\beta)\big)\big(c_1(\mathbf{L}_\eta, \mathbf{s}_\eta)[\mathfrak{M}(\mathfrak{W}, \Gamma)]^{\mathrm{virt}}\big)\Big).$$

---

[15]There is an exceptional case I should mention here. It is when there are $s \in S_\alpha$ so that $\bar\tau_\alpha(s) \in H_l$ while $f_\alpha^{-1}(\mathbf{D}_l) \cap \mathcal{X}_s = \emptyset$. Note that since $\mathcal{X}_s$ is connected, this is possible only when $l = 1$ or $n_\alpha + 1$. In either case, we agree that $\mathcal{X}_s$ decomposes into $C_1 \sqcup C_2$ with either $C_1 = \emptyset$ or $C_2 = \emptyset$, and the corresponding $\eta = (\Gamma_1, \emptyset, \emptyset)$ or $\eta = (\emptyset, \Gamma_2, \emptyset)$, which we agree is in $\tilde\Omega$. In this case we let $l_\eta$ be $1$ or $n_\alpha + 1$. With this agreement, this statement is true without exception.

*Proof.* First, by the commutativity of the diagram (3.6), we have $\Psi_\Gamma^{W_\xi}(\alpha(\xi), \beta) = \Psi_\Gamma^{W_0}(\alpha(0), \beta)$. In the later part of this paper, we will show that

$$\Psi_\Gamma^{W_0}(\alpha(0), \beta) = \mathbf{q}_{*0}\Big(\big(ev_0^*(\alpha(0)) \cup \pi_{g,k}^*(\beta)\big)\big(c_1(\mathbf{L}_0, \mathbf{r}_0)[\mathfrak{M}(\mathfrak{W}, \Gamma)]^{\mathrm{virt}}\big)\Big).$$

However, since $(\mathbf{L}_0, \mathbf{r}_0)$ is the tensor product of all $(\mathbf{L}_\eta, \mathbf{s}_\eta)$, the Chern class operations

$$c_1(\mathbf{L}_0, \mathbf{r}_0) = \sum_{\eta \in \Omega} c_1(\mathbf{L}_\eta, \mathbf{s}_\eta) : A_*\mathfrak{M}(\mathfrak{W}, \Gamma) \to A_*\mathfrak{M}(\mathfrak{W}_0, \Gamma).$$

The theorem then follows immediately.                                    ∎

### 3.2. Statement of the degeneration formula.

In this subsection, we will first construct the virtual moduli cycles of several moduli stacks related to the substack $\mathfrak{W}_0 \subset \mathfrak{W}$. After that, we will state the final version of the degeneration formula of the Gromov-Witten invariants of this paper. We will leave the proof of the key Lemmas to the next section.

Let $\eta \in \bar{\Omega}$ be any admissible triple. Associated to $\eta$ we have the substack $\mathfrak{M}(\mathfrak{Y}_1^{\mathrm{rel}} \sqcup \mathfrak{Y}_2^{\mathrm{rel}}, \eta) \subset \mathfrak{M}(\mathfrak{W}, \Gamma)$ that is the image stack of $\Phi_\eta$ in (3.1). Let $(\mathbf{L}_\eta, \mathbf{s}_\eta)$ be the C-divisor on $\mathfrak{M}(\mathfrak{W}, \Gamma)$ defined in the previous subsection. We define the substack

$$\mathfrak{M}(\mathfrak{W}_0, \eta) = \mathbf{s}_\eta^{-1}(0) \subset \mathfrak{M}(\mathfrak{W}, \Gamma).$$

Note that we have an increasing chain of closed substacks

$$\mathfrak{M}(\mathfrak{Y}_1^{\mathrm{rel}} \sqcup \mathfrak{Y}_2^{\mathrm{rel}}, \eta) \subset \mathfrak{M}(\mathfrak{W}_0, \eta) \subset \mathfrak{M}(\mathfrak{W}_0, \Gamma) \subset \mathfrak{M}(\mathfrak{W}, \Gamma),$$

where the first inclusion induces a homeomorphism on topological spaces. In the previous section, we have constructed the perfect-obstruction theories of $\mathfrak{M}(\mathfrak{W}, \Gamma)$ and of $\mathfrak{M}(\mathfrak{W}_0, \Gamma)$, and thus have constructed their virtual moduli cycles. In the first part of this subsection, we will show that the natural obstruction theories of $\mathfrak{M}(\mathfrak{Y}_1^{\mathrm{rel}} \sqcup \mathfrak{Y}_2^{\mathrm{rel}}, \eta)$ and of $\mathfrak{M}(\mathfrak{W}_0, \eta)$ are also perfect. Thus they have natural virtual moduli cycles.

We first investigate the obstruction theory of $\mathfrak{M}(\mathfrak{W}_0, \eta)$. The discussion is parallel to the obstruction theory of $\mathfrak{M}(\mathfrak{W}_0, \Gamma)$. We now present the details. Recall that for each $n$ we have the moduli $\mathfrak{M}(W[n], \Gamma)^{\mathrm{st}}$ of stable morphisms to $W[n]$ of topological type $\Gamma$ that are also stable as morphisms to the stack $\mathfrak{W}$. We let $\mathfrak{M}(W[n], \Gamma)^{\mathrm{st}} \to \mathfrak{M}(\mathfrak{W}, \Gamma)$ be the tautological morphism and let

$$\mathfrak{M}(W_0[n], \eta)^{\mathrm{st}} = \mathfrak{M}(W[n], \Gamma)^{\mathrm{st}} \times_{\mathfrak{M}(\mathfrak{W}, \Gamma)} \mathfrak{M}(\mathfrak{W}_0, \eta).$$

As in the case of $\mathfrak{M}(\mathfrak{W}, \Gamma)$, to study the obstruction theory of $\mathfrak{M}(\mathfrak{W}_0, \eta)$ it suffices to work out the obstruction theory of $\mathfrak{M}(W_0[n], \eta)^{\mathrm{st}}$, which we will do now. Of course for étale charts $S$ of $\mathfrak{M}(W[n], \Gamma)^{\mathrm{st}}$ we can define the notion of $\eta$-admissible as in Definition 3.3. Now let $S$ be an $\eta$-admissible chart of $\mathfrak{M}(W[n], \Gamma)^{\mathrm{st}}$ and let $S_\eta = S \times_{\mathfrak{M}(W[n],\Gamma)^{\mathrm{st}}} \mathfrak{M}(W_0[n], \eta)^{\mathrm{st}}$. As before, we let $f: \mathcal{X} \to W[n]$ be the universal family over $S$. We let $l$ be the integer so that $f^{-1}(\mathbf{D}_l)$ is where the $\eta$-decomposition of $f|_{S_\eta}$ takes place. We fix a covering of $f$ by charts $(\mathcal{U}_\alpha/\mathcal{V}_\alpha, f_\alpha)$ of the first and the second kinds, indexed by $\Lambda$. Following our convention, we denote by $\underline{\mathcal{U}}$ the étale covering $\{\mathcal{U}_\alpha\}_\Lambda$ of $\mathcal{X}$. We let $f_{\eta,\alpha}: \mathcal{U}_{\eta,\alpha} \to W[n]$ be the restriction of $f_\alpha$ to $\mathcal{U}_{\eta,\alpha} \triangleq \mathcal{U}_\alpha \times_S S_\eta$. We let $\mathcal{V}_{\eta,\alpha} = \mathcal{V}_\alpha \times_S S_\eta$ and let $\underline{\mathcal{V}_\eta}$ be the covering $\{\mathcal{V}_{\eta,\alpha}\}_\Lambda$ of $S_\eta$. We let $A = \Gamma(\mathcal{O}_S)$ and let $\mathbf{E}^\bullet \equiv \mathbf{E}(A)^\bullet$ be the complex constructed in section 2 associated to the covering $\Lambda$. Now let $L_\eta$ be the line bundle on $S_\eta$ that is the pull back of $\mathbf{L}_\eta$ via $S_\eta \to \mathfrak{M}(\mathfrak{W}, \Gamma)$. We form the ordinary Čech complex

$$(3.7) \qquad \mathbf{C}_\eta^\bullet \triangleq \mathbf{C}^\bullet(\underline{\mathcal{V}_\eta}, L_\eta)$$

of the invertible sheaf (line bundle) $L_\eta$ over $S_\eta$ associated to the covering $\underline{\mathcal{V}_\eta}$. Let $A_\eta = \Gamma(\mathcal{O}_{S_\eta})$. We now define a homomorphism of the complex:

$$(3.8) \qquad \delta : \mathbf{E}^\bullet \otimes_A A_\eta \Longrightarrow \mathbf{C}_\eta^{\bullet-1}$$

as follows: Let $(a, b) \in \mathbf{E}^1 \otimes_A A_\eta$ be any element. As argued in the proof of Proposition 1.16, $a$ defines, to each $\alpha \in \Lambda$, a flat extension $\tilde{\mathcal{U}}_{\eta,\alpha}/\tilde{\mathcal{V}}_{\eta,\alpha}$ of $\mathcal{U}_{\eta,\alpha}/\mathcal{V}_{\eta,\alpha}$ by the module $\Gamma(\mathcal{O}_{\mathcal{V}_{\eta,\alpha}})$ and pre-deformable extensions $\zeta_{\eta,\alpha} : \tilde{\mathcal{U}}_{\eta,\alpha} \to W[n]$ of $f_{\eta,\alpha} : \mathcal{U}_{\eta,\alpha} \to W[n]$ Let $\rho_{\eta,\alpha} : \mathcal{V}_{\eta,\alpha} \to \mathbf{A}^{n+1}$ be the tautological morphism (induced by $S \to \mathbf{A}^{n+1}$) and let $\tilde{\rho}_{\eta,\alpha} : \tilde{\mathcal{V}}_{\eta,\alpha} \to \mathbf{A}^{n+1}$ be the composite of $\zeta_{\eta,\alpha}$ and the projection $W[n]_\alpha \to \mathbf{A}^{n+1}$. Since $t_l \circ \rho_{\eta,\alpha} = 0$,

$$\xi_\alpha(a) \triangleq \mathbf{d}(t_l \circ \tilde{\rho}_{\eta,\alpha} - 0) \in \Gamma(\mathcal{V}_{\eta,\alpha}, \rho_{\eta,\alpha}^* N_{\mathbf{H}_l/\mathbf{A}^{n+1}}) \equiv \Gamma(\mathcal{V}_{\eta,\alpha}, L_\eta),$$

where $N_{\mathbf{H}_l/\mathbf{A}^{n+1}}$ is the normal bundle to $\mathbf{H}_l$ in $\mathbf{A}^{n+1}$. Here we have used the fact that $\rho_{\eta,\alpha}^* N_{\mathbf{H}_l/\mathbf{A}^{n+1}}$ is canonically isomorphic to $L_\eta$. On the other hand, for $b = \{b_\alpha\}$ with $b_\alpha \in \mathrm{Hom}_{\mathcal{U}_\alpha}(f^*\Omega_{W[n]}, A_\eta)^\dagger$, $b_\alpha(f^*(dt_l)) \in I$, where $\pi_n : W[n] \to \mathbf{A}^{n+1}$ is the projection. We define $\xi(b_\alpha) = b_\alpha(f^*(dt_l))$. This defines a homomorphism

$$(3.9) \qquad \delta_k : \mathbf{E}^k \otimes_A A_\eta \longrightarrow \mathbf{C}_\eta^{k-1} \triangleq \mathbf{C}^{k-1}(\underline{\mathcal{V}_\eta}, L_\eta)$$

via $\delta_1((a, b))_\alpha = \xi_\alpha(a) + \xi(b_\alpha)$ for $k = 1$. Note that for $k > 1$ elements in $\mathbf{E}^k$ are of the form $b = \{b_{\alpha_0 \cdots \alpha_k}\}$ and we can define $\delta_k(b_{\alpha_0 \cdots \alpha_k})$ similarly. This defines $\delta_k$ for $k > 1$. It is direct to check that the so defined map is a homomorphism of complexes. With this homomorphism of complexes, we can form a complex

$$(3.10) \qquad \mathbf{E}_\eta^k \triangleq \mathbf{C}_\eta^{k-2} \oplus \mathbf{E}^k \otimes_A A_\eta$$

with the induced differential. It is clear that the statement of Lemma 1.15 holds true to $\mathbf{E}_\eta^\bullet$. Namely, for sufficiently fine admissible covering $\Lambda$ of $f$, $\mathbf{E}_\eta^\bullet$ is a complex of flat $A_\eta$-modules. Further, to each $A_\eta$-module $I$ we have the following exact sequence of complexes

$$(3.11) \qquad 0 \Longrightarrow \mathbf{C}_\eta^{\bullet-2} \otimes_{A_\eta} I \Longrightarrow \mathbf{E}_\eta^\bullet \otimes_{A_\eta} I \Longrightarrow \mathbf{E}^\bullet \otimes_{A_\eta} I \Longrightarrow 0$$

which induces a long exact sequence of cohomologies

$$0 \longrightarrow h^1(\mathbf{E}_\eta^\bullet \otimes_{A_\eta} I) \longrightarrow h^1(\mathbf{E}^\bullet \otimes_A I) \longrightarrow I \otimes_{A_\eta} \mathcal{O}_{S_\eta}(L_\eta) \longrightarrow$$

$$(3.12) \qquad \longrightarrow h^2(\mathbf{E}_\eta^\bullet \otimes_{A_\eta} I) \longrightarrow h^2(\mathbf{E}^\bullet \otimes_A I) \longrightarrow 0.$$

In particular, $h^k(\mathbf{E}_\eta^\bullet \otimes_{A_\eta} I) = 0$ except $k = 1, 2$.

**Proposition 3.7.** *First, the functor of the first order deformations of $S_\eta$ is naturally isomorphic to the functor $\mathfrak{h}^1(\mathbf{E}_\eta^\bullet)$. Secondly, there is a natural obstruction theory to deformation of the families of pre-deformable morphisms $f|_{S_\eta} : \mathcal{X}|_{S_\eta} \to W[n] \times_{\mathbf{A}^{n+1}} \mathbf{H}_l$ taking values in $\mathfrak{h}^2(\mathbf{E}_\eta^\bullet)$. Finally, such obstruction theory is perfect.*

*Proof.* We will omit the proof here because it is parallel to the treatment of the obstruction theory of $\mathfrak{M}(\mathfrak{W}_0, \Gamma)$. $\qquad \square$

**Proposition 3.8.** *The perfect obstruction theory constructed in Proposition 3.7 induces a perfect obstruction theory of $\mathfrak{M}(\mathfrak{W}_0, \eta)$, which in turn defines a natural virtual moduli cycle $[\mathfrak{M}(\mathfrak{W}_0, \eta)]^{\mathrm{virt}}$.*

*Proof.* The proof is parallel to the construction of the perfect obstruction theory of $\mathfrak{M}(\mathfrak{W}, \Gamma)$ and of the virtual moduli cycle $[\mathfrak{M}(\mathfrak{W}, \Gamma)]^{\mathrm{virt}}$. We shall omit the details here. ∎

We next work out the obstruction theory of $\mathfrak{M}(\mathfrak{Y}_1^{\mathrm{rel}} \sqcup \mathfrak{Y}_2^{\mathrm{rel}}, \eta)$. We let $S = \operatorname{Spec} A \to \mathfrak{M}(W[n], \Gamma)^{\mathrm{st}}$ be an $\eta$-admissible chart and let $f : \mathcal{X} \to W[n]$ be the universal family over $S$. We then pick a sufficiently fine covering $(\mathcal{U}_\alpha / \mathcal{V}_\alpha, f_\alpha)$ indexed by $\Lambda$. Following section 1.2, we pick the complex $[F^0 \xrightarrow{d} F^1]$ in (1.5) and pick extensions $\zeta_\alpha$ as defined by (1.19). Based on these data we can form the complex $\mathbf{E}^\bullet$ as in (1.20) so that its cohomology is part of the obstruction theory of $S$. We let $S_0 = \operatorname{Spec} A_0 = S \times_{\mathfrak{M}(\mathfrak{W}, \Gamma)} \mathfrak{M}(\mathfrak{Y}_1^{\mathrm{rel}} \sqcup \mathfrak{Y}_2^{\mathrm{rel}}, \eta)$ and let $f_0 : \mathcal{X}_0 \to W[n]$ be the restriction of $f$ to fibers over $S_0$. We let $l \in [n+1]$ be the integer associated to $\eta$ defined in Definition 3.3. By definition, the family $f_0$ can be decomposed into two families of relative stable morphisms of types $\Gamma_1$ and $\Gamma_2$ respectively along a multi-section $\Sigma \subset \mathcal{X}_{0,\mathrm{node}}$ over $S_0$. We let $\mathcal{U}_{0,\alpha} = \mathcal{U}_\alpha \times_S S_0$, $\mathcal{V}_{0,\alpha} = \mathcal{V}_\alpha \times_S S_0$ and let $f_{0,\alpha} = f_\alpha|_{\mathcal{U}_{0,\alpha}}$. Note that $\Sigma$ is étale over $S_0$. Now let $\hat{\mathcal{X}}_0$ be the formal completion of $\mathcal{X}_0$ along $\Sigma \subset \mathcal{X}_0$. Then since $\Sigma \to S$ is finite and étale, and since $\Sigma \subset \mathcal{X}_0$ is a multiple section of the nodal points of the fibers of $\mathcal{X}_0/S_0$, the extension sheaf $\mathcal{E}xt^1_{\hat{\mathcal{X}}_0}(\Omega_{\hat{\mathcal{X}}_0/S_0}, \mathcal{O}_{\hat{\mathcal{X}}_0})$ is an invertible sheaf of $\mathcal{O}_\Sigma$-modules. We denote this sheaf by $\mathcal{M}_\Sigma$. Then we have a natural homomorphism

$$\operatorname{Ext}^1_{\mathcal{X}_0}(\Omega_{\mathcal{X}_0/S_0}, \mathcal{O}_{\mathcal{X}_0}) \longrightarrow \rho_*\big(\mathcal{E}xt^1_{\mathcal{X}_0}(\Omega_{\mathcal{X}_0/S_0}, \mathcal{O}_{\mathcal{X}_0})\big) \equiv \rho_*\mathcal{M}_\Sigma,$$

where $\rho : \mathcal{X}_0 \to S_0$ is the projection, which defines a canonical homomorphism

$$(3.13) \qquad F^1 \otimes_A A_0 \longrightarrow \operatorname{Ext}^1_{\mathcal{X}_0}(\Omega_{\mathcal{X}_0/S_0}, \mathcal{O}_{\mathcal{X}_0}) \longrightarrow \bar{\rho}_*\mathcal{M}_\Sigma.$$

Here $\bar{\rho} : \Sigma \to S_0$ is the tautological projection induced by $\rho$. Clearly, the composite of (3.13) is surjective. We let $F^0_{\eta,0} = F^0 \otimes_A A_0$ and let $F^1_{\eta,0}$ be the kernel of (3.13). Clearly, $F^1_{\eta,0}$ is a free $A_0$-module and $F^0 \otimes_A A_0 \to F^1 \otimes_A A_0$ factor through $F^0_{\eta,0} \to F^1_{\eta,0}$.

We now construct the complex that will give the obstruction theory of $S_0$ (which is a chart of $\mathfrak{M}(\mathfrak{Y}_1^{\mathrm{rel}} \sqcup \mathfrak{Y}_2^{\mathrm{rel}}, \eta)$). We let $\mathbf{E}^\bullet$ be the complex associated to the family $f$ and the covering $\Lambda$ mentioned before. We let $\mathbf{E}^k_{\eta,0} \triangleq \mathbf{E}^k \otimes_A A_0$ for $k \neq 1$ and let $\mathbf{E}^1_{\eta,0}$ be the kernel of the composite $\mathbf{E}^1 \otimes_A A_0 \xrightarrow{\mathrm{pr}} F^1 \otimes_A A_0 \to \bar{\rho}_*\mathcal{M}_\Sigma$. Clearly, the differentials in $\mathbf{E}^\bullet$ induce differentials in $\mathbf{E}^\bullet_{\eta,0}$.

**Proposition 3.9.** *The Lemma 1.15 holds true for the complex $\mathbf{E}^\bullet_{\eta,0}$. The Proposition 3.7 holds true for the charts $S_0$ with $\mathbf{E}^\bullet_\eta$ replaced by $\mathbf{E}^\bullet_{\eta,0}$. The Proposition 3.8 holds true for the moduli stack $\mathfrak{M}(\mathfrak{Y}_1^{rel} \sqcup \mathfrak{Y}_2^{rel}, \eta)$.*

*Proof.* We shall sketch the construction of the obstruction classes. The remainder part of the proof is similar to that of the Propositions 3.7 and 3.8 and will be omitted. Here we will follow closely the convention introduced in the proof of Proposition 1.18. Let $\xi = (B, I, \varphi)$ be an object in $\mathfrak{Tri}_{S_0}$. Let $f_T : \mathcal{X}_T \to W[n]$ be the pull back of $f$ via $T \to S_0$, where $T = \operatorname{Spec} B/I$, and let $\Sigma_T = \mathcal{X}_T \times_{\mathcal{X}_0} \Sigma$. By the definition of the subscheme $S_0 \subset S$, the formal completion of $\mathcal{X}_T$ along $\Sigma_T$ is isomorphic to $\operatorname{Spec} \Bbbk[\![z_1, z_2]\!]/(z_1 z_2) \times \Sigma_T$, at least after shrinking $S$ if necessary. Now let $\tilde{T} = \operatorname{Spec} B$. Then by the deformation theory of nodal curves, we can find a flat extension $\mathcal{X}_{\tilde{T}}/\tilde{T}$ of $\mathcal{X}_T/T$ so that the formal completion of $\mathcal{X}_{\tilde{T}}$ along $\Sigma_T$ is isomorphic to $\operatorname{Spec} \Bbbk[\![z_1, z_2]\!]/(z_1 z_2) \times \Sigma_{\tilde{T}}$, where $\Sigma_{\tilde{T}}$ is étale over $\tilde{T}$ so that $\Sigma_{\tilde{T}} \times_{\tilde{T}} T \equiv \Sigma_T$. In other words, the multiple section $\Sigma_T \subset \mathcal{X}$ extends to a multiple section $\Sigma_{\tilde{T}} \subset \mathcal{X}_{\tilde{T}}$ and the extended family $\mathcal{X}_{\tilde{T}}$ can be decomposed along $\Sigma_{\tilde{T}}$. Once we have chosen such an extension, to each $\alpha \in \Lambda$ we can pick a pre-deformable extension $\tilde{h}_\alpha : \tilde{\mathcal{U}}_{0,\alpha} \to W[n]$ of $f_{T,\alpha}$ so that the composite $\tilde{\mathcal{U}}_{0,\alpha} \to W[n] \to \mathbf{A}^{n+1}$ factor through the $l$-th coordinate hyperplane $\mathbf{H}_l \subset \mathbf{A}^{n+1}$. We let $b_{\alpha\beta}$ be as defined in (1.22). Then $b = \{b_{\alpha\beta}\} \in \mathbf{E}^2 \otimes_A I \equiv \mathbf{E}^2_{\eta,0} \otimes_{A_0} I$. We define

$$\mathfrak{ob}(\xi) = [b] \in h^2(\mathbf{E}^\bullet_{\eta,0} \otimes_{A_0} I).$$

It is direct to check that this defines an obstruction class. It is routine to check that such choices of obstruction classes satisfies the base change property. $\qquad\blacksquare$

Applying the basic construction of virtual moduli cycles formulated in the previous section, we obtain cycles $[\mathfrak{M}(\mathfrak{W}_0, \eta)]^{\mathrm{virt}}$ and $[\mathfrak{M}(\mathfrak{Y}^{\mathrm{rel}}_1 \sqcup \mathfrak{Y}^{\mathrm{rel}}_2, \eta)]^{\mathrm{virt}}$.

We are now ready to state the degeneration formula of the Gromov-Witten invariants of the family $W$. Let $\eta = (\Gamma_1, \Gamma_2, I)$ be an admissible triple in $\bar{\Omega}$. Let $\mu_i$ be the weight of the $i$-th root of $\Gamma_1$ and $\Gamma_2$, which are the same. We define the multiplicity of $\eta$ to be $\mathbf{m}(\eta) = \prod_{i=1}^r \mu_i$.

**Lemma 3.10.** *We have the identity*

$$c_1(\mathbf{L}_0, \mathbf{s}_0)[\mathfrak{M}(\mathfrak{W}, \Gamma)]^{\mathrm{virt}} = [\mathfrak{M}(\mathfrak{W}_0, \Gamma)]^{\mathrm{virt}} \in A_* \mathfrak{M}(\mathfrak{W}_0, \Gamma).$$

**Lemma 3.11.** *We have the identify*

$$[\mathfrak{M}(\mathfrak{W}_0, \eta)]^{\mathrm{virt}} = c_1(\mathbf{L}_\eta, \mathbf{s}_\eta)[\mathfrak{M}(\mathfrak{W}, \Gamma)]^{\mathrm{virt}} \in A_* \mathfrak{M}(\mathfrak{W}_0, \eta).$$

**Lemma 3.12.** *Under the natural isomorphism*

$$A_* \mathfrak{M}(\mathfrak{Y}^{rel}_1 \sqcup \mathfrak{Y}^{rel}_2, \eta) \cong A_* \mathfrak{M}(\mathfrak{W}_0, \eta)$$

*induced by the homeomorphism* $\mathfrak{M}(\mathfrak{Y}^{rel}_1 \sqcup \mathfrak{Y}^{rel}_2, \eta) \simeq \mathfrak{M}(\mathfrak{W}_0, \eta)$,

$$\mathbf{m}(\eta)[\mathfrak{M}(\mathfrak{Y}^{rel}_1 \sqcup \mathfrak{Y}^{rel}_2, \eta)]^{\mathrm{virt}} = [\mathfrak{M}(\mathfrak{W}_0, \eta)]^{\mathrm{virt}} \in A_* \mathfrak{M}(\mathfrak{W}_0, \eta).$$

Using the identity of C-divisors $(\mathbf{L}_0, \mathbf{s}_0) = \otimes_{\eta \in \bar{\Omega}}(\mathbf{L}_\eta, \mathbf{s}_\eta)$, we have

**Corollary 3.13.**

$$[\mathfrak{M}(\mathfrak{W}_0, \Gamma)]^{\mathrm{virt}} = \sum_{\eta \in \bar{\Omega}} \mathbf{m}(\eta)[\mathfrak{M}(\mathfrak{Y}^{rel}_1 \sqcup \mathfrak{Y}^{rel}_2, \eta)]^{\mathrm{virt}}.$$

We now state how virtual moduli cycles $[\mathfrak{M}(\mathfrak{Y}^{\mathrm{rel}}_1 \sqcup \mathfrak{Y}^{\mathrm{rel}}_2, \eta)]^{\mathrm{virt}}$ is related to $[\mathfrak{M}(\mathfrak{Y}^{\mathrm{rel}}_i, \Gamma_i)]^{\mathrm{virt}}$. As mentioned in [Li], we have the natural evaluation morphism $\mathfrak{M}(\mathfrak{Y}^{\mathrm{rel}}_i, \Gamma_i) \to D^r$ and the induced Cartesian diagram

$$(3.14) \quad \begin{array}{ccc} \mathfrak{M}(\mathfrak{Y}^{\mathrm{rel}}_1, \Gamma_1) \times_{D^r} \mathfrak{M}(\mathfrak{Y}^{\mathrm{rel}}_2, \Gamma_2) & \longrightarrow & \mathfrak{M}(\mathfrak{Y}^{\mathrm{rel}}_1, \Gamma_1) \times \mathfrak{M}(\mathfrak{Y}^{\mathrm{rel}}_2, \Gamma_2) \\ \downarrow & & \downarrow \\ D^r & \xrightarrow{\Delta} & D^r \times D^r \end{array}$$

Here the arrow $\Delta$ is the diagonal morphism. Let $\Phi_\eta$ be the finite étale morphism in (3.2), which has pure degree $|\operatorname{Eq}(\eta)|$ (see [Li, Section 4]).

**Lemma 3.14.** *We have the identity*

$$\frac{1}{|\operatorname{Eq}(\eta)|} \Phi_{\eta*} \Delta^! \big( [\mathfrak{M}(\mathfrak{Y}_1^{rel}, \Gamma_1)]^{\mathrm{virt}} \times [\mathfrak{M}(\mathfrak{Y}_2^{rel}, \Gamma_2)]^{\mathrm{virt}} \big) = [\mathfrak{M}(\mathfrak{Y}_1^{rel} \sqcup \mathfrak{Y}_2^{rel}, \eta)]^{\mathrm{virt}}.$$

We will prove these Lemmas in the next section.

The main degeneration formula of the Gromov-Witten invariants of $W$ follows immediately from these Lemmas and the first version of the degeneration formula proved in the previous subsection.

**Theorem 3.15.** *Let the notation be as before. Then as elements in $A_*\mathfrak{M}(\mathfrak{W}_0, \Gamma)$,*

$$[\mathfrak{M}(\mathfrak{W}_0, \Gamma)]^{\mathrm{virt}} = \sum_{\eta \in \bar{\Omega}} \frac{\mathbf{m}(\eta)}{|\operatorname{Eq}(\eta)|} \Phi_{\eta*} \Delta^! \big( [\mathfrak{M}(\mathfrak{Y}_1^{rel}, \Gamma_1)]^{\mathrm{virt}} \times [\mathfrak{M}(\mathfrak{Y}_2^{rel}, \Gamma_2)]^{\mathrm{virt}} \big).$$

Finally, we state the numerical corollary of this theorem. Let $\jmath_i: Y_i \to W$ be the inclusion and let

$$\jmath_i^*: H_C^0(R^*\pi_*\mathbb{Q}_W)^{\times k} \to H^*(Y_i, \mathbb{Q})^{\times k}$$

be the induced pull back homomorphism. Now let $\eta = (\Gamma_1, \Gamma_2, I) \in \bar{\Omega}$ be any admissible triple. For $i = 1$ or 2, we let $\mathfrak{M}_{\Gamma_i^\circ}$ be the moduli space of stable curves of topological type $\Gamma_i^\circ$ (See the definition before (2.9)). It is naturally a Deligne-Mumford stack. Further, we have a natural local immersion of stacks

$$\phi_\eta: \mathfrak{M}_{\Gamma_1^\circ} \times \mathfrak{M}_{\Gamma_2^\circ} \longrightarrow \mathfrak{M}_{g,k}$$

that associates to any pair of curves $(C_1, C_2) \in \mathfrak{M}_{\Gamma_1^\circ} \times \mathfrak{M}_{\Gamma_2^\circ}$ the gluing $C_1 \sqcup C_2$ by identifying the $i$-th distinguished marked point of $C_1$ with the $i$-th distinguished marked point of $C_2$ for all $i$. Now let $\beta \in H^*(\mathfrak{M}_{g,k})$. We assume $\beta$ has the following Kunneth type decomposition

$$\phi_\eta^*(\beta) = \sum_{j \in K_\eta} \beta_{\eta,1,j} \boxtimes \beta_{\eta,2,j}, \quad \beta_{i,j} \in H^*(\mathfrak{M}_{\Gamma_i^\circ}).$$

**Corollary 3.16.** *Let $W/C$ be the family and let $\Gamma = (g, b, k)$ be as before. Then for any closed point $\xi \neq 0 \in C$, $\alpha \in H_C^0(R^*\pi_*\mathbb{Q}_W)^{\times k}$ and $\beta \in H^*(\mathfrak{M}_{g,n})$ as before,*

$$\Phi_\Gamma^{W_\xi}(\alpha(\xi), \beta) = \sum_{\eta \in \bar{\Omega}} \frac{\mathbf{m}(\eta)}{|\operatorname{Eq}(\eta)|} \sum_{j \in K_\eta} \big[ \Psi_{\Gamma_1}^{Y_1^{rel}}(\jmath_1^*\alpha, \beta_{\eta,1,j}) \bullet \Psi_{\Gamma_2}^{Y_2^{rel}}(\jmath_2^*\alpha, \beta_{\eta,2,j}) \big]_0.$$

Here $\bullet$ is the intersection of the homology groups

$$H_*(D^r) \times H_*(D^r) \xrightarrow{\cap} H_*(D^r)$$

and $[\gamma]_0$ is the degree of the degree 0 part of the homology class $\gamma \in H_*(D^r)$.

## 4. Proof of the main theorem

The goal of this section is to prove Lemma 3.12-3.14. In essence, the proofs of these Lemmas (except Lemma 3.12) rely on the comparison of the virtual moduli cycles of stacks with the virtual moduli cycles of their substacks. This is precisely the situation studied in [LT2, Lemma 3.4], except that there we used the existence of certain locally free sheaves to construct the virtual moduli cycles. To prove Lemma 3.12, we need to study the situation more general than the one studied. In the first subsection we will revise [LT2, Lemma 3.4] to cover all the situations we need.

4.1. **Comparison of the virtual moduli cycles.** Let $\mathbf{M} \to \mathbf{N}$ be a representable morphism of stacks. In this subsection, we assume $\mathbf{M}$ is a DM-stack having a perfect-obstruction theory with the associated obstruction sheaves $\mathcal{O}b_{\mathbf{M}}$. For $\mathbf{N}$ we need to consider two possibilities: One is when $\mathbf{N}$ is a DM-stack having a perfect-obstruction theory with the obstruction sheaf $\mathcal{O}b_{\mathbf{N}}$. The other is when $\mathbf{N}$ is a closed substack of a smooth Artin stack $\mathbf{Q}$ defined by the vanishing of a section $s$ of a vector bundle $F$ on $\mathbf{Q}$.

Let $\mathcal{S} = \{S_\alpha\}_\Lambda$ be an atlas of $\mathbf{N}$ and let $\{(\mathcal{F}_\alpha^\bullet, \mathfrak{ob}_{S_\alpha})\}_\Lambda$ be the data associated to the perfect-obstruction theory of $\mathbf{N}$. In case $\mathbf{N}$ is a DM-stack this is specified in the Definition 2.1. In case $\mathbf{N}$ is an Artin stack, we assume that there is an atlas $\bar{\mathcal{S}} = \{\bar{S}_\alpha\}_\Lambda$ of $\mathbf{Q}$ so that $S_\alpha = \bar{S}_\alpha \times_{\mathbf{Q}} \mathbf{N}$. Then we simply take $\mathcal{F}_\alpha^\bullet = [\mathcal{F}_\alpha^1 \to \mathcal{F}_\alpha^2]$ to be $\mathcal{F}_\alpha^1 = \mathcal{O}_{S_\alpha}(T\bar{S}_\alpha)$, $\mathcal{F}_\alpha^2 = \mathcal{O}_{S_\alpha}(F|_{S_\alpha})$ and the arrow $\mathcal{F}_\alpha^1 \to \mathcal{F}_\alpha^2$ to be the one induced by the differential of the section $s \in H^0(F)$. The obstruction assignment $\mathfrak{ob}_{S_\alpha}$ taking values in $\mathcal{O}b_{S_\alpha} = h^2(\mathcal{F}_\alpha^\bullet)$ is the obvious one induced by the defining equation $s$. Since $\mathbf{M} \to \mathbf{N}$ is representable, $S_\alpha \times_{\mathbf{N}} \mathbf{M}$ is a scheme. For each $\alpha \in \Lambda$ we pick an affine étale universal open $R_\alpha \to S_\alpha \times_{\mathbf{N}} \mathbf{M}$. Without loss of generality, we can assume $\{R_\alpha\}_\Lambda$ is a covering of $\mathbf{M}$ in the sense that the image $\mathbf{R}_\alpha \triangleq \rho_\alpha(R_\alpha)$ of the tautological $p_\alpha : R_\alpha \to \mathbf{M}$ is an open substack and the collection $\{\mathbf{R}_\alpha\}_\Lambda$ forms an open covering of $\mathbf{M}$. Note that when $\mathbf{N}$ is a DM-stack, $R_\alpha \to \mathbf{M}$ are étale and then $\{R_\alpha\}_\Lambda$ forms an étale cover of $\mathbf{M}$. In case $\mathbf{N}$ is an Artin stack, then $R_\alpha \to \mathbf{M}$ are smooth morphisms. In this case we shall view $\{R_\alpha\}_\Lambda$ as an atlas in the smooth site[16] of $\mathbf{M}$.

We let $\{(\mathcal{E}_\alpha^\bullet, \mathfrak{ob}_{R_\alpha})\}_\Lambda$ be the date given by the perfect obstruction theory of $\mathbf{M}$ associated to the covering $\{R_\alpha\}_\Lambda$. Namely, each $(\mathcal{E}_\alpha^\bullet, \mathfrak{ob}_{R_\alpha})$ is a perfect obstruction theory of $R_\alpha$, the sheaves $\mathrm{Coker}\{\mathcal{E}_\alpha^1 \to \mathcal{E}_\alpha^2\}$ descends to the sheaf $\mathcal{O}b_{\mathbf{M}}$ of $\mathcal{O}_{\mathbf{M}}$-modules and the obstruction assignments $\{\mathfrak{ob}_{R_\alpha}\}$ are compatible over all $R_{\alpha\beta} = R_\alpha \times_{\mathbf{M}} R_\beta$.

We next assume $R_\alpha \to S_\alpha$ admits a perfect relative obstruction theory given by $(\mathcal{L}_\alpha^\bullet, \mathfrak{ob}_{R_\alpha/S_\alpha})$ as defined in Definition 1.11. We say $\mathbf{M} \to \mathbf{N}$ admits a perfect relative obstruction theory if we can choose $\{(\mathcal{L}_\alpha^\bullet, \mathfrak{ob}_{R_\alpha/S_\alpha})\}_\Lambda$ so that the relative obstruction sheaves $\mathcal{O}b_{R_\alpha/S_\alpha} \triangleq \mathrm{Coker}\{\mathcal{L}_\alpha^1 \to \mathcal{L}_\alpha^2\}$ descends to a global sheaf on $\mathbf{M}$ and the obstruction assignments $\mathfrak{ob}_{R_\alpha/S_\alpha}$ are compatible on the overlaps $R_{\alpha\beta}$.

**Definition 4.1.** *The perfect (relative) obstruction theories $\{\mathcal{E}_\alpha^\bullet, \mathfrak{ob}_{R_\alpha}\}$, $\{\mathcal{F}_\alpha^\bullet, \mathfrak{ob}_{S_\alpha}\}$ and $\{\mathcal{L}_\alpha^\bullet, \mathfrak{ob}_{R_\alpha/S_\alpha}\}$ are said to be compatible if to each $\alpha \in \Lambda$ there is an exact triangle of complexes*

$$\tag{4.1} \Longrightarrow \mathcal{L}_\alpha^\bullet \Longrightarrow \mathcal{E}_\alpha^\bullet \Longrightarrow \mathcal{F}_\alpha^\bullet \otimes_{\mathcal{O}_{S_\alpha}} \mathcal{O}_{R_\alpha} \Longrightarrow \mathcal{L}_\alpha^{\bullet+1} \Longrightarrow$$

*which induces a long exact sequence of cohomologies*

$$\tag{4.2} \to h^i(\mathcal{L}_\alpha^\bullet \otimes \mathcal{I}) \xrightarrow{\tau_{1,i}} h^i(\mathcal{E}_\alpha^\bullet \otimes \mathcal{I}) \xrightarrow{\tau_{2,i}} h^i(\mathcal{F}_\alpha^\bullet \otimes_{\mathcal{O}_{S_\alpha}} \mathcal{I}) \xrightarrow{\tau_{0,i+1}} h^{i+1}(\mathcal{L}_\alpha^\bullet \otimes \mathcal{I}) \to$$

*for any sheaf of $\mathcal{O}_{R_\alpha}$-modules $\mathcal{I}$ that satisfies the following properties:*
*(1). The first part of the above exact sequence is identical to the exact sequence*

$$0 \longrightarrow \mathfrak{Def}_{R_\alpha/S_\alpha}^1(\mathcal{I}) \longrightarrow \mathfrak{Def}_{R_\alpha}^1(\mathcal{I}) \longrightarrow \mathfrak{Def}_{S_\alpha}^1(\mathcal{I}) \longrightarrow$$

*under the canonical isomorphisms $\mathfrak{Def}^1 = h^1(\cdot \otimes \mathcal{I})$ given by the definition of the perfect obstruction theories.*
*(2). Let $\xi$ be any object in $\mathfrak{Tri}_{R_\alpha}$, which is also an object in $\mathfrak{Tri}_{S_\alpha}$. Then underr*

---

[16]Namely, the open covering are univeral open smooth morphisms.

the arrow $\tau_{2,2}$ in (4.2) we have $\mathfrak{ob}_{S_\alpha}(\xi) = \tau_{2,2}\big(\mathfrak{ob}_{R_\alpha}(\xi)\big)$.

(3). Let $\xi = (B, I, \varphi_0)$ be an object in $\mathfrak{Tri}_{R_\alpha/S_\alpha}$, which is also an object in $\mathfrak{Tri}_{R_\alpha}$. Then $\mathfrak{ob}_{R_\alpha}(\xi) = \tau_{2,1}\big(\mathfrak{ob}_{R_\alpha/S_\alpha}(\xi)\big)$. Further, suppose $\mathfrak{ob}_{R_\alpha}(\xi) = 0$ and $\varphi_0 \colon \operatorname{Spec} B/I \to R_\alpha$ extends to $\varphi \colon \operatorname{Spec} B \to R_\alpha$. Let $e \in h^1(\mathcal{F}_\alpha^\bullet \otimes_{B/I} I)$ be the difference of the tautological $\operatorname{Spec} B \to \operatorname{Spec} S_\alpha$ and the composite of $\operatorname{Spec} B \xrightarrow{\varphi} \operatorname{Spec} R_\alpha$ and $\operatorname{Spec} R_\alpha \xrightarrow{pr} \operatorname{Spec} S_\alpha$. Then $\mathfrak{ob}_{R_\alpha/S_\alpha}(\xi) = \tau_{0,2}(e)$.

Finally, the collection of exact sequences $h^2(\mathcal{L}_\alpha^\bullet) \to h^2(\mathcal{E}_\alpha^\bullet) \to h^2(\mathcal{F}_\alpha^\bullet) \longrightarrow 0$ descends to an exact sequence of sheaves

$$\mathcal{O}b_{\mathbf{M}/\mathbf{N}} \longrightarrow \mathcal{O}b_{\mathbf{M}} \longrightarrow \mathcal{O}b_{\mathbf{N}} \otimes_{\mathcal{O}_{\mathbf{N}}} \mathcal{O}_{\mathbf{M}} \longrightarrow 0.$$

The goal of this subsection is to show that with the data given in Definition 4.1, we can construct a class $[\mathbf{M}, \mathbf{N}]^{\mathrm{virt}} \in A_*\mathbf{M}$, called the relative virtual moduli cycle. We will then show that it is equal to $[\mathbf{M}]^{\mathrm{virt}}$ in $A_*\mathbf{M}$. This allows us to give a different interpretation of $[\mathbf{M}]^{\mathrm{virt}}$, useful in the proof of the key lemmas in the previous section.

We first study the local situation. We let $p \in \mathbf{M}$ be any closed point, $q \in \mathbf{N}$ be the image of $p$, $\bar{p} \in R_\alpha$ be a lift of $p \in \mathbf{M}$ and $\bar{q} \in S_\alpha$ be the image of $\bar{p}$. We let $T_1 = h^1(\mathcal{E}_\alpha^\bullet \otimes \Bbbk_{\bar{p}})$, $T_2 = h^1(\mathcal{F}_\alpha^\bullet \otimes \Bbbk_{\bar{q}})$ and let $T_{1/2} = h^1(\mathcal{L}_\alpha^\bullet \otimes \Bbbk_{\bar{p}})$. Similarly, we let $O_1 = h^2(\mathcal{E}_\alpha^\bullet \otimes \Bbbk_{\bar{p}})$, $O_2 = h^2(\mathcal{F}_\alpha^\bullet \otimes \Bbbk_{\bar{q}})$ and let $O_{1/2} = h^2(\mathcal{L}_\alpha^\bullet \otimes \Bbbk_{\bar{p}})$, all implicitly depending on the lift $\bar{p}$. Note that they fit into the exact sequence

$$(4.3) \qquad 0 \longrightarrow T_{1/2} \longrightarrow T_1 \longrightarrow T_2 \xrightarrow{\delta} O_{1/2} \longrightarrow O_1 \longrightarrow O_2 \longrightarrow 0$$

induced by (4.1). We now let $T = T_1 \oplus T_2/\ker(\delta)$ and pick a surjective homomorphism $T \to T_2$ that is an extension of $T_1 \to T_2$. Next we let $O = O_1 \oplus \operatorname{Im}(\delta)$. We then pick an injective homomorphism $\eta : O_{1/2} \to O$ so that the composite $O_{1/2} \to O \to O_{1/2}$ is the identity map. This way we have two exact sequences

$$(4.4) \qquad 0 \to T_{1/2} \to T \to T_2 \to 0 \quad \text{and} \quad 0 \to O_{1/2} \to O \to O_2 \to 0.$$

Further, with various homomorphisms chosen and the isomorphism $T_2/\ker(\delta) \cong \operatorname{Im}(\delta)$, we have the following induced exact sequence

$$0 \longrightarrow T_1 \longrightarrow T \longrightarrow O \longrightarrow O_1 \longrightarrow 0.$$

In the following for any vector space $V$ we denote by $\Bbbk[\![V]\!]$ the ring of formal power series $\lim_{\oplus} \oplus_{i=0}^n S^i(V)$. We now let

$$(4.5) \qquad f \in \Bbbk[\![T^\vee]\!] \otimes O, \ g \in \Bbbk[\![T_2^\vee]\!] \otimes O_2 \text{ and } h \in \Bbbk[\![T^\vee]\!]/(g) \otimes O_{1/2}$$

be the Kuranishi maps of the (relative) obstruction theories of $R_\alpha$, of $S_\alpha$ and of $R_\alpha/S_\alpha$ at $\bar{p}$ (or $\bar{q}$) respectively (see [LT2, Lemma 3.10]). By abuse of notation, we denote by $(g)$ the idea generated by the components of $g$ in $\Bbbk[\![T_2^\vee]\!]$ and in $\Bbbk[\![T^\vee]\!]$ via the inclusion $\Bbbk[\![T_2^\vee]\!] \to \Bbbk[\![T^\vee]\!]$.

**Lemma 4.2.** We can choose the Kuranishi maps $f$, $g$ and $h$ so that
(1) $\varphi_1(f_p) = \tau(g_p)$ under the naturally induced maps $\varphi_1 \colon \Bbbk[\![T^\vee]\!] \otimes O \to \Bbbk[\![T^\vee]\!] \otimes O_2$ and $\tau \colon \Bbbk[\![T_2^\vee]\!] \otimes O_2 \to \Bbbk[\![T^\vee]\!] \otimes O_2$;
(2) The differential $dh_p(0) \colon T \to O_{1/2}$ is identical to the composite $T \to T_2/\ker(\delta) \to O_{1/2}$ induced by $\delta$ in (4.3).
(3) Let $\tilde{h}$ and $\tilde{f}$ in $\Bbbk[\![T_{1/2}^\vee]\!]/(g) \otimes O$ be the images of $h$ and $f$ under the obvious maps induced by the arrows mentioned before, then $\tilde{f} = \tilde{h}$.

*Proof.* The proof follows from the construction of Kuranishi maps, as was demonstrated in [LT2, Lemma 3.10]. ∎

We now let $X = \operatorname{Spec} \Bbbk[\![T^\vee]\!]$, $V_{1/2} = O_{1/2} \times X$, $V = O \times X$ and let $V_2 = O_2 \times X$, all viewed as vector bundles (or their total spaces) over $X$. We let $\rho_1 : V_{1/2} \to V$ and $\rho_2 : V \to V_2$ be the vector bundle homomorphisms induced by arrows in (4.4). Then $V_2$ is the quotient vector bundle of $V$ by $V_{1/2}$. We consider the subscheme $\Gamma \subset V_1$ that is the graph of $f$, consider $Y = V_{1/2} \times_X \Gamma \subset V_{1/2} \times_X V$ and consider

$$\Theta_1 = \tilde{V}_{1/2} \times_{V_{1/2} \times_X V} Y \subset Y, \quad \Theta_2 = \tilde{V}_1 \times_{V_{1/2} \times_X V} Y \subset Y.$$

Here $\tilde{V}_{1/2}$ is the image scheme of the immersion $1_{V_{1/2}} \times \rho_1 : V_{1/2} \to V_{1/2} \times_X V$ and $\tilde{V} \triangleq 0_{V_{1/2}} \times V \subset V_{1/2} \times_X V$, where $0_{V_{1/2}}$ is the zero section of $V_{1/2}$. Note that $\tilde{V}_{1/2}$ and $\tilde{V}$ are isomorphic to $V_{1/2}$ and $V$, respectively. Following [LT2, Page 145], we denote the normal cone to $C_{\Theta_2/Y} \times_Y \Theta_1$ in $C_{\Theta_2/Y}$[17] by $\mathcal{B}(p)_1$ and denote the normal cone to $C_{\Theta_1/Y} \times_Y \Theta_2$ in $C_{\Theta_1/Y}$ by $\mathcal{B}(p)_2$. Both $\mathcal{B}(p)_1$ and $\mathcal{B}(p)_2$ are subcones in $\hat{V}_{1/2} \times_{\hat{X}} \hat{V}$, where $\hat{X} = \operatorname{Spec} \Bbbk[\![T_1^\vee]\!]/(f)$ and $\hat{V}_i = V_i \times_X \hat{X}$. As argued in [LT2, 146], based on the work of [Vis] (see also the recent [Kr1]) there is a canonical rational equivalence[18] $\mathcal{Q}(p) \in W_*(\hat{V}_{1/2} \times_{\hat{X}} \hat{V})$ so that

$$(4.6) \qquad \partial_0 \mathcal{Q}(p) = \mathcal{B}(p)_1 \quad \text{and} \quad \partial_\infty \mathcal{Q}(p) = \mathcal{B}(p)_2.$$

The cones $\mathcal{B}(p)_1$ and $\mathcal{B}(p)_2$ have the following interpretations as shown in [LT2, page 145]. Let $\mathcal{D}(p)_1 \subset \hat{V}$ be the normal cone to $\hat{X}$ in $X$, then

$$(4.7) \qquad\qquad \mathcal{B}(p)_1 = \phi_1^* \mathcal{D}(p)_1$$

where $\phi_1 : \hat{V}_{1/2} \times_{\hat{X}} \hat{V} \longrightarrow \hat{V}$ is the projection. Next, we let $W = \operatorname{Spec} \Bbbk[\![T_2^\vee]\!]$ and let $\hat{W} = \operatorname{Spec} \Bbbk[\![T_2^\vee]\!]/(g)$. We then form the normal cone to $C_{\hat{W}/W} \times_W \hat{X}$ in $C_{\hat{W}/W} \times_W X$, denoted by $\mathcal{D}(p)_2$. It is naturally a subcone in $\hat{V}_{1/2} \times_{\hat{X}} \hat{V}_2$. Then

$$(4.8) \qquad\qquad \mathcal{B}(p)_2 = \phi_2^* \mathcal{D}(p)_2,$$

where $\phi_2 = (1, \rho_2) : \hat{V}_{1/2} \times_{\hat{X}} \hat{V} \longrightarrow \hat{V}_{1/2} \times_{\hat{X}} \hat{V}_2$.

We caution that all the objects so far constructed depends on the lift $\bar{p}$, on the choices of arrows before (4.4) and on the Kuranishi maps. Later we will show that they are canonical in certain degree, up to the symmetry $\operatorname{Aut}(p)$.

In the following, we first construct the relative moduli cycle $[\mathbf{M}, \mathbf{N}]^{\mathrm{virt}}$ and several related cycles based on the collection of the infinitesimal models $\mathcal{D}(p)_2$, etc. First of all, Lemma 2.7 holds for the collection of cones $\{\mathcal{D}(p)_2\}_{p \in \mathbf{M}}$. Namely at each $p \in \mathbf{M}$ the cone

$$\mathcal{D}(p)_2 \times_{\hat{X}} 0 \subset (\hat{V}_{1/2} \times_{\hat{X}} \hat{V}_2) \times_X 0$$

is independent of the choices of the Kuranishi maps $f$, $g$ and $h$ and is invariant under the natural $\operatorname{Aut}(p)$ action. (Here $\hat{X}$ is the subscheme defined before that depends on $p$ implicitly.) Similar statements hold for the cycles $\mathcal{D}(p)_1 \in Z_*\hat{V}$, the cycles $\mathcal{B}(p)_i \in Z_*(\hat{V}_{1/2} \times_{\hat{X}} \hat{V})$ and the rational equivalence $\mathcal{Q}(p) \in W_*(\hat{V}_{1/2} \times_{\hat{X}} \hat{V})$. These were proved in [LT2, Section 3]. Secondly, for each $\alpha$ we pick a pair of locally free subsheaves[19] $\mathcal{L}_\alpha \subset \mathcal{E}_\alpha$ over $R_\alpha$, viewed as a complex $[\mathcal{L}_\alpha \to \mathcal{E}_\alpha]$, and a surjective homomorphism of complexes

$$(4.9) \qquad [\mathcal{L}_\alpha \to \mathcal{E}_\alpha] \Longrightarrow [\mathcal{O}b_{R_\alpha/S_\alpha} \to \mathcal{O}b_{R_\alpha}].$$

---

[17]For closed subscheme $A \subset B$, $C_{A/B}$ is the normal cone to $A$ in $B$.

[18] In this paper we use the convention that a rational equivalence $Q \in W_*Z$ is a cycle in $Z_*(Z \times \mathbf{P}^1)$ so that all its irreducible components are flat over $\mathbf{P}^1$. We then define $\partial_0 Q$ and $\partial_\infty Q$ to be $Q \cap Z \times \{0\}$ and $Q \cap Z \times \{\infty\}$ respectively. In case $\varphi : Z \to W$ is a flat morphism, we denote by $\varphi^* Q$ the flat pull back of the rational equivalence.

For simplicity, we denote by $W_\alpha$ the vector bundle $\mathrm{Vect}(\mathcal{L}_\alpha \oplus \mathcal{E}_\alpha)$. As shown in [LT2, Section 3] there is a unique cycle $B_{i,\alpha} \in Z_* W_\alpha$ so that the collection $\mathcal{B}_i = \{(R_\alpha, B_{i,\alpha}, \mathcal{L}_\alpha \oplus \mathcal{E}_\alpha)\}_\Lambda$ satisfies the cycle consistency criteria for the infinitesimal models $\mathcal{B}(p)_i$ over the pair $(\mathbf{M}, \mathcal{O}b_{\mathbf{M}/\mathbf{N}} \oplus \mathcal{O}b_{\mathbf{M}})$. We let $[\mathcal{B}_1]$ and $[\mathcal{B}_2]$ be the associated class in $A_* \mathbf{M}$ following the basic construction. Next we let $\mathcal{F}_\alpha \triangleq \mathcal{E}_\alpha/\mathcal{L}_\alpha$ be the quotient sheaf, which is locally free by our choice. Then $\mathcal{O}b_{S_\alpha} \otimes_{\mathcal{O}_{S_\alpha}} \mathcal{O}_{R_\alpha}$ is canonically a quotient sheaf of $\mathcal{F}_\alpha$. Again following [LT2, Section 3] we can find a unique cone cycle $\mathcal{D}_{2,\alpha} \subset \mathrm{Vect}(\mathcal{L}_\alpha \oplus \mathcal{F}_\alpha)$ so that the collection $\mathcal{D}_2 = \{(R_\alpha, \mathcal{D}_{2,\alpha}, \mathcal{L}_\alpha \oplus \mathcal{F}_\alpha)\}$ satisfies the cycle consistency criteria with the infinitesimal models $\mathcal{D}(p)_2 \subset \tilde{V}_{1/2} \times_{\tilde{X}} \tilde{V}_2$ over the pair $(\mathbf{M}, \mathcal{O}b_{\mathbf{M}/\mathbf{N}} \oplus \rho^* \mathcal{O}b_{\mathbf{N}})$, where $\rho : \mathbf{M} \to \mathbf{N}$ is the projection. Thus by applying the basic construction to this collection we obtain the class $[\mathcal{D}_2] \in A_* \mathbf{M}$. We will call the class $[\mathcal{D}_2]$ the relative cycle and denoted it by $[\mathbf{M}, \mathbf{N}]^{\mathrm{virt}}$.

**Lemma 4.3.** *We have* $[\mathbf{M}]^{\mathrm{virt}} = [\mathbf{M}, \mathbf{N}]^{\mathrm{virt}}$ *in* $A_* \mathbf{M}$.

To prove the Lemma we need to construct a rational equivalence $[\mathcal{Q}] \in W_* \mathbf{M}$ (or equivalently a class $[\mathcal{Q}] \in A_* \mathbf{M} \times \mathbf{P}^1$) so that $\partial[\mathcal{Q}]$ provides the identity in the Lemma. Here when $[\mathcal{Q}]$ is a class, we define $\partial[\mathcal{Q}] = \partial_0[\mathcal{Q}] - \partial_\infty[\mathcal{Q}]$ with $\partial_0[\mathcal{Q}]$ the image of the Gysin map $0^![\mathcal{Q}]$ associated to $0 \mapsto \mathbf{P}^1$. The $\partial_\infty$ is defined similarly. The proof is parallel to that of [LT2, Lemma 3.4] and will occupy the rest of this subsection. We need to provide a revised proof since the current construction does not rely on the existence of a global vector bundle as was assumed in [LT2].

We first set up the notation relevant to the construction of the cycles $[\mathcal{B}_i]$, etc., following the basic construction of the virtual cycles. First, in constructing the class $[\mathcal{B}_i]$ we first built the index set $\Xi(\mathcal{B}_i)$. Each $a \in \Xi(\mathcal{B}_i)$ associates to an integral closed substack $\mathbf{Y}_a \subset \mathbf{M}$, a multiplicity $m_a$, a maximal open substack $j_a : \mathbf{Y}_a^0 \to \mathbf{Y}_a$ so that $j_a^* \mathcal{O}b_{\mathbf{M}/\mathbf{N}} \oplus j_a^* \mathcal{O}b_{\mathbf{M}}$ is locally free. Over $\mathbf{Y}_a^0$ we have a canonical cone representative $\mathbf{N}_a^0 \subset \mathbf{F}_a^0$, where $\mathbf{F}_a^0$ is the vector bundle (stack) $\mathrm{Vect}(j_a^* \mathcal{O}b_{\mathbf{M}/\mathbf{N}} \oplus j^* \mathcal{O}b_{\mathbf{M}})$ over $\mathbf{Y}_a^0$. We then pick a projective variety $\varphi_a : Y_a \to \mathbf{Y}_a$ generically finite over $\mathbf{Y}_a$, a locally free sheaf $\mathcal{F}_a$ of $\mathcal{O}_{Y_a}$-modules and a quotient sheaf homomorphism $\mathcal{F}_a \to \varphi_a^* (\mathcal{O}b_{\mathbf{M}/\mathbf{N}} \oplus \mathcal{O}b_{\mathbf{M}})$. We form the vector bundle $F_a = \mathrm{Vect}(\mathcal{F}_a)$ with the induced flat morphism $F_a|_{\varphi_a^{-1}(\mathbf{Y}_a^0)} \to \mathbf{F}_a^0$. We let $N_a \subset F_a$ be the closure of the flat pull back of $\mathbf{N}_a^0$. Then $[\mathcal{B}_i]$ is the sum of $m_a \xi(a)$ over all $a \in \Xi(\mathcal{B}_i)$, where $\xi(a)$ is $\deg(\varphi_a)^{-1} \varphi_{a*} 0^!_{F_a}[N_a]$.

We set up the notation for the cycle $[\mathcal{D}_i]$ according to the same rule. We first prove $[\mathcal{B}_1] = [\mathbf{M}]^{\mathrm{virt}}$. This identity follows from the relation (4.7). The actual proof goes as follows: First the class $[\mathcal{D}_1] = [\mathbf{M}]^{\mathrm{virt}}$ (see [LT2, Section 3]). From (4.7) the index sets $\Xi(\mathcal{B}_1) \equiv \Xi(\mathcal{D}_1)$ naturally. For $a \in \Xi(\mathcal{B}_1)$ with $\bar{a} \in \Xi(\mathcal{D}_1)$ the corresponding element, their base substacks $\mathbf{Y}_a = \mathbf{Y}_{\bar{a}}$. Further the intrinsic representative $\mathbf{N}_a^0$ is the flat pull back of $\mathbf{N}_{\bar{a}}^0$ under an obvious homomorphism of vector bundles. Based on these, we can choose $Y_a = Y_{\bar{a}}$, choose locally free sheaves $\mathcal{F}_a = \mathcal{F}_{\bar{a}}$ and then the cone representatives $N_a \equiv N_{\bar{a}}$. This proves $\xi(a) = \xi(\bar{a})$ and hence $[\mathcal{B}_1] = [\mathcal{D}_1] = [\mathbf{M}]^{\text{virt}}$. Since the proof is straightforward, we will omit the details here. For the same reason, we prove $[\mathcal{B}_2] = [\mathbf{M}, \mathbf{N}]^{\text{virt}}$ based on the identity (4.8). Again we will omit the details here.

---

[19]In this paper by a pair of locally free sheaves we mean a locally free subsheaf of a locally free sheaf with locally free quotient sheaf, all of finite ranks.

We now construct the required cycle $[\mathcal{Q}]$ so that

$$(4.10) \qquad \partial_0[\mathcal{Q}] = [\mathcal{B}_1] \quad \text{and} \quad \partial_\infty[\mathcal{Q}] = [\mathcal{B}_2].$$

For each $\alpha \in \Lambda$ we fix the pair $\mathcal{L}_\alpha \subset \mathcal{E}_\alpha$ over $R_\alpha$ as in (4.9). As before we let $W_\alpha = \text{Vect}(\mathcal{L}_\alpha \oplus \mathcal{E}_\alpha)$. We claim that we can find a collection of rational equivalence $Q_\alpha \subset W_* W_\alpha$ indexed by $\alpha \in \Lambda$ that satisfy the following existence Lemma.

**Lemma 4.4.** *To each $\alpha \in \Lambda$ there is a unique rational equivalence $Q_\alpha \in W_* W_\alpha$ of which the following holds. Let $p \in \mathbf{M}$ be any point and let $\bar{p}$ be a lift of $p$ in some chart $R_{\bar{\beta}}$. We let $\hat{V}_{1/2}$, $\hat{V}$ and $\hat{V}_2$ be the vector bundles over $\hat{X}$ (associated to $\bar{p}$) defined before and after Lemma 4.2. Let $R_\alpha$ be the chart, let $R_{\alpha,p} = R_\alpha \times_{\mathbf{M}} p$ and let $\hat{R}_\alpha$ be the formal completion of $R_\alpha$ along $R_{\alpha,p}$. We also pick a morphism $\hat{R}_\alpha \to \hat{X}$ that commutes with $\hat{R}_\alpha \to \mathbf{M}$ and $\hat{X} \to \mathbf{M}$. Now let $O_{1/2} \hookrightarrow O$ be the pair in (4.4) associated to $\bar{p}$. Then up to $\text{Aut}(p)$ there are canonical induced homomorphisms of vector bundles over $R_{\alpha,p}$*

$$O_{1/2} \times R_{\alpha,p} \to h^2(\mathcal{L}_\alpha)|_{R_{\alpha,p}}, \ O \times R_{\alpha,p} \to h^2(\mathcal{E}_\alpha)|_{R_{\alpha,p}} \ \text{and} \ \mathcal{L}_\alpha|_{R_{\alpha,p}} \to O_{1/2} \times R_{\alpha,p}.$$

*We then pick a surjective homomorphism $\varphi_2 : \mathcal{E}_\alpha|_{R_{\alpha,p}} \to O \times R_{\alpha,p}$ so that the following diagram of the complexes of vector bundles over $R_{\alpha,p}$*

$$(4.11) \qquad \begin{array}{ccc} [\mathcal{L}_\alpha|_{R_{\alpha,p}} \to \mathcal{E}_\alpha|_{R_{\alpha,p}}] & \Longrightarrow & [O_{1/2} \times R_{\alpha,p} \to O \times R_{\alpha,p}] \\ & \searrow & \Downarrow \\ & & [h^2(\mathcal{L}_\alpha^\bullet)|_{R_{\alpha,p}} \to h^2(\mathcal{E}_\alpha^\bullet)|_{R_{\alpha,p}}] \end{array}$$

*is commutative. Then there is a vector bundle homomorphism*

$$\Phi_1 : W_\alpha \times_{R_\alpha} \hat{R}_\alpha \longrightarrow (\hat{V}_{1/2} \times_{\hat{X}} \hat{V}) \times_{\hat{X}} \hat{R}_\alpha$$

*extending the obvious homomorphism*

$$W_\alpha \times_{R_\alpha} R_{\alpha,p} \longrightarrow (V_{1/2} \times V) \times R_{\alpha,p}$$

*induced by the $\varphi_2$ before (4.11) so that*

$$\Phi_1^*(\psi^* \mathcal{Q}(p)) = \Phi_2^*(Q_\alpha),$$

*where $\psi : (\hat{V}_{1/2} \times_{\hat{X}} \hat{V}) \times_{\hat{X}} \hat{R}_\alpha \to \hat{V}_{1/2} \times_{\hat{X}} \hat{V}$ and $\Phi_2 : W_\alpha \times_{R_\alpha} \hat{R}_\alpha \to W_\alpha$ are the obvious (flat) projections.*

We will call this Lemma the *rational equivalence consistency criteria.*

This Lemma is proved in [LT2, Section 3] for the case where $\mathbf{M}$ is a quotient stack. The general case is exactly the same. We will not repeat the argument here.

We now fix the collection $\mathcal{Q} = \{(R_\alpha, Q_\alpha, \mathcal{L}_\alpha \oplus \mathcal{E}_\alpha)\}_\Lambda$. To each $\alpha \in \Lambda$ we let $\Xi(Q_\alpha)$ be the set of irreducible components of $Q_\alpha$. For each $a \in \Xi(Q_\alpha)$, we let $T_a \subset Q_\alpha$ be the corresponding irreducible component, let $m_a$ be the multiplicity of $Q_\alpha$ along $T_a$ and let $\mathbf{Y}_a \subset \mathbf{M}$ be the base substack of $a$, which is the closure of the image of $T_a \to \mathbf{M}$. We then let $\iota_a : \mathbf{Y}_a^0 \to \mathbf{Y}_a$ be an open substack so that $\iota_a^* \mathcal{O}b_{\mathbf{M}/\mathbf{N}}$, $\iota_a^* \mathcal{O}b_{\mathbf{M}}$ and $\mathcal{K}_a^0 = \ker\{\iota_a^* \mathcal{O}b_{\mathbf{M}/\mathbf{N}} \to \iota_a^* \mathcal{O}b_{\mathbf{M}}\}$ are locally free. We then pick a homomorphism

$$(4.12) \qquad \eta : \iota_a^* \mathcal{O}b_{\mathbf{M}/\mathbf{N}} \to \mathcal{K}_a^0$$

so that the composite $\mathcal{K}_a^0 \to \iota_*^* \mathcal{O}b_{\mathbf{M/N}} \to \mathcal{K}_a^0$ is the identity. (Such lift exists possibly after shrinking $\mathbf{Y}_a^0 \subset \mathbf{Y}$ if necessary.) We let $\mathbf{V}_{1/2,a}^0 = \mathrm{Vect}(\iota_a^* \mathcal{O}b_{\mathbf{M/N}})$, let $\mathbf{V}_a^0 = \mathrm{Vect}(\mathcal{K}_a^0 \oplus \iota_a^* \mathcal{O}b_{\mathbf{M}})$ and let $\bar{\eta} \colon \mathbf{V}_{1/2,a}^0 \to \mathbf{V}_a^0$ be the immersion induced by $\eta$ and the tautological $\mathcal{O}b_{\mathbf{M/N}} \to \mathcal{O}b_{\mathbf{M}}$.

We now construct the intrinsic representative of $a$ as a cycle in $\mathbf{F}_a^0 \times \mathbf{P}^1$, where $\mathbf{F}_a^0 = \mathbf{V}_{1/2,a}^0 \times_{\mathbf{Y}_a^0} \mathbf{V}_a^0$ is a vector bundle stack over $\mathbf{Y}_a^0$. Let $\mathrm{pr} \colon W_\alpha \to R_\alpha$ be the projection. It is easy to see that for some dense open subset $U$ of $\mathrm{pr}(T_a) \subset R_a$ we can lift the tautological $U \to \mathbf{Y}_a$ to $g \colon U \to \mathbf{Y}_a^0$ and lift the left vertical arrow (below) to surjective horizontal arrow (the top one) as shown in the commutative diagram

$$(4.13) \qquad \begin{array}{ccc} [\mathcal{L}_\alpha \to \mathcal{E}_\alpha] \otimes_{\mathcal{O}_{R_\alpha}} \mathcal{O}_U & \xrightarrow{\ lift\ } & g^*[\iota_a^* \mathcal{O}b_{\mathbf{M/N}} \to \mathcal{K}_a^0 \oplus \iota_a^* \mathcal{O}b_{\mathbf{M}}] \\ \downarrow & & \downarrow \\ g^*[\mathcal{O}b_{\mathbf{M/N}} \to \mathcal{O}b_{\mathbf{M/N}}] & \longrightarrow & g^*[\iota_a^* \mathcal{O}b_{\mathbf{M/N}} \to \iota_a^* \mathcal{O}b_{\mathbf{M/N}}]. \end{array}$$

Note that the other two arrows are tautological ones. We let

$$(4.14) \qquad (W_\alpha \times \mathbf{P}^1)|_{U \times \mathbf{P}^1} \longrightarrow \mathbf{F}_a^0 \times \mathbf{P}^1$$

be the induced projection. Then an easy argument shows that the statement in Lemma 4.4 implies that there is a reduced and irreducible cycle $\mathbf{Q}_{a,\eta}^0 \subset \mathbf{F}_a^0 \times \mathbf{P}^1$ so that the restriction of $T_a$ to fibers over $U \times \mathbf{P}^1$ is a dense open subset of the flat pull back of $\mathbf{Q}_{a,\eta}^0$ via the arrow in (4.14). We call $\mathbf{Q}_{a,\eta}^0$ an intrinsic representative of $a$.

Once we have constructed the intrinsic representative of $a \in \Xi(Q_\alpha)$, we can define an equivalence relation on $\cup \Xi(Q_\alpha)$ as we did before. Let $a \in \Xi(Q_\alpha)$ and $b \in \Xi(Q_\beta)$ be two elements. In case $\mathbf{Y}_a \neq \mathbf{Y}_b$, then $a \not\sim b$. In case $\mathbf{Y}_a = \mathbf{Y}_b$, we let $\mathbf{Q}_{a,\eta}^0$ and $\mathbf{Q}_{b,\eta}^0$ be the respective representatives of $a$ and $b$, based on the same $\eta$ in (4.12). Then we say $a \sim b$ if $\mathbf{Q}_{a,\eta}^0 \cap \mathbf{Q}_{b,\eta}^0$ is dense in both $\mathbf{Q}_{a,\eta}^0$ and $\mathbf{Q}_{b,\eta}^0$. This defines an equivalence relation. Further, whenever $a \sim b$ then $m_a = m_b$. Now let $\Xi(\mathcal{Q})$ be the set of equivalence classes. For any $a \in \Xi(\mathcal{Q})$, we pick a projective variety $Y_a$ and a generically finite morphism $\varphi_a \colon Y_a \to \mathbf{Y}_a$. We then pick a pair of locally free subsheaves $\mathcal{L}_a \hookrightarrow \mathcal{E}_a$ over $Y_a$, a surjective homomorphism of complexes

$$[\mathcal{L}_a \to \mathcal{E}_a] \Longrightarrow \varphi_a^*[\mathcal{O}b_{\mathbf{M/N}} \to \mathcal{O}b_{\mathbf{M}}]$$

and a lift (of the above arrow)

$$[\mathcal{L}_a \to \mathcal{E}_a] \Longrightarrow \varphi_a^*[\mathcal{O}b_{\mathbf{M/N}} \to \mathcal{K}_a^0 \oplus \mathcal{O}b_{\mathbf{M}}]$$

over a dense open subset $Y_a^0 \subset Y_a$. After that, we let $W_a = \mathrm{Vect}(\mathcal{L}_a \oplus \mathcal{E}_a)$ and let $j_a \colon W_a|_{Y_a^0} \to \mathbf{F}_a^0$ be the induced morphism. By shrinking $Y_a^0$ if necessary, we can assume that $W_a|_{Y_a^0} \to \mathbf{F}_a^0$ is flat. With all these chosen, we let $Q_a$ be the closure in $W_a \times \mathbf{P}^1$ of the flat pull back of $\mathbf{Q}_{a,\eta}^0$ and define $\xi(a) = \deg(\varphi_a)^{-1} \varphi_{a*} 0_{W_a}^![Q_a]$. We define

$$[\mathcal{Q}] = \sum_{a \in \Xi(\mathcal{Q})} m_a \xi(a) \in A_* \mathbf{M} \times \mathbf{P}^1.$$

It remains to prove the identity

$$(4.15) \qquad \partial_\infty[\mathcal{Q}] = \sum_{a \in \Xi(\mathcal{Q})} m_a \partial_\infty \xi(a) = [\mathcal{B}_2] \in A_* \mathbf{M}$$

and the similar identity with $\partial_\infty$ (resp. $[\mathcal{B}_2]$) replaced by $\partial_0$ (resp. $[\mathcal{B}_1]$). We will prove the identity (4.15). The proof of the other identity is similar. To achieve this, we will define a function $\mu : \Xi(\mathcal{Q}) \times \Xi(\mathcal{B}_2) \to \mathbb{Q}$ that satisfies

$$(4.16) \qquad \sum_{a \in \Xi(\mathcal{Q})} m_a \mu(a, b) = m_b, \qquad \forall\, b \in \Xi(\mathcal{B}_2)$$

and

$$(4.17) \qquad \partial_\infty \xi(a) = \sum_{b \in \Xi(\mathcal{B}_2)} \mu(a, b) \xi(b), \qquad \forall\, a \in \Xi(\mathcal{Q}).$$

Once these two are established, then

$$\partial_\infty [\mathcal{Q}] = \sum m_a \, \partial_\infty \xi(a) = \sum m_a \sum \mu(a, b)\, \xi(b) = \sum m_b\, \xi(b) = [\mathcal{B}_2],$$

which is (4.15).

We first construct the function $\mu$. We begin with any (smooth) chart $R_\alpha$. We let $\mathcal{L}_\alpha \oplus \mathcal{E}_\alpha$ be the locally free sheaves on $R_\alpha$ chosen before. For simplicity we denote by $\mathcal{W}_\alpha$ the sheaf $\mathcal{L}_\alpha \oplus \mathcal{E}_\alpha$ and continue to denote by $W_\alpha$ the vector bundle $\mathrm{Vect}(\mathcal{L}_\alpha \oplus \mathcal{E}_\alpha)$. We let $B_{2,\alpha} \subset W_\alpha$ and $Q_\alpha \subset W_\alpha \times \mathbf{P}^1$ be the cycles constructed before. For each $a \in \Xi(\mathcal{Q})$ we let $Q_{\alpha,a}$ be the closure in $W_\alpha \times \mathbf{P}^1$ of the pull back of $\mathbf{Q}^0_{a,\eta}$ via the arrow in (4.14). The cycle $Q_{\alpha,a}$ is also the union of irreducible components of $Q_\alpha$ in the equivalence class $a$. With this choice of $Q_{\alpha,a}$, we have the decomposition

$$(4.18) \qquad Q_\alpha = \sum_{a \in \Xi(\mathcal{Q})} m_a\, Q_{\alpha,a}.$$

Similarly, we have a canonical decomposition

$$B_{2,\alpha} = \sum_{b \in \Xi(\mathcal{B}_2)} m_b\, B_{2,\alpha,b}.$$

As before, we let $\Xi(B_{2,\alpha})$ be the index set of irreducible components of $B_{2,\alpha}$. For $c \in \Xi(B_{2,\alpha})$, we denote by $T_c$ the corresponding component. Then because $\partial_\infty Q_\alpha = B_{2,\alpha}$, there is a unique function $\mu_\alpha(a, \cdot) : \Xi(B_{2,\alpha}) \to \mathbb{Z}$ so that $\partial_\infty Q_{\alpha,a} = \sum \mu_\alpha(a, c) T_c$. Again, by cycle (rational equivalence) consistency criteria and the invariance of the cycle under $\mathrm{Aut}(p)$, we conclude that whenever $c_1 \sim c_2$ then $\mu_\alpha(a, c_1) = \mu_\alpha(a, c_2)$. Thus $\mu_\alpha(a, \cdot)$ descends to a function $\mu_\alpha(a, \cdot)$ from $\Xi(\mathcal{B}_2)$ to $\mathbb{Z}$. Finally, for $b \in \Xi(\mathcal{B}_2)$ we define $B_{2,\alpha,b}$ to be the union of $T_c$ for those $c \in \Xi(B_{2,\alpha})$ such that $c \sim b$. Then we have the identity

$$\partial_\infty Q_{\alpha,a} = \sum_{b \in \Xi(\mathcal{B}_2)} \mu_\alpha(a, b)\, B_{2,\alpha,b}.$$

For the same reason whenever $R_\alpha \times_{\mathbf{M}} R_\beta \times_{\mathbf{M}} \mathbf{Y}_b \neq \emptyset$, then $\mu_\alpha(a, b) = \mu_\beta(a, b)$. Hence, we can define $\mu(a, b)$ to be $\mu_\alpha(a, b)$ for those $\alpha$ so that $R_\alpha \times_{\mathbf{M}} \mathbf{Y}_b \neq \emptyset$. Then the identity (4.16) follows from collecting term $B_{2,\alpha,b}$ in

$$\sum_{b \in \Xi(\mathcal{B}_2)} m_b\, B_{2,\alpha,b} = B_{2,\alpha} = \partial_\infty Q_\alpha = \sum_{a \in \Xi(\mathcal{Q}} m_a \sum_{b \in \Xi(\mathcal{B}_2)} \mu(a, b)\, B_{2,\alpha,b}$$

for some $\alpha$ so that $R_\alpha \times_{\mathbf{M}} \mathbf{Y}_b \neq \emptyset$.

We now investigate (4.17). Let $a \in \Xi(\mathcal{Q})$ be any element, let $\varphi_a : Y_a \to \mathbf{Y}_a$ be the generically finite morphism and let $\mathcal{L}_a \subset \mathcal{E}_a$ be the pair of sheaves chosen before. We still denote by $W_a$ the vector bundle $\mathrm{Vect}(\mathcal{L}_a \oplus \mathcal{E}_a)$. We let $Q_a \in W_* W_a$ be the representative of $a$ and let $\partial_\infty Q_a = \sum m_i D_i$ be the decomposition into irreducible components. Since $Q_a \to \mathbf{M}$ factor through $\mathbf{Y}_a \subset \mathbf{M}$, all $D_i$ lie over $\mathbf{Y}_a$. We now let $\mathbf{A} \subset \mathbf{Y}_a$ be a closed integral substack. We let $\Upsilon_{\mathbf{A}}$ be (the index set of) those $D_i$ so that the image $D_i \to \mathbf{M}$ is exactly $\mathbf{A}$. Similarly we let $\Xi(\mathcal{B}_2)_{\mathbf{A}}$ be those $b \in \Xi(\mathcal{B}_2)$ so that $\mathbf{Y}_b = \mathbf{A}$. Clearly if we can show that

$$(4.19) \qquad \deg(\varphi_a)^{-1} \sum_{i \in \Upsilon_{\mathbf{A}}} m_i \, \varphi_{a*} 0^!_{W_a} [D_i] = \sum_{b \in \Xi(\mathcal{B}_2)_{\mathbf{A}}} \mu(a,b) \, \xi(b),$$

since $\mathbf{A} \subset \mathbf{Y}_a$ is arbitrary, (4.17) will follow immediately.

In the remainder of this section we will prove the identity (4.19) for any pair $(a, \mathbf{A})$. We pick an étale $\varphi_U : U \to \mathbf{Y}_a$ so that $U \times_{\mathbf{Y}_a} \mathbf{A} \to \mathbf{A}$ is dominant. We then pick a pair of locally free sheaves $\mathcal{L}_U \subset \mathcal{E}_U$ and surjective homomorphisms $\phi_1$ and $\phi_2$ of complexes as shown that make the diagram

$$
\begin{array}{ccc}
[\mathcal{L}_U \to \mathcal{E}_U] & \xrightarrow{\ \phi_1\ } & \varphi_U^*[\iota_a^* \mathcal{O}b_{\mathbf{M}/\mathbf{N}} \to \mathcal{K}_a^0 \oplus \iota_a^* \mathcal{O}b_{\mathbf{M}}] \\
\Big\downarrow {\scriptstyle\phi_2} & & \Big\downarrow \\
\varphi_U^*[\mathcal{O}b_{\mathbf{M}/\mathbf{N}} \to \mathcal{O}b_{\mathbf{M}}] & \longrightarrow & \varphi_U^*[\iota_a^* \mathcal{O}b_{\mathbf{M}/\mathbf{N}} \to \iota_a^* \mathcal{O}b_{\mathbf{M}}]
\end{array}
$$

commutative. Here the top-right corner is the complex in (4.13) and the right vertical arrow is the standard projection. As before, we let $W_U$ be the vector bundle $\mathrm{Vect}(\mathcal{L}_U \oplus \mathcal{E}_U)$. Then the homomorphism $\phi_1$ defines a flat morphism $W_U|_{\varphi_U^{-1}(\mathbf{Y}_a^0)} \to \mathbf{F}_a^0$. We let $Q_U \subset W_* W_U$ be the closure of the flat pull back of $\mathbf{Q}_{a,\eta}^0 \in W_* \mathbf{F}_a^0$.

We next consider the projections

$$q_1 : U \times_{\mathbf{Y}_a} Y_a \to U \quad \text{and} \quad \tilde{q}_1 : W_U \times_{\mathbf{Y}_a} Y_a \to W_U.$$

We let $\tilde{W} \subset W_U \times_{\mathbf{Y}_a} Y_a$ be a dense open subset so that the tautological $\tilde{W} \to \mathbf{Y}_a$ factor through $\mathbf{Y}_a^0$ and $\tilde{W} \to \mathbf{F}_a^0$ is quasi-finite. We then let $\tilde{Q}_U$ be the closure in $(W_U \times_{\mathbf{Y}_a} Y_a) \times \mathbf{P}^1$ of the flat pull back of $\mathbf{Q}_{a,\eta}^0$ via the obvious $\tilde{W} \to \mathbf{F}_a^0$. Then since $\mathbf{Q}_{a,\eta}^0$ is dominant over $\mathbf{Y}_a^0$ and since $\tilde{q}_1$ is proper, we have $\tilde{q}_{1*}(\tilde{Q}_U) = \deg(\varphi_a) Q_U$ and then

$$(4.20) \qquad \tilde{q}_{1*}(\partial_\infty \tilde{Q}_U) = \partial_\infty \tilde{q}_{1*}(\tilde{Q}_U) = \deg(\varphi_a) \partial_\infty Q_U.$$

For the convenience of the readers, we list the related rational equivalence relations constructed:

$$\mathbf{Q}_{a,\eta}^0 \in W_* \mathbf{F}_a^0, \ Q_U \in W_* W_U, \ \tilde{Q}_U \in W_*(W_U \times_{\mathbf{Y}_a} Y_a) \text{ and } Q' \in W_* W.$$

($Q'$ will be constructed in (4.23).) We now let $\Upsilon$ (resp. $\tilde{\Upsilon}$) be the (index) set of those irreducible components of $\partial_\infty Q_U$ (resp. $\partial_\infty \tilde{Q}_U$) that lie over and dominate $\mathbf{A}$. For $i \in \Upsilon$ or $\tilde{\Upsilon}$ we denote by $T_i$ the corresponding component and by $m_i$ the multiplicity of $Q_U$ or $\tilde{Q}_U$ along $T_i$. We need to divide $\tilde{\Upsilon}$ into two parts: One is $\tilde{\Upsilon}_1$ which consists of those $T_i$ so that $\tilde{q}_{1*}(T_i) \ne 0$. We let $\tilde{\Upsilon}_0 = \tilde{\Upsilon} - \tilde{\Upsilon}_1$. Because of (4.20), there is a map $\lambda : \tilde{\Upsilon}_1 \to \Upsilon$ so that under $\tilde{q}_{1*}$ the component $T_i$ is mapped onto $T_{\lambda(c)}$, say via $\psi_{i,\lambda(i)} : T_i \to T_{\lambda(i)}$. Then (4.20) implies

$$(4.21) \qquad \sum_{i \in \lambda^{-1}(j)} \deg(\psi_{i,j}) = \deg(\varphi_a), \quad j \in \Upsilon.$$

We now compare the collection $\tilde{\Upsilon}$ with the collection $\Upsilon_{\mathbf{A}}$. We let $p_1$ and $p_2$ be the projections of $U \times_{\mathbf{Y}_a} Y_a$ to $U$ and $Y_a$, respectively. We pick a pair of locally free sheaves $\mathcal{L} \hookrightarrow \mathcal{E}$ on $U \times_{\mathbf{Y}_a} Y_a$ so that there are two surjective homomorphisms $\phi_1$ and $\phi_2$ of complexes as shown in the commutative diagram

(4.22)
$$
\begin{array}{ccc}
[\mathcal{L} \to \mathcal{E}] & \xrightarrow{\ \phi_1\ } & p_2^*[\mathcal{L}_a \to \mathcal{E}_a] \\
\downarrow{\scriptstyle \phi_2} & & \downarrow \\
p_1^*[\mathcal{L}_U \to \mathcal{E}_U] & \longrightarrow & [\iota_a^* \mathcal{O}b_{\mathbf{M}/\mathbf{N}} \to \mathcal{K}_a^0 \oplus \iota_a^* \mathcal{O}b_{\mathbf{M}}].
\end{array}
$$

Here the two remaining arrows in the diagram are the ones chosen before. We let $W$ be the vector bundle $\mathrm{Vect}(\mathcal{L} \oplus \mathcal{E})$ over $U \times_{\mathbf{Y}_a} Y_a$ and let $\zeta_1 : W \to W_U \times_{\mathbf{Y}_a} Y_a$ and $\zeta_2 : W \to W_a$ be the morphisms to vector bundles over $U \times_{\mathbf{Y}_a} Y_a$ and $Y_a$ induced by $\phi_1$ and $\phi_2$, respectively. The map $\zeta_1$ is obviously smooth whose fibers are vector spaces. Since $U \times_{\mathbf{Y}_a} Y_a \to Y_a$ is étale, $\zeta_2$ is also smooth with affine fibers. Further, because of the commutative diagram (4.22) and the rational equivalence consistency criteria, the flat pull back of $Q_U$ via $\zeta_1$ and the flat pull back of $Q_a$ via $\zeta_2$ are identical over a dense open subset $\mathcal{U} \subset U \times_{\mathbf{Y}_a} Y_a$ that is flat over both $U$ and $Y_a$. Hence since both $Q_U$ and $Q_a$ are closed and since all their irreducible components dominate $\mathbf{Y}_a$,

(4.23)
$$
Q' \triangleq \zeta_1^*(Q_U) = \zeta_2^*(Q_a).
$$

We now let $\Upsilon'$ be the set of those irreducible components of $\partial_\infty Q'$ that lie over and dominate $\mathbf{A}$. Since the fibers of $\zeta_1$ are vector spaces, $\Upsilon'$ is naturally identical to $\tilde{\Upsilon}$, say via $\lambda' : \Upsilon' \to \tilde{\Upsilon}$. Further, for each $i \in \Upsilon'$ the component $T_i$ is a vector bundle over $T_{\lambda'(i)}$. On the other hand, since $p_2 : U \times_{\mathbf{Y}_a} Y_a \to Y_a$ is smooth, each irreducible component of $\partial_\infty Q'$ is a component of the flat pull back of a $T_i$ for some $i \in \Upsilon_{\mathbf{A}}$. Thus there is a map $\lambda_a : \Upsilon' \to \Upsilon_{\mathbf{A}}$ so that for each $i \in \Upsilon_{\mathbf{A}}$ the union $\sum_{j \in \lambda_a^{-1}(i)} T_j$ is exactly the pull back of $T_i$. For $i \in \Upsilon'$ and $j \in \lambda_a(i)$, we let $\varphi_{i,j} : T_i \to T_j$ be the tautological map.

We are now ready to prove the identity (4.19). For each $i \in \Upsilon'$, we let $b(T_i)$ be the image of $T_i \to U \times_{\mathbf{Y}_a} Y_a$. We pick a projective variety $Z_i$ so that a dense open subset $Z_i^0 \subset Z_i$ is a finite branched cover of $b(T_i)$ and the induced $Z_i^0 \to \mathbf{A}$ extends to $\varphi_i : Z_i \to \mathbf{A}$. We then pick a pair of locally free sheaves $\mathcal{L}_i \subset \mathcal{E}_i$ over $Z_i$ and a surjective homomorphism of complexes $\phi_1'$ that lifts to a surjective $\phi_2'$ as shown in the commutative diagram

$$
\begin{array}{ccc}
[\mathcal{L}_i \to \mathcal{E}_i] & \xrightarrow{\ \phi_2'\ } & [\mathcal{L}|_{Z_i^0} \to \mathcal{E}|_{Z_i^0}] \\
\downarrow{\scriptstyle \phi_1'} & & \downarrow \\
\varphi_i^*[\mathcal{O}b_{\mathbf{M}/\mathbf{N}} \to \mathcal{O}b_{\mathbf{M}}] & \longrightarrow & \varphi_i^*[\mathcal{O}b_{\mathbf{M}/\mathbf{N}}|_{Z_i^0} \to \mathcal{O}b_{\mathbf{M}}|_{Z_i^0}].
\end{array}
$$

Then following the definition of $\xi$ we let $W_i = \mathrm{Vect}(\mathcal{L}_i \oplus \mathcal{E}_i)$ and let $C_i \subset W_i \times \mathbf{P}^1$ be the closure of the flat pull back of $T_i \subset W|_{Z_i^0} \times \mathbf{P}^1$ via the induced $W_i|_{Z_i^0} \to W|_{Z_i^0}$. We then define (the non-normalized) $\bar{\xi}(i) = \varphi_{i*} 0^!_{W_i}[C_i] \in A_* \mathbf{M} \times \mathbf{P}^1$.

For $i \in \Upsilon$, $\tilde{\Upsilon}$ or $\Upsilon_{\mathbf{A}}$ we define the variety $Z_i^0$ to be the image of $T_i$ in $U$, in $U \times_{\mathbf{Y}_a} Y_a$ or $Y_a$ respectively, and then define $\varphi_i : Z_i \to \mathbf{A}$ and the class $\bar{\xi}(i)$ along the same line. We let $\psi_{U,\mathbf{A}}$ be the induced morphism $U \times_{\mathbf{Y}_a} \mathbf{A} \to \mathbf{A}$. We have the following Lemma concerning the relations $\lambda : \tilde{\Upsilon}_1 \to \Upsilon$, $\lambda' : \Upsilon' \to \tilde{\Upsilon}$ and $\lambda_a : \Upsilon' \to \Upsilon_{\mathbf{A}}$ defined before.

**Lemma 4.5.** *The following relations hold: 1. For each $i \in \Upsilon'$ we have $\bar{\xi}(i) = \bar{\xi}(\lambda'(i))$; 2. For each $i \in \tilde{\Upsilon}$ we have $\deg(\varphi_{i,\lambda_a(i)})\bar{\xi}(\lambda_a(i)) = \bar{\xi}(i)$ and for any $j \in \Upsilon_{\mathbf{A}}$ we have $\sum_{i \in \lambda_a^{-1}(j)} \deg(\varphi_{i,j}) = \deg(\psi_{U,\mathbf{A}})$; 3. For any $i \in \tilde{\Upsilon}_1$ we have $\deg(\psi_{i,\lambda(i)})\bar{\xi}(\lambda(i)) = \bar{\xi}(i)$; 4. For $i \in \tilde{\Upsilon}_0$ we have $\bar{\xi}(i) = 0$.*

*Proof.* The proof of (1) is parallel to the proof of Lemma 2.6, using the fact that the fibers of $\zeta_1$ are vector spaces. We will omit the proof here. The proof of the identity in (2) concerning $\bar{\xi}(\cdot)$ is similar and will be omitted too. The identity concerning the degrees in (2) follows from the fact that $U \times_{\mathbf{Y}_a} Y_a \to Y_a$ is étale. We now prove (3). Let $i \in \tilde{\Upsilon}$ and $j = \lambda(i) \in \Upsilon$. We let $Z_i^0 \subset Z_i$ and $Z_j^0 \subset Z_j$ be the pair of varieties constructed in defining the classes $\bar{\xi}(i)$ and $\bar{\xi}(j)$. By definition, there is a canonical dominant morphism $Z_i^0 \to Z_j^0$. Hence without loss of generality we can assume that it extends to $\rho \colon Z_i \to Z_j$. With this choice of $Z_i$ and $Z_j$, we can choose the pairs of sheaves $\mathcal{L}_i \subset \mathcal{E}_i$ and $\mathcal{L}_j \subset \mathcal{E}_j$ be so that the former is the pull back of the later via $\rho$. Then because of the relation (4.20), the cycle representatives $C_i \in Z_* \mathrm{Vect}(\mathcal{L}_i \oplus \mathcal{E}_i)$ and $C_j \in Z_* \mathrm{Vect}(\mathcal{L}_j \oplus \mathcal{E}_j)$ satisfies $\tilde{\rho}_* C_i = \deg(\psi_{i,j}) C_j$. Here we used the fact that $C_i$ is supported on a single variety. This relation implies (3) immediately. The proof of (4) is parallel and will be omitted. This proves the Lemma. $\qquad\square$

It follows from the Lemma that

$$\deg(\psi_{U,\mathbf{A}}) \sum_{i \in \Upsilon_{\mathbf{A}}} m_i \, 0^!_{W_a}[D_i] = \deg(\psi_{U,\mathbf{A}}) \sum_{i \in \Upsilon_{\mathbf{A}}} m_i \, \bar{\xi}(i) = \sum_{j \in \Upsilon'} m_j \, \bar{\xi}(j).$$

Here we used (2) to derive the last identity. Because $\bar{\xi}(i) = 0$ for $i \in \tilde{\Upsilon}_0$, by (1) of the Lemma the right hand side above is

$$\sum_{i \in \tilde{\Upsilon}_1} m_i \, \bar{\xi}(i) = \sum_{i \in \tilde{\Upsilon}_1} \deg(\psi_{i,\lambda(i)}) \, m_i \, \bar{\xi}(\lambda(i)) = \deg(\varphi_a) \sum_{i \in \Upsilon} m_i \, \bar{\xi}(i)$$

$$= \deg(\varphi_a) \deg(\psi_{U,\mathbf{A}}) \sum_{b \in \Xi(\mathcal{B}_2)_{\mathbf{A}}} \mu(a,b) \xi(b).$$

Here the first equality follows from (3) in the Lemma, the second equality follows from (4.21) and the definition of $\mu(a,b)$, using the fact that $U$ is étale over $\mathbf{Y}_a$. This proves the identity (4.19) and hence completes the proof of Lemma 4.3.

### 4.2. Three special cases.

In this subsection, we will study three cases of relative obstruction theories and derive some identities.

We first study the case where $\mathbf{N}$ is a DM-stack and $\mathbf{M} \to \mathbf{N}$ is a substack defined by the vanishing of a C-divisor $(\mathbf{L}, \mathbf{s})$. Namely if $S_\alpha \to \mathbf{N}$ is a chart (in an atlas $\Lambda$) and $(L_\alpha, s_\alpha)$ is the associated C-divisor on $S_\alpha$, then $\mathbf{M} \times_{\mathbf{N}} S_\alpha$ is defined by the vanishing of $s_\alpha$. We denote $\mathbf{M} \times_{\mathbf{N}} S_\alpha$ by $R_\alpha$. We claim that there is a canonical relative obstruction theory of $R_\alpha/S_\alpha$ taking values in the sheaf $\mathcal{O}_{R_\alpha}(L_\alpha)$ ($= \mathcal{O}_{R_\alpha}(L_\alpha|_{R_\alpha})$). Since $R_\alpha$ is a subscheme of $S_\alpha$, the relative first order deformations $\mathfrak{Def}^1_{R_\alpha/S_\alpha} = 0$. Now let $\xi = (B, I, \varphi)$ be any object in $\mathfrak{Tri}_{R_\alpha/S_\alpha}$. Since $B$ is a $\Gamma(\mathcal{O}_{S_\alpha})$-algebra, its associated morphism $\phi \colon \mathrm{Spec}\, B \to S_\alpha$ is an extension of $\mathrm{Spec}\, B/I \to S_\alpha$. We consider the section $s_\alpha \circ \phi$ of $\phi^* L_\alpha$. Since $s_\alpha \circ \phi|_{\mathrm{Spec}\, B/I} \equiv 0 \in \Gamma(\varphi^* L_\alpha)$, $s_\alpha \circ \phi$ is an element in $\Gamma(R_\alpha, L_\alpha) \otimes I$. We define this element to be the relative obstruction class $\mathfrak{ob}_{R_\alpha/S_\alpha}(\xi)$. Of course, $\mathfrak{ob}_{R_\alpha/S_\alpha}(\xi) = 0$ if and only if $\phi: \operatorname{Spec} B \to S_\alpha$ factor through $R_\alpha \subset S_\alpha$. Since $R_\alpha \to S_\alpha$ is an immersion, $\varphi$ extends to $\operatorname{Spec} B \to R_\alpha$ as $S_\alpha$-morphism if and only if $\phi$ already factor through $R_\alpha \subset S_\alpha$. This proves that $\mathfrak{ob}_{R_\alpha/S_\alpha}$ is an obstruction assignment. Since $(\mathbf{L}, \mathbf{s})$ is a C-divisor on $\mathbf{N}$, this defines a relative obstruction theory of $\mathbf{M}/\mathbf{N}$.

Now assume both $\mathbf{M}$ and $\mathbf{N}$ have perfect obstruction theories provided by the data $\{\mathcal{E}_\alpha^\bullet, \mathfrak{ob}_{R_\alpha}\}_\Lambda$ and $\{\mathcal{F}_\alpha^\bullet, \mathfrak{ob}_{R_\alpha/S_\alpha}\}_\Lambda$. We assume further that the obstruction theories of $\mathbf{M}$ and $\mathbf{N}$ are compatible to the relative obstruction theory of $\mathbf{M}/\mathbf{N}$ in the sense of Definition 4.1 with $\mathcal{L}_\alpha^\bullet = [0 \to \mathcal{O}_{R_\alpha}(L_\alpha)]$.

**Lemma 4.6.** *Let the notation be as before. Then* $[\mathbf{M}]^{\mathrm{virt}} = c_1(\mathbf{L}, \mathbf{s})[\mathbf{N}]^{\mathrm{virt}}$.

*Proof.* We will follow the notation developed before and after the Lemma 4.2. Let $p \in \mathbf{N}$ be any point, $S_\alpha$ be an étale chart of $\mathbf{N}$ with a lift $\bar{p} \in S_\alpha$ of $p$. We let $T_2 = h^1(\mathcal{F}_\alpha^\bullet \otimes \Bbbk_{\bar{p}})$ and let $O_2 = h^2(\mathcal{F}_\alpha^\bullet \otimes \Bbbk_{\bar{p}})$. We let $\hat{W}_p$ be the formal completion of $S_\alpha$ along $\bar{p}$. Here we use subscript $p$ instead $\bar{p}$ since $\hat{W}_p$ depends on $p$ up to $\mathrm{Aut}(p)$. We know $\hat{W}_p$ is $\operatorname{Spec} \Bbbk[\![T_2^\vee]\!]/(g)$, where $g$ is a Kuranishi map of the obstruction theory of $S_\alpha$ at $\bar{p}$. We denote by $C_p$ the normal cone to $\hat{W}_p$ in $W_p \triangleq \operatorname{Spec} \Bbbk[\![T_2^\vee]\!]$. The cone $C_p$ is naturally embedded in $O_2 \times \hat{W}_p$. Again, the pair $C_p \subset O_2 \times \hat{W}_p$ only depends on the point $p \in \mathbf{N}$, up to the symmetry $\mathrm{Aut}(p)$. By the construction of the virtual moduli cycles in [LT2] and in section 2.2, this collection of cone cycles $\{C_p\}_{p \in \mathbf{N}}$ can be algebraized and hence gives rise to a cycle, the virtual moduli cycle $[\mathbf{N}]^{\mathrm{virt}}$.

Now assume $p \in \mathbf{M}$. Then $\bar{p} \in R_\alpha$ where $R_\alpha = S_\alpha \times_{\mathbf{N}} \mathbf{M}$. As in (4.3), we let $T_1 = h^1(\mathcal{E}_\alpha \times \Bbbk_{\bar{p}})$ and $O_1 = h^2(\mathcal{E}_\alpha \otimes \Bbbk_{\bar{p}})$, etc. There are two cases to consider. One is when $\delta: T_2 \to O_{1/2}$ in (4.3) is 0. In this case $T = T_1 = T_2$, $O = O_1$ and $O_2 = O_1/O_{1/2}$. Further, we can choose the relative Kuranishi map $h \in \Bbbk[\![T^\vee]\!]/(g) \otimes O_{1/2}$ to be $\hat{s}_\alpha$, the pull back of $s_\alpha$ via $\operatorname{Spec} \Bbbk[\![T^\vee]\!]/(g) \cong \hat{W}_p \to S_\alpha$. Here we have used the canonical isomorphism $O_{1/2} = \mathcal{L}_\alpha|_{\bar{p}}$. Now let $X_p$ and $\hat{X}_p$ be the schemes $\operatorname{Spec} \Bbbk[\![T^\vee]\!]$ and $\operatorname{Spec}[\![T^\vee]\!]/(f)$, where $f$ is the Kuranishi map in (4.5). Then $X_p = W_p$ and the cone $\mathcal{D}(p)_2$ is the normal cone to $C_p \times_{\hat{W}_p} \{\hat{s}_\alpha = 0\}$ in $C_p$. The cone $\mathcal{D}(p)_2$ is naturally embedded in $\hat{V}_{1/2} \times_{\hat{X}_p} \hat{V}_2$. For the same reason, when $\delta$ is non-zero, thus surjective since $\dim O_{1/2} = 1$, the cone $\mathcal{D}(p)_2$ is also the normal cone to $C_p \times_{\hat{W}_p} \{\hat{s}_\alpha = 0\}$ in $C_p$.

We are now ready to prove the Lemma. First, from the discussion before Lemma 4.3 we know the collection $\{\mathcal{D}(p)_2\}_{p \in \mathbf{M}}$ can be algebraized. By definition, the cycle constructed based on $\{\mathcal{D}(p)_2\}$ following the basic construction in Section 2.2 is the relative virtual cycle $[\mathbf{M}, \mathbf{N}]^{\mathrm{virt}}$. However, since $\hat{s}_\alpha$ are pull back of the section $s_\alpha$, the cone $\mathcal{D}(p)_2$ is the normal cone to $C_p \times_{S_\alpha} \{s_\alpha = 0\}$ in $C_p$. Further, since the collection $\{(L_\alpha, s_\alpha)\}$ is the restriction of $(\mathbf{L}, \mathbf{s})$ to charts $S_\alpha$, a repetition of the proof of Lemma 4.3 shows that $[\mathbf{M}, \mathbf{N}]^{\mathrm{virt}} = c_1(\mathbf{L}, \mathbf{s})[\mathbf{N}]^{\mathrm{virt}}$. Then combined with the identity $[\mathbf{M}]^{\mathrm{virt}} = [\mathbf{M}, \mathbf{N}]^{\mathrm{virt}}$ in Lemma 4.3, we have $[\mathbf{M}]^{\mathrm{virt}} = c_1(\mathbf{L}, \mathbf{s})[\mathbf{N}]^{\mathrm{virt}}$. This completes the proof of the Lemma. □

We now investigate the second case. We let $\mathbf{N}$ be a DM-stack with a morphism $\mathbf{N} \to X$ to a scheme $X$. For simplicity we assume $X$ is smooth. Let $\xi: X_0 \longrightarrow X$ be a smooth subvariety and let $\mathbf{M}$ be defined by the Cartesian product

$$\mathbf{M} = \mathbf{N} \times_X X_0.$$

Then $\mathbf{M}$ is a substack of $\mathbf{N}$. Let $\mathbf{L}$ over $\mathbf{M}$ be the pull back of the normal bundle to $X_0$ in $X$. Similar to the case just studied, there is a canonical relative obstruction theory of $\mathbf{M}/\mathbf{N}$ taking values in the cohomology of the complex $\mathcal{L}^\bullet = [0 \to \mathcal{O}_\mathbf{M}(\mathbf{L})]$. We now assume $\mathbf{M}$ and $\mathbf{N}$ both have perfect obstruction theories and are compatible to the relative obstruction theory of $\mathbf{M}/\mathbf{N}$. By the intersection theory of DM-stacks [Vis], the Gysin map $\eta^!{[\mathbf{N}]}^{\mathrm{virt}} \in A_*\mathbf{M}$.

**Lemma 4.7.** *Let the notation be as above. Then* $[\mathbf{M}, \mathbf{N}]^{\mathrm{virt}} = \eta^!{[\mathbf{N}]}^{\mathrm{virt}}$.

*Proof.* The proof is similar to that of the previous Lemma with slight modification. We shall omit the proof here. ∎

We now investigate the third case. We let $\mathbf{Q}$ be a smooth Artin stack with $k$ C-divisor $(\mathbf{L}_i, \mathbf{s}_i)$ and $k$ positive integers $n_i$, $i = 1, \cdots, k$. We define $\mathbf{N} \subset \mathbf{Q}$ to be the substack defined by the vanishing of the sections $\mathbf{s}_1^{n_1}, \cdots, \mathbf{s}_k^{n_k}$. Then if we let $\bar{\rho}_\alpha : \bar{S}_\alpha \to \mathbf{Q}$ be a smooth chart and let $S_\alpha = \bar{S}_\alpha \times_\mathbf{Q} \mathbf{N}$, then $S_\alpha$ has a natural perfect obstruction theory taking values in the complex

$$(4.24) \qquad \mathcal{F}_\alpha^\bullet = [\mathcal{O}_{S_\alpha}(T\bar{S}_\alpha|_{S_\alpha}) \longrightarrow \oplus \mathcal{O}_{S_\alpha}(\rho_\alpha^* \mathbf{L}_i^{\otimes n_i})],$$

where $\rho_\alpha : S_\alpha \to \mathbf{N}$ is the tautological morphism. We now let $\mathbf{M}$ be a DM-stack with a morphism $\mathbf{M} \to \mathbf{N}$ and let $R_\alpha = S_\alpha \times_\mathbf{N} \mathbf{M}$. We assume $\mathbf{M}$ has a perfect obstruction theory given by $\{\mathcal{E}_\alpha^\bullet, \mathfrak{ob}_{R_\alpha}\}_\Lambda$, where $\Lambda$ is the atlas $\{R_\alpha\}_\Lambda$ of $\mathbf{M}$. We also assume $\mathbf{M}/\mathbf{N}$ has a perfect relative obstruction theory given by $\{\mathcal{L}_\alpha^\bullet, \mathfrak{ob}_{R_\alpha/S_\alpha}\}$. Finally, we assume all these obstruction theories are compatible. Hence by Lemma 4.3 we have $[\mathbf{M}]^{\mathrm{virt}} = [\mathbf{M}, \mathbf{N}]^{\mathrm{virt}}$.

What we are interested is to compare this cycle with the virtual moduli cycle of the substack $\mathbf{M}_0 \subset \mathbf{M}$ defined by $\mathbf{M}_0 = \mathbf{M} \times_\mathbf{N} \mathbf{N}_0$, where $\mathbf{N}_0 \subset \mathbf{Q}$ is defined by the vanishing of the sections $\mathbf{s}_1, \cdots, \mathbf{s}_k$. Note that $\mathbf{M}$ is homeomorphic to $\mathbf{M}_0$. Again we assume $\mathbf{M}_0$ has perfect obstruction theory so that it is compatible to the perfect obstruction theory of $\mathbf{N}_0$ and the perfect relative obstruction theory of $\mathbf{M}_0/\mathbf{N}_0$. For simplicity, we only consider the case when $\mathbf{N}_0$ is smooth and $\mathrm{Codim}(\mathbf{N}_0, \mathbf{Q}) = k$.

**Lemma 4.8.** *Suppose the relative obstruction theory of $\mathbf{M}_0/\mathbf{N}_0$ is induced from that of $\mathbf{M}/\mathbf{N}$. Then* $[\mathbf{M}, \mathbf{N}]^{\mathrm{virt}} = (\prod_{i=1}^k n_i)[\mathbf{M}_0, \mathbf{N}_0]^{\mathrm{virt}}$.

Here by the relative obstruction theory of $\mathbf{M}/\mathbf{N}$ inducing a relative obstruction theory of $\mathbf{M}_0/\mathbf{N}_0$ we mean that the relative obstruction sheaf $\mathcal{O}b_{\mathbf{M}_0/\mathbf{N}_0} = \mathcal{O}b_{\mathbf{M}/\mathbf{N}} \otimes_{\mathcal{O}_\mathbf{M}} \mathcal{O}_{\mathbf{M}_0}$ and the relative obstruction class assignment $\mathfrak{ob}_{\mathbf{M}_0/\mathbf{N}_0}(\xi) = \mathfrak{ob}_{\mathbf{M}/\mathbf{N}}(\xi)$ for any triple $\xi \in Ob(\mathfrak{Tri}_{\mathbf{M}_0/\mathbf{N}_0})$.

*Proof.* Without loss of generality, we can assume all $n_i \geq 2$ since otherwise we can replace $\mathbf{Q}$ by $\mathbf{Q} \cap \{s_i = 0 | n_i = 1\}$. Let $\bar{S}_\alpha \to \mathbf{Q}$, $S_\alpha = \bar{S}_\alpha \times_\mathbf{Q} \mathbf{N}$ and let $R_\alpha = S_\alpha \times_\mathbf{N} \mathbf{M}$ be as before. We let $(L_{i,\alpha}, s_{i,\alpha})$ be the restriction of $(\mathbf{L}_i, \mathbf{s}_i)$ to $\bar{S}_\alpha$. Then $S_\alpha = \bar{S}_\alpha \cap \{s_{i,\alpha}^{n_1} = \cdots = s_{i,\alpha}^{n_k} = 0\}$. We let $\mathcal{L}_\alpha^\bullet$ and $\mathcal{E}_\alpha^\bullet$ be the complexes of sheaves over $R_\alpha$ that are part of the (relative) obstruction theories of $\mathbf{M}/\mathbf{N}$ and $\mathbf{M}$ as stated in Definition 4.1. As before, we let $\mathcal{F}_\alpha^\bullet$ be the complex (4.24). Again, we assume that the exact sequences (4.1) hold. We now let $\bar{p} \in R_\alpha$ be any closed point and let $\bar{q} \in S_\alpha$ be the image of $\bar{p}$. We let $T_i$ and $O_i$ ($i = 1/2, 1, 2$ or $\emptyset$) be the vector spaces defined before (4.3) associated to $\bar{p}$ and the pair $R_\alpha/S_\alpha$. Since all $n_i \geq 2$, the vector space $O_2$ is $\oplus_{i=1}^k L_{i,\alpha}|_{\bar{p}}$. Also we can choose the Kuranishi map $g_p \in \Bbbk[\![T_2^\vee]\!] \otimes O_2$ to be the germ of $s_\alpha^{[n]} = (s_{1,\alpha}^{n_1}, \cdots, s_{k,\alpha}^{n_k})$. We let $W = \mathrm{Spec}\, \Bbbk[\![T_2]\!]$ and let $\hat{W}$, $X$ and $\hat{X}$ be the schemes defined after Lemma 4.2. Then since $\mathrm{Codim}(S_\alpha, \bar{S}_\alpha) = k$ and since $\bar{S}_\alpha$ is smooth, the cone $C_{\hat{W}/W}$ (which is defined after (4.7)) is the vector bundle $O_2 \times \hat{W}$ over $\hat{W}$. In this case the corresponding germ $\mathcal{D}(p)_2$, which is the normal cone to $C_{\hat{W}/W} \times_W \hat{X}$ in $C_{\hat{W}/W} \times_W \hat{W} X$, is $C_{\hat{W}/W} \times_W \hat{W} C_{\hat{X} \times_W \hat{W}/X \times_W \hat{W}}$.

Now we consider the parallel situation for the pair $\mathbf{M}_0 \to \mathbf{N}_0$. We let $S_{0,\alpha} = \bar{S}_\alpha \times_{\mathbf{Q}} \mathbf{N}_0$ and let $R_{0,\alpha} = R_\alpha \times_{\mathbf{Q}} \mathbf{N}_0$. We will denote by $\bar{p} \in R_{0,\alpha}$ and $\bar{q} \in S_{0,\alpha}$ the same points $\bar{p} \in R_\alpha$ and $\bar{q} \in S_\alpha$, via the inclusion $R_{0,\alpha} \subset R_\alpha$ and $S_{0,\alpha} \subset S_\alpha$. We then let $T_{0,i}$ and $O_{0,i}$ be the vector spaces defined before (4.3) associated to the pair $R_{0,\alpha}/S_{0,\alpha}$ over the point $\bar{p}$. By assumption, $T_{0,2} \subset T_2$, $T_{0,1/2} \equiv T_{1/2}$, $T_0 \subset T$ is a codimension $k$ linear subspace, $O_{0,2} = 0$ and $O_{0,1/2} = O_{1/2}$. Thus we have the following diagram

$$
\begin{array}{ccccccccc}
0 & \longrightarrow & O_{1/2} & \longrightarrow & O & \longrightarrow & O_2 & \longrightarrow & 0 \\
& & \| & & \downarrow & & \downarrow & & \\
0 & \longrightarrow & O_{0,1/2} & \longrightarrow & O_0 & \longrightarrow & 0 & &
\end{array}
$$

Further, we can choose the residue of $h \in \mathbb{k}[\![T^\vee]\!] \otimes O_{1/2}$ in $\mathbb{k}[\![T_0^\vee]\!] \otimes O_{1/2} \equiv \mathbb{k}[\![T_{0,1/2}^\vee]\!] \otimes O_{0,1/2}$, denoted by $h_0$, be the relative Kuranishi map of $R_{0,\alpha}/S_{0,\alpha}$ at $\bar{p}$. We let $W_0$, $\hat{W}_0$, $X_0$ and $\hat{X}_0$ be the similarly defined formal schemes associated to $\bar{q} \in S_{0,\alpha}$ and $\bar{p} \in R_{0,\alpha}$. Note that with this choice, the Kuranishi map $g_0 = 0$ and $\hat{W}_0 = W_0$. Hence the associated germ $\mathcal{D}(p)_{0,2}$ is the normal cone $C_{\hat{X}_0/X_0}$, which is a cycle in $O_{0,1/2} \times \hat{X}_0$.

Since $\hat{W}_0$ is $\hat{W}$ with the reduced scheme structure and since the relative Kuranishi map $h_0$ is the restriction of the Kuranishi map $h$ to $\hat{W}_0$, $Z_0 = Z \times_W W_0$ and $\hat{Z}_0 = \hat{Z} \times_{\hat{W}} W_0$. From this we see that the cycle $\mathcal{D}(p)_2$ is a multiple of the pullback of $\mathcal{D}(p)_{0,2}$ under the projection

$$(O_{1/2} \times O_2) \times \hat{Z} \xrightarrow{\text{proj}} O_{0,1/2} \times \hat{Z}$$

with the multiplier given by the multiplicity of $\hat{W}$ along $\hat{W}_0$, which is $\prod_{i=1}^k n_i$.

Based on this, we see that $\Xi(\mathbf{M}/\mathbf{N})$ is canonically isomorphic to $\Xi(\mathbf{M}_0/\mathbf{N}_0)$. Further, for any $a \in \Xi(\mathbf{M}/\mathbf{N})$ with the corresponding $\bar{a} \in \Xi(\mathbf{M}_0/\mathbf{N}_0)$, a representative $(A_a, F_a, \varphi_a)$ of $a$ is also a representative of $\bar{a}$. Of course their multiplicities obey $m_a = m_{\bar{a}} \prod n_i$. Therefore

$$[\mathbf{M}/\mathbf{N}]^{\text{virt}} = \sum_{a \in \Xi(\mathbf{M}/\mathbf{N})} m_a\, \xi(a) = \sum_{a \in \Xi(\mathbf{M}_0/\mathbf{N}_0)} m_a \prod_{i=1}^k n_i \cdot \xi(a) = (\prod_{i=1}^k n_i)[\mathbf{M}_0/\mathbf{N}_0]^{\text{virt}}.$$

This proves the Lemma. $\qquad\square$

### 4.3. Proof of Lemma 3.10 and 3.11.

The proof of Lemma 3.10 is similar to that of Lemma 3.11 while technically less involved. Hence we will prove Lemma 3.11 and omit the other. The strategy to prove Lemma 3.11 is to apply Lemma 4.6 to the case where $\mathbf{M} = \mathfrak{M}(\mathfrak{W}_0, \eta)$ and $\mathbf{N} = \mathfrak{M}(\mathfrak{W}, \Gamma)$. To this end, we need to work out the relative obstruction theory of $\mathfrak{M}(\mathfrak{W}_0, \eta)/\mathfrak{M}(\mathfrak{W}, \Gamma)$.

Following the argument in Section 2.1, we only need to look at the relative obstruction theory of $\mathfrak{M}(W_0[n], \eta)^{\text{st}}/\mathfrak{M}(W[n], \Gamma)^{\text{st}}$. We first cover $\mathfrak{M}(W[n], \Gamma)^{\text{st}}$ by affine étale charts $S_\alpha \to \mathfrak{M}(W[n], \Gamma)^{\text{st}}$. We let $(L_{\eta,\alpha}, s_{\eta,\alpha})$ be the restriction of $(\mathbf{L}_\eta, \mathbf{s}_\eta)$ to $S_\alpha$. By definition, $R_\alpha \triangleq S_\alpha \times_{\mathfrak{M}(W[n], \Gamma)^{\text{st}}} \mathfrak{M}(W_0[n], \eta)^{\text{st}}$ is the subscheme of $S_\alpha$ defined by the vanishing of $s_\alpha$. Hence $R_\alpha/S_\alpha$ admits an obvious relative obstruction theory induced by the pair $(L_{\eta,\alpha}, s_{\eta,\alpha})$, as defined in the first case in subsection 4.2. Namely, for any $\xi = (B, I, \varphi) \in \mathfrak{Tri}_{R_\alpha/S_\alpha}$ with $\varphi \colon \operatorname{Spec} B/I \to R_\alpha$, the relative obstruction class is

$$\mathfrak{ob}_{R_\alpha/S_\alpha}(\xi) = \mathbf{d}(s_\alpha \circ \varphi) \in \Gamma(R_\alpha, L_{\eta,\alpha}) \otimes_{\Gamma(R_\alpha)} I.$$

Hence the relative obstruction theory of $R_\alpha/S_\alpha$ takes values in the cohomology of the complex $[0 \to \mathcal{O}_{R_\alpha}(L_{\eta,\alpha})]$. Because $L_{\eta,\alpha}$ are the restriction of a global line bundle $\mathbf{L}_\eta$ on $\mathfrak{M}(W[n], \Gamma)^{\mathrm{st}}$ and the sections $s_\alpha$ over $S_\alpha$ are the restrictions of the global section $\mathbf{s}_\eta$, the collection of the relative obstruction assignments $\{\mathfrak{ob}_{R_\alpha/S_\alpha}\}$ are compatible over $R_{\alpha\beta}$ and thus defines a global relative obstruction theory of $\mathfrak{M}(W_0[n], \eta)^{\mathrm{st}}/\mathfrak{M}(W[n], \Gamma)^{\mathrm{st}}$ taking values in the cohomology of

$$[0 \to \mathcal{O}_{\mathfrak{M}(W_0[n],\eta)^{\mathrm{st}}}(\mathbf{L}_\eta)].$$

We now let $A_\alpha = \Gamma(\mathcal{O}_{S_\alpha})$ and $A_{\eta,\alpha} = \Gamma(\mathcal{O}_{R_\alpha})$. Without lose of generality, we can assume all $S_\alpha$ are $\eta$-admissible (cf. Definition 3.3). We then form the complex of $A_\alpha$-modules $\mathbf{E}_\alpha^\bullet$ as in (1.20) and the complex of $A_{\eta,\alpha}$-modules $\mathbf{E}_{\eta,\alpha}^\bullet$ as in (3.10), with the subscript $\alpha$ added to emphasize their dependence on $\alpha$. Recall that they are the respective complexes that are part of the perfect obstruction theories of $S_\alpha$ and $R_\alpha$. Further, we have the exact sequence of complexes (3.11). Thus to apply Lemma 4.6 we only need to show that the relative obstruction theory of $\mathfrak{M}(W_0[n], \eta)^{\mathrm{st}}/\mathfrak{M}(W[n], \Gamma)^{\mathrm{st}}$ is compatible to the obstruction theories of $\mathfrak{M}(W_0[n], \eta)^{\mathrm{st}}$ and of $\mathfrak{M}(W[n], \Gamma)^{\mathrm{st}}$, in the sense of Definition 4.1. Because of the exact sequence (3.11), we only need to show that for any triple $\xi = (B, I, \varphi)$ in $\mathfrak{Tri}_{R_\alpha/S_\alpha}$ we have

$$\zeta(\mathfrak{ob}_{R_\alpha/S_\alpha}(\xi)) = \mathfrak{ob}_{R_\alpha}(\xi).$$

Here $\zeta$ is the homomorphism $h^1(\mathcal{C}^{\bullet-1}(\mathcal{O}_{R_\alpha}(L_{\eta,\alpha}))) \to h^2(\mathbf{E}_{\eta,\alpha}^\bullet(\mathcal{O}_{R_\alpha}))$. But this follows directly from the construction of the respective (relative) obstruction theories. As argued before, the relative obstruction theory of $\mathfrak{M}(W_0[n], \eta)^{\mathrm{st}}/\mathfrak{M}(W[n], \Gamma)^{\mathrm{st}}$ descends to a relative obstruction theory of $\mathfrak{M}(\mathfrak{W}_0, \eta)/\mathfrak{M}(\mathfrak{W}, \Gamma)$ taking values in the cohomology of the complex $[0 \to \mathcal{O}_{\mathfrak{M}(\mathfrak{W}_0,\eta)}(\mathbf{L}_\eta)]$, and this relative obstruction theory is compatible to the obstruction theories of $\mathfrak{M}(\mathfrak{W}_0, \eta)$ and $\mathfrak{M}(\mathfrak{W}, \Gamma)$. Thus by applying the result proved in subsection 4.2, we conclude $[\mathfrak{M}(\mathfrak{W}_0, \eta)]^{\mathrm{virt}} = c_1(\mathbf{L}_\eta, \mathbf{s}_\eta)[\mathfrak{M}(\mathfrak{W}, \Gamma)]^{\mathrm{virt}}$. This proves Lemma 3.11.

### 4.4. Proof of Lemma 3.12.

We now prove Lemma 3.12. Here is our strategy. Let $\mathbf{M}_{g,n}$ be the (Artin) stack of $k$-pointed genus $g$ nodal curves and let

$$\mu \colon \mathfrak{M}(\mathfrak{W}, \Gamma) \longrightarrow \mathbf{M}_{g,n}$$

be the forgetful morphism. Let $\eta = (\Gamma_1, \Gamma_2, I) \in \Omega$ be as in Lemma 3.12 that has $r$ ordered roots of weights $\mu_1, \cdots, \mu_r$. We will show that in the formal neighborhood of $\mathfrak{M}(\mathfrak{W}_0, \eta)$ in $\mathfrak{M}(\mathfrak{W}, \Gamma)$ (possibly after an étale base change) there are divisors $(\mathbf{L}_i, \mathbf{s}_i)$ for $1 \leq i \leq r$ so that $\mathfrak{M}(\mathfrak{W}_0, \eta)$ is defined by the vanishing of $\mathbf{s}_1^{\mu_1}, \cdots, \mathbf{s}_r^{\mu_r}$ while $\mathfrak{M}(\mathfrak{Y}_1^{\mathrm{rel}} \sqcup \mathfrak{Y}_2^{\mathrm{rel}}, \eta)$ is defined by the vanishing of $\mathbf{s}_1, \cdots, \mathbf{s}_r$. This way we can reduce the proof of Lemma 3.12 to the situation studied in Lemma 4.8.

We now provide the detail of the proof. We first construct the desired base change of the formal neighborhood of $\mathfrak{M}(\mathfrak{W}_0, \eta)$ in $\mathfrak{M}(\mathfrak{W}, \Gamma)$. We let $\mathbf{M}_{\Gamma_i^\circ}$ be the moduli stack of pointed nodal curves (not necessary connected) of topological type $\Gamma_i^\circ$ (see the definition before (2.9)). Since we do not impose stability condition on such curves, $\mathbf{M}_{\Gamma_i^\circ}$ is an Artin stack. Then for any pair $(C_1, C_2) \in \mathbf{M}_{\Gamma_1^\circ} \times \mathbf{M}_{\Gamma_2^\circ}$ we can form a new curve $C_1 \sqcup C_2$ by gluing the pairs of the $i$-th distinguished marked point of $C_1$ and of $C_2$ for all $i = 1, \cdots, r$. This construction extends to families, thus defines a (local embedding) morphism $\mathbf{M}_{\Gamma_1^o} \times \mathbf{M}_{\Gamma_2^o} \to \mathbf{M}_{g,k}$.

**Definition 4.9.** *Let* $\mathbf{A}$ *be an Artin stack. We say* $\mathbf{B}$ *is a formal extension of* $\mathbf{A}$ *if* $\mathbf{A}$ *is a closed substack of* $\mathbf{B}$ *and the inclusion* $\mathbf{A} \to \mathbf{B}$ *is a homeomorphism.*

**Lemma 4.10.** *We can find a formal extension* $\mathbf{Q}$ *of* $\mathbf{M}_{\Gamma_1^o} \times \mathbf{M}_{\Gamma_2^o}$ *so that the morphism* $\mathbf{M}_{\Gamma_1^o} \times \mathbf{M}_{\Gamma_2^o} \to \mathbf{M}_{g,k}$ *extends to an étale morphism* $\mathbf{Q} \to \mathbf{M}_{g,k}$.

*Proof.* For schemes, this is the topological equivalence of étale morphisms [Mil]. Note that once such extensions exist, then they are canonical. The proof of the general case can be proved by applying this topological equivalence theorem to charts of the stacks. We will leave the details to the readers. $\square$

We next consider the gluing morphism

$$\Phi_\eta : \mathfrak{M}(\mathfrak{Y}_1^{\mathrm{rel}}, \Gamma_1) \times_{D^r} \mathfrak{M}(\mathfrak{Y}_1^{\mathrm{rel}}, \Gamma_1) \longrightarrow \mathfrak{M}(\mathfrak{Y}_1^{\mathrm{rel}} \sqcup \mathfrak{Y}_2^{\mathrm{rel}}, \eta) \subset \mathfrak{M}(\mathfrak{W}, \Gamma).$$

**Lemma 4.11.** *There is a formal extension* $\mathfrak{M}(\mathfrak{W}, \Gamma)\hat{\ }$ *of* $\mathfrak{M}(\mathfrak{Y}_1^{rel}, \Gamma_1) \times_{D^r} \mathfrak{M}(\mathfrak{Y}_1^{rel}, \Gamma_1)$ *(as DM-stack) so that the morphism* $\Phi_\eta$ *extends to an étale morphism*

$$\hat{\Phi}_\eta : \mathfrak{M}(\mathfrak{W}, \Gamma)\hat{\ } \longrightarrow \mathfrak{M}(\mathfrak{W}, \Gamma).$$

*Proof.* The proof exactly the same to Lemma 4.10, and will be omitted. $\square$

Because $\mathfrak{M}(\mathfrak{W}, \Gamma)\hat{\ }$ is homeomorphic to $\mathfrak{M}(\mathfrak{Y}_1^{\mathrm{rel}}, \Gamma_1) \times_{D^r} \mathfrak{M}(\mathfrak{Y}_1^{\mathrm{rel}}, \Gamma_1)$, the forgetful morphism

$$(4.25) \qquad \mathfrak{M}(\mathfrak{Y}_1^{\mathrm{rel}}, \Gamma_1) \times_{D^r} \mathfrak{M}(\mathfrak{Y}_2^{\mathrm{rel}}, \Gamma_2) \longrightarrow \mathbf{M}_{\Gamma_1^o} \times \mathbf{M}_{\Gamma_2^o}$$

extends to $\mathfrak{M}(\mathfrak{W}, \Gamma)\hat{\ } \to \mathbf{Q}$.

Our next task is to define the PD-divisors $(\mathbf{L}_i, \mathbf{s}_i)$ on $\mathbf{Q}$ for $i = 1, \cdots, r$ as mentioned. We let $\xi$ and $\zeta$ be the universal families over $\mathbf{M}_{g,k}$ and $\mathbf{M}_{\Gamma_1^o} \times \mathbf{M}_{\Gamma_2^o}$. We let $\xi_{\mathrm{node}}$ be the natural substack of all nodal points of the fibers of $\xi$. Then $\xi_{\mathrm{node}} \subset \xi$ is a smooth divisor. On the other hand, each fiber $\zeta_p$ (over $p \in \mathbf{M}_{\Gamma_1^o} \times \mathbf{M}_{\Gamma_2^o}$) contains $r$ ordered distinguished nodes. They define $r$ ordered distinguished sections

$$(4.26) \qquad \mathbf{n}_1, \cdots, \mathbf{n}_r : \mathbf{M}_{\Gamma_1^o} \times \mathbf{M}_{\Gamma_2^o} \longrightarrow \zeta.$$

The nodal locus $\zeta_{\mathrm{node}} \subset \zeta$ is the union of a smooth divisor with the images of these $r$ sections. We now consider the universal family $\bar{\zeta}$ over $\mathbf{Q}$. Since $\mathbf{Q} \to \mathbf{M}_{g,k}$ is étale, the nodal locus $\bar{\zeta}_{\mathrm{node}}$ is a smooth divisor in $\bar{\zeta}$. On the other hand, since $\bar{\zeta}|_{\mathbf{M}_{\Gamma_1^o} \times \mathbf{M}_{\Gamma_2^o}} \cong \zeta$, $\zeta_{\mathrm{node}} \subset \bar{\zeta}_{\mathrm{node}}$ and is a homeomorphism. Hence each $\mathrm{Im}(\mathbf{n}_i)$ is an open subset in $\bar{\zeta}_{\mathrm{node}}$. We let $\bar{\mathbf{B}}_i \subset \bar{\zeta}_{\mathrm{node}}$ be the open subscheme that contains and is homeomorphic to $\mathrm{Im}(\mathbf{n}_i)$. We let $\mathbf{B}_i \subset \mathbf{Q}$ be the image stack of $\bar{\mathbf{B}}_i$ under the projection $\bar{\zeta} \to \mathbf{Q}$. Clearly, $\mathbf{B}_i$ are smooth divisors of $\mathbf{Q}$. We define $(\mathbf{L}_i, \mathbf{s}_i)$ be the C-divisor on $\mathbf{Q}$ so that $\mathbf{B}_i = \mathbf{s}_i^{-1}(0)$.

For later application, we now give trivializations of $(\mathbf{L}_i, \mathbf{s}_i)$ on charts of $\mathbf{Q}$. let $T_\alpha \to \mathbf{Q}$ be any chart. Then $B_{i,\alpha} = T_\alpha \times_{\mathbf{Q}} \mathbf{B}_i$ is a smooth divisor in $T_\alpha$. Without loss of generality, we can assume that $B_{i,\alpha}$ is defined by the vanishing of a $u_{i,\alpha} \in \Gamma(\mathcal{O}_{T_\alpha})$. We then choose $L_{i,\alpha}$ be the line bundle over $T_\alpha$ so that $\mathcal{O}_{T_\alpha}(L_{i,\alpha}) = u_{i,\alpha}^{-1}\mathcal{O}_{T_\alpha}$ and let $s_{i,\alpha} \in \Gamma(L_{i,\alpha})$ be the constant 1. Put it differently, $e_{i,\alpha} = u_{i,\alpha}^{-1}1$ is a global holomorphic basis of $L_{i,\alpha}$ while the section $s_{i,\alpha} = 1 = u_{i,\alpha}e_{i,\alpha}$ vanishes on $B_{i,\alpha}$. In case $\tilde{u}_{i,\alpha}$ is another defining equation of $B_{i,\alpha}$, we define $(\tilde{L}_{i,\alpha}, \tilde{s}_{i,\alpha})$ similarly via a basis $\tilde{e}_{i,\alpha} = \tilde{u}_{i,\alpha}^{-1} 1$ and $\tilde{s}_{i,\alpha} = \tilde{u}_{i,\alpha} \tilde{e}_{i,\alpha}$. The transition function is via $\tilde{e}_{i,\alpha} = u_{i,\alpha} / \tilde{u}_{i,\alpha} \, e_{i,\alpha}$.

Now let $\mu_i$ be the weight of the $i$-th root of $\eta$. We define

(4.27) $\quad \mathbf{N}_0 = \{\mathbf{s}_1 = \cdots = \mathbf{s}_r = 0\} \quad \text{and} \quad \mathbf{N} = \{\mathbf{s}_1^{\mu_1} = \cdots = \mathbf{s}_r^{\mu_r} = 0\} \subset \mathbf{Q},$

both are substacks of $\mathbf{Q}$. Clearly, $\mathbf{N}_0 = \mathbf{M}_{\Gamma_1^\circ} \times \mathbf{M}_{\Gamma_2^\circ}$. Hence

$$\mathfrak{M}(\mathfrak{Y}_1^{\mathrm{rel}} \sqcup \mathfrak{Y}_2^{\mathrm{rel}}, \eta) = \mathfrak{M}(\mathfrak{W}, \Gamma)\hat{} \times_{\mathbf{Q}} \mathbf{N}_0.$$

We define

$$\mathfrak{M}(\mathfrak{W}_0, \eta)^{et} = \mathfrak{M}(\mathfrak{W}, \Gamma)\hat{} \times_{\mathbf{Q}} \mathbf{N}.$$

**Lemma 4.12.** *Let* $\pi_Q : \mathfrak{M}(\mathfrak{W}, \Gamma)\hat{} \to \mathbf{Q}$ *and* $\pi_M : \mathfrak{M}(\mathfrak{W}, \Gamma)\hat{} \to \mathfrak{M}(\mathfrak{W}, \Gamma)$ *be the tautological projections. Then we have isomorphisms of C-divisors*

$$\pi_Q^*(\mathbf{L}_i, \mathbf{s}_i)^{\otimes \mu_i} \cong \pi_M^*(\mathbf{L}_\eta, \mathbf{s}_\eta).$$

*Further,* $\mathfrak{M}(\mathfrak{W}_0, \eta)^{et} \to \mathfrak{M}(\mathfrak{W}, \Gamma)$ *factor through an étale* $\mathfrak{M}(\mathfrak{W}_0, \eta)^{et} \to \mathfrak{M}(\mathfrak{W}_0, \eta)$.

*Proof.* We cover $\mathfrak{M}(\mathfrak{W}, \Gamma)\hat{}$ by an atlas of étale charts $S_\alpha \to \mathfrak{M}(\mathfrak{W}, \Gamma)\hat{}$ indexed by $\Lambda$. To each $\alpha \in \Lambda$, we let $f_\alpha : \mathcal{X}_\alpha \to W[n_\alpha]$ be the pull back of the universal family over $\mathfrak{M}(\mathfrak{W}, \Gamma)$. We then cover $\mathbf{Q}$ by charts $T_\alpha$ indexed by the same set $\Lambda$. Without lose of generality, we can assume that the composite of $S_\alpha \to \mathfrak{M}(\mathfrak{W}, \Gamma)\hat{} \to \mathbf{Q}$ factor through $h_\alpha : S_\alpha \to T_\alpha \subset \mathbf{Q}$ and the family $\mathcal{X}_\alpha$ is the pull back of the universal family $\bar{\zeta}_\alpha$ on $T_\alpha$. As before, we let $B_{i,\alpha} = T_\alpha \times_{\mathbf{Q}} \mathbf{B}_i \subset T_\alpha$ and let $u_{i,\alpha}$ be a defining equation of $B_{i,\alpha}$. We let $(L_{i,\alpha}, s_{i,\alpha})$ be the restriction of $(\mathbf{L}_i, \mathbf{s}_i)$ to the chart $T_\alpha$. As mentioned, we can choose $e_{i,\alpha} = u_{i,\alpha}^{-1} 1$ to be a holomorphic basis of $L_{i,\alpha}$, and hence $s_{i,\alpha} = u_{i,\alpha} e_{i,\alpha}$. Now we let $n_{i,\alpha} : B_{i,\alpha} \to \bar{\zeta}_\alpha|_{B_{i,\alpha}}$ be the lift of the section $\mathbf{n}_i$ in (4.26) to $B_{i,\alpha}$. Since $\mathbf{Q}$ is homeomorphic to $\mathbf{M}_{\Gamma_1^\circ} \times \mathbf{M}_{\Gamma_2^\circ}$, such a lift exists and is unique. We then let $N_{i,\alpha} = n_{i,\alpha}(B_{i,\alpha})$ and let $\hat{N}_{i,\alpha}$ be the formal completion of $\bar{\zeta}_\alpha$ along $N_{i,\alpha}$. Without loss of generality, we can assume

$$\hat{N}_{i,\alpha} \cong \mathrm{Spec}\,\Bbbk[\![z_1, z_2]\!] \times_{\mathrm{Spec}\,\Bbbk[\![t]\!]} T_\alpha.$$

Here $\mathrm{Spec}\,\Bbbk[\![z_1, z_2]\!] \to \mathrm{Spec}\,\Bbbk[\![t]\!]$ is defined by $t \mapsto z_1 z_2$ and $T_\alpha \to \mathrm{Spec}\,\Bbbk[\![t]\!]$ is defined by $t \mapsto u_{i,\alpha}$.

We now back to the family $f_\alpha : \mathcal{X}_\alpha \to W[n_\alpha]$. Recall $\mathcal{X}_\alpha = \bar{\zeta}_\alpha \times_{\mathbf{Q}} S_\alpha$. We let $\hat{\mathcal{X}}_{i,\alpha}$ be the formal completion of $\mathcal{X}_\alpha$ along $N_{i,\alpha} \times_{\mathbf{Q}} S_\alpha$, which is isomorphic to $\mathcal{X}_\alpha \times_{\bar{\zeta}_\alpha} \hat{N}_{i,\alpha}$. By shrinking $S_\alpha$ if necessary, we can assume that there is a parameterization of a neighborhood of nodes $\mathcal{W}_\alpha \subset W[n_\alpha]$ given by

$$\psi_\alpha : \mathcal{W}_\alpha \longrightarrow \mathrm{Spec}\,\Bbbk[w_1, w_2] \otimes_{\Bbbk[t_{l_\alpha}]} \Gamma(\mathbf{A}^{n_\alpha + 1})$$

for some $l_\alpha \in [n_\alpha + 1]$, as in (1.1), so that the induced morphism $\hat{f}_\alpha : \hat{\mathcal{X}}_{i,\alpha} \to W[n_\alpha]$ factor through

(4.28) $\qquad \tilde{f}_\alpha : \hat{\mathcal{X}}_{i,\alpha} \to \mathcal{W}_\alpha \quad \text{via} \quad \tilde{f}_\alpha^*(w_j) = \beta_{j,\alpha} \cdot z_j^{\mu_j}, \quad j = 1, 2.$

Now let $(L_{\eta,\alpha}, s_{\eta,\alpha})$ be the restriction to $S_\alpha$ of $(\mathbf{L}_\eta, \mathbf{s}_\eta)$. By definition, a trivialization of $(L_{\eta,\alpha}, s_{\eta,\alpha})$ is given by $\mathcal{O}(L_{\eta,\alpha}) = t_{l_\alpha}^{-1} \mathcal{O}_{T_\alpha}$ with the basis $\epsilon_{\eta,\alpha} = t_{l_\alpha}^{-1} 1$ and $s_{\eta,\alpha} \equiv t_{l_\alpha} \epsilon_{\eta,\alpha}$. As mentioned, a trivialization of $(L_{i,\alpha}, s_{i,\alpha})$ is given by $\mathcal{O}(L_{i,\alpha}) = u_{i,\alpha}^{-1} \mathcal{O}_{S_\alpha}$ with $e_{i,\alpha} = u_{i,\alpha}^{-1} 1$ and $s_{i,\alpha} = u_{i,\alpha} e_{i,\alpha}$. We then define an isomorphism $L_{\eta,\alpha} \cong L_{i,\alpha}^{\otimes \mu_i}$ via

$$\epsilon_{\eta,\alpha} = (\beta_{\alpha,1} \beta_{\alpha,2}) e_{i,\alpha}^{\otimes \mu_i}.$$

Note that since $\beta_{\alpha,1}\beta_{\alpha,2} \in \Gamma(\mathcal{O}_{S_\alpha}^\times)$, the above identity defines an isomorphism $L_{\eta,\alpha} \cong L_{i,\alpha}^{\otimes\mu_i}$. Further, because of the relations (4.28) and $w_1 w_2 = t_{l_\alpha}$,

$$(4.29) \qquad s_{\eta,\alpha} = s_{i,\alpha}^{\otimes\mu_i}$$

under this isomorphism. Hence it defines an isomorphism of the corresponding C-divisors. It is routine to check that the so defined isomorphisms extend to an isomorphism of C-divisors

$$\pi_Q^*(\mathbf{L}_\eta, \mathbf{s}_\eta) \cong \pi_M^*(\mathbf{L}_i, \mathbf{s}_i)^{\otimes\mu_i}.$$

This proves the isomorphisms of the C-divisors.

The last statement follows directly from the construction of $\mathfrak{M}(\mathfrak{W}_0, \eta)^{et}$. $\qquad\square$

We now consider the pair $\rho : \mathfrak{M}(\mathfrak{W}_0, \eta)^{et} \to \mathbf{N}$. Since $\mathfrak{M}(\mathfrak{W}_0, \eta)^{et}$ is étale over $\mathfrak{M}(\mathfrak{W}_0, \eta)$, we can take the obstruction theory of $\mathfrak{M}(\mathfrak{W}_0, \eta)^{et}$ to be the pull back of that of $\mathfrak{M}(\mathfrak{W}_0, \eta)$. We now workout the relative obstruction theory of $\mathfrak{M}(\mathfrak{W}_0, \eta)^{et}/\mathbf{N}$ and to show that it satisfies the set up in Lemma 4.8. We will follow closely the notation developed in the beginning of Section 4.1. For convenience, we denote $\mathfrak{M}(\mathfrak{W}_0, \eta)^{et}$ by $\mathbf{M}$.

We begin with a smooth chart $\bar{S}$ of $\mathbf{Q}$ and the associated chart $S = \bar{S} \times_{\mathbf{Q}} \mathbf{N}$. We let $(L_i, s_i)$ be the restriction of $(\mathbf{L}_i, \mathbf{s}_i)$ to $\bar{S}$. Since $\mathbf{Q}$ and hence $\bar{S}$ are smooth, the pairs $(L_i, s_i)$ for $i = 1, \cdots, r$ define a natural obstruction theory of $S$ taking values in the cohomology of the complex

$$(4.30) \qquad \mathcal{F}^\bullet = [\mathcal{O}_S(T\bar{S}) \xrightarrow{ds^{[\mu]}} \oplus_{i=1}^r \mathcal{O}_S(L_i^{\otimes\mu_i})],$$

as defined in the beginning of the section 4.1. Here $ds^{[\mu]}$ is the abbreviation of $(ds_1^{\mu_1}, \cdots, ds_r^{\mu_r})$. Now let $\mathcal{X}$ be the universal family over $S$ and let $\mathcal{D} \subset \mathcal{X}$ be the divisor of the marked sections of $\mathcal{X}$. Then by the deformation theory of nodal curves, there is a canonical homomorphism of sheaves (the Kodaira map)

$$\mathcal{O}_S(T\bar{S}) \longrightarrow \mathcal{E}xt^1_{\mathcal{X}/S}(\Omega_{\mathcal{X}/S}(\mathcal{D}), \mathcal{O}_\mathcal{X}).$$

Next we pick an affine étale universal open $R \to S \times_{\mathbf{N}} \mathbf{M}$ with $\rho : R \to S$ the projection. Without loss of generality we can assume that the universal family of $R$ is of the form $f : \rho^*\mathcal{X} \to W[n]$ for some integer $n$. We now construct the standard obstruction theory of $R$. Since $R$ is a smooth chart (not necessary étale) of $\mathbf{M}$, its obstruction theory is slightly different from that defined in Section 3.1, which is for étale charts of $\mathbf{M}$.

We begin with the construction of the complex $\mathbf{E}_\eta^\bullet$ that will be part of the obstruction theory of $R$. We let $B = \Gamma(\mathcal{O}_R)$ with $A = \Gamma(\mathcal{O}_S)$ as before. First of all, we choose the complex $F^\bullet$ in (1.15) to be

$$F^0 = 0 \quad \text{and} \quad F^1 = \Gamma(R, \rho^*T\bar{S}) = \Gamma(\mathcal{O}_S(T\bar{S})) \otimes_A B.$$

We pick a collection of charts $(\mathcal{U}_\alpha/\mathcal{V}_\alpha, f_\alpha)$ of $f$ indexed by $\Lambda$, as mentioned in the paragraph after (1.17). By shrinking $R$ and/or make an étale base change of $\bar{S}$, we can assume that there is a collection of charts $\{U_\alpha/V_\alpha\}$ of $\mathcal{X}/S$ indexed by the same $\Lambda$ so that $\mathcal{V}_\alpha = R \times_S V_\alpha$ and $\mathcal{U}_\alpha = R \times_S U_\alpha$. We form the modules $\mathrm{Hom}_{\mathcal{U}_\alpha}(f^*\Omega_{W[n]}, B)^\dagger$, after picking the necessary data as shown in the paragraphs after equation (1.17). Since we have chosen $F^0 = 0$, the homomorphism $\zeta_\alpha(\cdot)$ in (1.19) is zero. The homomorphism $\zeta_{\alpha\beta}(\cdot)$ is exactly the one defined in (1.18). We then form the complex $\mathbf{D}^\bullet$, the homomorphism $\delta$ and the complex $\mathbf{E}^\bullet$, following the definitions after Lemma 1.14, line by line except that we replace $F^0$ by 0. To form the complex $\mathbf{E}_\eta^\bullet$, we shall follow the discussion after (3.8). We let $\mathbf{C}_\eta^\bullet$ be the complex defined in (3.7) and let $\delta\colon \mathbf{E}^\bullet \Rightarrow \mathbf{C}_\eta^{\bullet-1}$ be the homomorphism defined exactly as in (3.8). We then let $\mathbf{E}_\eta^\bullet$ be the complex defined in (3.10). It follows from the discussion before that the obstruction theory of $\mathbf{M}$ induces a natural obstruction theory of $R$ taking values in the cohomology of the complex $\mathbf{E}_\eta^\bullet$.

We now construct a complex $\mathbf{F}^\bullet$ of $A$-modules that is quasi-isomorphic to $\mathcal{F}^\bullet$ (in (4.30)) and a homomorphism of complexes $\mathbf{E}_\eta^\bullet \Rightarrow \mathbf{F}^\bullet$. Let $i$ be any integer between 1 and $r$. We let $\Lambda_i \subset \Lambda$ be the subset of the indices $\alpha$ so that $U_\alpha \to \mathcal{X}$ contains the $i$-th distinguished nodes of some fibers of $\mathcal{X}/S$. Then the collection $\{U_\alpha\}_{\alpha\in\Lambda_i}$ forms a covering of a neighborhood of the $i$-th distinguished nodes of the fibers of $\mathcal{X}/S$ and the collection $\{V_\alpha\}_{\Lambda_i}$ forms an étale covering of $S$. We let $\mathbf{C}_i^k = \mathcal{C}^k(\Lambda_i, L_i^{\otimes\mu_i})$ be the Čech complex of $k$-cochains of the line bundle $L_i^{\otimes\mu_i}$ associated to the covering $\{V_\alpha\}_{\Lambda_i}$. The complex $\mathbf{C}_i^k$ comes with the standard coboundary operator $\mathbf{C}_i^k \to \mathbf{C}_i^{k+1}$. We define $\mathbf{F}^1 = \Gamma(\mathcal{O}_S(T\bar{S}))$ and define $\mathbf{F}^k = \oplus_{i=1}^r \mathbf{C}_i^{k-2}$ for $k \geq 2$.

We next define the differential $\partial^k\colon \mathbf{F}^k \to \mathbf{F}^{k+1}$. Let $\mathcal{I}$ be the ideal sheaf of the zero section of the total space $T\bar{S}$ and let $\bar{S}^{(2)}$ be the subscheme of $T\bar{S}$ defined by the ideal sheaf $\mathcal{I}^2$. Then there is a tautological morphism $\kappa\colon \bar{S}^{(2)} \to \bar{S}$ characterized by the following property: Let $v \in T_p\bar{S}$ with $\mathrm{Spec}\,\Bbbk[t]/(t^2) \to T_p\bar{S}$ its associated morphism that lifts to the morphism $[v]\colon \mathrm{Spec}\,\Bbbk[t]/(t^2) \to \bar{S}^{(2)}$. Then the composite $\kappa \circ [v]\colon \mathrm{Spec}\,\Bbbk[t]/(t^2) \to \bar{S}$ is the tangent vector $v \in T_p\bar{S}$. Now let $\xi \in \Gamma(\mathcal{O}_S(T\bar{S}))$ be any element. We let $\tilde{S} = \mathrm{Spec}\,\Gamma(\mathcal{O}_S) * A$ be the trivial extension of $S$ by $A$ and let $\theta_\xi\colon \tilde{S} \to \bar{S}$ be the composite $\tilde{S} \to \bar{S}^{(2)} \xrightarrow{\kappa} \bar{S}$, where $\tilde{S} \to \bar{S}^{(2)}$ is defined by $\xi$. Since $\theta_\xi^*(s_i^{\mu_i}) \in \Gamma(\theta_\xi^* L_i^{\otimes\mu_i})$ and $\theta_\xi^*(s_i^{\mu_i})|_S \equiv 0$,

$$\mathbf{d}[\theta_\xi^*(s_i^{\mu_i}) - 0] \in \Gamma(S, \mathcal{O}_S(L_i^{\otimes\mu_i}) \otimes \mathcal{I}_{S\subset\bar{S}}) = \Gamma(S, \mathcal{O}_S(L_i^{\otimes\mu_i})).$$

Here $I_{S\subset\bar{S}}$ is the ideal of $S \subset \tilde{S}$, which is isomorphic to $A$. Now let $\rho_\alpha^* \theta_\xi^*(s_i^{\mu_i}) \in \Gamma(\mathcal{V}_\alpha, \mathcal{O}_S(L_i^{\otimes\mu_i}))$ be the pull back under $\rho_\alpha\colon \mathcal{V}_\alpha \to S$. We then define

$$\partial_i^1(\xi)_\alpha = \mathbf{d}[\rho_\alpha^* \theta_\xi^*(s_i^{\mu_i}) - 0] \in \Gamma(V_\alpha, \mathcal{O}_S(L_i^{\otimes\mu_i})), \quad \alpha \in \Lambda_i,$$

and define $\partial^1 = \oplus_i \partial_i^1$. For $k \geq 1$ we let $\partial^k\colon \mathbf{F}^k \to \mathbf{F}^{k+1}$ to be the direct sum of the coboundary operators of $\mathbf{C}_i^\bullet$. Clearly, the so defined operator $\partial^\bullet$ satisfies $\partial^k \circ \partial^{k+1} = 0$, and hence defines a complex $\mathbf{F}^\bullet = (\mathbf{F}^k, \partial^k)$. Further, it follows from our construction that the complex $\mathbf{F}^\bullet$ is quasi-isomorphic to the complex $\mathcal{F}^\bullet$ in (4.30).

We now define the promised homomorphism $\varphi^\bullet\colon \mathbf{E}_\eta^\bullet \Rightarrow \mathbf{F}^\bullet \otimes_A B$. Recall that $\mathbf{E}_\eta^1 = \Gamma(\mathcal{O}_S(T\bar{S})) \otimes_A B \oplus \mathbf{D}^0$ and $\mathbf{F}^1 = \Gamma(\mathcal{O}_S(T\bar{S}))$. The homomorphism $\varphi^1\colon \mathbf{E}_\eta^1 \to \mathbf{F}^1 \otimes_A B$ is the one induced by the identity of $\Gamma(\mathcal{O}_S(T\bar{S})) \otimes_A B$. For $k \geq 2$, we notice that $\mathbf{E}_\eta^k = \mathbf{E}^k \oplus \mathbf{C}_\eta^{k-2}$ and $\mathbf{F}^k = \oplus_{i=1}^r \mathbf{C}_i^{k-2}$. The homomorphism $\varphi^k$ will be induced by $\varphi_i^k\colon \mathbf{C}_\eta^k \to \mathbf{C}_i^k \otimes_A B$, which we define now. Let $\xi \in \mathbf{C}_\eta^k$ be any element and let $(\alpha_0 \cdots \alpha_k)$ be a $(k+1)$-tuple in $\Lambda_i$. Then $\xi_{\alpha_0 \cdots \alpha_k}$ is an element in $\Gamma(\mathcal{V}_{\alpha_0 \cdots \alpha_k}, L_\eta)$. Here $L_\eta$ is the restriction of $\mathbf{L}_\eta$ to $R$. Using the canonical isomorphism $L_\eta \cong \rho^* L_i^{\otimes\mu_i}$, where $\rho\colon R \to S$ is the tautological projection,

$$\xi_{\alpha_0 \cdots \alpha_k} \in \Gamma(\mathcal{V}_{\alpha_0 \cdots \alpha_k}, \rho^* L_i^{\otimes\mu_i}) = \Gamma(V_{\alpha_0 \cdots \alpha_k}, L_i^{\otimes\mu_i}) \otimes_A B.$$

We denote this element by $\tilde{\xi}_{\alpha_0 \cdots \alpha_k}$. We define

$$\varphi_i^k(\xi)_{\alpha_0 \cdots \alpha_k} = \tilde{\xi}_{\alpha_0 \cdots \alpha_k} \in \Gamma(V_{\alpha_0 \cdots \alpha_k}, L_i^{\otimes\mu_i}) \otimes_A B.$$

This defines a homomorphism $\mathbf{C}_\eta^k \to \oplus_{i=1}^r \mathbf{C}_i^k \otimes_A B$.

We claim that the so defined homomorphisms form a homomorphism of complexes $\mathbf{E}_\eta^\bullet \Rightarrow \mathbf{F}^\bullet \otimes_A B$. For this, we need to check the commutativity of the following diagram

(4.31)
$$
\begin{array}{ccc}
\mathbf{E}_\eta^k & \xrightarrow{d^k} & \mathbf{E}_\eta^{k+1} \\
{\scriptstyle \varphi^k} \downarrow & & \downarrow {\scriptstyle \varphi^{k+1}} \\
\mathbf{F}^k \otimes_A B & \xrightarrow{\partial^k} & \mathbf{F}^{k+1} \otimes_A B.
\end{array}
$$

We will prove the commutativity of (4.31) for the case $k = 1$. The other cases are similar and will be omitted. First, we show that (4.31) commutes on $\mathbf{D}^0 \subset \mathbf{E}_\eta^1$. Let $\xi \in \mathbf{D}^0 \subset \mathbf{E}_\eta^1$ be any element given by a collection $\{\xi_\alpha\}$ of $\xi_\alpha \in \mathrm{Hom}_{\mathcal{U}_\alpha}(f^*\Omega_{W[n]}, \mathcal{O}_S)^\dagger$, as defined in (1.11). Because of the relation (1.12), the composite $\mathbf{D}^0 \subset \mathbf{E}_\eta^1 \to \mathbf{E}_\eta^2 \xrightarrow{\varphi^2} \mathbf{F}^2 \otimes_A B$ is zero. On the other hand, by definition we have $\varphi^1(\xi) = 0$. Hence $\varphi^2(\partial^1(\xi)) = d^1(\varphi^1(\xi))$ for all $\xi \in \mathbf{D}^0$.

We next check that (4.31) commutes on $\Gamma(\mathcal{O}_S(T\bar{S})) \subset \mathbf{E}_\eta^1$. We first recall the definition of $\varphi^2 \circ d^1$. Let $\xi \in \Gamma(\mathcal{O}_S(T\bar{S}))$ be any element. It determines a morphism $\tilde{\kappa}_\xi : \tilde{R} = \mathrm{Spec}\,\Gamma(\mathcal{O}_R) * B \to \bar{S}$ that is the pull back of $\kappa_\xi : \tilde{S} \to \bar{S}$ via the tautological extension $\tilde{\rho} : \tilde{R} \to \tilde{S}$ of $\rho : R \to S$. We let $\alpha \in \Lambda_i$ be any index with $U_\alpha/V_\alpha$ the associated chart of the universal family $\mathcal{X}/S$ and with $\tilde{U}_\alpha/\tilde{V}_\alpha$ its minimal extension to the family $\tilde{\mathcal{X}}/\tilde{S}$. We let $\mathcal{U}_\alpha/\mathcal{V}_\alpha$ and $\tilde{\mathcal{U}}_\alpha/\tilde{\mathcal{V}}_\alpha$ be the corresponding pull back charts of $\rho^*\mathcal{X}/R$ and $\tilde{\rho}^*\tilde{\mathcal{X}}/\tilde{R}$. We then pick a local parameterization of the nodes of $U_\alpha/V_\alpha$ and its extension to $\tilde{U}_\alpha/\tilde{V}_\alpha$. We let $(z_{\alpha,i}, s_\alpha)$ and $(\tilde{z}_{\alpha,i}, \tilde{s}_\alpha)$ be the relevant functions associated to these parameterizations[20]. We let the parameterizations of the nodes of $\mathcal{U}_\alpha/\mathcal{V}_\alpha$ and $\tilde{\mathcal{U}}_\alpha/\tilde{\mathcal{V}}_\alpha$ be the pull back of those from $U_\alpha/V_\alpha$ and $\tilde{U}_\alpha/\tilde{V}_\alpha$. We let $f_\alpha : \mathcal{U}_\alpha \to W[n]$ be the restriction of $f$ to $\mathcal{U}_\alpha$ and let $\tilde{f}_\alpha : \tilde{\mathcal{U}}_\alpha \to W[n]$ be a pre-deformable extension of $f_\alpha$. We let $\iota_\alpha : \mathcal{U}_\alpha \to U_\alpha$, $\tilde{\iota} : \tilde{\mathcal{U}}_\alpha \to \tilde{U}_\alpha$, $j_\alpha : \mathcal{V}_\alpha \to V_\alpha$, $\tilde{j}_\alpha : \tilde{\mathcal{V}}_\alpha \to \tilde{V}_\alpha$, $\rho_\alpha : \mathcal{U}_\alpha \to V_\alpha$ and $\tilde{\rho}_\alpha : \tilde{\mathcal{U}}_\alpha \to V_\alpha$ be the tautological projections. Then after picking a local parameter of $f_\alpha(\mathcal{U}_\alpha) \subset W[n]$, say $(w_1, w_2)$ with $w_1 w_2 = t_{l_\alpha}$, we have

$$
f_\alpha^*(w_{\alpha,i}) = \iota_\alpha^*(z_{\alpha,i}^{\mu_i}) \cdot h_{\alpha,i} \quad \text{and} \quad \tilde{f}_\alpha^*(w_{\alpha,i}) = \tilde{\iota}_\alpha^*(\tilde{z}_{\alpha,i}^{\mu_i}) \cdot \tilde{h}_{\alpha,i}
$$

for some $h_{\alpha,i} \in \Gamma(\mathcal{O}_{\mathcal{U}_\alpha}^\times)$ and their extensions $\tilde{h}_{\alpha,i} \in \Gamma(\mathcal{O}_{\tilde{\mathcal{U}}_\alpha}^\times)$ that satisfy $h_{\alpha,1} h_{\alpha,2} \in \Gamma(\mathcal{O}_{\mathcal{V}_\alpha})$ and $\tilde{h}_{\alpha,1} \tilde{h}_{\alpha,2} \in \Gamma(\mathcal{O}_{\tilde{\mathcal{V}}_\alpha})$. Since $\tilde{z}_{\alpha,1} \tilde{z}_{\alpha,2} = \tilde{s}_\alpha$, we have

$$
\tilde{f}_\alpha^*(t_{l_\alpha}) = (\tilde{h}_{\alpha,1} \tilde{h}_{\alpha,2}) \tilde{\iota}_\alpha^*(\tilde{s}_\alpha^{\mu_i}).
$$

Then

$$
\varphi^2(d^1(\xi))_\alpha = \mathbf{d}[\tilde{f}_\alpha^*(t_{l_\alpha}) - 0] = (h_{\alpha,1} h_{\alpha,2}) \mathbf{d}[j_\alpha^*(\tilde{s}_\alpha^{\mu_i}) - 0] \in \Gamma(\mathcal{V}_\alpha, \rho^* L_{\eta,\alpha} \otimes I_{R \subset \tilde{R}}).
$$

Here the last relation follows because $\tilde{f}_\alpha$ is pre-deformable.

As to $d^1 \circ \varphi^1$, by definition

$$
d^1(\varphi^1(\xi))_\alpha = \mathbf{d}[\tilde{s}_\alpha^{\mu_i} - 0] \in \Gamma(V_\alpha, \rho^* L_i^{\otimes \mu_i} \otimes I_{S \subset \tilde{S}}).
$$

Therefore we have $\varphi^2(\partial^1(\xi))_\alpha = d^1(\varphi^1(\xi))_\alpha$ because of the isomorphism $L_\eta \cong L_i^{\otimes \mu_i}$, (4.31), the relation $t_{l_\alpha} = (h_{\alpha,1} h_{\alpha,2}) s_\alpha^{\mu_i}$ and $I_{R \subset \tilde{R}} \cong B$ and $T_{S \subset \tilde{S}} \cong A$. This proves the commutativity of (4.31) for $k = 1$. The case $k \geq 2$ is similar and will be omitted.

---

[20] By abuse of notation we will view $s_\alpha$ as functions on $\mathcal{V}_\alpha$ and on $\mathcal{U}_\alpha$ via the pull back $\mathcal{O}_{\mathcal{V}_\alpha} \to \mathcal{O}_{\mathcal{U}_\alpha}$. The same convention applies to $\tilde{s}_\alpha$ as well.

We define $\mathbf{G}^k = \ker\{\mathbf{E}_\eta^k \to \mathbf{F}^k \otimes_A B\}$. Since $\mathbf{E}_\eta^k \to \mathbf{F}^k \otimes_A B$ is surjective, $\mathbf{G}^k$ is a flat $\mathcal{O}_S$-module. The differentials of $\mathbf{E}_\eta^\bullet$ induces differentials of $\mathbf{G}^\bullet$ and the resulting complex fits into the following exact sequence

$$0 \Longrightarrow \mathbf{G}^\bullet \Longrightarrow \mathbf{E}_\eta^\bullet \Longrightarrow \mathbf{F}^\bullet \otimes_A B \Longrightarrow 0.$$

It is routine to check that for each $\xi = (B', I, \varphi_0) \in \mathfrak{Tri}_{S/R}$ there is a canonical obstruction class $\mathfrak{ob}_{R/S}(\xi) \in h^2(\mathbf{G}^\bullet \otimes I)$ to extending $\varphi_0 \colon \operatorname{Spec} B'/I \to R$ to an $S$-morphism $\operatorname{Spec} B' \to R$, and further such assignment satisfies the requirement in Definition 4.1. Finally, we remark that though the complexes $\mathbf{G}^\bullet$, $\mathbf{E}_\eta^\bullet$ and $\mathbf{F}^\bullet$ depend on the choice of the covering $\{\mathcal{U}_\alpha/\mathcal{V}_\alpha\}$ of $f$, they as elements in the derived category are unique. In particular the modules (sheaves) $\mathcal{O}b_{R/S} \triangleq h^2(\mathbf{G}^\bullet)$, $\mathcal{O}b_R \triangleq h^2(\mathbf{E}_\eta^\bullet)$ and $\mathcal{O}b_S \triangleq h^2(\mathbf{F}^\bullet)$ are independent of the choice of the coverings.

We now cover $\mathbf{N}$ be an atlas $\{S_\alpha\}_\Xi$ and for each $\alpha \in \Xi$ we pick an open étale $R_\alpha \to \mathbf{M} \times_\mathbf{N} S_\alpha$ so that $\{R_\alpha\}_\Xi$ forms an atlas of $\mathbf{M}$. For each $\alpha$ we pick a sufficiently fine covering of its universal family and then form the associated complexes $\mathbf{G}_\alpha^\bullet$, $\mathbf{E}_{\eta,\alpha}^\bullet$ and $\mathbf{E}_\alpha^\bullet$. Here we added the subscript $\alpha$ to indicate the dependence on the chart $R_\alpha/S_\alpha$. To be consistent with the notation in Definition 4.1, we let $\mathcal{L}_\alpha^\bullet = \mathbf{G}_\alpha^\bullet$, $\mathcal{E}_\alpha^\bullet = \mathbf{E}_{\eta,\alpha}^\bullet$ and $\mathcal{F}_\alpha^\bullet = \mathbf{F}_\alpha^\bullet$, viewed as complexes of sheaves of $\mathcal{O}_{R_\alpha}$ or $\mathcal{O}_{S_\alpha}$-modules accordingly.

**Lemma 4.13.** *There are standard relative obstruction theories of $R_\alpha/S_\alpha$ for $\alpha \in \Xi$ taking values in the complexes $\mathcal{L}_\alpha^\bullet$ such that the (relative) obstruction theories $\{\mathcal{E}_\alpha^\bullet, \mathfrak{ob}_{R_\alpha}\}$, $\{\mathcal{F}_\alpha^\bullet, \mathfrak{ob}_{S_\alpha}\}$ and $\{\mathcal{L}_\alpha^\bullet, \mathfrak{ob}_{R_\alpha/S_\alpha}\}$ are compatible in the sense of definition 4.1. Further, the so defined obstruction theories of $R_\alpha/S_\alpha$ define a relative obstruction theory of $\mathbf{M}/\mathbf{N}$ that is compatible to the obstruction theories of $\mathbf{M}$ and $\mathbf{N}$.*

*Proof.* The proof is routine and will be omitted. $\qquad\square$

Let $\mathbf{N}_0 \subset \mathbf{N}$ be as in (4.27) and let $\mathbf{M}_0 = \mathbf{M} \times_\mathbf{N} \mathbf{N}_0$. By definition $\mathbf{M}_0$ is an étale cover of $\mathfrak{M}(\mathfrak{Y}_1^{\mathrm{rel}} \sqcup \mathfrak{Y}_2^{\mathrm{rel}}, \eta)$. Note that $\mathbf{N}_0$ is smooth. As to $\mathbf{M}_0$, we endow it with the induced obstruction theory of $\mathfrak{M}(\mathfrak{Y}_1^{\mathrm{rel}} \sqcup \mathfrak{Y}_2^{\mathrm{rel}}, \eta)$, which is perfect since that of $\mathfrak{M}(\mathfrak{Y}_1^{\mathrm{rel}} \sqcup \mathfrak{Y}_2^{\mathrm{rel}}, \eta)$ is. We let $[\mathbf{M}_0]^{\mathrm{virt}}$ be the virtual moduli cycle of $\mathbf{M}_0$. Then Lemma 3.12 is equivalent to

$$(4.32) \qquad\qquad \mathfrak{m}(\eta)[\mathbf{M}_0]^{\mathrm{virt}} = [\mathbf{M}]^{\mathrm{virt}} \in A_*\mathbf{M}.$$

By Lemma 4.8, to prove this identity it suffices to show that *1) there is a relative obstruction theory of $\mathbf{M}_0/\mathbf{N}_0$ that is compatible to the obstruction theory of $\mathbf{M}_0$ and $\mathbf{N}_0$ and 2) the relative obstruction theory $\mathbf{M}_0/\mathbf{N}_0$ is compatible to the relative obstruction theory of $\mathbf{M}/\mathbf{N}$.*

The proof of 1) is parallel to the construction of the relative obstruction theory of $\mathbf{M}/\mathbf{N}$. The proof of 2) is immediate once the relative obstruction theory was constructed. Since the proof is routine, we will leave it to the readers. We state it as a Lemma.

**Lemma 4.14.** *The standard relative obstruction theory of $\mathbf{M}_0/\mathbf{N}_0$ is compatible to the obstruction theories of $\mathbf{M}_0$ and $\mathbf{N}_0$. Further, the relative obstruction theory of $\mathbf{M}_0/\mathbf{N}_0$ is induced from the obstruction theory of $\mathbf{M}/\mathbf{N}$.*

In the end, we apply Lemma 4.8 to the pairs $\mathbf{M}_0/\mathbf{N}_0 \subset \mathbf{M}/\mathbf{N}$ to conclude (4.32). This completes the proof of Lemma 3.12.

4.5. **Proof of Lemma 3.14.** It remains to prove Lemma 3.14. Let

$$(4.33) \qquad \Phi_\eta : \mathfrak{M}(\mathfrak{Y}_1^{\mathrm{rel}}, \Gamma_1) \times_{D^r} \mathfrak{M}(\mathfrak{Y}_2^{\mathrm{rel}}, \Gamma_2) \longrightarrow \mathfrak{M}(\mathfrak{Y}_1^{\mathrm{rel}} \sqcup \mathfrak{Y}_2^{\mathrm{rel}}, \eta)$$

be the étale morphism in (3.2). Using the Cartesian product (3.14), we can give $\mathfrak{M}(\mathfrak{Y}_1^{\mathrm{rel}}, \Gamma_1) \times_{D^r} \mathfrak{M}(\mathfrak{Y}_2^{\mathrm{rel}}, \Gamma_2)$ a canonical obstruction theory. We call such obstruction theory the obstruction theory induced by the Cartesian product. On the other hand, since $\Phi_\eta$ is étale, the obstruction theory of $\mathfrak{M}(\mathfrak{Y}_1^{\mathrm{rel}} \sqcup \mathfrak{Y}_2^{\mathrm{rel}}, \eta)$ induces an obstruction theory of $\mathfrak{M}(\mathfrak{Y}_1^{\mathrm{rel}}, \Gamma_1) \times_{D^r} \mathfrak{M}(\mathfrak{Y}_2^{\mathrm{rel}}, \Gamma_2)$ as well.

**Lemma 4.15.** *The two obstruction theories of $\mathfrak{M}(\mathfrak{Y}_1^{rel}, \Gamma_1) \times_{D^r} \mathfrak{M}(\mathfrak{Y}_2^{rel}, \Gamma_2)$, one defined by the Cartesian product and the other induced by that the $\mathfrak{M}(\mathfrak{Y}_1^{rel} \sqcup \mathfrak{Y}_2^{rel}, \eta)$, are identical.*

*Proof.* The proof is similar to that in [LT2], and will be omitted. ∎

We now consider the virtual moduli cycle of $\mathfrak{M}(\mathfrak{Y}_1^{\mathrm{rel}}, \Gamma_1) \times_{D^r} \mathfrak{M}(\mathfrak{Y}_2^{\mathrm{rel}}, \Gamma_2)$. Since the obstruction theory of $\mathfrak{M}(\mathfrak{Y}_1^{\mathrm{rel}}, \Gamma_1) \times_{D^r} \mathfrak{M}(\mathfrak{Y}_2^{\mathrm{rel}}, \Gamma_2)$ is induced by the Cartesian product,

$$[\mathfrak{M}(\mathfrak{Y}_1^{\mathrm{rel}}, \Gamma_1) \times_{D^r} \mathfrak{M}(\mathfrak{Y}_2^{\mathrm{rel}}, \Gamma_2)]^{\mathrm{virt}} = \Delta^! \big( [\mathfrak{M}(\mathfrak{Y}_1^{\mathrm{rel}}, \Gamma_1)]^{\mathrm{virt}} \times [\mathfrak{M}(\mathfrak{Y}_2^{\mathrm{rel}}, \Gamma_2)]^{\mathrm{virt}} \big).$$

On the other hand, because of Lemma 4.15 and because $\Phi_\eta$ is étale of pure degree $|\operatorname{Eq}(\eta)|$,

$$\frac{1}{|\operatorname{Eq}(\eta)|} \Phi_{\eta*} \Delta^! \big( [\mathfrak{M}(\mathfrak{Y}_1^{\mathrm{rel}}, \Gamma_1)]^{\mathrm{virt}} \times [\mathfrak{M}(\mathfrak{Y}_2^{\mathrm{rel}}, \Gamma_2)]^{\mathrm{virt}} \big) = [\mathfrak{M}(\mathfrak{Y}_1^{\mathrm{rel}} \sqcup \mathfrak{Y}_2^{\mathrm{rel}}, \eta)]^{\mathrm{virt}}.$$

This is exactly Lemma 3.14.

This completes the proof of the degeneration formula of the Gromov-Witten invariants stated in the beginning of this paper.

## 5. Appendix

5.1. **The tangent and the obstruction spaces.** In this appendix, we will express the first order deformation and the obstruction spaces of an $\mathfrak{M}(\mathfrak{W}, \Gamma)$ and $\mathfrak{M}(\mathfrak{Z}^{\mathrm{rel}}, \Gamma)$ in terms of some known cohomology groups. As a corollary, we will show that the obstruction theories we constructed in this paper are all perfect.

Let $S \to \mathfrak{M}(W[n], \Gamma)^{\mathrm{st}}$ be an affine étale chart. As before, we denote by $f : \mathcal{X} \to W[n]$ the universal family with $\mathcal{D} \subset \mathcal{X}$ the divisor of the ordinary marked sections. As before, we let $\pi : \mathcal{X} \to S$ be the projection and let $\rho : S \to \mathbf{A}^{n+1}$ be the morphism under $f$. By shrinking and making an étale base change, we can assume the following holds for $S$: For each $l$ the projection induced by $\pi$

$$f^{-1}(\mathbf{D}_l)_{\mathrm{red}} \longrightarrow (S \times_{\mathbf{A}^{n+1}} \mathbf{H}^l)_{\mathrm{red}}$$

is a union of $r_l$ disjoint sections $\sigma_{l,i} : (S \times_{\mathbf{A}^{n+1}} \mathbf{H}^l)_{\mathrm{red}} \to f^{-1}(\mathbf{D}_l)_{\mathrm{red}}$ for $i = 1, \cdots, r_l$. (Note that $r_l$ could be zero.) We pick an atlas $(\mathcal{U}_\alpha/\mathcal{V}_\alpha, f_\alpha)$ of $f$ indexed by $\Lambda$ so that for each $l$ and $i \in [1, r_l]$ there is exactly one and only one $\alpha$ so that $\mathcal{U}_\alpha \cap \operatorname{Im}(\sigma_{l,i}) \neq \emptyset$, and hence $\operatorname{Im}(\sigma_{l,i})$ is covered by $\mathcal{U}_\alpha$. As before, we let $\Lambda_l$ be the collection of those $\alpha$ so that $\mathcal{U}_\alpha/\mathcal{V}_\alpha$ covers $\operatorname{Im}(\sigma_{l,i})$ for some $i$. When $\alpha$ is of the second kind, we let $(z_{\alpha,1}, z_{\alpha,2}, s_\alpha)$ be the parameterization of the distinguished nodes of $\mathcal{U}_\alpha/\mathcal{V}_\alpha$. We require $s_\alpha^{m_\alpha} = g^*(t_l)$ in case $\alpha \in \Lambda_l$. With such assumptions and choices made, the standard log structure of $S$ is given by the pre-log structure $\mathcal{N}_S = \oplus N_l \to \mathcal{O}_S$ given in (1.7).

We now let $\mathbf{E}^\bullet$ and $\mathbf{D}^\bullet$ be the complexes associated to the perfect obstruction theory of $S$ constructed in section 1. Let $A = \Gamma(\mathcal{O}_S)$. For simplicity, we give an ordering of $\Lambda_l$ and thus the $r_l$ charts in $\Lambda_l$ are $\mathcal{U}_{l,1}/\mathcal{V}_{l,1}, \cdots, \mathcal{U}_{l,r_l}/\mathcal{V}_{l,r_l}$. We let $\mathbf{R}_l^\bullet$ be the complex

$$(s_{l,1}^{m_{l,1}}, \cdots, s_{l,r_l}^{m_{l,r_l}}) : \mathcal{O}_S^{\oplus r_l} \longrightarrow (\mathcal{O}_S^{\oplus r_l}/\mathcal{O}_S) \otimes_{\mathcal{O}_{\mathbf{A}^{n+1}}} \mathcal{O}_{\mathbf{A}^{n+1}}(\mathbf{H}^l),$$

where $\mathcal{O}_S^{\oplus r_l}/\mathcal{O}_S$ is the quotient of $\mathcal{O}_S^{\oplus r_l}$ by the diagonal $\mathcal{O}_S \hookrightarrow \mathcal{O}_S^{\oplus r_l}$. In case $r_l = 0$, we agree $\mathbf{R}_l^\bullet = [A \to 0]$. Here the complex $\mathbf{R}_l^\bullet$ is indexed at $[0,1]$.

**Proposition 5.1.** *For any $A$-module $I$, we have the following two exact sequences*

$$0 \longrightarrow \mathrm{Ext}_{\mathcal{X}}^0(\Omega_{\mathcal{X}/S}(\mathcal{D}), \mathcal{I}) \longrightarrow H^0(\mathbf{D}^\bullet \otimes_A I) \longrightarrow H^1(\mathbf{E}^\bullet \otimes_A I) \longrightarrow$$
$$\longrightarrow \mathrm{Ext}_{\mathcal{X}}^1(\Omega_{\mathcal{X}/S}(\mathcal{D}), \mathcal{I}) \longrightarrow H^1(\mathbf{D}^\bullet \otimes_A I) \longrightarrow H^2(\mathbf{E}^\bullet \otimes_A I) \longrightarrow 0$$

$$0 \longrightarrow H^0(\mathcal{H}om(f^*\Omega_{W[n]^\dagger/\mathbf{A}^{n+1\dagger}}, \mathcal{I})) \longrightarrow H^0(\mathbf{D}^\bullet \otimes_A I) \xrightarrow{b_0} \oplus_{l=1}^{n+1} H^0(\mathbf{R}_l^\bullet \otimes_{\mathcal{O}_S} I) \longrightarrow$$
$$\xrightarrow{\delta} H^1(\mathcal{H}om(f^*\Omega_{W[n]_0^\dagger/\mathbf{A}^{n+1\dagger}}, \mathcal{I})) \longrightarrow H^1(\mathbf{D}^\bullet \otimes_A I) \xrightarrow{b_1} \oplus_{l=1}^{n+1} H^0(\mathbf{R}_l^\bullet \otimes_{\mathcal{O}_S} I) \longrightarrow 0$$

*and the vanishing $H^i(\mathbf{E}^\bullet) = 0$ for $i > 2$ and $H^i(\mathbf{D}^\bullet) = 0$ for $i \geq 2$. Here $\mathcal{I} = I \otimes_A \mathcal{O}_{\mathcal{X}}$.*

*Proof.* The first exact sequence follows directly from the construction of the complexes $\mathbf{D}^\bullet$ and $\mathbf{E}^\bullet$. We now prove the second exact sequence. We first construct the arrow $b_0$. Let $\xi \in H^0(\mathbf{D}^\bullet \otimes I)$ be any element. Then $\xi$ is represented by $\{\xi_\alpha\}$ where $\xi_\alpha \in \mathrm{Hom}_{\mathcal{U}_\alpha}(f^*\Omega_{W[n]}, I)^\dagger$. In case $\alpha \in \Lambda_l$ then

$$(5.1) \qquad \xi_\alpha = (\varphi_\alpha, \eta_{\alpha,1}, \eta_{\alpha,2}) \in \mathrm{Hom}_{\mathcal{U}_\alpha}(f^*\Omega_{W[n]}, \mathcal{I}_\alpha) \oplus \mathcal{I}_\alpha^{\oplus 2}, \ \mathcal{I}_\alpha = \mathcal{I} \otimes_{\mathcal{O}_{\mathcal{X}}} \mathcal{O}_{\mathcal{U}_\alpha},$$

following the convention in (1.10). Specifically, the $\eta$. means $dw_i/w_i \mapsto \eta_{\alpha,i}$, under the appropriate parameterization $(w_1, w_2)$ of $f_\alpha(\mathcal{U}_\alpha)$. Recall $\Lambda_l$ was given an ordering, thus for $\alpha = i \in \Lambda_l$ we denote $\eta_{(l,i)} = \eta_{\alpha,1} + \eta_{\alpha,2}$. Then we define

$$b_0(\xi)_l = (\eta_{(l,1)}, \cdots, \eta_{(l,r_l)}) \in A^{\oplus r_l}$$

and then the arrow $b_0$ is

$$b_0(\xi) = (b_0(\xi)_1, \cdots, b_0(\xi)_{r_l}).$$

In case $r_l = 0$, then we pick an $\alpha \in \Lambda$ and define $\eta_l = \varphi_\alpha(dt_l)$, where $\varphi_\alpha$ is part of the $\xi_\alpha$ as in (5.1). It is direct to check that if $\xi$ is a cohomology class then $b_0(\xi)_l \in H^0(\mathbf{R}_l^\bullet \otimes I)$. This defines the corresponding exact sequence in the second exact sequence.

We next construct the arrow $H^1(\mathbf{D}^\bullet) \to H^1(\mathbf{R}^\bullet)$. Let $\Sigma_l = f^{-1}(\mathbf{D}_l) \subset \mathcal{X}$ and let $I_{\Sigma_l \subset \mathcal{X}}$ be the relative locally constant ideal sheaf defined by $\Gamma(\mathcal{U}_\alpha, I_{\Sigma_l \subset \mathcal{X}}) = \Gamma(\mathcal{V}_\alpha, \mathcal{O}_{\mathcal{V}_\alpha})$ in case $\mathcal{U}_\alpha \cap f^{-1}(\mathbf{D}_l) = \emptyset$ and is $s_\alpha^{m_\alpha} \Gamma(\mathcal{V}_\alpha, \mathcal{O}_{\mathcal{V}_\alpha})$ otherwise. Similarly, we let $\pi^{-1}\mathcal{O}_S$ be the pull back sheaf, namely $\Gamma(\mathcal{U}_\alpha, \pi^{-1}\mathcal{O}_S) = \Gamma(\mathcal{V}_\alpha, \mathcal{O}_{\mathcal{V}_\alpha})$. Then we have the exact sequence in étale site

$$0 \longrightarrow I_{\Sigma_l \subset \mathcal{X}} \longrightarrow \pi^{-1}\mathcal{O}_S \longrightarrow \oplus_{\alpha \in \Lambda_l}(\mathcal{O}_S/s_\alpha^{m_\alpha}\mathcal{O}_S) \longrightarrow 0$$

and its induced exact sequence in cohomologies (of $\mathcal{O}_S$-modules)

$$0 \longrightarrow g^*(\mathcal{O}_S) \to \mathcal{O}_S \longrightarrow \oplus_{\alpha \in \Lambda_l}(\mathcal{O}_S/s_\alpha^{m_\alpha}\mathcal{O}_S) \longrightarrow H^1(\mathbf{R}_l^\bullet) \longrightarrow 0.$$

Now we back to the complex $\mathbf{D}^\bullet$. It follows from the construction of $\mathbf{D}^\bullet$ that there is a homomorphism of complexes $\mathbf{D}^\bullet \to \mathbf{C}^\bullet(\Lambda, I_{\Sigma_l \subset \mathcal{X}})$. Thus we have a homomorphism of cohomologies $b_{1,l} : H^1(\mathbf{D}^\bullet \otimes_A I) \to H^1(\mathbf{R}_l^\bullet \otimes_A I)$ for each $l$. This forms the arrow in the second exact sequence of the Proposition. The arrow $\delta$ in the second sequence is the ordinary connecting homomorphism. We leave it to readers to check that with these arrows the second sequence is exact.

The vanishing result stated in Lemma 1.19 follows from these exact sequences. $\qquad\square$

The tangent and the obstruction to deformation of $\mathfrak{M}(\mathfrak{Z}^{\mathrm{rel}}, \Gamma)$ is similar. We let $S \to \mathfrak{M}(Z[n]^{\mathrm{rel}}, \Gamma)^{st}$ be an affine étale chart, satisfying the similar property as in the case just studied. We let $f \colon \mathcal{X} \to Z[n]$ be the universal family with over $S$ with $\mathcal{D} \subset \mathcal{X}$ be the divisor of the union of all ordinary and the distinguished marked sections. We let $(z_{\alpha,1}, z_{\alpha,2}, s_\alpha)$ be the parameterization of charts $\mathcal{U}_\alpha/\mathcal{V}_\alpha$ of $f$ as in the previous case. We let $\Lambda_l$ be those $\alpha$ such that $\mathcal{U}_\alpha \cap f^{-1}(\mathbf{B}_l) \neq \emptyset$. Let $\mathbf{E}^\bullet$ and $\mathbf{D}^\bullet$ be the complexes constructed in the section 1 that is part of the perfect obstruction theory of $S$.

**Proposition 5.2.** *The two exact sequences in Proposition 5.1 still hold with the sheaf $\Omega_{W[n]^\dagger/\mathbf{A}^{n+1\dagger}}$ replaced by $\Omega_{Z[n]^\dagger/\mathbf{A}^{n\dagger}}$. The same vanishing results hold as well.*

We will close this section by working out the obstruction sheaf of an example suggested to us by E. Ionel.

Let $Z_{\mathrm{rel}} = (Z, D)$ be a pair of smooth variety and a smooth divisor. We let $\Gamma$ be the graph consisting of one vertex and one leg. We assign the weights of the vertex to be $g = 1$ and $d = 0$. Thus $\mathfrak{M}(\mathfrak{Z}^{\mathrm{rel}}, \Gamma)$ is the moduli of relative stable morphisms to $Y$ from 1-pointed genus 1 curves to $Z$ of degree 0. Since $d = 0$, all $f : X \to Z$ in $\mathfrak{M}(\mathfrak{Z}^{\mathrm{rel}}, \Gamma)$ are constant maps. Hence $\mathfrak{M}(\mathfrak{Z}^{\mathrm{rel}}, \Gamma)$ is isomorphic to $\mathfrak{M}_{1,1} \times Z$. We now show that its obstruction sheaf is

$$(5.2) \qquad \mathcal{O}b = \pi_2^* \Omega_Z(\log D)^\vee,$$

where $\pi_2 : \mathfrak{M}_{1,1} \times Z \to Z$ is the second projection.

Let $f_0 \in \mathfrak{M}(\mathfrak{Z}^{\mathrm{rel}}, \Gamma)$ be a relative stable morphism. Since $d = 0$, we can always represent $f_0$ by a morphism $f_0 : X \to Z[1]^o$, where $Z[1]^o = Z[1] - D[1] \cup Z[1]_{0,\mathrm{sing}}$ with $Z[1]_{0,\mathrm{sing}}$ is the singular locus of $Z[1]_0$. Then $Z[1]^o/\mathbb{C}^* \equiv Z$. The obstruction to deforming $f_0$ as morphism to $Z[1]^o$ is $H^1(f_0^* T_{Z[1]^o/\mathbf{A}^1}) \equiv T_{Z[1]^o/\mathbf{A}^1}|_{f_0(X)}$, where $Z[1]^o \to \mathbf{A}^1$ is the tautological projection and $T_{Z[1]^o/\mathbf{A}^1}$ is the relative tangent bundle. We let $f : \mathcal{X} \to Z[1]^o$ be the family over $\mathfrak{M}_{1,1} \times Z[1]^o$ so that $\mathcal{X}$ is the pull back of the universal family over $\mathfrak{M}_{1,1}$ while the morphism $f$ is the composite of the projection $\mathcal{X} \to \mathfrak{M}_{1,1} \times Z[1]^o$ with the second projection $\mathfrak{M}_{1,1} \times Z[1]^o \to Z[1]^o$. Clearly, $f$ is the universal family of $\mathfrak{M}(Z[1]^{\mathrm{rel}}, \Gamma)^{st}$. The obstruction bundle to the moduli space $\mathfrak{M}(Z[1]^{\mathrm{rel}}, \Gamma)^{st}$ over $\mathfrak{M}_{1,1} \times Z[1]^o$ is $p_2^* T_{Z[1]^o/\mathbf{A}^1}$, where $p_2$ is the second projection of $\mathfrak{M}_{1,1} \times Z[1]^o$. The $\mathbb{C}^*$-action lifts canonically to $p_2^* T_{Z[1]^o/\mathbf{A}^1}$ and the obstruction sheaf of $\mathfrak{M}(\mathfrak{Z}^{\mathrm{rel}}, \Gamma)$ is the descent of $p_2^* T_{Z[1]^o/\mathbf{A}^1}$. It is direct to check that under the quotient map $\mathfrak{M}_{1,1} \times Z[1]^o/\mathbb{C}^* \cong \mathfrak{M}_{1,1} \times Z$, equivariant part $(p_2^* T_{Z[1]^o/\mathbf{A}^1})^{\mathbb{C}^*}$ is canonically isomorphic to $\pi_2^* \Omega_Z(\log D)^\vee$. This proves the identity (5.2).

## 5.2. Local and formal pre-deformable morphisms.

In this part we prove that formal pre-deformable morphisms are automatically local pre-deformable.

**Lemma 5.3.** *Let $\hat{L} = (\mathbb{k}[\![z_1, z_2]\!] \otimes_{\mathbb{k}[\![s]\!]} A)^\hat{}$, be as in Definition 1.2. Suppose there are units $f_1, f_2, g_1, g_2 \in \hat{L}$ and an integer $n \geq 1$ so that $f_1 f_2$ and $g_1 g_2 \in \hat{A}$ and that $z_1^m f_i = z_1^m g_i$ for $i = 1, 2$. Then $f_1 = g_1$ and $f_2 = g_2$.*

*Proof.* We let $r_i = f_i/g_i$. Then $z_i^m(r_i - 1) = 0$ for $i = 1$ and $2$. Recall that elements in $\hat{L}$ have unique normal form $a_0 + \sum_{i \geq 1}(a_i z_1^i + b_i z_2^i)$ for $a_i, b_i \in \hat{A}$ (cf. [Li]). By the uniqueness of the normal form, the normal form of $r_1$ (resp. $r_2$) must be of the form $r_1 = 1 + \sum_{j \geq 1} a_j z_2^j$ (resp. $r_2 = 1 + \sum_{j \geq 1} b_j z_1^j$). Then $f_1 f_2/g_1 g_2 = r_1 r_2 \in \hat{A}$ implies that $1 + \sum a_j z_2^j = \epsilon(1 + \sum b_j z_1^j)^{-1}$ for some unit $\epsilon \in \hat{A}$, which is impossible unless all $a_j$ and all $b_j$ are zero. This proves the uniqueness Lemma. $\square$

**Lemma 5.4.** *The notion of pure contact is independent of the choice of the charts of the nodes of $\mathcal{U}/\mathcal{V}$.*

*Proof.* Let $\phi$ in (1.4) and $\tilde{\phi} : \mathbb{k}[\tilde{z}_1, \tilde{z}_2] \otimes_{\mathbb{k}[s]} A \to R$ be two charts of the nodes of $\mathcal{U}/\mathcal{V}$. Without loss of generality we can assume that the vanishing locus of $\phi(z_1)$ and $\tilde{\phi}(\tilde{z}_1)$ are identical in $\mathrm{Spec}\, R$. Then by the proof of [Li, Lemma 2.9], $c_i = z_i/\tilde{z}_i \in R_{\mathcal{S}}$ and hence $c_1 c_2 \in A_{\mathcal{T}}$. Now let $\varphi : \mathbb{k}[w_1, w_2] \to R$ be of pure contact with respect to the chart $\phi$. If $\varphi(w_i) = z_i^m h_i$ in $R_{\mathcal{S}}$ for $h_1, h_2 \in R_{\mathcal{S}}$ satisfying $h_1 h_2 \in A_{\mathcal{T}}$, then $\varphi(w_i) = \tilde{z}_i^m(c_i^m h_i)$ in $R_{\mathcal{S}}$ and $c_1^m h_1 c_2^m h_2 \in A_{\mathcal{T}}$ as well. Thus $\varphi$ is of pure contact with respect to $\tilde{\phi}$ as well. This proves the Lemma. $\square$

We now state and prove the following equivalence result.

**Lemma 5.5.** *Let the notation be as in Definition 1.4. Then $\varphi$ is of pure contact if and only if it is formally of pure contact.*

*Proof.* Clearly, if $\varphi$ is of pure $m$-contact, then it is so formally. We now prove the other direction. Assume $\varphi$ is formally of pure $m$-contact. Then there are $\beta_1$ and $\beta_2 \in \hat{R}$ such that $\hat{\varphi}(w_i) = z_i^m \beta_i$, where as usual $\hat{\varphi} \colon \Bbbk[w_1, w_2] \to \hat{R}$ is the homomorphism induced by $\varphi$. We first show that there are $f_1$ and $f_2 \in R_{\mathcal{S}}$ so that $\varphi(w_i) = z_i^m f_i$ in $R_{\mathcal{S}}$ for $i = 1$ and 2. Let $x = \varphi(w_1) \in R$, then $\hat{x} = \hat{\varphi}(w_1) \in z_1^m \hat{R}$. Thus by [Ma, Thm 8.1], the residue class of $\bar{x} \in R_{\mathcal{S}}/(z_1^m)$ is $(\bar{x})^{\wedge} = 0 \in (R/(z_1^m))^{\wedge}$. Hence $\bar{x} \in \cap_{m \geq 1} I^m R/(z_1^m)$, and by [Ma, Thm 8.9] there is an $a \in R/(z_1^m)$ satisfying $a \equiv 1 \mod I$ such that $a\bar{x} = 0$. Then $a \in \mathcal{S}$ and by our assumption $a$ is a unit in $R_{\mathcal{S}}/(z_1^m)$. Hence $\bar{x} = 0$ in $R_{\mathcal{S}}/(z_1^m)$. This proves that $\varphi(w_1) \in z_1^m R_{\mathcal{S}}$ and hence there is an $f_1 \in R_{\mathcal{S}}$ such that $\varphi(w_1) = z_1^m f_1$ in $R_{\mathcal{S}}$. For the same reason, $\varphi(w_2) = z_2^m f_2$ in $R_{\mathcal{S}}$ for some $f_2 \in R_{\mathcal{S}}$.

We next prove the following induction hypothesis: *For any non-negative integer $k$, there are $g_1$ and $g_2 \in R_{\mathcal{S}}$ such that*

$$(5.3) \qquad z_1^m f_1 - z_1^m g_1, \quad z_2^m f_2 - z_2^m g_2 \in s^k R_{\mathcal{S}} \quad \text{and} \quad g_1 g_2 \in A_{\mathcal{T}} + s^k R_{\mathcal{S}}.$$

Clearly, this statement is true for $k = 0$. We now show that this statement is true for $k + 1$ if it is true for $k$. Let $g_1$ and $g_2$ be elements in $R_{\mathcal{S}}$ satisfying (5.3) for an integer $k$. To carry out the induction we need to find $r_1$ and $r_2 \in R_{\mathcal{S}}$ so that

$$(5.4) \qquad z_i^m f_i - z_i^m (g_i + r_i) \in s^{k+1} R_{\mathcal{S}} \quad \text{for} \quad i = 1, 2$$

and

$$(5.5) \qquad (g_1 + r_1)(g_2 + r_2) \in A_{\mathcal{T}} + s^{k+1} R_{\mathcal{S}}.$$

Since $\varphi$ formally is of pure $m$-contact, there are units $\eta_1$ and $\eta_2 \in \hat{R}_{\mathcal{S}}$ such that $z_i^m \hat{f}_i = z_i^m \eta_i$ for $i = 1$ and 2 and $\eta_1 \eta_2 \in \hat{A}$. Because $\hat{g}_i$ and $\eta_i$ satisfy the relation

$$z_i^m \hat{g}_i = z_i^m \eta_i \mod s^k \quad \text{and} \quad \hat{g}_1 \hat{g}_2, \eta_1 \eta_2 \in \hat{A} + s^k \hat{R},$$

by the uniqueness Lemma 5.3, $\hat{g}_i \equiv \eta_i \mod s^k$. Hence

$$(z_i^m f_i - z_i^m g_i)^{\wedge} = z_i^m \hat{f}_i - z_i^m \hat{g}_i = z_i^m (\eta_i - \hat{g}_i) \in z_i^m s^k \hat{R}.$$

As we argued in the existence of $f_i$, this implies that $z_i^m (f_i - g_i) \in z_i^m s^k R_{\mathcal{S}}$.

We now let $l_i \in R_{\mathcal{S}}$ be such that $z_i^m f_i - z_i^m g_i = z_i^m s^k l_i$. We let $r_1 = s^k(l_1 + z_2 h_1)$ and $r_2 = s^k(l_2 + z_1 h_2)$ with $h_1$ and $h_2 \in R_{\mathcal{S}}$ to be determined. Clearly, (5.4) will hold with this choice of $r_1$ and $r_2$. As for (5.5), we need

$$(g_1 + r_1)(g_2 + r_2) = g_1 g_2 + s^k(g_1 l_2 + g_1 z_1 h_1 + g_2 l_1 + g_2 z_2 h_2) \in A_{\mathcal{T}} + s^{k+1} R_{\mathcal{S}}.$$

Since $g_1 g_2 \in A_{\mathcal{T}} + s^k R_{\mathcal{S}}$, we can find $c \in R_{\mathcal{S}}$ so that $g_1 g_2 - cs^k \in A_{\mathcal{T}} + s^{k+1} R_{\mathcal{S}}$. Hence we need to find $h_1$ and $h_2$ in $R_{\mathcal{S}}$ so that

$$(c + g_1 l_2 + g_2 l_1) + g_1 h_1 z_1 + g_2 h_2 z_2 \in A_{\mathcal{T}} + s R_{\mathcal{S}}.$$

Let $d = c + g_1 l_2 + g_2 l_1$. Clearly, there are $\gamma_1$ and $\gamma_2 \in \hat{R}$ and $\alpha \in A_{\mathcal{T}}$ so that $\hat{d} - (\alpha + \gamma_1 z_1 + \gamma_2 z_2) \in s R_{\mathcal{S}}$. Hence following the argument for the existence of $f_i$, there are $c_1$ and $c_2 \in R_{\mathcal{S}}$ so that $d - (c_1 z_1 + c_2 z_2) \in A_{\mathcal{T}} + s R_{\mathcal{S}}$. Hence the choice $h_1 = -c_1 g_1^{-1}$ and $h_2 = -c_2 g_2^{-1}$ will do the job for (5.5). Here $g_1$ and $g_2$ are units since $\beta_1$ and $\beta_2$ are units. This proves that for each $k$ we can find $g_1$ and $g_2$ that satisfy (5.3).

Now we show that there are $h_1$ and $h_2 \in R_{\mathcal{S}}$ as required by the Lemma. We first let $M \subset R_{\mathcal{S}}$ be the set of those elements that are annihilated by some power of $s$. It is an ideal, and since $R_{\mathcal{S}}$ is Noetherian, there is an $N$ so that $s^N M = 0$. We let $g_1$ and $g_2$ be the pair satisfying (5.3) for $k = N + 1$. In case $R/M = 0$ then the Lemma is already proved. Now assume $R/M \neq 0$. We consider the ring $R_{\mathcal{S}}/M_{\mathcal{S}}$. Since $\hat{\varphi}(t) = \hat{\varphi}(w_1)\hat{\varphi}(w_2) = \epsilon s^m$ in $\hat{R}$ for some unit $\epsilon \in \hat{R}$, $R_{\mathcal{S}}/M_{\mathcal{S}}$ is flat over $\Bbbk[t]$. Now consider the homomorphism $\bar{\varphi} \colon \Bbbk[w_1, w_2] \longrightarrow R_{\mathcal{S}}/M_{\mathcal{S}}$ induced by $\varphi$. By [Li, Prop. 2.2], it is formally of pure contact. Hence there are $\bar{h}_1$ and $\bar{h}_2 \in R_{\mathcal{S}}/M_{\mathcal{S}}$ so that $\bar{\varphi}(w_i) = z_i^m \bar{h}_i$ in $R_{\mathcal{S}}/M_{\mathcal{S}}$ for $i = 1$ and 2. Because $R_{\mathcal{S}}/M_{\mathcal{S}}$ is flat over $\Bbbk[t]$, $\bar{\varphi}_1(w_1)\bar{\varphi}_2(w_2) \in R_{\mathcal{S}}/M_{\mathcal{S}}$ implies that $\bar{h}_1 \bar{h}_2 \in R_{\mathcal{S}}/M_{\mathcal{S}}$. We then apply the uniqueness Lemma 5.3 to conclude that the residue classes of $\bar{h}_i$ and $g_i$ in $R_S/(M_S, s^k)$ are identical. Hence we can find $h_i \in R_S$ so that its residue class in $R_S/M_S$ and $R_S/(s^k)$ are $\bar{h}_i$ and $g_i$ respectively. Therefore we have $\varphi(w_i) = z_i^m h_i$ in $R_S$ and $h_1 h_2 \in A_T$. This completes the proof of the Lemma. $\qquad\square$

Department of Mathematics, Stanford University, Stanford, CA 94305
    *E-mail address*: `jli@math.stanford.edu`